\documentclass[twoside,english]{elsarticle}
\usepackage[LGR,T1]{fontenc}
\usepackage[latin9]{inputenc}
\usepackage{array}
\usepackage{verbatim}
\usepackage{rotating}
\usepackage{float}
\usepackage{calc}
\usepackage{units}
\usepackage{bbding}
\usepackage{multirow}
\usepackage{amsthm}
\usepackage{amsbsy}
\usepackage{amstext}
\usepackage{amssymb}
\usepackage{arcs}
\usepackage{graphicx}
\usepackage{xcolor}
\usepackage{caption}
\makeatletter

\DeclareRobustCommand{\greektext}{%
  \fontencoding{LGR}\selectfont\def\encodingdefault{LGR}}
\DeclareRobustCommand{\textgreek}[1]{\leavevmode{\greektext #1}}
\DeclareFontEncoding{LGR}{}{}
\DeclareTextSymbol{\~}{LGR}{126}
\newcommand{\binom}[2]{{#1 \choose #2}}

\providecommand{\tabularnewline}{\\}

\theoremstyle{plain}
\newtheorem{thm}{\protect\theoremname}
\theoremstyle{plain}
\newtheorem{conjecture}[thm]{\protect\conjecturename}

\journal{Example: Nuclear Physics B}

\makeatother

\usepackage{babel}
\providecommand{\conjecturename}{Conjecture}
\providecommand{\theoremname}{Theorem}

\begin{document}
\journal{{\tt arXiv}\it}

\title{On the nature of Mersenne fluctuations\tnoteref{t1,t2}}

\tnotetext[t1]{This document deepens aspects addressed in a previous article titled
``Parafermi algebra and interordinality'' (see \cite{Merkel}).}

\tnotetext[t2]{In ``Parafermi algebra and interordinality'', the central theme
was the implications raised by the special case that two parafermi
algebras are of Mersenne-wise neighboring orders. The present document
is meant to be largely self-contained, but an in-depth study of the
previous work is helpful and therefore recommended.}

\author{U.~Merkel\corref{cor1}}

\ead{merkel.u8@googlemail.com}

\cortext[cor1]{Corresponding author}

\address[rvt]{Universit{\"a}tsstr. 38, 70569 Stuttgart, Germany}

\address{}
\begin{abstract}
In Part I, crotons are introduced, multifaceted pre-geometric objects
that occur both as labels encoded on the boundary of a ``volume''
and as complementary aspects of geometric fluctuations within that
volume. If you think of crotons as linear combinations, then the scalars
used are croton base numbers. Croton base numbers can be combined
to form the amplitudes and phases of Mersenne fluctuations which,
in turn, form qphyla. Volume normally requires space or space-time
as a prerequisite; in a pregeometric setting, however, ``volume''
is represented by a qphyletic assembly. Various stages of pre-geometric
refinement, expressed through the aspects crotonic amplitude or phase,
combine to eventually form and/or dissolve sphere-packed chunks of
Euclidean space. A time-like crotonic refinement is a rough analog
of temporal resolution in tenacious time, whereas space-like crotonic
refinement is analogous to spatial resolution in sustained space.
The analogy suggests a conceptual link between the ever-expanding
scope of Mersenne fluctuations and the creation and lifetime patterns
of massive elementary particles. A three-stage process of ideation,
organization and intraworldly action is introduced to back this up.
In Part II, the intrawordly aspect is analyzed first, including our
preon model of subnuclear structure, and the organizer aspect thereafter,
based on three types of Mersenne numbers, $M_{\textrm{reg}}$, $M_{\nicefrac{5}{8}}$,
$M_{\nicefrac{9}{8}}$, and two formal principles: juxtaposition $x$
vs. $\! x\!\pm\!1(2)$ and the interordinal application of functional
{\normalsize{$1\!\divideontimes(f^{(a)}\!\circ(f^{(b)}\!\divideontimes f^{(c)}))$.}}{\normalsize \par}\end{abstract}
\maketitle

\part{\ }

\section{Introduction}

Crotons are pregeometric objects that emerge both as labels encoded
on the boundary of a ``volume'' and as complementary aspects of
geometric fluctuations within that volume. To express their multifacetedness,
the name \emph{croton} was chosen, after Crotos, son of Pan and Eupheme,
who, once a mortal $3D$ being, was put in sky by Muses as the celestial
fixture Sagittarius. The term volume is normally linked to the categories
space or space-time. In a pre-geometric setting, more basic categories
are needed $-$ Mersenne fluctuations and qphyla. Both require various
stages of \emph{pregeometric refinement} which, expressed through
the complementary aspects croton amplitude and phase, combine to eventually
form -- or dissolve -- real geometric objects. Advancing from mark
$n$ to $n+1$ thus, in what follows, means a time-like refinement
$2^{-n}c\;\mapsto\;2^{-n-1}c$ (roughly the analog of an exponential
increase of temporal resolution in tenacious time), and an increase
from $\alpha$ to mark $\alpha+1$ a space-like refinement ${\displaystyle {\displaystyle {\textstyle \frac{a}{b_{\alpha}}}}}\;\mapsto\;{\displaystyle {\displaystyle {\textstyle \frac{a}{b_{\alpha}+\frac{1}{b_{\alpha+1}}}}}}$
$(b_{\alpha},b_{\alpha+1}>0)$ (analoguous to increase of spatial
resolution in sustained space). On the boundary, these increases find
expression in additionally encoded labels.

\noindent In a previous work \cite{Merkel}, basic croton components
have been identfied, though at the time the name croton was not yet
used. The starting point was the equivalence between a Mersennian
identity, destilled from the special case that two parafermi algebras
\cite{Green} are of neighboring orders $p=2^{i}-1$, $p'=p^{i+1}\!-1$
(order marked by parenthesized superscript):\vspace{-0.25cm}

\begin{equation}
\frac{1}{2}\textrm{\large\{}\boldsymbol{b}^{(p')},1\!\!\!\boldsymbol{1}^{\otimes i}\otimes\boldsymbol{b}^{(1)}\textrm{\large\}}=\boldsymbol{b}^{(p)}\otimes1\!\!\!\boldsymbol{1},\label{eq:inter-b}
\end{equation}
\vspace{-0.5cm}
and the identities\negthinspace{} %
\footnote{\ \ where $\boldsymbol{f}^{(1)}\equiv\boldsymbol{b}^{(1)}=\left(\begin{array}{cc}
0 & 0\\
1 & 0
\end{array}\right)$,$\; c_{3}=\left(\begin{array}{cc}
0 & 1\\
-1 & 0
\end{array}\right)$, $c_{2}=\left(\begin{array}{cc}
0 & 1\\
1 & 0
\end{array}\right),1\!\!\!\boldsymbol{1}=\left(\begin{array}{cc}
1 & 0\\
0 & 1
\end{array}\right)$: \\
$\boldsymbol{\qquad f}^{(p)}\equiv1\!\!\!\boldsymbol{1}^{\otimes i-1}\otimes\boldsymbol{b}^{(1)}+(G_{\mu\nu}^{(p)})\otimes c_{3}$,$\quad\boldsymbol{h}^{(p)}\equiv1\mathbf{\!\!\!1}^{\otimes i-1}\otimes\boldsymbol{b}^{(1)}+(J_{\mu\nu}^{(p)})\otimes c_{2}$,
$\quad(i=2,3,\ldots)$%
}

\begin{equation}
(\boldsymbol{f}^{(p')})^{2}=\boldsymbol{f}^{(p)}\otimes1\mathbf{\!\!\!1},\label{eq:inter-f}
\end{equation}
\vspace{-0.9cm}

\begin{equation}
(\boldsymbol{h}^{(p')})^{2}=\boldsymbol{h}^{(p)}\otimes1\mathbf{\!\!\!1}.\label{eq:inter-f-1}
\end{equation}

\noindent Leaving the details to \ref{sec:Crotons-as-boundary},
the way croton base numbers are derived and how they are subdivided
into bases pop out naturally when the matrix elements of $\boldsymbol{f}^{(p)}$
and $\boldsymbol{h}^{(p)}$ are constructed. Crotons, conceived of
as linear combinations, use the following croton base numbers as scalars
(underlining explained later): for $i=2$, $G^{(3)}=1$, $J^{(3)}=1$;
for $i=3$, $G^{(7)}=\underline{1}$, $(J_{\varrho}^{(7)})=(-\underline{1},3)$;
for $i=4$, $(G_{\rho}^{(15)})=(3,\underline{5},11,17,41,113)$, $(J_{\varrho}^{(15)})=(-\underline{5},15,-43,149),$
to name only the first few (singletons and bases). They are instructive
enough to show how label encoding works on the boundary.

\section{\label{sec:Crotons-on-the}Crotons on the boundary}

We first concentrate on order $p=15$, dropping the parenthesized
superscript and just asking the reader to bear in mind that the crotons
examined belong to $i=\log_{2}(p+1)=4$. Our boundary is then defined
by the $3^{T}-1$ outer nodes of a $T$-cube complex, $T$ being the
number of croton base numbers to handle: $T=6$ for $(G_{\rho})=(3,\underline{5},11,17,41,113)$,
and $T=4$ for $(J_{\varrho})=(-\underline{5},15,-43,149)$. Let the
$x$-th node out of the $728=3^{6}-1$ of the first boundary bear
the label $\Gamma_{x}=E_{x}^{\rho}G_{\rho}$, and, correspondingly,
the $y$-th node out of the $80=3^{4}-1$ of the second boundary the
label $\chi_{y}=E_{y}^{\varrho}J_{\varrho}$(summation convention,
and $E$ denoting all non-null $T$-tuples out of $3^{T}$ possible
from $-1,0,1$). It's easy to see that the total of labels form a
croton field in either case: $\Gamma$ and $\chi$. The fact aside
that nodes can be grouped into pairs bearing values of opposite sign,
field values may occur multiply, for instance $6=(0,-1,1,0,0,0)\cdot G^{t}=(0,0,-1,1,0,0)\cdot G^{t}$.
With each field defined on its own boundary, it's far from obvious
they should have anything in common. Yet, as we assume either one
deals with a distinct crotonic aspect $-$ $\Gamma$ with the global
perspective, $\chi$ with  $T$-cube complexes $\Lambda_{\alpha}^{(n)}$
to be introduced in \ref{sub:Crotonic-implementation-locally} $-$
we have to find ways of considering them side by side.

\subsection{\label{sub:Croton-field-duality}Croton field duality and complementarity}

We may, for instance, ask how many distinct labels there can be expressed
\emph{potentially}, neglecting mere sign reversals. Counting from
1 on and taking as the highest conceivable value the sum of croton
base numbers in absolute terms, we arrive at the number $190$ of
potential labels from $G$. Out of these, 170 are realized as node
labels $\Gamma_{x}$. Those not realizable are 20 in number: $7,34,48,\ldots,189$.
The converse holds true for the $J$ case. Of 212 potentially attainable
labels, 40 are realized by $\chi_{y}$ (sign-reversals included, that's
the stock of nodes), leaving 172 labels in potential status. 

\noindent A comparable situation arises when we bunch together croton
base numbers that are rooted in neighboring Mersenne orders, a process
we have previously termed \emph{interordinal} to express this kind
of hybridization. We now have $T=7$ for $(G_{\rho}^{(7,15)})$ $=(\underline{1},3,\underline{5},11,17,41,113)$,
and $T=6$ for $(J_{\rho}^{(7,15)})=(-\underline{1},3,-\underline{5},15,-43,$
$149)$. Neglecting sign reversals and counting again from 1 on, we
get 191 potential labels from the enlarged $G$ and 216 from the enlarged
$J$. All of the 191s' bunch are realized as $\Gamma_{x}$ on the
expanded boundary's nodes; but a singularity also springs up, $0=(1,0,1,1,-1,0,0)\cdot G^{t}$.
By contrast, 202 out of the 216s' bunch are realized as $\chi_{y}$,
on another expanded boundary's nodes and with no singularity popping
up, leaving 14 in potential status: $68,69,\ldots,81$. The conclusion
is that the fields are dual to each other with respect to realizability
of labels on the boundary. The duality is controlled by two quantities,
Catalan number $C_{q\pm1}$ and the number $5\cdot2^{i-r}$$\:(q\in\{1,3\},r\in\{2,3\})$: 

\noindent \emph{Intraordinal case}:\vspace{-0.25cm}
\begin{equation}
\begin{array}{ccc}
 & _{C_{2}}\\
\#\:\Gamma_{x}=170 & \longleftrightarrow & \#\:_{{\textstyle \urcorner}}\chi_{y}=172\\
 & _{5\cdot2^{2}}\\
\#\:\chi_{y}=40 & \longleftrightarrow & \#\:_{{\textstyle \urcorner}}\Gamma_{x}=20.
\end{array}\label{eq:control1}
\end{equation}
\noindent \emph{Interordinal case}:

\begin{equation}
\begin{array}{ccc}
 & _{5\cdot2}\\
\#\:\Gamma_{x}=192^{^{*}} & \longleftrightarrow & \#\:\chi_{y}=202\\
 & _{C_{4}}\\
\#\:_{{\textstyle \urcorner}}\chi_{y}=14 & \longleftrightarrow & \#\:_{{\textstyle \urcorner}}\Gamma_{x}=0.
\end{array}\label{eq:control2}
\end{equation}
($^{{\scriptscriptstyle *}}$The singularity assignment included.)
The key role in that duality is taken by the quantity $C_{q}$ $(q=(p-3)/4)$
around which the croton base numbers for a specific basis of order
$p$ are built (hence the underlining of $C_{q}$ in Sect. \ref{sub:Connecting-boundary-and}
where the bases of order 31 are presented):

\emph{\vspace{1cm}
}\\
\emph{Intraordinal case}:\emph{ }
\begin{equation}
\begin{array}{ccc}
 & _{C_{3}\textrm{\,\ sign\,\ reversal}}\\
G_{\rho}^{(15)} & \longleftrightarrow & J_{\varrho}^{(15)}
\end{array}
\end{equation}
\emph{Interordinal case}:

\begin{equation}
\begin{array}{ccc}
 & {\scriptstyle C_{1},C_{3}\textrm{\,\ sign\,\ reversals}}\\
G_{\delta}^{(7,15)} & \longleftrightarrow & J_{\vartheta}^{(7,15)}.
\end{array}
\end{equation}

\section{Crotons in the volume}

For ``volume'' as the term is used here, a multitude of Mersenne
fluctuations are constituive. They assume a descriptive $\bigwedge$
shape when amplitude is plotted versus ``time''. Nodes on legs of
a `$\bigwedge$' each bear a croton amplitude $\varphi_{\alpha_{n\mp r}}^{(n\mp r)}$
$\in\mathbb{N}$ that emerges with a specific time-like and space-like
refinement $-$ on the left leg $n-r,\alpha_{n-r}$, on the right
$n+r,\alpha_{n+r}$ $-$ and the peak amplitude is reached at $n,\alpha_{n}$.
The left-leg structure is given by
\begin{equation}
\varphi_{\alpha_{n-r+1}}^{(n-r+1)}=2\varphi_{\alpha_{n-r}}^{(n-r)}+1+\epsilon,\qquad(\epsilon\in\{-1,0,1\})\label{eq:left-leg}
\end{equation}
the right-leg structure by 
\begin{equation}
\varphi_{\alpha_{n+r+1}}^{(n+r+1)}=\left\lfloor \varphi_{\alpha_{n+r}}^{(n+r)}/2\right\rfloor -\overline{\delta}\qquad\left(\overline{\delta}=\left\{ \begin{array}{cc}
0\textrm{ or }1 & \varphi_{\alpha_{n+r}}^{(n+r)}\textrm{ even}\\
0 & \textrm{else}
\end{array}\right.\right)\label{eq:right-leg}
\end{equation}
 under the constraint
\[
\left|\varphi_{\alpha_{n-r}}^{(n-r)}-\varphi_{\alpha_{n+r}}^{(n+r)}\right|=\left\{ \begin{array}{cc}
0\textrm{ or }1 &0<r< h-1\\
0 & r=h-1
\end{array}\right.(h\textrm{ the fluctuation's height})
\]
 A typical Mersenne fluctuation is shown in Fig. \ref{fig:A-prototype-geometric}:
\begin{figure}[H]
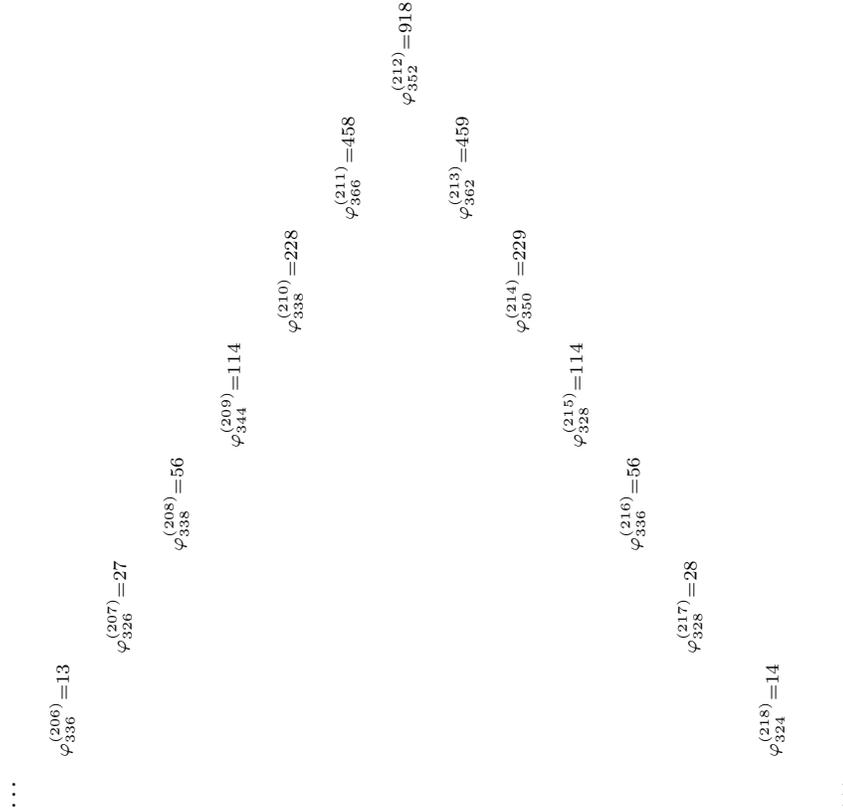

\caption{A geometric fluctuation of Mersenne type\label{fig:A-prototype-geometric}}
\bigskip{}
\begin{tabular*}{12cm}{@{\extracolsep{\fill}}cccccccccccccc>{\centering}p{1cm}>{\centering}p{1cm}}
 &  &  &  &  &  &  &  &  &  &  &  &  &  &  & \tabularnewline
 &  &  &  &  &  &  &  &  &  &  &  &  &  &  & \tabularnewline
 &  &  &  &  &  &  &  & \begin{sideways}
${\scriptstyle \varphi_{352}^{(212)}}{\scriptstyle =918}$
\end{sideways} &  &  &  &  &  &  & \tabularnewline
 &  &  &  &  &  &  & \begin{sideways}
${\scriptstyle \varphi_{366}^{(211)}=458}$
\end{sideways} &  & \begin{sideways}
${\scriptstyle \varphi_{362}^{(213)}=459}$
\end{sideways} &  &  &  &  &  & \tabularnewline
 &  &  &  &  &  & \begin{sideways}
${\scriptstyle \varphi_{338}^{(210)}=228}$
\end{sideways} &  &  &  & \begin{sideways}
${\scriptstyle \varphi_{350}^{(214)}=229}$
\end{sideways} &  &  &  &  & \tabularnewline
 &  &  &  &  & \begin{sideways}
${\scriptstyle \varphi_{344}^{(209)}=114}$
\end{sideways} &  &  &  &  &  & \begin{sideways}
${\scriptstyle \varphi_{328}^{(215)}=114}$
\end{sideways} &  &  &  & \tabularnewline
 &  &  &  & \begin{sideways}
${\scriptstyle \varphi_{338}^{(208)}=56}$
\end{sideways} &  &  &  &  &  &  &  & \begin{sideways}
${\scriptstyle \varphi_{336}^{(216)}=56}$
\end{sideways} &  &  & \tabularnewline
 &  &  & \begin{sideways}
${\scriptstyle \varphi_{326}^{(207)}=27}$
\end{sideways} &  &  &  &  &  &  &  &  &  & \begin{sideways}
${\scriptstyle \varphi_{328}^{(217)}=28}$
\end{sideways} &  & \tabularnewline
\begin{sideways}
\end{sideways} &  & \begin{sideways}
${\scriptstyle \varphi_{336}^{(206)}=13}$
\end{sideways} &  &  &  &  &  &  &  &  &  &  &  & \begin{sideways}
${\scriptstyle \varphi_{324}^{(218)}=14}$
\end{sideways} & \tabularnewline
 & $\vdots$ &  &  &  &  &  &  &  &  &  &  &  &  &  & $\vdots\qquad\qquad$\tabularnewline
 &  &  &  &  &  &  &  &  &  &  &  &  &  &  & \tabularnewline
\end{tabular*}
\end{figure}
\noindent We can stay in the (``time'',amplitude) coordinate system
and observe how fluctuations which share amplitudes that differ maximally
by $\delta$ at each node but peak at different heights, grow to what
we have previously termed \emph{qphylum.}%
\footnote{\emph{\label{fn:One-such-qphylum}}~One such qphylum would for instance
house (peaks in boldface) the Mersenne fluctuations\\
($\ldots,$17,35,72,145,291,584,$\mathbf{1170}$,585,292,145,72,35,17,$\ldots$),\\
($\ldots,$18,36,72,145,291,584,1169,$\mathbf{2340}$,1170,585,292,145,72,36,18,$\ldots$),\\
($\ldots,$17,35,72,146,292,584,1169,2340,$\mathbf{4681}$,2340,1169,584,292,145,72,36,17,$\ldots$)
etc. However, the Mersenne fluctuation ($\ldots,$584,1168,2337,4675,$\mathbf{9350}$,4674,2336,1168,583,$\ldots$)
definitely belongs to a different qphylum.%
}~ 
\begin{figure}[H]
\caption{A prototype qphylum\label{fig:A-prototype-qphylum}}
\vspace{1cm}
\qquad{}\qquad{}\qquad{}\includegraphics[scale=0.5]{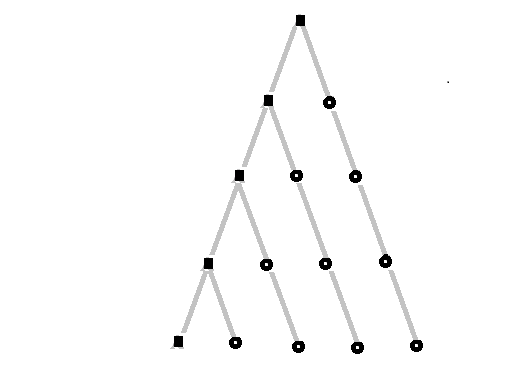}
\end{figure}
\bigskip{}

\begin{figure}[th]
\caption{\label{fig:1-sphere-packing-with}1-sphere packing with or without
centerpiece}

\vspace{2cm}

\qquad{}\qquad{}\qquad{}\qquad{}\qquad{}\qquad{}\includegraphics[scale=0.5]{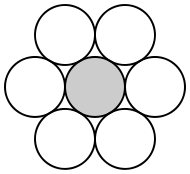}
\end{figure}
\noindent Seen top-down, a qphylum is a left-complete binary tree,
that is: a rooted tree whose root node and left child nodes have left
and right child nodes, while right child nodes have only right child
nodes, as shown in Fig. \ref{fig:A-prototype-qphylum}. Typically,
qphyleticly related amplitudes are rooted in different time-like and
space-like refinements; nodes of a qphylum thus are associated with
a set of frozen-in pregeometric ``time'' and ``space'' signatures.%
\footnote{As for unfreezing, see piece in \ref{sub:Link-between-local} titled
`The role of organizers.'%
} ``Volume'' then becomes the assembly of all distinct qphyla. But
let us go back one step and ask what it means when an amplitude in
a given fluctuation reaches a certain level. If that level coincides
with $L_{m}+1$ or $L_{m}$, where $L_{m}$ denotes the kissing number
of $m$-dimensional Euclidean space, it could mean that a chunk of
space containing an $(m-1)$-sphere packing \emph{with} or \emph{without}
centerpiece was created in that fluctuation $-$ or dissolved if the
amplitude did not peak: crotons which wax and wane. As an example,
Fig. \ref{fig:A-prototype-geometric} shows a fluctuation that peaks
at 918, a quantity considered to be the proper kissing number of 13-dimensional
Euclidean space. A chunk of 13$D$-space containing 918 12-spheres
certainly is hard to visualize, so a 2$D$ version may suffice to
give a first impression (see Fig. \ref{fig:1-sphere-packing-with}).

\subsection{\label{sub:Connecting-boundary-and}Connecting boundary and volume
croton data}

We may ask if and how the peak amplitude 918 is related to a label
$\Gamma_{x}$ on the boundary. Certainly it is a realizable label,
and one realizable \emph{intra}ordinally: All kissing numbers lying
$-$ on the basis of $(G_{\rho}^{(15)})$ $-$ in the range of potentially
attainable labels can be seen to be realizable intraordinally, and
this holds true too for the current basis, the 18-tuple%
\footnote{~in full length, the tuple reads\vspace{-0.2cm}
\[
\begin{array}{ll}
(G_{\rho}^{(31)})= & (19,43,115,155,\underline{429},1275,1595,1633,4819,4905,15067.\\
 & \:15297,18627,58781,189371,227089,737953,2430289);
\end{array}
\]
its origin and the origin of the tuple 
\[
\begin{array}{ll}
(J_{\rho}^{(31)})= & (13,-41,117,143,-\underline{429},1319,1343,1547,-4823,-4903,15547,\\
 & \:17989,18269,-58791,\,194993,\,223573,-747765,2886235);
\end{array}
\]
are elucidated in \ref{sec:Crotons-as-boundary}; in \ref{sec:Crotons-in-the}
(see Table \ref{tab:Specific-fractions-in-1} and \ref{tab:Specific-fractions-in})
various kissing numbers and kissing number-related croton amplitudes
are tabularized, among them also the peak amplitude 918 from Fig.\ref{fig:A-prototype-geometric}.%
} $(G_{\rho}^{(31)})=(19,43,115,155,\underline{429},\ldots,1275,\ldots,$
$4819$, $4905,\ldots)$ where 918 belongs: $918=(0,1,1,0,-1,0,{\scriptstyle \ldots},1,0,{\scriptstyle \ldots},1,-1,0,{\scriptstyle \ldots})\cdot(G^{(31)})^{t}$.
The crucial question is, Do we require \emph{all} croton amplitudes
from a given Mersenne fluctuation with one of them ``geometrizing''
to have counterparts in intraordinally realizable labels on the boundary,
in a narrow interpretation of the holographic principle? Amplitudes
``on the way to/therefrom'' may at least in principle be amenable
to an answer. And, what does this mean for Mersenne fluctuations ``making
detours'' which presumably are by far in the majority? If one of
the croton amplitudes, call it \emph{pivotal}, comes only close and
does not ``geometrize'', it is because some residual Mersenne fluctuations
co-evolve in \emph{different} qphyla. Yet, with sufficiently tight
space-like and time-like refinement constraints, fluctuations that
are inter-qphyleticly linked to the pivotal fluctuation can be identified
and examined. See the example below where one of the residual partial
amplitudes is $102$ as $n-2=1556$, and the pivotal amplitude $5219$,
together with a second residual partial amplitude $24$, is closing
in on $L_{17}(=5346)$ as $n=1558$:
\begin{figure}[H]
\caption{\label{fig:Pivotal-amplitude-closing}Pivotal amplitude closing in
on $L_{17}(=5346)$ plus two residual partial amplitudes}
\vspace{2cm}

\qquad{}\quad{}\includegraphics[scale=0.9]{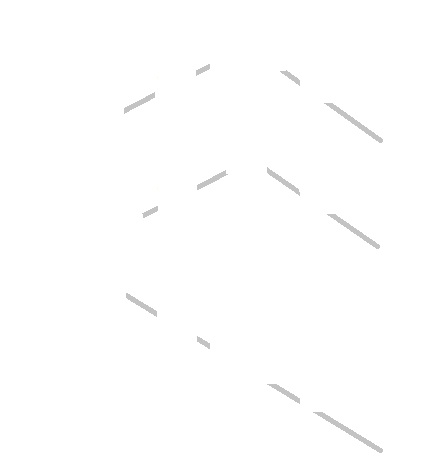}\vspace{-10.8cm}

\qquad{}\qquad{}\qquad{}\qquad{}%
\begin{tabular}{cccccc}
 &  &  &  &  & \tabularnewline
\cline{3-3} 
 & \multicolumn{1}{c|}{} & \multicolumn{1}{c|}{$\begin{array}{c}
{\scriptstyle \varphi_{401}^{(1558)}}\\
{\scriptstyle =5219}
\end{array}$} &  &  & \tabularnewline
\cline{3-3} 
 & $\begin{array}{c}
{\scriptstyle \varphi_{442}^{(1557)}}\\
{\scriptstyle =2609}
\end{array}$ &  & $\begin{array}{c}
{\scriptstyle \varphi_{407}^{(1559)}}\\
{\scriptstyle =2609}
\end{array}$ &  & \tabularnewline
$\begin{array}{c}
{\scriptstyle \varphi_{448}^{(1556)}}\\
{\scriptstyle =1304}
\end{array}$ &  &  &  & $\ddots$ & \tabularnewline
 &  & $\begin{array}{c}
{\scriptstyle \varphi_{390}^{(1558)}}\\
{\scriptstyle =411}
\end{array}$ &  &  & \tabularnewline
 & $\begin{array}{c}
{\scriptstyle \varphi_{380}^{(1557)}}\\
{\scriptstyle =205}
\end{array}$ &  & $\begin{array}{c}
{\scriptstyle \varphi_{398}^{(1559)}}\\
{\scriptstyle =205}
\end{array}$ &  & \tabularnewline
\cline{1-1} 
\multicolumn{1}{|c|}{$\begin{array}{c}
{\scriptstyle \varphi_{404}^{(1556)}}\\
{\scriptstyle =102}
\end{array}$} &  &  &  & $\ddots$ & \tabularnewline
\cline{1-1} 
 &  &  &  &  & \tabularnewline
$\begin{array}{c}
{\scriptstyle \varphi_{406}^{(1556)}}\\
{\scriptstyle =100}
\end{array}$ &  &  &  &  & \tabularnewline
 & $\begin{array}{c}
{\scriptstyle \varphi_{382}^{(1557)}}\\
{\scriptstyle =49}
\end{array}$ &  &  &  & \tabularnewline
\cline{3-3} 
 & \multicolumn{1}{c|}{} & \multicolumn{1}{c|}{$\begin{array}{c}
{\scriptstyle \varphi_{392}^{(1558)}}\\
{\scriptstyle =24}
\end{array}$} &  &  & \tabularnewline
\cline{3-3} 
 &  &  & $\begin{array}{c}
{\scriptstyle \varphi_{400}^{(1559)}}\\
{\scriptstyle =12}
\end{array}$ &  & \tabularnewline
 &  &  &  & $\ddots$ & \tabularnewline
\end{tabular}
\end{figure}
\noindent Another example is shown in Fig. \ref{fig:Pivotal-amplitude-closing-1}
where two residual partial amplitudes attain the levels $58$ and
$193$ respectively as $n-2=1003$, allowing a pivotal amplitude $207679$
to close in on $L_{29}(=207930)$ as $n=1005$: 
\begin{figure}[H]
\caption{\label{fig:Pivotal-amplitude-closing-1}Pivotal amplitude closing
in on $L_{29}(=207930)$ plus two residual partial amplitudes}
\vspace{2cm}

\qquad{}\qquad{}\includegraphics[scale=0.9]{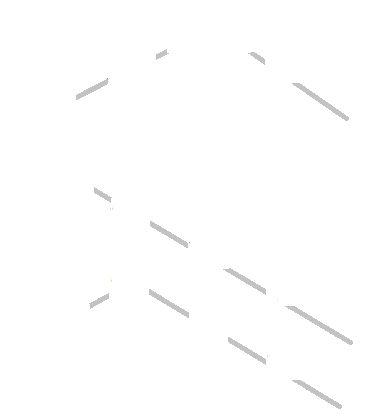}\vspace{-10.1cm}

\qquad{}\qquad{}\quad{}%
\begin{tabular}{cccccc}
 &  &  &  &  & \tabularnewline
\cline{3-3} 
 & \multicolumn{1}{c|}{} & \multicolumn{1}{c|}{$\begin{array}{c}
{\scriptstyle \varphi_{438}^{(1005)}}\\
{\scriptstyle =207679}
\end{array}$} &  &  & \tabularnewline
\cline{3-3} 
 & $\begin{array}{c}
{\scriptstyle \varphi_{442}^{(1004)}}\\
{\scriptstyle =103839}
\end{array}$ &  & $\begin{array}{c}
{\scriptstyle \varphi_{442}^{(1006)}}\\
{\scriptstyle =103839}
\end{array}$ &  & \tabularnewline
$\begin{array}{c}
{\scriptstyle \varphi_{448}^{(1003)}}\\
{\scriptstyle =51919}
\end{array}$ &  &  &  & $\ddots$ & \tabularnewline
 &  &  &  &  & \tabularnewline
 &  &  &  &  & \tabularnewline
\cline{1-1} 
\multicolumn{1}{|c|}{$\begin{array}{c}
{\scriptstyle \varphi_{449}^{(1003)}}\\
{\scriptstyle =193}
\end{array}$} &  &  &  &  & \tabularnewline
\cline{1-1} 
 & $\begin{array}{c}
{\scriptstyle \varphi_{443}^{(1004)}}\\
{\scriptstyle =96}
\end{array}$ &  &  &  & \tabularnewline
 &  & $\begin{array}{c}
{\scriptstyle \varphi_{439}^{(1005)}}\\
{\scriptstyle =48}
\end{array}$ &  &  & \tabularnewline
 & $\begin{array}{c}
{\scriptstyle \varphi_{477}^{(1004)}}\\
{\scriptstyle =117}
\end{array}$ &  & $\begin{array}{c}
{\scriptstyle \varphi_{427}^{(1006)}}\\
{\scriptstyle =23}
\end{array}$ &  & \tabularnewline
\cline{1-1} 
\multicolumn{1}{|c|}{$\begin{array}{c}
{\scriptstyle \varphi_{479}^{(1003)}}\\
{\scriptstyle =58}
\end{array}$} &  & $\begin{array}{c}
{\scriptstyle \varphi_{465}^{(1005)}}\\
{\scriptstyle =58}
\end{array}$ &  & $\ddots$ & \tabularnewline
\cline{1-1} 
 &  &  & $\begin{array}{c}
{\scriptstyle \varphi_{449}^{(1006)}}\\
{\scriptstyle =29}
\end{array}$ &  & \tabularnewline
 &  &  &  & $\ddots$ & \tabularnewline
 &  &  &  &  & \tabularnewline
\end{tabular}
\end{figure}
\noindent The pivotal amplitudes in Figs. \ref{fig:Pivotal-amplitude-closing}
and \ref{fig:Pivotal-amplitude-closing-1} each coincide with the
peak of their parental fluctuation, but peak amplitude is not a necessary
condition. Fig. \ref{fig:Pivotal-edge-amplitudes} describes a situation
where an $(n\mp r)$-pair of pivotal amplitudes on a fluctuation's\emph{
}legs are about to close in on $L_{29}$; here, since only one time-like
refinement lies between each candidate and the peak, one further time-like
refinement also suffices to determine the residual partial amplitudes,
where they originate and which of the two `leggy' candidates $207646$
and $207647$ would have succeeded in filling the bill had it peaked:
\begin{figure}[H]
\caption{\label{fig:Pivotal-edge-amplitudes}Leggy pivot closing in on $L_{29}(=207930)$
plus two residual partial amplitudes}

\includegraphics[scale=0.9]{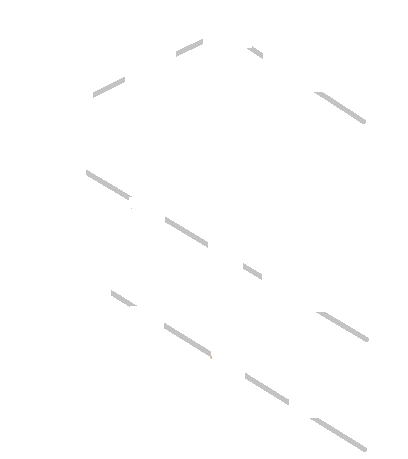}\vspace{-10.9cm}

\hspace{0.8cm}%
\begin{tabular}{cccccc}
 &  &  &  &  & \tabularnewline
 &  & $\begin{array}{c}
{\scriptstyle \varphi_{421}^{(992)}}\\
{\scriptstyle =415294}
\end{array}$ &  &  & \tabularnewline
\cline{4-4} 
 & $\begin{array}{c}
{\scriptstyle \varphi_{417}^{(991)}}\\
{\scriptstyle =207646}
\end{array}$ & \multicolumn{1}{c|}{} & \multicolumn{1}{c|}{$\begin{array}{c}
{\scriptstyle \varphi_{423}^{(993)}}\\
{\scriptstyle =207647}
\end{array}$} &  & \tabularnewline
\cline{4-4} 
$\begin{array}{c}
{\scriptstyle \varphi_{433}^{(990)}}\\
{\scriptstyle =103823}
\end{array}$ &  &  &  & $\ddots$ & \tabularnewline
 &  &  &  &  & \tabularnewline
 &  &  &  &  & \tabularnewline
$\begin{array}{c}
{\scriptstyle \varphi_{421}^{(990)}}\\
{\scriptstyle =1178}
\end{array}$ &  &  &  &  & \tabularnewline
 & $\begin{array}{c}
{\scriptstyle \varphi_{409}^{(991)}}\\
{\scriptstyle =589}
\end{array}$ &  &  &  & \tabularnewline
 &  & $\begin{array}{c}
{\scriptstyle \varphi_{413}^{(992)}}\\
{\scriptstyle =294}
\end{array}$ &  &  & \tabularnewline
\cline{1-1} \cline{4-4} 
\multicolumn{1}{|c|}{$\begin{array}{c}
{\scriptstyle \varphi_{485}^{(990)}}\\
{\scriptstyle =136}
\end{array}$} &  & \multicolumn{1}{c|}{} & \multicolumn{1}{c|}{$\begin{array}{c}
{\scriptstyle \varphi_{415}^{(993)}}\\
{\scriptstyle =147}
\end{array}$} &  & \tabularnewline
\cline{1-1} \cline{4-4} 
 & $\begin{array}{c}
{\scriptstyle \varphi_{465}^{(991)}}\\
{\scriptstyle =68}
\end{array}$ &  &  & $\ddots$ & \tabularnewline
 &  & $\begin{array}{c}
{\scriptstyle \varphi_{473}^{(992)}}\\
{\scriptstyle =33}
\end{array}$ &  &  & \tabularnewline
 &  &  & $\begin{array}{c}
{\scriptstyle \varphi_{471}^{(993)}}\\
{\scriptstyle =16}
\end{array}$ &  & \tabularnewline
 &  &  &  & $\ddots$ & \tabularnewline
\end{tabular}
\end{figure}

\subsection{\label{sec:croton-phase}Croton phase and its inter-qphyletic role}

``Phase'' in the pre-geometric setting assumed here just means `having
an ordinate value fluctuate between positions above and below an imaginary
baseline,' with consecutive marks on that line corresponding to stepwise
increases of space-like refinement. That ordinate value, let us call
it $\psi_{\alpha_{\psi}}^{(n)}$ for a given ``time'' level $n$,
is linked to the croton amplitude $\varphi_{\alpha_{\varphi}}^{(n)}$
by the condition:
\begin{equation}
\textrm{if}\;\left|\psi_{\alpha_{\psi}}^{(n)}\right|=\varphi_{\alpha_{\varphi}}^{(n)}+\delta\;\textrm{then}\enskip\psi_{\alpha_{\psi}}^{(n)}=\left\{ \begin{array}{cl}
\!\!\!-\varphi_{\alpha_{\varphi}}^{(n)}-1 & \!\!\!(\alpha_{\psi}\:\textrm{even})\\
\!\!\!\varphi_{\alpha_{\varphi}}^{(n)}+\delta & \!\!\!(\alpha_{\psi}\:\textrm{odd})
\end{array}\right.(\delta\in\{0,1\}).
\end{equation}
Whilst introducing croton phase in the volume now, the discussion
will be limited to the fluctuations already considered in order to
keep things as coherent as possible. Let us first follow two Mersenne
fluctuations' amplitudes (Fig. \ref{fig:Pivotal-edge-amplitudes})
and their associated $\psi^{(n)}$, one steering a pivotal, one a
selected residual's course:\vspace{-1cm}

\begin{table}[H]
\caption{\label{tab:Pivotal-and-residual-fig5}Pivotal and inter-qphyleticly
accompanying residual fluctuation for $987\leq n\leq997$}
\medskip{}

\begin{tabular}{|c|cc|cc|c|cc|cc|c|}
\hline 
 &  & \multicolumn{1}{c}{} &  & \multicolumn{1}{c}{} &  &  & \multicolumn{1}{c}{} &  & \multicolumn{1}{c}{} & \tabularnewline
 & \multicolumn{5}{c|}{Pivot} & \multicolumn{5}{c|}{Residue}\tabularnewline
$n$ & $\varphi_{P}$ & $\alpha_{\varphi}$ & $\psi_{P}$ & $\alpha_{\psi}$ & $\triangle\alpha$ & $\varphi_{R}$ & $\alpha_{\varphi}$ & $\psi_{R}$ & $\alpha_{\psi}$ & $\triangle\alpha$\tabularnewline
 &  &  &  &  &  &  &  &  &  & \tabularnewline
\hline 
\hline 
987 & 12977 & \negthinspace{}407 & $1$2977 & \negthinspace{}437 & 30 & 9434 & \negthinspace{}397 & -9435 & \negthinspace{}426 & 29\tabularnewline
\hline 
988 & 25955 & \negthinspace{}411 & -25956 & \negthinspace{}418 & 7 & 4716 & \negthinspace{}399 & -4717 & \negthinspace{}410 & 11\tabularnewline
\hline 
989 & 51911 & \negthinspace{}397 & 51911 & \negthinspace{}449 & 52 & 2357 & \negthinspace{}385 & -2358 & \negthinspace{}438 & 53\tabularnewline
\hline 
990 & 103823 & \negthinspace{}433 & 103824 & \negthinspace{}441 & 8 & 1178 & \negthinspace{}421 & -1179 & \negthinspace{}430 & 9\tabularnewline
\hline 
\noalign{\vskip0.5mm}
991 & 207646 & \negthinspace{}417 & 207647 & \negthinspace{}457 & %
\framebox{\begin{minipage}[t][0.6\totalheight]{0.03\columnwidth}%
40%
\end{minipage}} & 589 & \negthinspace{}409 & 589 & \negthinspace{}449 & %
\framebox{\begin{minipage}[t][0.6\totalheight]{0.03\columnwidth}%
40%
\end{minipage}}\tabularnewline[0.5mm]
\hline 
\noalign{\vskip0.5mm}
992 & 415294 & \negthinspace{}421 & 415294 & \negthinspace{}469 & %
\framebox{\begin{minipage}[t][0.6\totalheight]{0.03\columnwidth}%
48%
\end{minipage}} & 294 & \negthinspace{}413 & 294 & \negthinspace{}461 & %
\framebox{\begin{minipage}[t][0.6\totalheight]{0.03\columnwidth}%
48%
\end{minipage}}\tabularnewline[0.5mm]
\hline 
993 & 207647 & \negthinspace{}423 & 207647 & \negthinspace{}473 & 50 & 147 & \negthinspace{}405 & 147 & \negthinspace{}461 & 56\tabularnewline
\hline 
994 & 103823 & \negthinspace{}443 & 103823 & \negthinspace{}457 & 14 & 73 & \negthinspace{}435 & -74 & \negthinspace{}446 & 11\tabularnewline
\hline 
995 & 51911 & \negthinspace{}425 & 51912 & \negthinspace{}453 & 28 & 36 & \negthinspace{}413 & 36 & \negthinspace{}445 & 32\tabularnewline
\hline 
996 & 25955 & \negthinspace{}431 & 25955 & \negthinspace{}477 & 46 & 17 & \negthinspace{}417 & 17 & \negthinspace{}461 & 44\tabularnewline
\hline 
997 & 12977 & \negthinspace{}439 & 12978 & \negthinspace{}469 & 30 & 8 & \negthinspace{}427 & 8 & \negthinspace{}457 & 30\tabularnewline
\hline 
\end{tabular}
\end{table}
\noindent From the table one can glean that, as the $991${\small -th}
time-like refinement level is reached, the offsets $\triangle\alpha(\equiv\alpha_{\psi}-\alpha_{\varphi})$
under consideration get correlated for pivot and residual $-$ first
in $\triangle\alpha=40$, then in $\triangle\alpha=48$ $-$ signalling
the residual amplitude's share 147 to the target bill and concomitant
decorrelation of $\triangle\alpha$ as $n=993$. The same holds true
for the pivot's and the largest residual's amplitudes from Fig. \ref{fig:Pivotal-amplitude-closing}
and their associated $\psi^{(n)}$: (1) correlation in $\triangle\alpha=20$
as $n=1003$; (2) correlation in $\triangle\alpha=22$ as $n=1004$;
(3) residual amplitude's (belated) contribution 193 and decorrelation
of $\triangle\alpha$ as $n=1005$ (see Table \ref{tab:Pivotal-and-residual-fig4}).
So a first conclusion is that, from a volume point of view, target-seeking
implies phase correlation irrespective of a pivotal amplitude's coincidence
with a peak or not. \vspace{-0.5cm}

\begin{table}[H]
\caption{\label{tab:Pivotal-and-residual-fig4}Pivotal and inter-qphyleticly
accompanying residual fluctuation for $1003\leq n\leq1005$}

\medskip{}
\ %
\begin{tabular}{|c|cc|cc|c|cc|cc|c|}
\hline 
 &  & \multicolumn{1}{c}{} &  & \multicolumn{1}{c}{} &  & \multicolumn{5}{c|}{}\tabularnewline
 & \multicolumn{5}{c|}{Pivot} & \multicolumn{5}{c|}{Residue}\tabularnewline
$n$ & $\varphi_{P}$ & $\alpha_{\varphi}$ & $\psi_{P}$ & $\alpha_{\psi}$ & $\triangle\alpha$ & $\varphi_{R}$ & $\alpha_{\varphi}$ & $\psi_{R}$ & $\alpha_{\psi}$ & $\triangle\alpha$\tabularnewline
 &  &  &  &  &  &  &  &  &  & \tabularnewline
\hline 
\hline 
\noalign{\vskip0.5mm}
1003 & 51919 & 448 & -51920 & 468 & %
\framebox{\begin{minipage}[t][0.6\totalheight]{0.03\columnwidth}%
20%
\end{minipage}} & 193 & 449 & 194 & 469 & %
\framebox{\begin{minipage}[t][0.6\totalheight]{0.03\columnwidth}%
20%
\end{minipage}}\tabularnewline[0.5mm]
\hline 
\noalign{\vskip0.5mm}
1004 & 103839 & 442 & 103839 & 464 & %
\framebox{\begin{minipage}[t][0.6\totalheight]{0.03\columnwidth}%
22%
\end{minipage}} & 96 & 443 & 96 & 465 & %
\framebox{\begin{minipage}[t][0.6\totalheight]{0.03\columnwidth}%
22%
\end{minipage}}\tabularnewline[0.5mm]
\hline 
\noalign{\vskip0.5mm}
1005 & 207679 & 438 & 207679 & 443 & 5 & 48 & 439 & 48 & 445 & 6\tabularnewline[0.5mm]
\hline 
\end{tabular}
\end{table}

\noindent  There is more to croton phase than just that. Let us once
more go back one step and consider the target-matching case first.
If a croton amplitude reaches a level $L_{m}+1$ or $L_{m}$, it was
assumed a chunk of $m$-dimensional Euclidean space containing an
$(m-1)$-sphere packing \emph{with} or \emph{without} centerpiece
was created in that fluctuation (or dissolved if the amplitude did
not peak). Whether that creation succeeded depends on the quantity
$\delta\in\{0,1\}$: only if the amplitude $\varphi$ peaks on $L_{m}$
and the phase $\psi$, in absolute terms, on $L_{m}+\delta\:(\delta\in\{0,1\})$
can we be sure of successful creation; if $\varphi<L_{m}$, or if
$\left|\psi\right|>L_{m}+1$, we'd be uncertain whether to settle
on success or state failure. There are situations where that criterion
applies to more than one Mersenne fluctuation. See Table \ref{tab:Pivotal-and-residual-fig4-1}
which illuminates the stance of three detouring fluctuations:

\vspace{-0.5cm}
\begin{table}[H]
\caption{\label{tab:Pivotal-and-residual-fig4-1}Co-occurrent fluctuations
targeted at $L_{29}(=207930)$, $L_{10}(=336)$ and $L_{8}(=240)$}

\medskip{}
\quad{}%
\begin{tabular}{|c|ccccc|cc|cc|c|}
\hline 
 & \multicolumn{4}{c}{} &  & \multicolumn{5}{c|}{}\tabularnewline
 & \multicolumn{5}{c|}{Pivot} & \multicolumn{5}{c|}{Residue 1}\tabularnewline
$n$ & $\varphi_{P}$ & \multicolumn{1}{c|}{$\alpha_{\varphi}$} & $\psi_{P}$ & \multicolumn{1}{c|}{$\alpha_{\psi}$} & $\triangle\alpha$ & $\varphi_{R_{1}}$ & $\alpha_{\varphi}$ & $\psi_{R_{1}}$ & $\alpha_{\psi}$ & $\triangle\alpha$\tabularnewline
 &  & \multicolumn{1}{c|}{} &  & \multicolumn{1}{c|}{} &  &  &  &  &  & \tabularnewline
\hline 
\hline 
\noalign{\vskip0.5mm}
609 & 208430 & \multicolumn{1}{c|}{72} & -208431 & \multicolumn{1}{c|}{72} & 0 & 10 & 58 & -11 & 52 & %
\framebox{\begin{minipage}[t][0.6\totalheight]{0.03\columnwidth}%
6%
\end{minipage}}\tabularnewline[0.5mm]
\hline 
\hline 
 &  &  &  &  &  & \multicolumn{5}{c|}{}\tabularnewline
 &  &  &  &  &  & \multicolumn{5}{c|}{Residue 2}\tabularnewline
 &  &  &  &  &  & $\varphi_{R_{2}}$ & $\alpha_{\varphi}$ & $\psi_{R_{2}}$ & $\alpha_{\psi}$ & $\triangle\alpha$\tabularnewline
 &  &  &  &  &  &  &  &  &  & \tabularnewline
\hline 
\hline 
\noalign{\vskip0.5mm}
609 & \multicolumn{5}{c|}{} & 66 & 78 & 66 & 83 & %
\framebox{\begin{minipage}[t][0.6\totalheight]{0.03\columnwidth}%
5%
\end{minipage}}\tabularnewline[0.5mm]
\hline 
\end{tabular}
\end{table}
\noindent Contrary to the former examples of detouring, in the above
there is only one time-like refinement that counts because the largest
amplitude (still called pivot) overshoots as $n=609$. The three fluctuations
could, in a ``covert conspiracy'', strive concurrently after three
targets, $\varphi_{P}+\varphi_{R_{1}}+\varphi_{R_{2}}=L_{29}+L_{10}+L_{8}$,
where $L_{29}=207930,$ $L_{10}=336$ and $L_{8}=240$. Residuality
not only assumes a different meaning here, the space-like refinements
get also symmetrized, one residual's being lower than the pivot's,
the other one's higher, and the offsets in question get correlated
as $n=609$. Offset equality obviously is uncertain by a factor $\delta=\left|\triangle\alpha_{R_{2}}-\triangle\alpha_{R_{1}}\right|\;(\delta\in\{0,1\})$,
and $\delta=1$ above since the phase inversions that enter at $n=609$
are followed only by two of the three contributors. A very similar
example is shown in Table \ref{tab:Pivotal-and-residual-fig4-1-1}:\vspace{-0.5cm}
 
\begin{table}[H]
\caption{\label{tab:Pivotal-and-residual-fig4-1-1}Co-occurrent fluctuations
targeted at $L_{12}(=756)$, $L_{10}(=336)$ and $L_{2}(=6)$}

\medskip{}
\quad{}\quad{}%
\begin{tabular}{|c|cc|ccc|cc|cc|c|}
\hline 
 & \multicolumn{4}{c}{} &  & \multicolumn{5}{c|}{}\tabularnewline
 & \multicolumn{5}{c|}{Pivot} & \multicolumn{5}{c|}{Residue 1}\tabularnewline
$n$ & $\varphi_{P}$ & $\alpha_{\varphi}$ & $\psi_{P}$ & \multicolumn{1}{c|}{$\alpha_{\psi}$} & $\triangle\alpha$ & $\varphi_{R_{1}}$ & $\alpha_{\varphi}$ & $\psi_{R_{1}}$ & $\alpha_{\psi}$ & $\triangle\alpha$\tabularnewline
 &  &  &  & \multicolumn{1}{c|}{} &  &  &  &  &  & \tabularnewline
\hline 
\hline 
\noalign{\vskip0.5mm}
1000 & 758 & 239 & 758 & \multicolumn{1}{c|}{269} & 30 & 335 & 356 & -336 & 388 & %
\framebox{\begin{minipage}[t][0.6\totalheight]{0.03\columnwidth}%
32%
\end{minipage}}\tabularnewline[0.5mm]
\hline 
\hline 
 &  & \multicolumn{1}{c}{} &  &  &  & \multicolumn{5}{c|}{}\tabularnewline
 & \multicolumn{5}{c|}{} & \multicolumn{5}{c|}{Residue 2}\tabularnewline
 &  & \multicolumn{1}{c}{} &  &  &  & $\varphi_{R_{2}}$ & $\alpha_{\varphi}$ & $\psi_{R_{2}}$ & $\alpha_{\psi}$ & $\triangle\alpha$\tabularnewline
 &  & \multicolumn{1}{c}{} &  &  &  &  &  &  &  & \tabularnewline
\hline 
\noalign{\vskip0.5mm}
1000 & \multicolumn{5}{c|}{} & 5 & 135 & -6 & 168 & %
\framebox{\begin{minipage}[t][0.6\totalheight]{0.03\columnwidth}%
33%
\end{minipage}}\tabularnewline[0.5mm]
\hline 
\end{tabular}
\end{table}

\noindent A natural question to ask is if the aforementioned amplitudes
and phases and the conclusions drawn from them stand a boundary check.

\subsection{Boundary check}

Clearly, the boundary must be checked at this stage because it has
yet to be decided if crotons from Mersenne fluctuations making a detour
around a kissing number get encoded \emph{intra}ordinally or \emph{inter}ordinally.
We may put together the relevant facts here by starting with a recollection
and extrapolating therefrom: 

(1) Out of all $\Gamma^{(15)}$ labels realizable on the basis of
$(G_{\rho}^{(15)})$ (the intraordinal case), one subset of labels
can be extracted that encode croton amplitudes coincident with $\pm L_{m},\pm(L_{m}+1)$
$(m=1,2,\ldots7)$. Not $\pm(L_{2}+1)=\pm7$, however. The complete
realization (interordinal case) demands an enlarged basis $(G_{\varrho}^{(7,15)})$
that brings singular labels in its wake. On the basis of $(J_{\rho}^{(15)})$,
no subset of $\chi^{(15)}$ labels is able to encode $L_{m},L_{m}+1$
$(m=1,2,\ldots7)$ or sign reversed versions thereof; that encoding
only catches up when the basis is enlarged to $(J_{\varrho}^{(7,15)})$,
facing us with two adamant cases yet: $\pm(L_{7}+1)=\pm127$ (realizable
before) and $\pm L_{6}$, $\pm(L_{6}+1)$ (unrealizable after); no
singular labels are entailed. Altogether a complex picture. 

(2) In contrast, singular labels spring up directly on the basis of
either $(G_{\rho}^{(31)})$ or $(J_{\rho}^{(31)})$ (see \ref{sec:Crotons-as-boundary});
a simplification that, in turn, pays off with \emph{twofold}-realizable
labels all the way up for croton amplitudes coincident with $\pm L_{m}$,
$\pm(L_{m}+1)$ $(m=8,9,\ldots,31)$, also making up leeway to the
former special cases $m=2,6,7$.

\noindent So what can, on top of that, be checked is the realizability
of our example \emph{pivots} from Figs.\ref{fig:Pivotal-amplitude-closing}
to \ref{fig:Pivotal-edge-amplitudes}, including their associated
$\varphi_{P}$. Because of what they are targeted at, it is decidable
unequivocally where their images are to be sought: $\Gamma^{(31)}$
and $\chi^{(31)}$, based on the 18-tuples $(G_{\rho}^{(31)})$ and
$(J_{\rho}^{(31)})$. The targets under discussion behave as one would
expect: $\pm L_{17}(=\pm5346)$, $\pm(L_{17}+1)$ as well as $\pm L_{29}(=\pm207930)$,
$\pm(L_{29}+1)$ are realizable altogether by $\Gamma^{(31)}$ and
$\chi^{(31)}$. The same statement holds true for our example pivots
$-$ 5219, 5220. And, peaking or not, amplitudes 207646, 207647 and
207679 are perfectly twofold-realizable either. Surprisingly, twofold-realizability
holds out for the whole pivotal and residual stopovers and co-occurrent
targets mentioned in the discussion of detouring fluctuations. The
holographic principle, according to which all volume quantities $\varphi_{P}$,
$\psi_{P}$ and $\varphi_{R}$, $\psi_{R}$ from Tables \ref{tab:Pivotal-and-residual-fig5}
to \ref{tab:Pivotal-and-residual-fig4-1-1} should have boundary counterparts,
$\varphi$ in $\Gamma^{(31)}$, $\psi$ in $\chi^{(31)}$, is exceedingly
satisfied $-$ the aforesaid quantities are invariably twofold-realizable
(see \ref{sec:Crotons-as-boundary}). Since the same encompassing
holographicity also obtains for the stopovers in the target matching
Mersenne fluctuation of Fig. \ref{fig:A-prototype-geometric}, one
can in summary say that amplitude and phase data from target-seeking
Mersenne fluctuations in the volume have a perfect image on the boundary.\newpage{}

\section{Sources of Mersenne fluctuations}

Thus far, examples of Mersenne fluctuations have been alleged without
specifying their sources. What we expect from actual sources is that
they reveal the conditions under which Mersenne fluctuations (1) develop
and (2) grow into qphyla that in turn define ``volume'' in a pregeometric
context. The apparatus employed here is continued fractions
\begin{equation}
b_{0}+\frac{\left.a_{1}\right|}{\left|b_{1}\right.}+\frac{\left.a_{2}\right|}{\left|b_{2}\right.}+\frac{\left.a_{3}\right|}{\left|b_{3}\right.}+\cdots
\end{equation}
where the shorthand $\left[b_{0};b_{1},b_{2},\ldots\right]$ is used
for the regular case $(a_{\alpha}=1)$; a shorthand for the case $a_{0}=a_{2\mu-1}=1,$
$a_{2\mu}=-1$ will be given soon.

\subsection{The role of continued fractions in refinement}

\begin{figure}[H]
\caption{\label{fig:Sqare-spiral-of}Square spiral representation of $\mathbb{N}_{0}$}

\hspace{-1.6cm}\includegraphics[scale=0.33]{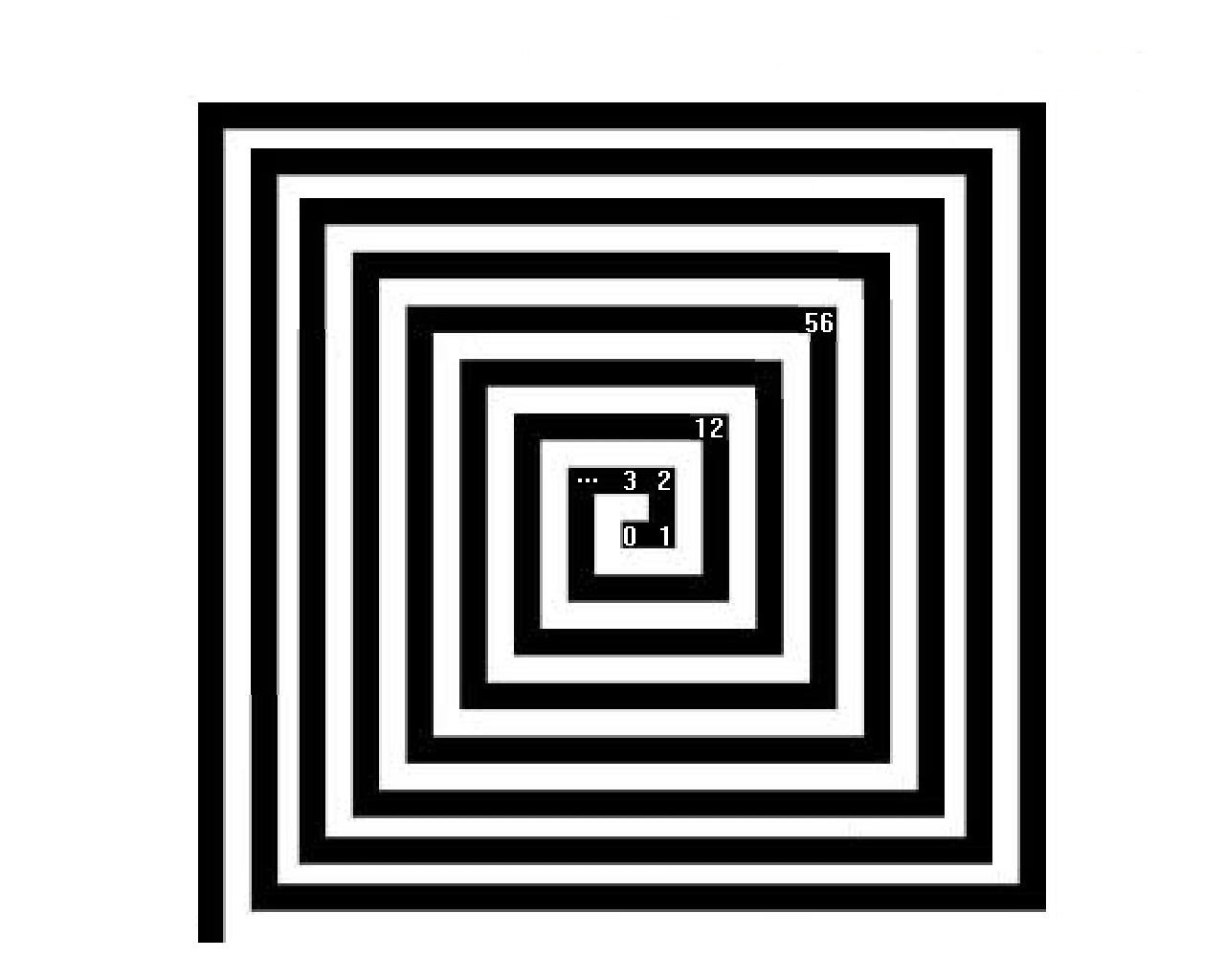}\vspace{-11.7cm}

\begin{tabular}{ccccccccccc}
 &  &  &  &  &  &  &  &  &  & \tabularnewline
 &  &  &  &  &  &  &  & $\longleftarrow\qquad(p_{n}+1)/\sqrt{2}\qquad\longrightarrow$ &  & $\qquad\qquad p_{n}(p_{n}+1)$\tabularnewline
\end{tabular}\vspace{11cm}
\end{figure}

To illustrate the role of contined fractions in refinement, let us
start with a time-honored example, the square spiral formed by the
numbers $\mathbb{N}_{0}$ accompanied by a generalization of $p=2^{i}-1$
$(i=1,2,\ldots5)$ to Mersenne numbers $p_{n}\equiv2^{n}-1$ $(n\in\mathbb{N})$.
As indicated in Fig. \ref{fig:Sqare-spiral-of}, with $\mathrm{B(,)}$
the Beta function, the terms $(C_{p_{n}}\mathrm{B}(p_{n},p_{n}+1))^{-1}=p_{n}(p_{n}+1)$
figure as marks on a subset of corners along the number pattern's
diagonal: For $p_{1}$, this is one corner away from the origin, for
$p_{2}$ two corners, and for $p_{n}$, $p_{n-1}+1$ corners generally.
Taking the number of corners as a measure, we can say the square spiral
is endowed with an \emph{expansion} parameter: $(p_{n-1}+1)/\sqrt{2}$,
the radius of an inscribed circle of a square with side length $(p_{n}+1)/\sqrt{2}$.
That in turn is equivalent to saying a fixed irrational quantity $\sqrt{2}$
gets refined in steps of powers of two, $\left(\frac{2^{n}}{\sqrt{2}}\right)^{-1}$.
The denominators from a convergent's regular continued fraction representation
$\left(\frac{2^{n}}{\sqrt{2}}\right)^{-1}\rightarrow\left[b_{0}^{(n)};b_{\alpha}^{(n)}\right]$
then unveil the time-like and space-like aspects of refinement: One
just proceeds from $n$ to $n+1$ in the superscript of the denominators
to follow the convergent's \emph{time}-like refinement, and follows
its respective \emph{space}-like refinements by proceeding from $\alpha=1$
to $\alpha=2$ to further increments of $\alpha$ in the cf term subscripts
ad infinitum.

\subsection{Mersenne fluctuations and randomness}

The continued fraction representations $\left(\frac{2^{n}}{\sqrt{2}}\right)^{-1}\rightarrow\left[b_{0}^{(n)};b_{\alpha}^{(n)}\right]$
yield orderly Mersenne fluctuations. Moreover, allowing for (CFR)
$\left(\frac{2^{n}}{\sqrt{2}}\right)^{-1}\pm\mathbb{N}\rightarrow\left[\tilde{b}_{0}^{(n)};\tilde{b}_{\alpha}^{(n)}\right]$,
amplitudes can conveniently be generalized to trajectories across
lattice points via the identities $b_{\alpha}^{(n)}=\tilde{b}{}_{\alpha}^{(n)}\;(\alpha>2)$.
In both representations, however, terms are confined to a period after
which they repeat (like in a Wilczek time crystal) and lead to relatively
modest target matching: Results are, at least for online CFR calculators
with the typical limits $n\leq3324,\:\alpha\leq499$ restricted to
the kissing numbers $L_{m}$, $m=1,\ldots,13$ (see Table \ref{tab:Specific-fractions-in-1}).

\noindent It was mentioned in the introduction that the bases $G^{(p_{n})}$,
$J^{(p_{n})}$ are rooted in the matrix representations of the operators
$\boldsymbol{f}^{(p_{n})}$, $\boldsymbol{h}^{(p_{n})}$. In the same
previous work that introduced them it was further noted that $\left(\boldsymbol{f}^{(p_{n})}\right)^{p_{n}+1}\!\!=0$,$\left(\boldsymbol{h}^{(p_{n})}\right)^{p_{n}+1}\!\!=0$,
and that the length of an arc on a cardioid parametrized by $p_{n}$
shows similar behavior: After $n$ steps taken in reverse, the cardiod's
arclength too becomes zero: 
\begin{equation}
(A_{n},A_{n})\equiv2\,\frac{A_{n}\cdot\textrm{co-}A_{n}}{A_{n}+\bar{A}_{n}}=A_{n-1}\qquad\textrm{etc.}\label{eq:f-squared-1}
\end{equation}
where

\begin{equation}
\begin{array}{c}
\begin{array}{cr}
A_{n}=4c\sin\frac{\pi}{p_{n}+1},\\
\textrm{co-}A_{n}=4c\cos\frac{\pi}{p_{n}+1}, & \:(c\;\textrm{a}\:\textrm{parameter})\\
\bar{A}_{n}=4c-A_{n}.
\end{array}\end{array}\:\label{eq:arcs-per-se-1}
\end{equation}
Normalizing to $c=\frac{1}{4}$ in Eqs. (\ref{eq:arcs-per-se-1})
and associating a regular (CFR) $A_{n}\rightarrow\left[\tau_{0}^{(n)};\tau_{1}^{(n)},\tau_{2}^{(n)},\ldots\right]$
to $A_{n}$, one notes that the denominator $\tau_{1}$ satisfies
\begin{equation}
\tau_{1}^{(n+1)}=2\tau_{1}^{(n)}+\delta_{\tau}^{(n)}\qquad(\delta_{\tau}^{(n)}\in\{0,1\},n>2).\label{eq:generalized-neigh-1}
\end{equation}
 In particular, $\tau_{1}^{(n)}$ is the integer part of ${\scriptstyle {\displaystyle \frac{2^{n}}{\pi}}}$:

\begin{equation}
\begin{array}{cccccc}
 & \tau_{1}^{(4)}=5 & \tau_{1}^{(5)}=10\quad & \tau_{1}^{(6)}=20\quad & \tau_{1}^{(7)}=40\quad & \tau_{1}^{(8)}=81\quad\\
\\
 & {\scriptstyle {\scriptstyle {\displaystyle \frac{2^{4}}{\pi}}}}=5.09\quad & {\scriptstyle {\displaystyle \frac{2^{5}}{\pi}}}=10.18\quad & {\scriptstyle {\displaystyle \frac{2^{6}}{\pi}}}=20.37\quad & {\scriptstyle {\displaystyle \frac{2^{7}}{\pi}}}=40.74\quad & {\scriptstyle {\displaystyle \frac{2^{8}}{\pi}}}=81.48\quad
\end{array}\cdots\label{eq:rows-1}
\end{equation}
with ${\scriptstyle {\displaystyle \frac{2^{n}}{\pi}}}$ as a new
expansion parameter. Again, we have an irrational quantity that allows
for time-like refinement in steps of powers of two of the denominators
from a convergent's regular (CFR) $\left(\frac{2^{n}}{\pi}\right)^{-1}\rightarrow\left[b_{0}^{(n)};b_{\alpha}^{(n)}\right]$
and space-like refinement in the cf terms constitutive for Mersenne
flucuations.\label{with--as}%
\footnote{With the option of treating terms as croton amplitudes and generalizing
them to trajectories via (CFR) $\left(\frac{2^{n}}{\pi}\right)^{-1}\pm\mathbb{N}_{0}\rightarrow\left[\tilde{b}_{0}^{(n)};\tilde{b}_{\alpha}^{(n)}\right]$;%
}$^{,}$%
\footnote{the clamp which connects the sources $\left(\frac{2^{n}}{\pi}\right)^{-1}$
and $\left(\frac{2^{n}}{\sqrt{2}}\right)^{-1}$ is the quantity $\delta\in\{0,1\}$
which links a specific cf term at $n$ to the first cf term at $n+1$;
in Eq. (\ref{eq:generalized-neigh-1}), $\delta_{\tau}^{(n)}=1-\left\lceil {\displaystyle \frac{\tau_{2}^{(n)}-1}{\tau_{2}^{(n)}}}\right\rceil $,
while for (CFR) $\left(\frac{2^{n}}{\sqrt{2}}\right)^{-1}\rightarrow\left[b_{0}^{(n)};b_{\alpha}^{(n)}\right]$
the analog to Eq. (\ref{eq:generalized-neigh-1}) is $b_{1}^{(n+1)}=2b_{\textrm{last}}^{(n)}+\delta_{b}^{(n)}$
and $\delta_{b}^{(n)}=1-\left\lceil {\displaystyle \frac{b_{\textrm{ntl}}^{(n)}-1}{b_{\textrm{ntl}}^{(n)}}}\right\rceil $,
with $b_{\textrm{ntl}}^{(n)}$ and $b_{\textrm{last}}^{(n)}$ respectively
denoting the next-to-last and last cf term terminating a period ($\sqrt{2}$
is an algebraic number). %
} But contrary to the case of the square spiral, the cardiod has a
\emph{second} expansion parameter, $c$. It was previously shown \cite{Merkel}
that the croton base number $C_{q_{n}}$ $(q_{n}\equiv(p_{n}-3)/4)$
comes with the identity
\begin{equation}
{\scriptscriptstyle -}\left\lceil {\scriptscriptstyle {\textstyle n/2}}\right\rceil +\Sigma_{i=1}^{n-2}p_{i}=\left\lfloor \log_{2}C_{q_{n}}\right\rfloor \;(n>3),\label{eq:identity-1}
\end{equation}
which makes for an ideal candidate regarding second parametrizing
via $c$. We can therefore conceive of the irrationals
\begin{equation}
\begin{array}{lc}
\textrm{Type}\:\textrm{I}: & \left\lfloor \log_{2}C_{q_{s}}\right\rfloor \end{array}\left(\frac{2^{n}}{\pi}\right)^{-1},\label{eq:full-mers-1-1}
\end{equation}

\begin{equation}
\begin{array}{lc}
\textrm{Type}\:\textrm{II}: & \log_{2}(C_{q_{s}})\end{array}\left(\frac{2^{n}}{\pi}\right)^{-1},\label{eq:full-mers-2-1}
\end{equation}
and

\begin{equation}
\begin{array}{lc}
\textrm{Type}\:\textrm{III}: & \log_{2}(C_{q_{s}})\left(\left\lfloor {\displaystyle \frac{2^{n}}{\pi}}\right\rfloor \right)^{-1}(=\log_{2}(C_{q_{s}})/\tau_{1}^{(n)}).\end{array}\label{eq:trunc-mers-1}
\end{equation}
The symmetry is not a perfect one: Mersenne fluctuations of type I
or II (a shorthand saying they have their habitat in the CFR of irrationals
of type I and II) are fully traceable $-$ as are those based on $\frac{2^{n}}{\sqrt{2}}$.
In contrast, rounding $\frac{2^{n}}{\pi}$ to $\left\lfloor \frac{2^{n}}{\pi}\right\rfloor (=\tau_{1}^{(n)})$
in the CFR of irrationals of type III leads to truncated Mersenne
fluctuations, referred to as Mersenne fluctuations of type III. Truncation
occurs whenever a sequence of like $\delta$'s that determine the
upper row of Eq.(\ref{eq:rows-1}) above breaks: $\delta_{\tau}^{(n)}=\delta_{\tau}^{(n+1)}=\cdots=\delta_{\tau}^{(n+i-1)}\neq\delta_{\tau}^{(n+i)}$.
Conversely, a Mersenne fluctuation of type III is given birth when
a like delta sequence is initiated $-$ to stay in the picture, when
$\delta_{\tau}^{(n-1)}\neq\delta_{\tau}^{(n)}$. Truncated fluctuations
can be hard to assign to a qphylum. The longer the fragment that coincides
with a connected path in a qphylum the lesser the risk of misassignment;
if only few predecessor and successor nodes are available to escort
the insertion, assignment is fraught with uncertainty. Mersenne fluctuations
of type III thus occupy a middle position between randomness and qphyleticly
founded ``volume'' definition.

\subsection{CFR aspect of the examples shown in Sect. 3}

In agreement with the symbol choices of Sect. 3, a regular CFR associated
with irrationals of type I, II or III will be denoted $\left[\varphi_{0};\varphi_{1},\varphi_{2},\ldots\right]$,
while the alternating case corresponding to $b_{0}+\frac{\left.1\right|}{\left|b_{1}\right.}+\frac{\left.-1\right|}{\left|b_{2}\right.}+\frac{\left.1\right|}{\left|b_{3}\right.}+\cdots$
is denoted $\left[\psi_{0};\psi_{1},\psi_{2},\ldots\right]$. With
$\log_{2}C_{q_{s}}\:(q_{s}\equiv(2^{s}-1))$ modulating the outcome
of refinements, the start value of $s$ must be 2 to guarantee a nonvanishing
$\log_{2}C_{q_{s}}$ term; otherwise, the natural numbers $s$ and
$n$ can be chosen freely due to a remarkable property: As was shown
in Sect. 3, the amplitudes and phases produced by Mersenne fluctuations
are by the holographic principle (hgp) linked to the croton bases
underlying the boundary definition. The bases $G^{(15)}$, $J^{(15)}$,
$G^{(31)}$, $J^{(31)},\ldots$ are in turn rooted in the square-matrix
reprentations of $\boldsymbol{f}^{(15)}$, $\boldsymbol{h}^{(15)}$,
$\boldsymbol{f}^{(31)}$, $\boldsymbol{h}^{(31)},\ldots$ If peak
amplitudes fit with boundary $\Gamma^{(p)}$ at $\log_{2}(p+1)=n_{\textrm{hgp}}$
(the same reasoning applies to $\chi^{(p)}$), they can nevertheless
originate in Mersenne fluctuations for which the margins $n_{\textrm{min}}$,$n_{\textrm{max}}\gg n_{\textrm{hgp}}$
because of an in-built recursivity. Using the shorthands LL and UR
for lower left and upper right square-matrix quadrants, the quadrant
$\textrm{LL}(\boldsymbol{f}^{(p)})$ can be shown to coincide with
the subquadrant $\textrm{UR}(\textrm{LL}(\boldsymbol{f}^{(p')}))$,
a recursivity%
\footnote{there are subsubquadrantal (SSQ) relationships that owe their existence
to recursivity too and can be classified in accordance with the equations'
(\ref{eq:arcs-per-se-1}) leading/next-to-leading cf term behavior:
\emph{interordinal} SSQ identities forming structural analogs of $p'=2p+1$
can be classified as sine-like, and \emph{intraordinal} SSQ identities
forming structural analogs of $p''=2p'+1=2(2p+1)+1$ as cosine-like
(see \cite{Merkel} for details)%
} that allows for arbitrarily large assignments $n_{\textrm{min}}$,$n_{\textrm{max}}$
without impeding the amplitude's membership to $\Gamma^{(p)}$. The
assignment of $s$ is similarly open-ended; for technical reasons,%
\footnote{~online CF calculators command a scope of 499 denominators; with
a value of $\pi$ accurate to 1000 decimals, that makes for a limit
$n_{\textrm{max}}\approx1024$ at full coverage, and $n_{\textrm{max}}\approx3030$
at lesser and lesser denominator production%
} however, only fluctuations with $n_{\textrm{max}}\leq3330$ and modulations
$\log_{2}C_{q_{s}}$ with a domain $2\leq s\leq9$ could be considered.

\noindent In what follows, Mersenne fluctuations are listed in the
order they occurred in the text. The fluctuation shown in Fig. \ref{fig:A-prototype-geometric}
is taken from the CFR of $\left\lfloor \log_{2}C_{3}\right\rfloor \cdot$
$\left(\frac{2^{n}}{\pi}\right)^{-1}(=\frac{\pi}{2^{n-1}})$. The
values $\varphi_{\alpha}^{(n)}$ at $n<206$ and $n>218$ have been
omitted. A principal limitation of the CFR approach becomes apparent
at this point: The more frequent and closer to unity $\varphi$ gets,
the less can we tell its affinity.%
\footnote{\,In \emph{Sect}.\ref{sec:Application-to-subatomic}, we will postulate
a mere 2-$3D$-space relatedness of $\varphi_{\alpha}\leq13$ $-$
quasi as the \emph{conditio sine qua non} of continuum illusion.%
} The fact that $\alpha\,\textrm{mod}\,2$ must be invariant and offsets
in $\alpha$ not become too decorrelated from one stopover $n$ to
the next $n+1$ (see Fig. \ref{fig:A-prototype-geometric}) is a help
in telling right from wrong candidates, but that criterion fails if
candidates satisfying $\alpha\,\textrm{mod}\,2$ equivalence come
close to one another. Worst are instances of $\varphi_{\alpha}=1$
with like $\alpha\,\textrm{mod}\,2$ $-$ they are truly legion. 

\noindent The rest of the examples mentioned are, in the right order,
based on the CFR of:\\
\begin{tabular}{ll}
\hspace{-0.2cm}Footnote \ref{fn:One-such-qphylum}: & $\left\lfloor \log_{2}C_{31}\right\rfloor $ $\left(\frac{2^{n}}{\pi}\right)^{-1}$$(1488\leq n\leq1500)$;\tabularnewline
 & $\log_{2}C_{31}\left(\frac{2^{n}}{\pi}\right)^{-1}$$(1068\leq n\leq1082)$;\tabularnewline
 & $\log_{2}C_{31}\left(\frac{2^{n}}{\pi}\right)^{-1}$$(2012\leq n\leq2028)$; \tabularnewline
 & $\left\lfloor \log_{2}C_{31}\right\rfloor $ $\left(\frac{2^{n}}{\pi}\right)^{-1}$$(1951\leq n\leq1959)$; \tabularnewline
\hspace{-0.2cm}Fig. \ref{fig:Pivotal-amplitude-closing}:  & $\left\lfloor \log_{2}C_{127}\right\rfloor $ $\left(\frac{2^{n}}{\pi}\right)^{-1}$$(1556\leq n\leq1559)$; \tabularnewline
\hspace{-0.2cm}Fig. \ref{fig:Pivotal-amplitude-closing-1}, Table
\ref{tab:Pivotal-and-residual-fig4}:  & $\log_{2}C_{3}\left(\frac{2^{n}}{\pi}\right)^{-1}$$(1003\leq n\leq1006)$; \tabularnewline
\hspace{-0.2cm}Fig. \ref{fig:Pivotal-edge-amplitudes}, Table \ref{tab:Pivotal-and-residual-fig5}: & $\log_{2}C_{3}\left(\frac{2^{n}}{\pi}\right)^{-1}$$(987\leq n\leq997)$;\tabularnewline
\hspace{-0.25cm}Table \ref{tab:Pivotal-and-residual-fig4-1}:  & $\log_{2}C_{511}\left\lfloor \left(\frac{2^{n}}{\pi}\right)\right\rfloor ^{-1}$$(n=609)$; \tabularnewline
\hspace{-0.25cm}Table \ref{tab:Pivotal-and-residual-fig4-1-1}: & $\log_{2}C_{127}\left\lfloor \left(\frac{2^{n}}{\pi}\right)\right\rfloor ^{-1}$$(n=1000)$. \tabularnewline
\end{tabular}

\section{Conclusions}

The notable thing about CFR-based Mersenne fluctuations is that whether
the underlying irrational quantity to be refined is an algebraic or
a transcendental number does not matter, as long as there exists an
\emph{interordinal} connection $b_{1}^{(n+1)}=2b_{x}^{(n)}+\delta^{(n)}\:(\delta^{(n)}\in\{0,1\}$.%
\footnote{For (CFR) $\left(\frac{2^{n}}{\sqrt{2}}\right)^{-1}\rightarrow\left[b_{0}^{(n)};b_{\alpha}^{(n)}\right]$,
$b{}_{x}$ is the last cf term in a finite period, while the CFR of
the cardiod-arclength-function argument $\left(\frac{2^{n}}{\pi}\right)^{-1}$
has an infinite period, hence $\tau{}_{x}=\tau_{1}$ in Eq. (\ref{eq:generalized-neigh-1}).%
} The tables shown in \ref{sec:Crotons-in-the} are the outcome of
an in-depth study of denominators emerging with CFR-based Mersenne
fluctuations. Table \ref{tab:Specific-fractions-in-1} summarizes
the results of scanning the CFRs $\left(\frac{2^{n}}{\sqrt{2}}\right)^{-1}\rightarrow\left[b_{0}^{(n)};b_{\alpha}^{(n)}\right]$
for kissing-number matches and hit frequencies, while Table \ref{tab:Specific-fractions-in}
summarizes the corresponding figures for Eqs. (\ref{eq:full-mers-1-1})-(\ref{eq:trunc-mers-1}).
Results for closest pivots and largest peaks are set in parentheses.
A note on largest peaks: While those given in Table \ref{tab:Specific-fractions-in}
for type I and III, 12\,986\,152 and 9\,996\,953, are definitely
beyond $\Gamma^{(31)}$, $\chi^{(31)}$ realizabilities,%
\footnote{~the largest number realizable in $\Gamma^{(31)}$ being 3\,707\,462,
the largest in $\chi^{(31)}$, 4\,177\,840,%
} the largest peak in Table \ref{tab:Specific-fractions-in-1}, 2\,445\,930,
and the largest type-II peak in Table \ref{tab:Specific-fractions-in},
3\,614\,855, are twofold-realizable in $\Gamma^{(31)}$, $\chi^{(31)}$.
That does not mean that a kissing number $L_{m}$ they may have as
target must have an image in both boundaries $\Gamma^{(31)}$ and
$\chi^{(31)}$. For example, many numbers realizable in $\Gamma^{(15)}$,
$\chi^{(15)}$ may in theory become pivotal with respect to $L_{6}(=72)$;
but, although the numbers $L_{m}$ $(m=1,\ldots,7)$ are realizable
in $\Gamma^{(15)}$, $L_{6}$ is not in $\chi^{(15)}$ $-$ it becomes
(twofold-)realizable in $\Gamma^{(31)}$, $\chi^{(31)}$ at last.
By the same token, it may be that peaks such as 2\,445\,930 and
3\,614\,855 do not become true pivots until a target $L_{m}$ for
them realizable in $\Gamma^{(63)}$, $\chi^{(63)}$ is found. Unfortunately,
kissing number candidates $L_{m}\:(m>31)$ are, with the exception
$L_{48}(=52\,416\,000)$, notoriously uncertain. For the time being,
more matches with $L_{m}\:(m\leq31)$ and interesting pivots may only
be obtained by enlarging the scope of $s$. 

\noindent At the very beginning, however, there is a master Mersenne
fluctuation that pauses at $n=1$: If we take the geometrizing hypothesis
via Mersenne fluctuations literally, the continued fraction $[0;1,\bar{2}]$
$-$ a special case of (CFR) $\left(\frac{2^{n}}{\sqrt{2}}\right)^{-1}\!\!\rightarrow\!\left[b_{0}^{(n)};b_{\alpha}^{(n)}\right]$
for $n=1$ $-$ means ``recursive geometrization into a centerpiece-free
pair of a 0-spheres'' or, a self-similar laminar pattern ``dash-space-dash''
for all space-like refinements. \newpage{}

\part{\ }

\section{\label{sec:Pregeometric-categories-relevant}Pregeometric categories
relevant to physics}

Mersenne fluctuations and qphyla may be viewed as \emph{global} categories
that mediate between infinite expansion and infinite refinement and
thereby provide the center stage for physics. This is envisioned here
in form of a three-step process. 

\noindent (1) Co-occurrent coincidences 
\[
\varphi_{\alpha}^{(n)}+\Delta\varphi^{(n)}=L_{m}+\delta\:(\delta\in\left\{ 0,1\right\} )
\]
and, to a lesser degree, coincidences within  time-like refinements $n_{1},n_{2},\scriptstyle{\ldots}$ lying close to each other,
 
\[
\varphi_{\alpha_{1}}^{(n_{1})}+\varphi_{\alpha_{2}}^{(n_{2})}+\ldots=L_{m}+\delta\:(\delta\in\left\{ 0,1\right\} ),
\]

\noindent are seen as\emph{ ideations} of space and matter. \\
(2) Materialization is accomplished via organizers,\emph{ }in organization
centers
\[
\begin{array}{c}
\varkappa_{\alpha}^{(n)}+\Delta\varkappa^{(n)}\rightarrow L_{m}+\delta,\\
\left(\varkappa'\right)_{\alpha}^{(n)}+\Delta\left(\varkappa'\right)^{(n)}\rightarrow L_{m}+\delta,
\end{array}\:(\delta\in\left\{ 0,1\right\} )
\]

\noindent or, to a lesser degree, dispersed in time and space,
\[
\begin{array}{c}
\varkappa{}_{\alpha_{1}}^{(n_{1})}+\varkappa_{\alpha_{2}}^{(n_{2})}+\ldots\rightarrow L_{m}+\delta,\\
\left(\varkappa'\right)_{\alpha_{1}}^{(n_{1})}+\left(\varkappa'\right)_{\alpha_{2}}^{(n_{2})}+\ldots\rightarrow L_{m}+\delta,
\end{array}\:(\delta\in\left\{ 0,1\right\} )
\]

\noindent where $2^{-n}\kappa\!\rightarrow\!\left[\varkappa_{0}^{(n)}\!;\varkappa_{\alpha}^{(n)}\right]$
and $2^{-n}\kappa'\!\rightarrow\!\left[\left(\varkappa'\right)_{0}^{(n)}\!;\left(\varkappa'\right)_{\alpha}^{(n)}\right]$
are CFRs of $\kappa$ and $\kappa'$, irrational quantities which
we will introduce later and use in approximations of the electromagnetic
and weak coupling constant, respectively. \\
(3) The result of the organization is a `world container' whose expanses
$3^{2^{n}}$ are `filled' with \emph{local} categories. Two kinds
of crotons making up the local categories can be distinguished. Crotons
of a first kind spring from $3^{2^{n}}$ divided by its first `derivative,'
\begin{equation}
\frac{3^{2^{n}}}{\frac{\partial(3^{2^{x}})}{\partial x}\left|_{x=n}\right.}=\left(2^{n}\log(2)\log(3)\right)^{-1},\label{eq:velocity}
\end{equation}
and form the time- and space-like refinement scheme of a charge in
(uniform) motion
\begin{equation}
(\textrm{CFR})\:\left(2^{n}\log(2)\log(3)\right)^{-1}\!\rightarrow\!\left[\between_{\,0}^{(n)}\!;\between_{\,\alpha}^{(n)}\right].\label{eq:internal-qphyla-1}
\end{equation}
Crotons of a second kind spring from $3^{2^{n}}$ over its second
`derivative,' 
\begin{equation}
\frac{3^{2^{n}}}{\frac{\partial^{2}(3^{2^{x}})}{\partial x^{2}}\left|_{x=n}\right.}=\left(2^{n}\log^{2}(2)\log(3)(2^{n}\log(3)+1)\right)^{-1},\label{eq:accelerate}
\end{equation}
and form 
\begin{equation}
(\textrm{CFR})\:\left(2^{n}\log^{2}(2)\log(3)(2^{n}\log(3)+1)\right)^{-1}\!\rightarrow\!\left[\gamma_{0}^{(n)}\!;\gamma_{\alpha}^{(n)}\right].\label{eq:radiation}
\end{equation}
Either kind is more than an ideation via plain inclusion of charge
and motion. We postpone the discussion of organizers that enable those
inclusions and deal with the `world container' aspect first. Crotons
of the second kind (Eq. (\ref{eq:radiation})) show a pecularity:
Increasing or decreasing $n$ does \emph{not} help them trace a particular
Mersenne fluctuation, much less a qphyletic embedding. We associate
this category with the time-like and space-like refinement scheme
of radiation emitted by an accelerated charge. Everytime radiation
is created this marks an instant of time-like refinement $n_{0}$,
and the exacting task is to determine its characteristics for observation
points $n>n_{0}$. (Classical electrodynamics' incomplete description
of radiation is reverberating here.) Crotons of the first kind (Eq.
(\ref{eq:internal-qphyla-1})), by contrast, entail Mersenne fluctuations
and are embeddable in qphyla in the usual way. So we tackle them first.

\subsection{Crotons of the first kind}

They form a kind of Charles Dodgson universe whose underlying principles
are:

\bigskip{}
\begin{minipage}[t]{0.9\columnwidth}%
{\small{1. It is represented by a qphyletic assembly of Mersenne fluctuations.}}{\small \par}
\smallskip{}
{\small{2. Time-like refinements determine the course of time via
continued halving of the time unit (time order).}}{\small \par}
\smallskip{}
{\small{3. As necessitated by the notion of particles (constituents)
in motion, the Mersenne fluctuations'  terms represent `lengths' $-$
replacing what was previously called `amplitudes.' The `length order'
comes in two forms: \mbox{(near-)}dou\-bling of the `intrinsic length 1' bottom-up,
(near-)halving of the `peak length' top-down.}}{\small \par}
\smallskip{}
{\small{4a) On a Mersenne fluctuation's left leg, time order and `length order' run antiparallel.
Depending on the direction of time, `lengths' inflate or deflate;
the co-represented particles (constituents) are {\em virtual} ones.
At time one, a term owes its existence to a seed at time zero which,
depending on its space-like refinement, can be almost arbitrarily
large; among the largest terms attainable, the one that turns out to
continue in left-leg mode from time one on and reach its peak only
after a huge number of time-like refinement steps, determines the
size and age of the universe.}}{\small \par}

{\small{4b) On a Mersenne fluctuation's right leg (including lower-lying peaks), time order
and `length order' run parallel and constitute uniform motion; the
co-represented particles or,
depending on time direction, antiparticles $-$ or constituents thereof $-$ are {\em real}.}}{\small \par}
\smallskip{}
{\small{5. A `low-expansion' regime $-$ one (low-lying) left-leg step up, one (low-lying) right-leg
step down$-$ mirrors time-zero conditions.}}%
\end{minipage}

\newpage

\noindent Points 1 and 2 follow from our definitions in Part I; point
3 does as well once the motion concept is incorporated. Next are points
4a and 4b. The continued fraction of $1/\log(2)$ has been computed
by E. Weisstein to 9\,702\,699\,208 terms \cite{Weisstein}  $-$ the
largest that he found, 53\,155\,160\,769, at, using our locution,
space-like refinement $\alpha=$\,2\,565\,310\,827. No one has
since computed that many terms of the continued fraction of $1/[(\log(2)\log(3)]$,
but we may safely assume that similarly large terms arise for large
enough $\alpha$. As regards the two thousand five hundred or so terms
the Wolfram$\mid$Alpha app provides, the two largest ones are, using our denotation, $\between_{2398}^{(0)}=1752$
and $\between_{2308}^{(0)}=838$. The first is seed to an (uninteresting)
right-leg term $\between_{2402}^{(1)}=876$, the second a seed to
an inflaton $\between_{2268}^{(1)}=1677$, or miniscule universe,
indeed. The highest Mersenne fluctuation provided by the app, of $h=24$,
has a left leg that runs from $\between_{\,\textrm{?}}^{6104}=1$ via $\between_{218}^{(6107)}=13$
to $\between_{200}^{(6127)}=14\,571\,717$. If seeds comparable in
size to or larger than those found in the Weisstein quest do indeed exist, they
may combine with a giant left leg to puff up a Charles Dodgson universe
of respectable size and age. On right legs, a halving of the time
unit $2^{-n}\rightarrow2^{-n-1}$ is accompanied by a \mbox{(near-)}\-halving
of   `length' $\between_{\alpha}^{(n)}\rightarrow\left\lfloor \between_{\alpha}^{(n)}\!\!/2\right\rfloor -{\bar{\delta}}$.
For the real particle or constituent co-represented, it indicates
their uniform motion. The reverse holds for antiparticles: the reverse
time direction implies a  doubling of the time unit $2^{-n}\rightarrow2^{-n+1}$
and is accompanied by a \mbox{(near-)}dou\-bling of  `length' $\between_{\alpha}^{(n)}\rightarrow2\!\between_{\alpha}^{(n)}\!+1+\epsilon$.
Such congruences are absent on the left leg. To time halvings $2^{-n}\rightarrow2^{-n-1}$
there $-$ our normal course of time $-$ correspond relative `length \mbox{(near-)}quadruplings' and signal inflation,
and to $2^{-n}\rightarrow2^{-n+1}$ correspond relative `length \mbox{(near-)}quarterings'
or deflation, and, because such behavior is unobserved in real particles,
co-represented particles (constituents) must be virtual ones. The
tag `virtual' understates their part, though. When imposing normal time
order over the entire course of the Mersenne fluctuation, a virtual particle on
its way up the left leg and endowed with `intrinsic length 1'  ($\leftarrow$ `length
order' bottom-up) is by the relationship
\[
\left|\between_{\alpha_{n-r}}^{(n-r)}-\between_{\alpha_{n+r}}^{(n+r)}\right|=\left\{ \begin{array}{cc}
0\textrm{ or }1 & \;0\!<r\!<h\!-\!1\\
0 &\; r=h\!-\!1
\end{array}\right.\!\!\!\!\,\,\,\ (h\textrm{ the fluctuation's height}),
\]
tied to its counterpart, a real particle on its way down the right
leg and endowed with `peak length' ($\leftarrow$ `length order' top-down).
For some $r$, one of the two, or both, may reveal their $L_{m}$-based
(particulate) identity. This cannot be a Moir{\'e} pattern, since
there is no sustained space involved. Instead, we would suspect that
the qphyletic assembly is responsible for the occasional coincidence
of `current length' and $L_{m}$-based identity.

\noindent As to point 5, even though not rigorously established, Khinchin's constant $\mathcal{K}_{0}\approx2.685452\scriptstyle{\ldots}$
seems to apply to the continued
fraction of $1/\log(2)$ \cite{Weisstein}. Maybe the continued
fraction of $1/[\log(2)\log(3)]$ is a further candidate. Successive
values $(\between_{1}^{(0)})^{1/1}$, $(\between_{1}^{(0)}\between_{2}^{(0)})^{1/2}$,
$(\between_{1}^{(0)}\between_{2}^{(0)}\cdots\between_{k}^{(0)})^{1/k}$,
if true, would then converge to $\mathcal{K}_{0}$. It's natural to
expect that a `low-expansion' regime based on one-step creations of
the terms 2 and 3 (the roundings off and up of $\mathcal{K}_{0}$)
would, of their own, approach $\mathcal{K}_{0}$ as closely as possible.
This appears to be the case. Starting from a left-leg term $\between_{\alpha_{n-h+1}}^{(n-h+1)}=1$
($h$ the relevant Mersenne fluctuation's height, and $n\geq h$), 
the outcome of an upward transformation $\between_{\alpha_{n-h+1}}^{(n-h+1)}\rightarrow2\!\between_{\alpha_{n-h+1}}^{(n-h+1)}\!+1+\epsilon$
is three possible terms: $1\rightarrow2,3,4$. Sequences of terms containing
2 and 3 are obtained from the right-leg term $\between_{\alpha_{n+h-3}}^{(n+h-3)}\in\{4,5,6,7,8\}$
via the downward transformation $\between_{\alpha_{n+h-3}}^{(n+h-3)}\rightarrow\left\lfloor \between_{\alpha_{n+h-3}}^{(n+h-3)}\!\!/2\right\rfloor -{\bar{\delta}}$.
Seven possible oucomes can be distinguished: $8\rightarrow4,3$; $7\rightarrow3$;
$6\rightarrow3,2$; $5\rightarrow2$; $4\rightarrow2$. Thus, the
associated geometric mean $(2\cdot3\cdot4\cdot4\cdot3\cdot3\cdot3\cdot2\cdot2\cdot2)^{1/10}\approx2.701920$
mirrors time-zero conditions provided they are characterized by an
average `length' unit equal to Khinchin's constant.\footnote{\,According to Barcus {\em et al.}\,\cite {Barcus}, a physically justifiable bound eleventh order polynomial regression of  the most recent proton radius measurements yields a value of 0.854 fm. This  would make the fm equal  $\frac{\pi}{\mathcal{K}_{0}}$ proton radii.} 
\newline
 To illustrate how a pre-geometric quantity like $\between_{\alpha}^{(n)}$
can be akin to `length' $-$ more specifically, `wavelength' $-$, let us start
with a nontrivial example: two motion analogs, a) of $L_{2}+1$ and
b) of $L_{2}$; they are elongated in one direction and represent
a squashed and a stretched 2$D$ space chunk, respectively:

\smallskip{}
{\label{honey-squash-str}\includegraphics[scale=0.08]{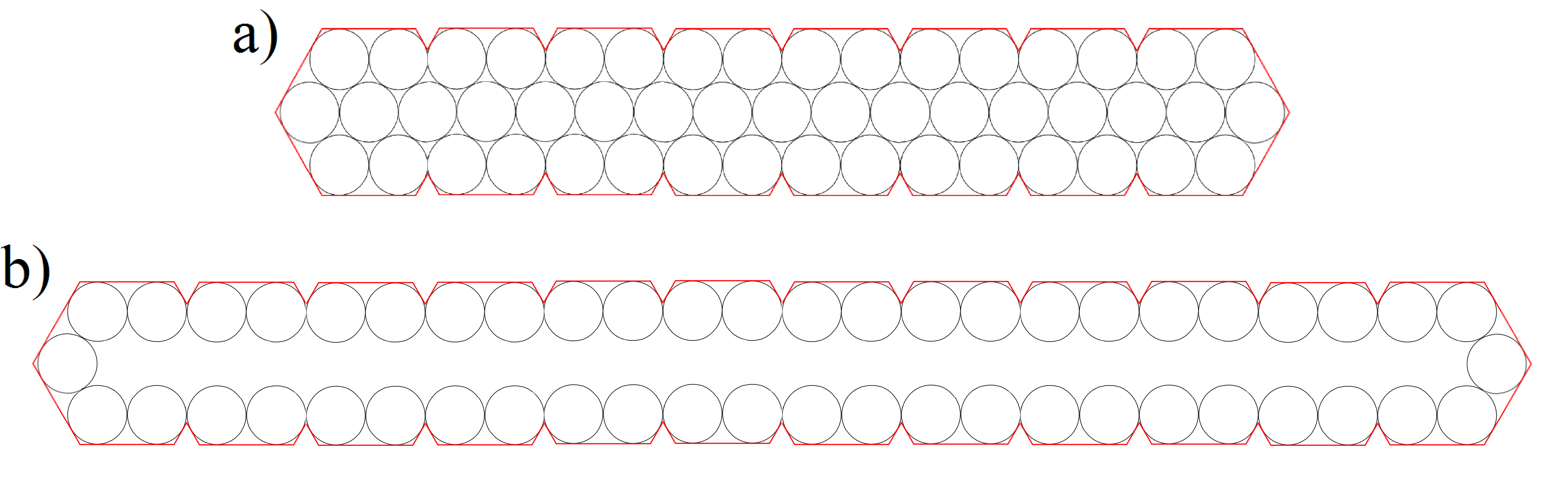}
\smallskip{}

\noindent The graph exhibits two pecularities: (1) 1-spheres, potential chargeholders,
are enclosed in interlaced honeycombs which may arise thus: ordinary
lattice honeycombs are picked up, then tilted (by 30 degrees or a multiple
thereof) and interlaced with each other until they overlap in one
1-sphere position; configuration a) has forty-nine 1-spheres in eight
interlaced honeycombs, configuration b) fifty of them in twelve interlaced
honeycombs. (2) In a), there is a homogeneous array of 3-star interstices
$-$ adjacent lattice honeycombs, by contrast, have also 4-star interstices
at their interfaces. In b), we find just one $(6+4k)$-star interstice, namely 
for $k=11$. \\
 Thus far, only crotons coinciding with $L_{m}+\delta\:(\delta\in\left\{ 0,1\right\} )$
have been considered successful sphere-packing creations. When they
are lumped together in larger chunks, however, `fails' equaling $L_{m}-1$
are perfectly admissible, too. Let us apply the sphere-packing notion,
$L_{m}$, together with the supplements 
\[
L_{m}^{+}=L_{m}+1,\quad L_{m}^{-}=L_{m}-1,
\]
to the 2$D$ lattice honeycombs \\

 \qquad{}\includegraphics[scale=0.42]{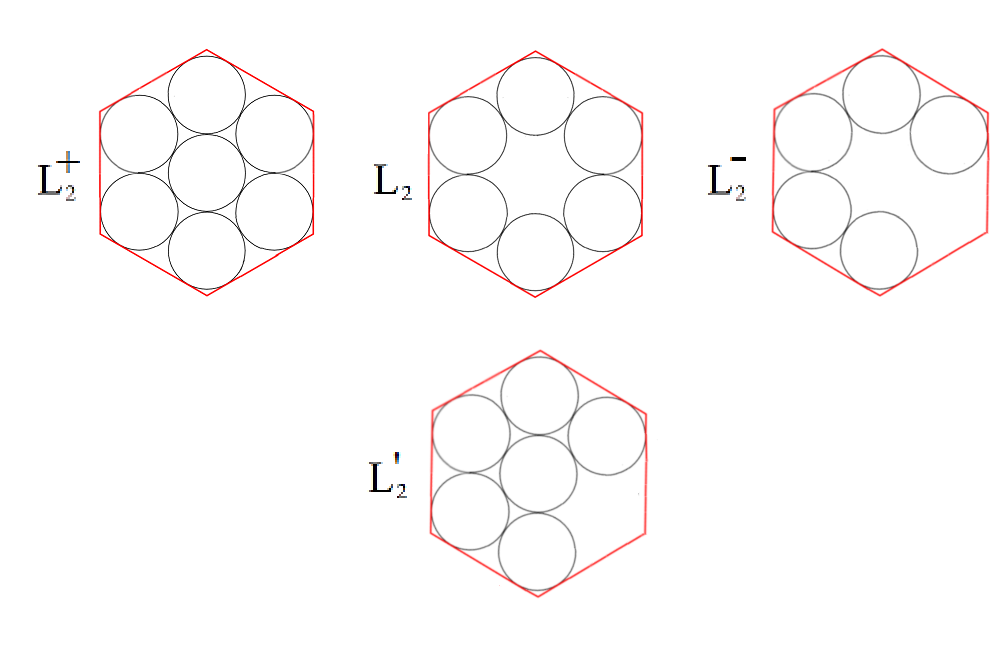}
 \\
 and also assume, with the context charge-in-motion in
mind, that 1-spheres get movable $-$ the 1-sphere at the  near-southeast
position in the $L_{2}$ honeycomb, for instance, may wander to the
center as shown in $L'_{2}$. It is easy to see that if you take seven
copies of the $L'_{2}$ honeycomb, tilt each by -30 degrees and
interlace the first of them with the likewise tilted $L_{2}$ honeycomb,
then interlace the duo with the remaining copies one at a time and
always in the same direction, you'll get the squashed chunk a). Regarding
b), tilting/interlacing an $L_{2}^{-}$ honeycomb with a mirrored
copy of itself would  only yield a shortened b) version $-$ one containing
ten 1-spheres; all 1-spheres have to get movable for a version containing
fifty of them to arise, a relaxation that opens up the possibility
to do statistics on the various ways a) and b) can be formed from
$L_{2}^{+},L_{2}$ and $L_{2}^{-}$ (here, eight different ones each).
\bigskip{}

\noindent Formation of a) :

\bigskip{}

\qquad{}%
\begin{tabular}{|c|c|c|c|c|}
\hline 
$aL_{2}^{+}$  & $bL_{2}$  & $cL_{2}^{-}$  & \multirow{1}{*}{$a+b+c$}  & \# of 1-spheres\tabularnewline
\hline 
\hline 
$7L_{2}^{+}$  &  &  & \multirow{1}{*}{7}  & 49\tabularnewline
\hline 
$4L_{2}^{+}$  & $1L_{2}$  & $3L_{2}^{-}$  & 8  & \textquotedbl{}\tabularnewline
\hline 
$3L_{2}^{+}$  & $3L_{2}$  & $2L_{2}^{-}$  & 8  & \textquotedbl{}\tabularnewline
\hline 
$2L_{2}^{+}$  & $5L_{2}$  & $1L_{2}^{-}$  & 8  & \textquotedbl{}\tabularnewline
\hline 
$1L_{2}^{+}$  & $7L_{2}$  &  & 8  & \textquotedbl{}\tabularnewline
\hline 
$2L_{2}^{+}$  &  & $7L_{2}^{-}$  & 9  & \textquotedbl{}\tabularnewline
\hline 
$1L_{2}^{+}$  & $2L_{2}$  & $6L_{2}^{-}$  & 9  & \textquotedbl{}\tabularnewline
\hline 
 & $4L_{2}$  & $5L_{2}^{-}$  & 9  & \textquotedbl{}\tabularnewline
\hline 
$\Sigma a=20$  & $\Sigma b=22$  & $\Sigma c=24$  & $\Sigma=66$  & \tabularnewline
\hline 
\end{tabular}\bigskip{}

\newpage
\noindent Formation of b) :\smallskip{}

\qquad{}%
\begin{tabular}{|c|c|c|c|c|}
\hline 
$aL_{2}^{+}$  & $bL_{2}$  & $cL_{2}^{-}$  & \multirow{1}{*}{$a+b+c$}  & \# of 1-spheres\tabularnewline
\hline 
\hline 
$5L_{2}^{+}$  &  & $3L_{2}^{-}$  & \multirow{1}{*}{8}  & 50\tabularnewline
\hline 
$4L_{2}^{+}$  & $2L_{2}$  & $2L_{2}^{-}$  & 8  & \textquotedbl{}\tabularnewline
\hline 
$3L_{2}^{+}$  & $4L_{2}$  & $1L_{2}^{-}$  & 8  & \textquotedbl{}\tabularnewline
\hline 
$2L_{2}^{+}$  & $6L_{2}$  &  & 8  & \textquotedbl{}\tabularnewline
\hline 
$2L_{2}^{+}$  & $1L_{2}$  & $6L_{2}^{-}$  & 9  & \textquotedbl{}\tabularnewline
\hline 
$1L_{2}^{+}$  & $3L_{2}$  & $5L_{2}^{-}$  & 9  & \textquotedbl{}\tabularnewline
\hline 
 & $5L_{2}$  & $4L_{2}^{-}$  & 9  & \textquotedbl{}\tabularnewline
\hline 
 &  & $10L_{2}^{-}$  & 10  & \textquotedbl{}\tabularnewline
\hline 
$\Sigma a=17$  & $\Sigma b=21$  & $\Sigma c=31$  & $\Sigma=69$  & \tabularnewline
\hline 
\end{tabular}\bigskip{}

\noindent The combined proportions $\Sigma a$ :$\,\Sigma b$ :$\,\Sigma c$,
37:43:55, can be compared to those ari\-sing in the incidences $I\left(\between_{\,\alpha}^{(n)}=x\mid x\in\left\{ L_{2}^{+},L_{2},L_{2}^{-}\right\} ;n\in P_{30},\alpha\!<\!500\right)$,
where $P_{30}$ denotes the first thirty primes. For this comparison
to make sense, we draw on the `length order' bottom-up which endows
virtual particles (constituents) with `intrinsic length 1,' corresponding
to one honeycomb, here. Then $I$ encompasses those cases where a honeycomb
identity of the virtual particles and/or their right-leg real counterpart
is being revealed; the corresponding proportions coming out as 359:452:569
(or roughly ten times the former), they underpin the role of uniform
motion as a common underlying context. \vspace*{2mm}

\noindent It is obvious that there is a link between the number of
lattice honeycombs contiguously aligned in one direction, $\lambda=a+b+c$,
and the intrinsic number of interlaced honeycombs, $\underline{\lambda}$
for a) and $\overline{\lambda}$ for b). On the `lattice wavelength'
$\lambda$, then,\emph{ a}spect-\emph{o}riented \emph{t}ransformations
may be applied to match the `interlaced wavelengths' $\underline{\lambda}$,
$\overline{\lambda}$. The transformations turn out to be a combination
of a $x$ vs.$\,(x+2)$ or $x$ vs.$\,(x+1)$ juxtaposition and `Doppler'
terms $\mp1$:

\[
\begin{array}{c}
\textrm{tao}_{1}(x):=x+2-1,\\
\textrm{tao}_{3}(x):=x+2+1,\\
\textrm{tao}_{0}(x):=x+1-1,\\
\textrm{tao}_{2}(x):=x+1+1.
\end{array}
\]
Thus, 
\[
\begin{array}{c}
\textrm{tao}_{1}(7)=\textrm{tao}_{0}(8)=\underline{\lambda}\,(=8)\qquad\quad\left[\textrm{case a})\right],\\
\\
\textrm{tao}_{3}(9)=\textrm{tao}_{2}(10)=\overline{\lambda}\,(=12)\qquad\left[\textrm{case b})\right].
\end{array}
\]

\noindent One of the basic tenets of uniform motion is that configurations
in relative motion to one another are equally entitled to stationary
status. Under $\underline{\lambda}=8$, $\lambda=8$ would connote
a stationary status and $\lambda=7$, $\lambda=9$ a Doppler-shifted
pair in relative motion to it; under $\overline{\lambda}=12$, $\lambda=9$
would connote a stationary status and $\lambda=8$, $\lambda=10$
Doppler-shifted twins. To grant all of them equal right, we have to
invoke the inverse of at least one of the transformations $\textrm{tao}_{1},\textrm{tao}_{3}$: 
$
\begin{array}{c}
\textrm{tao}_{1}^{-1}(x):=x-2+1,\end{array}
\!\!$
say, in order to get 
\[
\begin{array}{c}
\textrm{tao}_{1}(7)=\textrm{tao}_{1}^{-1}(9)=8,\\
\\
\textrm{tao}_{1}(8)=\textrm{tao}_{1}^{-1}(10)=9.
\end{array}
\]
If the number of 1-spheres is uncertain, \emph{i.e}., if one can't
discern between case a) and case b), the status of motion is indeterminate
in general. In accordance with the above equations, $\lambda=8$ would
connote a stationary status with respect to the triple $\lambda=7,8,9$
and a nonstationary one with the triple $\lambda=8,9,10$. Vice versa
for $\lambda=9$. Of course, a wavelength $\lambda=7$ on one hand
and $\lambda=10$ on the other would still point the way to case a)
and case b), respectively. Alas, this is but a unique property of
$2D$ space chunks $-$ perhaps because they are exempt from the loss
of interior-defining interstices that follows for analogs of a) and
b) with $m>2$. Writing 
\[
\begin{array}{c}
\lambda=7:=\left(\lambda_{\textrm{min}}\right)_{2},\qquad\underline{\lambda}=8:=\underline{\lambda}_{2},\\
\\
\lambda=10:=\left(\lambda_{\textrm{max}}\right)_{2},\qquad\overline{\lambda}=12:=\overline{\lambda}_{2},
\end{array}
\]
a natural generalization would be 
\begin{equation}
\begin{array}{c}
L_{m}^{-}\left(\lambda_{\textrm{max}}\right)_{m}-1=L_{m}^{+}+L_{m}\left(\lambda_{\textrm{min}}\right)_{m},\end{array}\qquad(m=2,3,\ldots)\label{eq:natgen}
\end{equation}
which becomes an identity for 
\begin{equation}
\left(\lambda_{\textrm{max}}\right)_{m}:=2L_{m}^{-},\qquad\left(\lambda_{\textrm{min}}\right)_{m}:=2L_{m}-5.\label{eq:iden-natgen}
\end{equation}

\noindent To a) then corresponds a general squashed space chunk with
$2\left(L_{m}^{-}\right)^{2}-1$ $(m\!-\!1)$-spheres and characteristic
`wavelength' $\underline{\lambda}_{m}=2L_{m}-3$; and to b), a general
stretched space chunk with $2\left(L_{\ensuremath{m}}^{-}\right)^{2}$
$(m\!-\!1)$-spheres and a characteristic `wavelength' $\overline{\lambda}_{m}=2L_{m}$.
And the number $w_{m}$ of different alignments of lattice polytopes
to be interlaced to produce a) and b) analogs would recursively be
given by 
\begin{equation}
\begin{array}{c}
w_{m}=w_{m-1}+L_{m}\qquad(w_{1}\equiv2).\end{array}\label{eq:recurways}
\end{equation}
We may now check, for $m=3$ for instance, whether combined proportions
and incidences $I(\,\between\,)$ remain in unison. The squashed $3D$
space chunk, which contains 241 spheres, would result from the 20
different alignments\vspace{3mm}

\noindent $\lambda=19\left(=(\lambda_{\textrm{min}})_{3}\right)$:
\[
(13+k)L_{3}^{+}+(6-2k)L_{3}+kL_{3}^{-}=241,\quad k=0,1,2,3;
\]
$\lambda=20$: 
\[
(1+k)L_{3}^{+}+(19-2k)L_{3}+kL_{3}^{-}=241,\quad k=0,1,\ldots,9;
\]
$\lambda=21$: 
\[
kL_{3}^{+}+(10-2k)L_{3}+(11+k)L_{3}^{-}=241,\quad k=0,1,\ldots,5;
\]

\noindent and the stretched $3D$ space chunk from the 20 alignments\vspace{3mm}

\noindent $\lambda=19\left(=(\lambda_{\textrm{min}})_{3}\right)$:
\[
(16-k)L_{3}^{+}+(1+2k)L_{3}+(2-k)L_{3}^{-}=242,\quad k=0,1,2;
\]
$\lambda=20$: 
\[
(2+k)L_{3}^{+}+(18-2k)L_{3}+kL_{3}^{-}=242,\quad k=0,1,\ldots,9;
\]
$\lambda=21$: 
\[
kL_{3}^{+}+(11-2k)L_{3}+(10+k)kL_{3}^{-}=242,\quad k=0,1,\ldots,5;
\]
$\lambda=22\left(=(\lambda_{\textrm{max}})_{3}\right)$: 
\[
22L_{3}^{-}=242;
\]
While the `wavelength' $(\lambda_{\textrm{min}})_{3}=19$ turns out
to be minimal in either 3\emph{D} space chunk form, its counterpart
$(\lambda_{\textrm{max}})_{3}=22$ unswervingly points towards the
stretched $3D$ space chunk $-$ and continues this way all along
$m=4,5,\scriptstyle{\ldots}$. In case the chunk's number of spheres is uncertain,
this is clearly an asymmetry. And, whilst transformations $\textrm{tao}{}_{v}^{\pm1}$
$(v=0,1,2,3)$ still suffice to convert $\lambda$s into one another,
the agreement between the combined proportions $\Sigma a$ :$\,\Sigma b$
:$\,\Sigma c$ and those from, say, the incidences $I\left(\between_{\,\alpha}^{(n)}\!=\! x\mid x\!\in\!\left\{ L_{3}^{+},\! L_{3},\! L_{3}^{-}\right\} \!;n\in P_{60},\alpha\!<\!500\right)$ $-$
all $x$ taken either from virtual particles (constituents) in `length'
order bottom-up or from their right-leg real counterparts $-$ begins
to deteriorate: 253\,:\,277\,:\,277 vs. 184\,:\,249\,:\,278.
(We'll later see that, from $m=4$ on, a different order emerges that
modulates $I$, as will
be illustrated for $8\leq m\leq11$ in Sect. \ref{sub:The-plateau-effect}.)

\noindent We now change the perspective and consider representation
of  `length' relative to the present, $n$, asking what the infinity of $\between_{\,\alpha}^{(n)}\notin\!\left\{ L_{m}^{+},L_{m},L_{m}^{-}\right\} $
$(m>3)$ might mean from that vantage point. To this end, we narrow the range
of answers to the cases 
\begin{equation}
\between_{\,\alpha}^{(n)}=x_{m}\pm\Delta_{m}\quad x_{m}\!\in\!\left\{ L_{m}^{+},L_{m},L_{m}^{-}\right\} \label{eq:close to L_m}
\end{equation}
under the constraint 
\begin{equation}
\begin{array}{c}
x_{m}\pm\Delta_{m}\lessgtr x_{m\pm1}\mp\Delta_{m}.\end{array}\label{eq:close L_m constraint}
\end{equation}
The first thing that catches the eye is the indeterminacy of $x_{m}$.
Its identification with one of the elements $\left\{ L_{m}^{+},L_{m},L_{m}^{-}\right\} $
depends  on the type of correlation between the element in question
and one or more $\between_{\,\beta}^{(n)}$ or $\between_{\,\beta'}^{(n')}$
that make up $\Delta_{m}$. It works like with entanglement: given
$\between_{\,\alpha}^{(n)}$, determine $\Delta_{m}$ in terms of $\between_{\,\beta}^{(n)}$
or $\between_{\,\beta'}^{(n')}$ and you know the outcome of $x_{m}$.
We have proposed diverse forms of correlation in Part I, and more
of them are proposed in Sect.\ref{sub:The-role-of-centralizers} below
when dealing with organizers. Here, the tool we use to analyze $\Delta_{m}$,
once the correlation is known, is an accumulator 
\[
\Delta_{m}=\textrm{Pre}(\mathfrak{q}).
\]
Remotely reminiscent of a Fourier expression, $\mathfrak{q}=\sum_{s \in \mathfrak{S}}\mathfrak{q}_{s}$
$(\mathfrak{S}\textrm{ a subset of }$ $\{1,{\scriptstyle{\ldots}},\textstyle m-2\})$ sums up quaternions
of the form $\mathfrak{q}_{s}\!=u_{s}(a_{0,s}+a_{1,s}\boldsymbol{i}+a_{2,s}\boldsymbol{j}+a_{3,s}\boldsymbol{k})$.
They obey Hamilton's equations for quaternion  units $\boldsymbol{i}^{2}=\boldsymbol{j}^{2}=\boldsymbol{i}\boldsymbol{j}\boldsymbol{k}=-1$
and have  `length' terms as coefficients,
\[
\begin{array}{c}
a_{0,s}=\textrm{sgn}(y)\left\lfloor \left|y\right|\right\rfloor \:\textrm{with}\: y=L_{m-s}\cos(2\pi \nu/(\lambda_{\textrm{min}})_{m-s}),\\
a_{1,s}=\textrm{sgn}(y)\left\lfloor \left|y\right|\right\rfloor \:\textrm{with}\: y=L_{m-s}\sin(2\pi\nu/(\lambda_{\textrm{min}})_{m-s}),\\
a_{2,s}=\textrm{sgn}(y)\left\lceil \left|y\right|\right\rceil \:\textrm{with}\: y=L_{m-s}\cos(2\pi\nu/(\lambda_{\textrm{max}})_{m-s}),\\
a_{3,s}=\textrm{sgn}(y)\left\lceil \left|y\right|\right\rceil \:\textrm{with}\: y=L_{m-s}\sin(2\pi \nu/(\lambda_{\textrm{max}})_{m-s}),
\end{array}\quad(s\in \mathfrak{S})
\]
 where  (1) $a_{0,s},a_{1,s}$ entail rounding toward zero, $a_{2,s},a_{3,s}$
rounding away from zero; (2) each $u_{s}$ is a quaternion
with unit components in $\left\{ \pm1,\pm\boldsymbol{i},\pm\boldsymbol{j},\pm\boldsymbol{k}\right\} $ and (3) $\nu:=n$ in case of co-occurrence, $\nu:=n'$ else.
The function Pre() is  equivalent to applying a function $f()$  to the real part of the sum of  quaternions
in the subset: $\textrm{Pre}(\mathfrak{q})=f(\textrm{Re}(\sum_{s \in \mathfrak{S}}\mathfrak{q}_{s}))$. $f()$ takes account of the time order of the decomposition
of $\Delta_{m}$ in terms of those $a_{\mu,s}$ that do not vanish. For a co-occurrent $\between_{\beta}^{(n)},$ a term
$a_{\mu,s}$ becomes part of a $\Delta_m$ that originates  in the present by definition and
$f$ is just the identiy; for $n-n'=1$, $a_{\mu,s}$ becomes part of a  recall that  installs a `stored'
item of the past  in the present;  thusly,
on a left-leg ascent $\between_{\alpha_{n-1}}^{(n-1)}$ $\rightarrow$ $\between_{\alpha}^{(n)}$, $f$ becomes $f(x)=2x+1+\epsilon$ and on a right-leg descent
 $\between_{\alpha_{n-1}}^{(n-1)}$ $\rightarrow$ $\between_{\alpha}^{(n)}$, $f(x)=\left\lfloor x/2\right\rfloor -\bar{\delta}$.
For $n-n'=-1$, instead, $a_{\mu,s}$ becomes part of an anticipation that installs a future item in the present;
here, $f(x)=\left\lfloor x/2\right\rfloor -\bar{\delta}$
on  a left-leg ascent $\between_{\alpha}^{(n)}$ $\rightarrow$ $\between_{\alpha_{n+1}}^{(n+1)}$ and  $f(x)=2x+1\!+\epsilon$ on  a right-leg descent $\between_{\alpha}^{(n)}$ $\rightarrow$ $\between_{\alpha_{n+1}}^{(n+1)}$. 

\noindent The admissible number
of terms  per quaternion component in the quaternion sum is   a rule-based pick out of $\,2, 4$. When   the design of programming languages was in its infancy, the pros and cons of COMEFROM vs. $\!$GOTO were hotly debated, notions that can be used to advantage here with the co-occurrence problem. The rules for picking the admissible number of terms, first with  co-occurrence and  then for $n-n'=\pm1$, with `$\rightarrow$' denoting the forward motion along the respective Mersenne fluctuation, are as follows:
\\

\noindent COMEFROM:  \hspace*{1mm} 2 \hspace*{5mm}if $\between_{\alpha_{n-1}}^{(n-1)}$ $\rightarrow$ $\between_{\alpha}^{(n)}$ and $\between_{\beta_{n-1}}^{(n-1)}$
$\rightarrow$ $\between_{\beta}^{(n)}$ both describe a left-leg ascent or a right-leg descent,

\noindent \hspace*{15mm}4 \hspace*{5mm}if   $\between_{\alpha_{n-1}}^{(n-1)}$ $\rightarrow$ $\between_{\alpha}^{(n)}$
is a left-leg ascent and  $\between_{\beta_{n-1)}}^{(n-1)}$
$\rightarrow$ $\between_{\beta}^{(n)}$ a right-leg descent
or vice versa; 

\noindent GOTO: 
 \hspace*{1mm} 2 \hspace*{5mm}if $\between_{\alpha}^{(n)}$ $\rightarrow$ $\between_{\alpha_{n+1}}^{(n+1)}$ and $\between_{\beta}^{(n)}$
$\rightarrow$ $\between_{\beta_{n+1}}^{(n+1)}$ both describe a left-leg ascent or a right-leg descent,

\noindent \hspace*{15mm}4 \hspace*{4mm} if there's ascent/descent opposition;

\noindent $n-n'=1$: 

\noindent \hspace*{14mm} 2\hspace*{5mm} if $\between_{\alpha_{n-1}}^{(n-1)}$ $\rightarrow$ $\between_{\alpha}^{(n)}$ and $\between_{\beta'}^{(n')}$
$\rightarrow$ $\between_{\beta_{n'+1}}^{(n'+1)}$ both describe a left-leg ascent or a right-leg descent,

\noindent \hspace*{15mm}4\hspace*{5mm}if there's ascent/descent opposition;

\noindent $n-n'=-1$: 

\noindent \hspace*{14mm} 2\hspace*{5mm} if $\between_{\alpha}^{(n)}$ $\rightarrow$ $\between_{\alpha_{n+1}}^{(n+1)}$ and $\between_{\beta_{n'-1}}^{(n'-1)}$
$\rightarrow$ $\between_{\beta'}^{(n')}$ both describe a left-leg ascent  or a right-leg descent,

\noindent \hspace*{15mm}4\hspace*{5mm} if there's  ascent/descent opposition.

\newpage
\noindent Let us begin with a simple example%
\footnote{\ free online tools do not do such demanding calculations in general;
it requires the WolframAlpha app to get the present terms computed%
} 
 $-$ co-occurrence, with   $f$   just the identity and $\between$'s  presumed to be correlated such that $\Delta_{20}=100$
and $x_{20}=L_{20}^{-}$: 
\[
\begin{array}{c}
\between_{\,566}^{(9736)}=17299,\qquad\between_{\,1251}^{(9736)}=100.\end{array}
\]
Tracking  the relevant Mersenne fluctuation segments, $({\scriptstyle{\ldots}},8649,\boldsymbol{17299},8649,\scriptstyle{\ldots})$ and $({\scriptstyle{\ldots}},50,100,200,\scriptstyle{\ldots})$, one registers two distinct correlative situations. By COMEFROM, `twice   left-leg ascent' and by GOTO, a `right-leg descent' vs. a  `left-leg ascent,'   indicating two terms being used per quaternion component in the quaternion sum in the former setting and four such terms in the latter:
\vspace{6mm}

\noindent\hspace*{0.5mm}%
\begin{tabular}{|c|c|c|c|c|c|}
\hline 
$s$  & $L_{20-s}$  & $(\lambda_{\textrm{max}})_{20-s}$  & $a_{2,s}$  & $a_{3,s}$  & $\Delta_{20}\;\textrm{  (at }n=9736)$\tabularnewline
\hline 
\hline 
$\vdots$  &  &  &  &  & \tabularnewline
\hline 
9  &438  & 874  & $281$  & $337$  & \tabularnewline
\hline 
10  & 336  & 670  & $-330$  & $-66$  &$\!$ COMEFROM:\hspace*{25mm}\tabularnewline
\hline 
11  & 272  & 542  & $265$  & $-63$  &$\!\!\!\!\!\!\!\!\!100=63+37$ \tabularnewline
\hline 
12  & 240  & 478  & $-163$  & 177  & \tabularnewline
\hline 
13  & 126  & 250  & 119  & $-44$  & \tabularnewline
\hline 
14  & 72  & 142  & $-67$  & $-28$  & GOTO:\hspace*{35mm}\tabularnewline
\hline 
15  & 40  & 78  & 18  & $-37$  & $100=(-281+265+\!119-\!18)+$\tabularnewline
\hline 
16  & 24  & 46  & $-14$  & $-20$  &  $(44-37+4+4)$\tabularnewline
\hline 
17  & 12  & 22  & $-12$  & $-4$  & \tabularnewline
\hline 
18  & 6  & 10  & $-5$  & $-4$  &$\{100=119-14-5\}$\tabularnewline
\hline 
\end{tabular}\bigskip{}

\noindent True to interpreting them as `lengths,' only terms of same
sign are considered physical $-$ unphysical decompositions, as indicated
by curly braces in the last row above, are dismissed. What remains are two  possible decompositions: 
\[\!\!\!\!\!\!\!\!\!\!\!\!\!\!\!\!\!\!\!\!\!\!\!\!\!\!\!\!\!\!\!\!\!\!\!\!\!\!\!\!\!\!\!\!\!\!\!\!\!\!\!\!\textrm{COMEFROM:} \qquad\qquad\Delta_{20} =\textrm{\textrm{Pre}(\ensuremath{\mathfrak{q}})}=f(\textrm{Re}(\mathfrak{q}_{11}+\mathfrak{q}_{15}))=
\]
\[ \quad\boldsymbol{k}a_{3,11}\boldsymbol{k}+\boldsymbol{k}a_{3,15}\boldsymbol{k}=63+37=100\quad\textrm{ with } u_{11}=u_{15}=\boldsymbol{k};
\]
\[\!\!\!\!\!\!\!\!\!\!\!\!\!\!\!\!\!\!\!\!\!\!\!\!\textrm{GOTO:}\qquad\Delta_{20}=\textrm{\textrm{Pre}(\ensuremath{\mathfrak{q}})}=f(\textrm{Re}(\mathfrak{q}_{9}+\mathfrak{q}_{11}+\mathfrak{q}_{13}+\mathfrak{q}_{15}+\mathfrak{q}_{17}+\mathfrak{q}_{18}))=
\]
\[(\boldsymbol{j}a_{2,9}\boldsymbol{j}+(-\boldsymbol{j})a_{2,11}\boldsymbol{j}+(-\boldsymbol{j})a_{2,13}\boldsymbol{j}+\boldsymbol{j}a_{2,15}\boldsymbol{j})+
\]
\[(\boldsymbol{k}a_{3,13}\boldsymbol{k}+(-\boldsymbol{k})a_{3,15}\boldsymbol{k}+\boldsymbol{k}a_{3,17}\boldsymbol{k}+\boldsymbol{k}a_{3,18}\boldsymbol{k})=
\]
\[(-281+265+119-18)+(44-37+4+4)=100
\]
\[\textrm{ with }u_9=-u_{11}=\boldsymbol{j}, u_{17}=u_{18}=\boldsymbol{k}, u_{13}=-u_{15}=-\boldsymbol{j}+\boldsymbol{k}. 
\]
No terms $a_{0,s},a_{1,s}$ are taken into account; they are set to zero  because even-numbered
$\Delta_{m}$'s are interpreted as `lengths of stretched chunks of
space' which ask  for decomposition in terms of $-$ non-vanishing $-$ $a_{2,s},a_{3,s}$.

\noindent Another even-numbered candidate is $\Delta_{20}=102$, as found in the  ($n-n'=1$) example 
\[
\begin{array}{c}
\between_{\,566}^{(9736)}=17299,\qquad\between_{\,563}^{(9735)}=102.\end{array}
\]
Here, the correlation is rather obvious by the proximity of the respective space-like refinements
 $-$ 566 vs. 563 $-$  and the `entanglement' becomes $x_{20}=L_{20}^{+}$.
From the perspective of $\between_{\,566}^{(9736)}=17299$, $\between_{\,563}^{(9735)}=102$
is a `stored' item which has to be installed in the present by a 
recall as $\Delta_{20}=102$. Which $f$ has to be chosen is  already  determined by the Mersenne fluctuation segment $({\scriptstyle{\ldots}},8649,\boldsymbol{17299})$, namely   $f(x)=2x+1+\epsilon$, and how
many terms of same sign per quaternion component in the quaternion sum shall be used depends
on the new Mersenne fluctuation segment in question. That being
 $\!\!$ $({\scriptstyle{\ldots}},102,\boldsymbol{206})$,
the correlative situation is `twice left-leg ascent,' asking for
 two terms per quaternion component in the quaternion sum:
\bigskip{}

\hspace*{1.5mm}%
\begin{tabular}{|c|c|c|c|c|c|}
\hline 
$s$  & $L_{20-s}$  & $(\lambda_{\textrm{max}})_{20-s}$  & $a_{2,s}$  & $a_{3,s}$  & $\Delta_{20}\;\textrm{ (at }n'=9735)$\tabularnewline
\hline 
\hline 
$\vdots$  &  &  &  &  & \tabularnewline
\hline 
10  & 336  & 670  & $-331$  & $-63$  & \tabularnewline
\hline 
11  & 272  & 542  & 264  & $-66$  & \tabularnewline
\hline 
12  & 240  & 478  & $-160$  & 179  & \tabularnewline
\hline 
13  & 126  & 250  & 118  & $-47$  & \tabularnewline
\hline 
14  & 72  & 142  & $-68$  & $-25$  & \tabularnewline
\hline 
15  & 40  & 78  & 15  & $-38$  & \tabularnewline
\hline 
16  & 24  & 46  & $-17$  & $-18$  & $102=2(68-17)+1-1$\tabularnewline
\hline 
17  & 12  & 22  & $-12$  & $(0)$  & \tabularnewline
\hline 
18  & 6  & 10  & $-5$  & $(0)$  & \tabularnewline
\hline 
\end{tabular}\bigskip{}

The decomposition reads in extenso: 
\[\Delta_{20}=\textrm{\textrm{Pre}(\ensuremath{\mathfrak{q}})}=
\]
\[f(\textrm{Re}(\mathfrak{q}_{14}+\mathfrak{q}_{16}))=
\]
\[2(\boldsymbol{j}a_{2,14}\boldsymbol{j}+(-\boldsymbol{j)}a_{2,16}\boldsymbol{j})+1+\epsilon=
\]
\[2(68-17)+1+\epsilon=102
\]
\[\textrm{where }\epsilon=-1 \textrm{ and } u_{14}=-u_{16}=\boldsymbol{j}.
\]
\noindent (Had there  been a  $\between_{\,?}^{(9735)}=101$ correlated
with $\between_{\,566}^{(9736)}=17299$ (with `entanglement' $x_{m}=L_{20}$),
its oddness would have constituted a case of `length of a squashed
chunk of space,' with $a_{2,s},a_{3,s}$  set to zero, instead of $a_{0,s},a_{1,s}$, and an instalment as $\Delta'_{20}=101$ by a 
recall using, once again,
an $f$ of the form $f(x)=2x+1+\epsilon$;   moreover, a hypothetical left-leg ascent $(\ldots,101,202+\!1+\!\epsilon)$ would have  led  to two relevant terms $a_{0,s},a_{1,s}$  per  quaternion component in the quaternionic sum that decomposes $\Delta_{20}=\textrm{\textrm{Pre}(\ensuremath{\mathfrak{q}})}$,
whereas  a  hypothetical right-leg descent $(\ldots,101,50)$ would have asked for four such terms.)

\noindent \noindent A more interesting ($n-n'=1$) example is 
\[
\between_{\,718}^{(13883)}=195798,\qquad\between_{61}^{(13882)}=762,
\]
 \smallskip{}

\noindent
\begin{tabular}{|c|c|c|c|c|c|}
\hline 
$s$  & $L_{24-s}$  & $(\lambda_{\textrm{max}})_{24-s}$  & $a_{2,s}$  & $a_{3,s}$  & $\Delta_{24}\;\textrm{ (at }n'=13882)$\tabularnewline
\hline 
\hline 
$\vdots$  &  &  &  &  & \tabularnewline
\hline 
3  & 27720  & 55438  & $-71$  & 27720  & \tabularnewline
\hline 
$\vdots$  &  &  &  &  & \tabularnewline
\hline 
9  & 2340  & 4678  & 2292  & -475  & \tabularnewline
\hline 
10  & 1422  & 2842  & $1065$  & $-944$  & \tabularnewline
\hline 
11  & 918  & 1834  & $-833$  & $-387$  & \tabularnewline
\hline 
12  & 756  & 1510  & 264  & 709  & \tabularnewline
\hline 
13  & 438  & 874  & 326  & $-294$  & \tabularnewline
\hline 
14  & 336  & 670  & $-65$  & $-330$  & \tabularnewline
\hline 
15  & 272  & 542  & $-207$  & $-177$  & \tabularnewline
\hline 
16  & 240  & 478  & 232  & 63  & \tabularnewline
\hline 
17  & 126  & 250  & $-125$  & $-23$  & $762=2((833\!-\!65+207\!-\!125)+$\tabularnewline
\hline 
18  & 72  & 142  & 5  & $-72$  & $(-387\!-\!177+23+72))+1\!-\!1$\tabularnewline
\hline 
19  & 40  & 78  & 40  & $-7$  & \tabularnewline
\hline 
20  & 24  & 46  & 5  & $-24$  & \tabularnewline
 &  &  &  &  & \tabularnewline
\hline 
\end{tabular}\bigskip{}

\noindent in that  the correlation, here $\Delta_{24}=\!762$, $x_{24}=\!L_{24}=\!196560$,
 may be attributed to the factorizations $195798=2\times3\times32633$,
and $762=2\times3\times\!127$  which share the prefix $2\times3$.%
\footnote{\, A second aspect is that when the correlation is written in the form
$L_{2}+L_{12}+195798=L_{24}$, the index on the RHS becomes the product
of the indexes on the LHS.%
}  %
$\!\!\!^,$\footnote{\, There exists an alternative choice yet: for $n'=13888$, we have
$\between_{699}^{(13888)}=762$, too. Note that the factorization
$n'=13888=2^{6}\times7\times31$, possibly coincidentally, completes
the quartet of the first four Mersenne primes 3,7,31,127; and the relative proximity of $\alpha$ and $\beta'$, 718 vs. 699,  makes a correlation even likelier. Unfortunately, though, decomposition over a span $n-n'=-5$    is  out of reach with the present methodology.}%
 $\;$
The relevant Mersenne fluctuation segments being
 $({\scriptstyle \ldots},97899,$ $\boldsymbol{195798})$
 vs. $\!\!({\scriptstyle \ldots},\boldsymbol{762},381)$, 
we are dealing with two peaks, but peak 762 precedes peak 195798: from the perspective
of $\between_{\,718}^{(13883)}=195798$, $\between_{\,61}^{(13882)}=762$ is a `stored' item requiring instalment
in the present by a  recall as $\Delta_{24}=762$, using
$f(x)=2x+1+\epsilon$ as well as four terms per quaternion component in the quaternion sum since the correlative situation is `left-leg ascent' vs. $\!\!$ `right-leg descent.'
 Hence, we decompose: $\Delta_{24}=\textrm{\textrm{Pre}(\ensuremath{\mathfrak{q}})}=f(\textrm{Re}(\mathfrak{q}_{11}+\mathfrak{q}_{14}+\mathfrak{q}_{15}+\mathfrak{q}_{17}+\mathfrak{q}_{18}))=2((\boldsymbol{j}a_{2,11}\boldsymbol{j}+(-\boldsymbol{j})a_{2,14}\boldsymbol{j}+\boldsymbol{j}a_{2,15}\boldsymbol{j}+(-\boldsymbol{j})a_{2,17}\boldsymbol{j})+((-\boldsymbol{k})a_{3,11}\boldsymbol{k}+(-\boldsymbol{k)}a_{3,15}\boldsymbol{k}+(-\boldsymbol{k})a_{3,17}\boldsymbol{k}+(-\boldsymbol{k})a_{3,18}\boldsymbol{k}))+1+\epsilon=2((833-65-207-125)+(-387-177+23+72))+1+\epsilon=762$
where $\epsilon=-1$, $u_{14}=\boldsymbol{j}$, $u_{18}=\boldsymbol{k}$, $u_{11}=u_{15}=-u_{17}=\boldsymbol{j}-\boldsymbol{k}$. \bigskip{}

\noindent What has been missing, thus far, is the counterpart to 
recall $-$ anticipation.

\noindent A case in point is the `overshoot' example 
\[
\between_{\,761}^{(14079)}=5534,\quad\between_{\,586}^{(14080)}=188:
\]

\vspace*{\smallskipamount}
 \hspace{2mm} %
\begin{tabular}{|c|c|c|c|c|c|}
\hline 
$s$  & $L_{17-s}$  & $(\lambda_{\textrm{max}})_{17-s}$  & $a_{2,s}$  & $a_{3,s}$  & $\Delta_{17}\;\textrm{ (at }n'=14080)$\tabularnewline
\hline 
\hline 
$\vdots$  &  &  &  &  & \tabularnewline
\hline 
5  & 756  & 1510  & $-342$  & 675  & \tabularnewline
\hline 
6  &438  &874  & $338$  & 279  & \tabularnewline
\hline 
7  & 336  & 670  & 335  & 32  & \tabularnewline
\hline 
8  & 272  & 542  & 270  & $-38$  & \tabularnewline
\hline 
10  & 240  & 478  & $-231$  & $66$  & \tabularnewline
\hline 
11  & 126  & 250  & $-54$  & 115  &  \tabularnewline
\hline 
12  & 72  & 142  & $41$  & $60$  &$188=\lfloor((342-231)+$ \tabularnewline
\hline 
13  & 40  & 78  & $-40$  & $-4$  & \tabularnewline
\hline 
14  & 24  & 46  & $21$  & 13  & $(279-13))/2\rfloor$ \tabularnewline
\hline 
15  & 12  & 22  & 12  & (0)  & \tabularnewline
\hline 
16  & 6  & 10  & 6  & (0)  & \tabularnewline
\hline 
\end{tabular}

\vspace*{5mm}
\noindent From the perspective of $\between_{\,586}^{(14079)}=5534$,
$\between_{\,586}^{(14080)}=188$ is a future item (with `entanglement'
$x_{17}=L_{17}=5346$) which has to be installed in the present by anti\-cipation
as $\Delta_{17}=188$. Even though the tracking of the  Mersenne fluctuation for 5534 is incomplete
due to computational limits, $({\scriptstyle \ldots},5534,{\scriptstyle \ldots},5533,{\scriptstyle \ldots})$,
one recognizes a left-leg ascent, hence $f(x)=\lfloor x/2 \rfloor-\bar{\delta}$. Also that, with
 $({\scriptstyle \ldots},93,188,{\scriptstyle \ldots})$, the correlative situation  is `twice left-leg' and implies using two terms per quaternion component in the quaternion sum,
resulting in the decomposition $\Delta_{17}=\textrm{Pre}(\ensuremath{\mathfrak{q}})=f(\textrm{Re}(\mathfrak{q}_{5}+\mathfrak{q}_{6}+\mathfrak{q}_{10}+\mathfrak{q}_{14}))=\lfloor((\boldsymbol{j}a_{2,5}\boldsymbol{j}+(-\boldsymbol{j)}a_{2,10}\boldsymbol{j})+((-\boldsymbol{k)}a_{3,6}\boldsymbol{k}+\boldsymbol{k}a_{3,14}\boldsymbol{k}))/2\rfloor-\bar{\delta}=\lfloor((342-231)+(279-13))/2\rfloor-\bar{\delta}=188$
where $\bar{\delta}=0,u_{5}=-u_{10}=\boldsymbol{j},-u_6=u_{14}=\boldsymbol{k}$.

\bigskip{}
\noindent Does the observed even-(odd-)ness of $\,\alpha\,$ vs. $\!\!$ odd-(even-)ness
of $\beta$ (or $\beta'$) depend on whether \\
 a) $\Delta_{m}$ represents the `length of a stretched space chunk' or the `length
of a squashed space chunk,' or whether\\
 b) $\Delta_{m}$ has to be added or subtracted? 
 (Our example
\[
\between_{\,566}^{(9736)}=17299,\quad\between_{\,1251}^{(9736)}=100,
\]
after all, has a counterpart in
\[
\between_{\,323}^{(14108)}=17704,\quad\between_{\,558}^{(14109)}=303.)
\]

\noindent The following example discards possibility a). It consists
of the four terms 
\[
\between_{\,761}^{(14079)}=5534,\between_{\,214}^{(14079)}=187,\between_{\,300}^{(14079)}=189,\between_{\,586}^{(14080)}=188
\]
which, under $L_{17}=5346$, allow for the choices 
\[
\between_{\,761}^{(14079)}=5534=L_{17}^{+}+187=L_{17}^{-}+189=L_{17}+188.
\]
And a further example, 
\[
\between_{\,21}^{(14125)}=4041,\qquad\between_{\,240}^{(14126)}=280,
\]
which allows for $x_{16}+\Delta_{16}=L_{16}^{+}=4321$, discards b).
Thus, regardless of whether $\Delta_{m}$ stands for the `length
of a squashed chunk of space' or the `length of a stretched chunk of space' and/or is to be added
or subtracted, the underlying rule, as a necessary condition, seems to simply read: 
\[
\alpha=2a+\delta\Rightarrow\beta(\beta')=2b+\delta-1>1\qquad(\delta\in\{0,1\}).
\]
The latter example has been picked as it constitutes a case of anticipation  with descent/ascent opposition:
From $n=14125$ on, the relevant Mersenne
fluctuation segments are $(4041,2020,{\scriptstyle \ldots})$
vs. $\!\!(140,280,{\scriptstyle \ldots})$ so that,
this time, the instalment $-$ of a future item
$\between_{\,240}^{(14126)}=280$ in the present as $\Delta_{16}$ $-$ works with $f(x)=2x+1+\epsilon$ and {\em  four}
terms per quaternion component in the quaternion sum. The  decomposition reads in full
\[
 \Delta_{16}=
\textrm{\textrm{Pre}(\ensuremath{\mathfrak{q}})}=f(\textrm{Re}(\mathfrak{q}_{3}+\mathfrak{q}_{4}+\mathfrak{q}_{5}+\mathfrak{q}_{6}+\mathfrak{q}_{7}+\mathfrak{q}_{9}+\mathfrak{q}_{10}+\mathfrak{q}_{13}))=
\]
\[
2(((-\boldsymbol{j)}a_{2,3}\boldsymbol{j}+\boldsymbol{j}a_{2,4}\boldsymbol{j}+\boldsymbol{j}a_{2,9}\boldsymbol{j}+(-\boldsymbol{j})a_{2,10}\boldsymbol{j})+
\]
\[
(\boldsymbol{k}a_{3,5}\boldsymbol{k}+(-\boldsymbol{k})a_{3,6}\boldsymbol{k}+
(-\boldsymbol{k})a_{3,7}\boldsymbol{k}+\boldsymbol{k}a_{3,13}\boldsymbol{k}))+1+\epsilon=
\]
\[
2((-272+464+126-72)+(-374+169+105-7))+1+\epsilon=280
\]
\[ \textrm{ where }\epsilon=1,-u_{3}=u_{4}=u_{9}=-u_{10}=\boldsymbol{j}, u_{5}=-u_{6}=-u_{7}=u_{13}=\boldsymbol{k} :
\]
\noindent \vspace{6mm}

\noindent %
\noindent\hspace*{0.5mm}
\begin{tabular}{|c|c|c|c|c|c|}
\hline 
$s$  & $L_{16-s}$  & $(\lambda_{\textrm{max}})_{16-s}$  & $a_{2,s}$  & $a_{3,s}$  & $\Delta_{16}\;\textrm{ (at }n'=14126)$\tabularnewline
\hline 
\hline 
$\vdots$  &  &  &  &  & \tabularnewline
\hline 
3  & 918  & 1834  & $-272$  & $-878$  & $280=\quad\quad\quad\quad\qquad\qquad$ \tabularnewline
\hline 
4  & 756  & 1510  & $-464$  & 598  &  $2((-272+464+\!126\!-\!72)+$\tabularnewline
\hline 
5  & 438  & 874  & 229  & 374  & $(-374+169+\!105\!-\!7))+\!1+\!1$ \tabularnewline
\hline 
6  & 336  & 670  & 291  & 169  & \tabularnewline
\hline 
7  & 272  & 542  & $252$  & 105  & \tabularnewline
\hline 
8  & 240  & 478  & $-228$  & $-78$  & \tabularnewline
\hline 
9  & 126  & 250  & $-126$  & -4  &\tabularnewline
\hline 
10  & 72  & 142  & $-72$  & $10$  & \tabularnewline
\hline 
11  & 40  & 78  & $32$  & $25$  & \tabularnewline
\hline 
12  & 24  & 46  & $21$  & 13  &\tabularnewline
\hline 
13  & 12  & 22  & 11  & 7  &   \tabularnewline
\hline 
14  & 6  & 10  & $-5$  & -4  & \tabularnewline
\hline 
\end{tabular}

\newpage
\noindent For $n=14079$ (see above), there are two co-occurrent squashed cases, $\Delta^{\textrm{\tiny nd}}_{17}=187$, $\Delta^{\textrm{\tiny d}}_{17}=189$, so we take
the opportunity to  decompose them using
 $a_{0,s},a_{1,s}$:

\vspace*{4mm}
\noindent\hspace*{1.5mm}
\begin{tabular}{|c|c|c|c|c|c|}
\hline 
$s$  & $L_{17-s}$  & $(\lambda_{\textrm{min}})_{17-s}$  & $\!a_{0,s}$  & $\!a_{1,s}$  & $\Delta^{\textrm{\tiny nd}}_{17},\Delta^{\textrm{\tiny d}}_{17}\;\textrm{ (at }n=14079)$\tabularnewline
\hline 
\hline 
$\vdots$  &  &  &  &  & $\!\!\!\!\!\!\!\!\!\!\!\!\!\!\!\!\!\!\!$COMEFROM:\tabularnewline
\hline 
2  & 2340  & 4675  & $2333$  & $169$ &$\!\!\!\!\!\!\!\!\!\!\!\!\!\!\!\!\!\!\!\!\!\!\!\!\!\!\!\!\!\!\!\!\!\!\!\!\!\!\!\!\!\!187=$ \tabularnewline
\hline 
3  & 1422  & 2839  & $1375$  & $-361$  &$(-7+20)+(169+5)$ \tabularnewline
\hline 
4  & 918  & 1831  & $-342$  & $-851$  & \tabularnewline
\hline 
5  & 756  & 1507  & $-414$  & 632  & $\!\!\!\!\!\!\!\!\!\!\!\!\!\!\!\!\!\!\!\!\!\!\!\!\!\!\!\!\!\!\!\!\!\!\!\!\!\!\!\!$GOTO:\tabularnewline
\hline 
6  & 438  & 871  & 224  & 375  & $\!\!\!\!\!\!\!\!\!\!187=\qquad\qquad\qquad$\tabularnewline
\hline 
7  & 336  & 667  & 261  & 210  &  $(-342+\!414\!-\!152+\!16)+$\tabularnewline
\hline 
8  & 272  & 539  & $197$  & 186  &$\!\!\!\!\!\!\!\!\!\!\!(632-375-11+5)$ \tabularnewline
\hline 
9  & 240  & 475  & $-152$  & $-184$  & \tabularnewline
\hline 
10  & 126  & 247  & 126  & (0)  &\tabularnewline
\hline 
11  & 72  & 139  & $-16$  & $69$  & COMEFROM \& GOTO: \tabularnewline
\hline 
12  & 40  & 75  & $-7$  & $-39$  &$189=152\!+16\!+20\!+1$ \tabularnewline
\hline 
13  & 24  & 43  & $-20$  & 11  & \tabularnewline
\hline 
14  & 12  & 19  & 12  & (0)  & \tabularnewline
\hline 
15  & 6  & 7  & $-1$  & 5  & \tabularnewline
 &  &  &  &  &\tabularnewline
 &  &  &  &  &\tabularnewline
\hline 
\end{tabular}

\vspace*{4mm}

\noindent Regardless of the poor tracking of one of the Mersenne
fluctuation segments  $-$ $(\scriptstyle{?}\textstyle{,5534,}\scriptstyle{?}\textstyle{,}\scriptstyle{?}\textstyle{,}\scriptstyle{?}\textstyle{,5534,}{\scriptstyle \ldots})$
 vs. $\!\!(93,187,93,{\scriptstyle \ldots})$ or $(378,189,94,{\scriptstyle \ldots})$ $-$,
 one can distinguish  three distinct correlative situations. For  $\Delta^{\textrm{\tiny nd}}_{17}\,=187$: `twice left-leg ascent' under COMEFROM and   `left-leg ascent' vs. $\!\!$`right-leg descent' under GOTO. For $\Delta^{\textrm{\tiny d}}_{17}\,=189\,$:  `left-leg ascent' vs. $\!\!$`right-leg descent'  under  both COMEFROM and GOTO. With   $f$  the identity, we thus get the decompositions:  \vspace*{1mm}

\noindent $\Delta^{\textrm{\tiny nd}}_{17}\,=\textrm{\textrm{Pre}(\ensuremath{\mathfrak{q}})}=f(\textrm{Re}(\mathfrak{q}_{2}+\mathfrak{q}_{12}+\mathfrak{q}_{13}+\mathfrak{q}_{15}))$
 $=$ $(a_{0,12}+(-1)a_{0,13})+((-\boldsymbol{i})a_{1,2}\boldsymbol{i}+(-\boldsymbol{i})a_{1,15}\boldsymbol{i})=(-7+20)+(169+5)=187$ for $u_2=u_{15}=-\boldsymbol{i},u_{12}=-u_{13}=1$; 

\noindent$\Delta^{\textrm{\tiny nd}}_{17}\,=\textrm{\textrm{Pre}(\ensuremath{\mathfrak{q}})}=f(\textrm{Re}(\mathfrak{q}_{4}+\mathfrak{q}_{5}+\mathfrak{q}_{6}+\mathfrak{q}_{9}+\mathfrak{q}_{11}+\mathfrak{q}_{13}+\mathfrak{q}_{15}))$ $=$ $(a_{0,4}+(-1)a_{0,5}+a_{0,9}+(-1)a_{0,11})+((-\boldsymbol{i})a_{1,5}\boldsymbol{i}+\boldsymbol{i}a_{1,6}\boldsymbol{i}+\boldsymbol{i}a_{1,13}\boldsymbol{i}+(-\boldsymbol{i})a_{1,15}\boldsymbol{i})=(-342+414-152+16)+(632-375-11+5)=187$ 
for $u_{4}=u_{9}=-u_{11}=1,u_{6}=u_{13}=-u_{15}=\boldsymbol{i},u_{5}=-1-\boldsymbol{i}$;

\noindent $\Delta^{\textrm{\tiny d}}_{17}$ $=$ $\textrm{\textrm{Pre}(\ensuremath{\mathfrak{q}})}\!=\!f(\textrm{Re}(\mathfrak{q}_{9}\!+\mathfrak{q}_{11}\!+\mathfrak{q}_{13}\!+\mathfrak{q}_{15}))\!=(-1)a_{0,9}+(-1)a_{0,11}+(-1)a_{0,13}+(-1)a_{0,15}=152+16+20+1=189$
for $u_{9}=u_{11}=u_{13}=u_{15}=-1$.
 \vspace*{1.5mm}

\noindent The fact that  5534 accepts only one  candidate as partner may seem to amount to a question of principle: nondeterminism `$=$   COMEFROM \& GOTO offering decompositions (of  187) that differ' $\!$ vs. $\!\!$ determinism `$=$ COMEFROM \& GOTO with  one and the same  decomposition (of 189).'  (Which in retrospect motivates the use of upper tags in $\Delta^{\textrm{\tiny nd}}_{17},\Delta^{\textrm{\tiny d}}_{17}$.) But, as a never endangered first option, `entangling' (here, $5534=L^+_{17}+187\,$ vs.   $5534=L^-_{17}+189,$ and in  our previous example, $5534=L_{17}+188$)  is the act that sets the course.

\subsection{Crotons of the second kind: radiation}
Let us now discuss $\gamma_{\alpha}^{(\nu)}$. Here, we will have
a d{\'e}j{\`a}-vu with $x$ vs.$\,(x+2)$ and $x$ vs.$\,(x+1)$
juxtaposition in connection with two important extensions of the regular
Mersenne numbers $M_{\textrm{reg}}:=p_{i}=1,3,7,{\scriptstyle \ldots}$.
The  first such extension was  the Mersenne fluctuations $-$ whereas regular Mersenne numbers are formed by $p_{i+1}= 2p_i+1$, the term to follow $\between_\alpha^{(n)}$ (on the left leg) is formed by  $2\between_\alpha^{(n)}+1+\epsilon$ (with $\epsilon$ fluctuating). The next extension is

\begin{equation}
M_{\nicefrac{5}{8}}:=o_{i}=4,9,19,\ldots,\label{eq:M5/8}
\end{equation}
induced by $\frac{5}{8}p_{i}$ in that each member of $\frac{5}{8}M_{\textrm{reg}}\:(i\geq3)$
exceeds its $M_{\nicefrac{5}{8}}\:(i\geq1)$ counterpart by $\frac{3}{8}$:
$\frac{5}{8}p_{3}=o_{1}+\frac{3}{8}$, $\frac{5}{8}p_{4}=o_{2}+\frac{3}{8},\ldots$
, which is a case of $i$ vs.$\,(i+2)$ juxtaposition in the index
of $o\,$ vs. $\!$the index of $p$.

\medskip{}
\noindent And a third extension is
\begin{equation}
M_{\nicefrac{9}{8}}:=\sigma_{i}=\frac{7}{2},8,17,\ldots\label{eq:M9/8}
\end{equation}
and is induced by $\frac{9}{8}p_{i}$ in that each member of $\frac{9}{8}M_{\textrm{reg}}\:(i\geq2)$
differs from its $M_{\nicefrac{9}{8}}\:(i\geq1)$ counterpart by $-\frac{1}{8}$:
$\frac{9}{8}p_{2}=\sigma_{1}\!-\!\frac{1}{8}$, $\frac{9}{8}p_{3}=\sigma_{2}\!-\!\frac{1}{8},\ldots$
, a case of $i$ vs.$\,(i+1)$ juxtaposition in the index of $\sigma$
vs. index of $p$. \\
These extensions have complements,
\begin{equation}
M_{\nicefrac{5}{8}}^{+}:=o_{i}+1\equiv o_{i}^{+}=5,10,20,\ldots,\label{eq:M5/8+1}
\end{equation}
and 
\begin{equation}
M_{\nicefrac{9}{8}}^{-}:=\sigma_{i}-1\equiv\sigma_{i}^{-}=\frac{5}{2},7,16,\ldots\,.\label{eq:M9/8-1}
\end{equation}
The new Mersenne numbers set the scene for quasi-supersymmetric `boson'
and `fermion' lists where juxtaposition once more takes effect: that
of $l\,$ vs. $l-2$ as applied to lower summation bounds, with $\frac{5}{8}$
and $\frac{9}{8}$ as factors common to the sum members. Replacing
$\omega$ in the energy equation 
\[
\mathcal{E}=\hbar\omega
\]
with sums $\frac{5}{8}\sum_{l=3}^{\ell_{i}}p_{l}$, we respectively
get a list of `fermionic' excitations

\bigskip{}

\noindent $\frac{5}{8}\hbar\sum_{l=3}^{6}p_{l}=\frac{145}{2}\hbar,\:$%
\framebox{\begin{minipage}[t][0.9\totalheight]{0.29\columnwidth}%
$\frac{5}{8}\hbar\sum_{l=3}^{14}p_{l}=\frac{40,935}{2}\hbar,\:$%
\end{minipage}}$\quad\frac{5}{8}\hbar\sum_{l=3}^{22}p_{l}=\frac{10,485,725}{2}\hbar,$\medskip{}

\noindent %
\framebox{\begin{minipage}[t][0.9\totalheight]{42mm}%
$\frac{5}{8}\hbar\!\sum_{n=3}^{30}p_{n}\!=\!\frac{2,684,354,515}{2}\hbar,\:$%
\end{minipage}}$\quad\ldots$

\bigskip{}

\noindent  and, with replacements $\omega\rightarrow\frac{5}{8}\sum_{l=3}^{\ell_{i}+4}p_{l}$,
`bosonic' ones 

\bigskip{}

\noindent  $\frac{5}{8}\hbar\sum_{l=3}^{10}p_{l}=1,\!270\,\hbar,\:\frac{5}{8}\hbar\sum_{l=3}^{18}p_{n}=327,\!665\,\hbar,\:\frac{5}{8}\hbar\sum_{l=3}^{26}p_{l}=83,\!886,\!060\,\hbar,$\medskip{}

\noindent $\frac{5}{8}\hbar\sum_{l=3}^{34}p_{l}=21,\!474,\!836,\!455\,\hbar,\:\ldots$\bigskip{}

\noindent  which, with replacements $\omega\rightarrow\frac{5}{8}\sum_{l=1}^{\ell_{i}}p_{l}$,
find their juxtaposition in `bosonic' excitations

\bigskip{}

\noindent $\frac{5}{8}\hbar\sum_{l=1}^{6}p_{l}=75\hbar,\quad$%
\framebox{\begin{minipage}[t][0.9\totalheight]{38mm}%
$\frac{5}{8}\hbar\sum_{l=1}^{14}p_{l}=20,\!470\,\hbar,$%
\end{minipage}} $\:\frac{5}{8}\hbar\sum_{l=1}^{22}p_{l}=5,\!242,\!865\,\hbar,$\medskip{}

\noindent  %
\framebox{\begin{minipage}[t][0.9\totalheight]{50mm}%
$\frac{5}{8}\hbar\sum_{l=1}^{30}p_{l}=1,\!342,\!177,\!260\,\hbar,$%
\end{minipage}}$\quad\ldots$ \bigskip{}

\noindent and, with replacements $\omega\rightarrow\frac{5}{8}\sum_{l=1}^{\ell_{i}-4}p_{l}$,
`fermionic' ones \bigskip{}

\noindent  $\frac{5}{8}\hbar\sum_{l=1}^{2}p_{l}=\frac{5}{2}\hbar,\:\frac{5}{8}\hbar\sum_{l=1}^{10}p_{l}=\frac{2,545}{2}\hbar,\:\frac{5}{8}\hbar\sum_{l=1}^{18}p_{l}=\frac{655,455}{2}\hbar,$\medskip{}

\noindent  $\frac{5}{8}\hbar\sum_{l=1}^{26}p_{l}=\frac{167,772,125}{2}\hbar,\:\ldots$
.\bigskip{}

\noindent The upper summation bounds $\ell_{i}=6,14,\ldots$ 
obviously mark instances of quasi-supersymmetry (framed terms). They show up 
once again in juxtapositions with replacements $\;\omega\rightarrow\frac{9}{8}\sum_{l=3}^{\ell_{i}}p_{l}$,$\quad\omega\rightarrow\frac{9}{8}\sum_{l=1}^{\ell_{i}\pm4}p_{l}\;$,
respectively $-$ \medskip{}

\noindent `fermion' list:\bigskip{}

\noindent  $\frac{9}{8}\hbar\sum_{l=3}^{6}p_{l}=\frac{261}{2}\,\hbar,\quad$%
\framebox{\begin{minipage}[t][0.9\totalheight]{37mm}%
$\frac{9}{8}\hbar\sum_{l=3}^{14}p_{l}=\frac{73,683}{2}\,\hbar,\:$%
\end{minipage}}$\quad\frac{9}{8}\hbar\sum_{l=3}^{22}p_{l}=\frac{18,874,305}{2}\,\hbar,\:\ldots;$\bigskip{}

\noindent `boson' list: \bigskip{}

\noindent  $\frac{9}{8}\hbar\sum_{l=3}^{10}p_{l}=2,\!286\,\hbar,\:\frac{9}{8}\hbar\sum_{l=3}^{18}p_{l}=589,\!797\,\hbar,\:\frac{9}{8}\hbar\sum_{l=3}^{26}p_{l}=150,\!994,\!908\,\hbar,\:\ldots\,.$\bigskip{}

\noindent and juxtaposed `boson' list: \bigskip{}

\noindent $\frac{9}{8}\hbar\sum_{l=1}^{6}p_{l}=135\,\hbar,\quad$
\framebox{\begin{minipage}[t][0.9\totalheight]{38mm}%
$\frac{9}{8}\hbar\sum_{l=1}^{14}p_{l}=36,\!846\,\hbar,\:$%
\end{minipage}}$\quad\frac{9}{8}\hbar\sum_{l=1}^{22}p_{l}=9,\!437,\!157\,\hbar,\:\ldots;$\bigskip{}

\noindent juxtaposed `fermion' list:\bigskip{}

\noindent  $\frac{9}{8}\hbar\sum_{l=1}^{2}p_{l}=\frac{9}{2}\,\hbar,\:\frac{9}{8}\hbar\sum_{l=1}^{10}p_{l}=\frac{4,581}{2}\,\hbar,\:\frac{9}{8}\hbar\sum_{l=1}^{18}p_{l}=\frac{1,179,603}{2}\,\hbar,\:\ldots$
.\bigskip{}

\noindent The upper summation bounds $\ell_{i}=6,14,\ldots$ correspond to `length' differences of $\between$-type Mersenne fluctuations, one homgeneous in $\epsilon=1$, the other homogeneous in  $\epsilon=-1$.  While for the former the  $\between$'s proceed like $1,4,10,22,{\scriptstyle\ldots}$, for the latter they proceed like $1,2,4,8,{\scriptstyle\ldots}$, so that differences  $\ell_{i}$ will form. Therefore interaction between particles  switching from a Mersenne fluctuation of given  $\epsilon$-homogeneity to an alloqphyletic  one of opposite  $\epsilon$-homogeneity is involved $-$ generally of different space-like and time-like refinements. The  radiation exchanged (photons in the $\frac{5}{8}$ case, massive bosons in the $\frac{9}{8}$ case)  can
accordingly be expressed in terms of  $M_{\nicefrac{5}{8}},M_{\nicefrac{5}{8}}^{+}$
and $\left\lceil M_{\nicefrac{9}{8}}\right\rceil,\left\lceil M_{\nicefrac{9}{8}}^{-}\right\rceil $.
We state without proof that every natural number $>2$ except 6, depending
on the situation and allowing for empty sums, comes with one or with
either decomposition
\begin{equation}
\sum_{i}o_{i}\!+\!\sum_{j}o_{j}^{+},\quad\sum_{k}\left\lceil \sigma_{k}\right\rceil \!+\!\sum_{l}\left\lceil \sigma_{l}^{-}\right\rceil .\label{eq:decomp-5/8or9/8}
\end{equation}
It turns out the juxtaposed `bosons' behave complementarily in that
the set that corresponds to $\frac{5}{8}\sum_{l=1}^{\ell_{i}}p_{l}$
obeys the rule 
\begin{equation}
\begin{array}{lr}
\#\,\,\textrm{of}\,\textrm{\,\ distinct}\: o_{k}: & 5(\frac{\ell_{i}+2}{8}-1),\\
\#\,\textrm{\,\ of}\,\textrm{\,\ distinct}\: o_{k}^{+}: & 4+3(\frac{\ell_{i}+2}{8}-1),
\end{array}\label{eq:comp5/8}
\end{equation}
while the set that corresponds to $\frac{9}{8}\sum_{l=1}^{\ell_{i}}p_{l}$
obeys

\begin{equation}
\begin{array}{lr}
\#\,\textrm{\,\ of}\,\textrm{\,\ distinct}\:\left\lceil \sigma_{k}\right\rceil : & 5+7(\frac{\ell_{i}+2}{8}-1),\\
\#\,\textrm{\,\ of}\,\textrm{\,\ distinct}\:\left\lceil \sigma_{k}^{-}\right\rceil : & 1(\frac{\ell_{i}+2}{8}-1).
\end{array}\label{eq:comp9/8}
\end{equation}
Samples of the first kind are: $\frac{5}{8}\sum_{l=1}^{6}p_{l}\cong5+10+20+40=75$,
\\
$\frac{5}{8}\sum_{l=1}^{14}p_{l}\cong(5\!+\!10\!+\!20+\ldots+320)+(639+\!1279+\ldots+\!10239)=20470$,
etc.\\
And samples of the second kind are: $\frac{9}{8}\sum_{l=1}^{6}p_{l}\cong4+8+17+35+71=135$,
\\
$\frac{9}{8}\sum_{l=1}^{14}p_{l}\cong(4+8+17+\ldots+2303+4607+9215)+18430=36846$,
etc.\\
\vspace{-3mm}

\noindent  Of course, the results remain unchanged for pairwise role
exchange, such as $(9215\mapsto9214,18430\mapsto18431)$ etc. \\
The above juxtapostion should be familiar to readers of \cite{Merkel}.
In Sect.4.2 of that reference, it was shown that Catalan numbers find
a natural place in the framework of the secondary trace structure
of $\textrm{LL}(G_{\alpha\beta}^{(p_{i})})$: The total `secondary
traces,' secondary trace proper and adjacents, of $\textrm{LL}(G_{\alpha\beta}^{(p_{i})})$
determine the Catalan numbers $C_{q_{i}},C_{q_{i}+1},{\scriptstyle \ldots},C_{2q_{i}}$,
which by definition can be written as $C_{p_{i-2}},\! C_{p_{i-2}+1},{\scriptstyle \ldots},\! C_{2p_{i-2}}$
via $q_{i}=(p_{i}\!-\!3)/4$. We note in passing that for $i=4,5$,
the major building blocks of $\textrm{LL}(G_{\alpha\beta}^{(p_{i})})$
modulo 8 are $\left(_{3\;\:5}^{5\;\:3}\right)$, in nice correspondence
with the rule expressed in (\ref{eq:comp5/8}). Blocks $\left(_{1\;\:7}^{7\;\:1}\right)$,
complying to rule (\ref{eq:comp9/8}), do not emerge until at $i=6$,
$\left(\begin{array}{c}
\\
\\
\end{array}\right.$\hspace{-5mm}%
\begin{tabular}{cc}
${\scriptstyle G_{17,5}^{(63)}(=58791)}$ & ${\scriptstyle G_{17,6}^{(63)}(=18633)}$\tabularnewline
${\scriptstyle G_{18,5}^{(63)}(=189393)}$ & ${\scriptstyle G_{18,6}^{(63)}(=58791)}$\tabularnewline
\end{tabular}\hspace{-5mm}$\left.\begin{array}{c}
\\
\\
\end{array}\right)$$\equiv\left(\begin{array}{cc}
{\scriptstyle 7} & {\scriptstyle 1}\\
{\scriptstyle 1} & {\scriptstyle 7}
\end{array}\right)\,\left(\textrm{mod}\,8\right)$ being a party concerned, and oscillations $\left(_{1\;\:7}^{7\;\:1}\right)\leftrightarrow\left(_{3\;\:5}^{5\;\:3}\right)$
do not occur until from then on. Thus, rule (\ref{eq:comp9/8}) would
underpin a physical meaning of quasi-supersymmetry beginning with
$\log_{2}(p_{6}+1)$, whereas rule (\ref{eq:comp5/8}) would do this
beginning with $\log_{2}(p_{4}+1)$. \\
Our next example of juxtapostion $x\,$ vs. $x-2$ directly affects
the CFRs $\gamma_{\alpha}^{(n)}$. 

\noindent Let $a,b\in\mathbb{Z}$, $c,m,n\in\mathbb{N}$ and the
greatest lower bound $\textrm{GLP}(n)$ and least upper bound $\textrm{LUP}(n)$
of $3^{2^{n}}$ be respectively given by its adjacent prime numbers
\[
\begin{array}{c}
\textrm{GLP}(n)=3^{2^{n}}-\vartheta_{\textrm{glp}}(n),\\
\textrm{LUP}(n)=3^{2^{n}}+\vartheta_{\textrm{lup}}(n).
\end{array}
\]
\noindent Then we observe relationships $a+bm=c\:$ with 
\[
\begin{array}{c}
\left|a\right|,\left|b\right|\in\left\{ \vartheta_{\textrm{glp}}(n)\right\} \cup\left\{ \vartheta_{\textrm{lup}}(n)\right\} \:\left(\left|a\right|\neq\left|b\right|\right),\\
c\in\left\{ T(n)\right\} \cup\left\{ (\xi+\zeta)(n)\right\} ,\\
m\in M_{\nicefrac{5}{8}}\cup M_{\nicefrac{5}{8}}^{+},
\end{array}
\]
where $\: T(n)$ denotes number of basis elements of $(G_{\rho}^{(p)})$
and $\:(\xi+\zeta)(n)$ a sum of characteristic numbers (whose elucidation
has to wait until Sect. \ref{sec:Supplementary-remarks-and}): 
\[
\begin{array}{c}
T(n)=A(n)\frac{\partial(3^{2^{x}})}{\partial x}\left|_{x=n}\right.,\;(\xi+\zeta)(n)=\Omega(n)\frac{\partial(3^{2^{x}})}{\partial x}\left|_{x=n}\right.\\
\\
\left(A(n)=\frac{2^{-n+1}}{\log(2)}\!\cdot\!\frac{3^{-2^{n}+n-3}}{\log(3)},\;\Omega(n)=\frac{2^{n-6}+2^{-3}-2^{-n+1+(n-5)\log(5)/\log(2)}}{\log(2)}\!\cdot\!\frac{3^{-2^{n}}}{\log(3)}\right).
\end{array}
\]

\noindent With growing $n$, the reciprocal actions of $\frac{\partial(3^{2^{x}})}{\partial x}\left|_{x=n}\right.$
outdistancing $3^{2^{n}}$ and $A(n)$, $\Omega(n)$ turning infinitesimal
conspire to make recording of $T(n)$ and $(\xi+\zeta)$ $(n)$ in
the form $a+bm$ possible via encounter of $3^{2^{n}}$with distances
$\vartheta_{\textrm{glp}}(n)$, $\vartheta_{\textrm{lup}}(n)$ which
are mutually exclusively used as single units $\left|a\right|$ or
units $\left|b\right|\,$ in bins $m=4,5,9,10,{\scriptstyle \ldots}\,$.
To get this manifested, we recur to juxtaposition $n\,$ vs. $n-2$
in the pair of tables \bigskip{}

\qquad{}\qquad{}\qquad{}%
\begin{tabular}{|c|}
\hline 
$n$\tabularnewline
\hline 
\hline 
5\tabularnewline
\hline 
6\tabularnewline
\hline 
7\tabularnewline
\hline 
8\tabularnewline
\tabularnewline
\tabularnewline
\hline 
\end{tabular}\,%
\begin{tabular}{|c|c|}
\hline 
$T$ & $\xi+\zeta$\tabularnewline
\hline 
\hline 
18 & 18\tabularnewline
\hline 
54 & 62\tabularnewline
\hline 
162 & 222\tabularnewline
\hline 
486 & 806\tabularnewline
 & \tabularnewline
 & \tabularnewline
\hline 
\end{tabular}\quad{}$\longleftrightarrow$\quad{}%
\begin{tabular}{|c|}
\hline 
$n-2$\tabularnewline
\hline 
\hline 
3\tabularnewline
\hline 
4\tabularnewline
\hline 
5\tabularnewline
\hline 
6\tabularnewline
\hline 
7\tabularnewline
\hline 
8\tabularnewline
\hline 
\end{tabular}\,%
\begin{tabular}{|c|c|}
\hline 
$\vartheta_{\textrm{glp}}$ & $\vartheta_{\textrm{lup}}$\tabularnewline
\hline 
\hline 
8 & 2\tabularnewline
\hline 
102 & 26\tabularnewline
\hline 
34 & 70\tabularnewline
\hline 
104 & 92\tabularnewline
\hline 
62 & 190\tabularnewline
\hline 
398 & 788\tabularnewline
\hline 
\end{tabular}\enskip{}.\bigskip{}

\noindent By setting the limits on $n$ and $n-2$ in the respective
halves to $n_{\textrm{max}}=8$, the outcome is $30$ solutions for
$m\in M_{\nicefrac{5}{8}}\cup M_{\nicefrac{5}{8}}^{+}$ and a subset
of 14 solutions for $m\in M_{\nicefrac{5}{8}}$: \bigskip{}

\qquad{}\qquad{}\enskip{}%
\begin{tabular}{|c|}
\hline 
$m$\tabularnewline
\hline 
\hline 
\tabularnewline
4\tabularnewline
\tabularnewline
\tabularnewline
\tabularnewline
\tabularnewline
\tabularnewline
\tabularnewline
\tabularnewline
\tabularnewline
\hline 
\tabularnewline
9\tabularnewline
\tabularnewline
\tabularnewline
\hline 
\tabularnewline
19\tabularnewline
\tabularnewline
\tabularnewline
\hline 
\end{tabular}\,%
\begin{tabular}{|c|c|c|}
\hline 
$a$ & $b$ & $c$\tabularnewline
\hline 
\hline 
26 & -2 & 18\tabularnewline
26 & 34 & 162\tabularnewline
-26 & 62 & 222\tabularnewline
70 & -2 & 62\tabularnewline
70 & 104 & 486\tabularnewline
62 & -2 & 54\tabularnewline
190 & 8 & 222\tabularnewline
190 & -34 & 54\tabularnewline
398 & 102 & 806\tabularnewline
-398 & 104 & 18\tabularnewline
\hline 
 &  & \tabularnewline
788 & 2 & 806\tabularnewline
 &  & \tabularnewline
 &  & \tabularnewline
\hline 
-8 & 26 & 486\tabularnewline
70 & 8 & 222\tabularnewline
92 & -2 & 54\tabularnewline
 &  & \tabularnewline
\hline 
\end{tabular}\,%
\begin{tabular}{|c|}
\hline 
$m$\tabularnewline
\hline 
\hline 
\tabularnewline
5\tabularnewline
\tabularnewline
\tabularnewline
\tabularnewline
\tabularnewline
\tabularnewline
\tabularnewline
\tabularnewline
\tabularnewline
\hline 
\tabularnewline
10\tabularnewline
\tabularnewline
\tabularnewline
\hline 
\tabularnewline
20\tabularnewline
\tabularnewline
\tabularnewline
\hline 
\end{tabular}\,%
\begin{tabular}{|c|c|c|}
\hline 
$a$ & $b$ & $c$\tabularnewline
\hline 
\hline 
8 & 2 & 18\tabularnewline
-8 & 34 & 162\tabularnewline
102 & -8 & 62\tabularnewline
26 & 92 & 486\tabularnewline
-34 & 104 & 486\tabularnewline
92 & 26 & 222\tabularnewline
-398 & 92 & 62\tabularnewline
-788 & 190 & 162\tabularnewline
 &  & \tabularnewline
 &  & \tabularnewline
\hline 
-26 & 8 & 54\tabularnewline
34 & 2 & 54\tabularnewline
-62 & 8 & 18\tabularnewline
-398 & 62 & 222\tabularnewline
\hline 
2 & 8 & 162\tabularnewline
102 & -2 & 62\tabularnewline
-34 & 26 & 486\tabularnewline
62 & 8 & 222\tabularnewline
\hline 
\end{tabular} 

\bigskip{}

\noindent Obviously, the number of solutions, 14 and 30, form a subset
of the upper bounds of summation found in our quasi-supersymmetric
`boson' and `fermion' lists. Although all excitations listed therein
might seem meaningful (not hypo\textcompwordmark{}thetical particles,
but as-of-yet-unidentified atomic or subatomic features), the scopes
of the results $a+bm=c$ found mean that, for the subset of upper
bounds of summation in question, a `boson' that is equivalent to its
juxtaposed `fermion' refers to radiation. The subset must obey the
condition 
\begin{equation}
\sum_{\nu=1}^{\ell_{i}}p_{\nu}=2^{\ell_{i}+1}-2^{i+1}.\label{eq:cond6-14}
\end{equation}
The right-hand side is an even number, which matches the left-hand
side $\sum_{\nu=1}^{\ell_{i}}p_{\nu}$ only if the exponent $\ell_{i}+1=p_{i+1}=3,7,15,31,\ldots$.
We note however that the selection rule (\ref{eq:cond6-14}) may equivalently
be expressed as
\begin{equation}
\frac{\sum_{\nu=1}^{\ell_{i}}p_{\nu}}{2^{\ell_{i}+1}}=\frac{2^{{\scriptscriptstyle \sum}_{\nu=1}^{i}p_{\nu}}-1}{2^{{\scriptscriptstyle \sum}_{\nu=1}^{i}p_{\nu}}}.\label{eq:cond1-4}
\end{equation}

\noindent  Setting $n_{0}=0$, the conditions (\ref{eq:cond6-14})-(\ref{eq:cond1-4})
are compatible with either assignment $n\!:=\ell_{i}+1$, $n\!:={\textstyle \Sigma}_{\nu=1}^{i}p_{\nu}$,
in $\gamma_{\alpha}^{(n)}$. %
\footnote{\ If the restriction to fermions and bosons is given up, all $n=1,2,\ldots$
are of course admissible. In this case, the analysis of $\Delta_{m}$
in $\between_{\,\alpha}^{(n)}=x_{m}\pm\Delta_{m}\quad x_{m}\!\in\!\left\{ L_{m}^{+},L_{m},L_{m}^{-}\right\} $
would ask for an \emph{anyonic} quaternion accumulator; anyonic quaternion
algebra rules have been studied by the author in \cite{Merkel1}.%
}\medskip{}

\noindent The solution to this paradox lies in the identity 
\begin{equation}
p_{r+1}-\sum_{\nu=1}^{r}p_{\nu}\;(=\ell_{r}+1-\sum_{\nu=1}^{r}p_{\nu})=r+1\label{eq:old-ident}
\end{equation}
previously introduced in \cite{Merkel}. We may interpret $r+1$ as
the time dilation that corresponds to a redshift of radiation in a
gravity potential: assignment in question $n:=\sum_{\nu=1}^{r}p_{\nu}+n_{0}$;
radiation unaffected by gravity then corresponds to the assignment
$n':=\ell_{r}+1+n_{0}$ so that $r+1=n'-n$. We recognize that $n'$
is context-sensitive: If $r+1=10$ and $n_{0}=0$, then $n=1013,n'=1023$.
$1023$ may in its turn be classified as red-shifted with respect
to $r+1=10$ and $n_{0}'=10$, with $n=1023,n'=1033$. How the situation
can be disambiguated tells us the equivalence principle which allows
to partition terms $\gamma_{\alpha}^{(n)}$ into components (Unruh
effect): \\
 `observer in accelerated reference frame,' 

\quad{}$\,\begin{array}{c}
\mathcal{O}=(\gamma_{\alpha_{\xi}}^{(n)};\:\gamma_{\alpha_{\xi}}^{(n)}\!\in\mathcal{\! L}^{(n)}),\end{array}$ 

\noindent and components of the `heat bath felt,'

\quad{}\negthinspace{} $\begin{array}{c}
\mathcal{H}=(\gamma_{\alpha_{\zeta}}^{(n)};\:\gamma_{\alpha_{\zeta}}^{(n)}\!\notin\mathcal{\! L}^{(n)}),\end{array}$

\noindent where $\mathcal{L}^{(n)}$ \negthinspace{}is a subset
of generators, $G_{\mu\nu}^{(p)}$ and auxiliaries, of $\{R(L_{m}^{,}\!\!\!\!\,^{\mp});m\!\textrm{=}\!{\scriptstyle {\textstyle 2,3}},{\scriptstyle \ldots}\}$,
$R(\cdot)$ being representations of $L_{m}^{,}\!\!\,^{\mp}$ (see
\cite{Merkel}). For $2\leq m\leq7$, LL($G_{\mu\nu}^{(15)}$) yields
\[
\begin{array}{rl}
\!\!\! R(L_{2}^{,}\!\!\,^{\mp})=\!\! & \!\!\! G_{5,1}+G_{5,3}\left(\mp G_{5,4}\right)=5+1(\mp1)=6(\mp1),\\
\!\!\! R(L_{3}^{,}\!\!\,^{\mp})=\!\! & \!\!\! G_{7,2}-G_{7,3}\left(\mp G_{5,4}\right)=17-5(\mp1)=12(\mp1),\\
\!\!\! R(L_{4}^{,}\!\!\,^{\mp})=\!\! & \!\!\! G_{7,1}-G_{7,2}\left(\mp G_{5,4}\right)=41-17(\mp1)=24(\mp1),\\
\!\!\! R(L_{5}^{,}\!\!\,^{\mp})=\!\! & \!\!\! G_{7,1}\!-\! G_{7,2}\!+\! G_{8,3}\!+\! G_{8,4}\left(\mp G_{5,4}\right)=41\!-\!17\!+\!11\!+5(\mp1)=40(\mp1),\\
\!\!\! R(L_{6}^{,}\!\!\,^{\mp})=\!\! & \!\!\! G_{8,1}-G_{8,2}\left(\mp G_{5,4}\right)=113-41(\mp1)=72(\mp1),\\
\!\!\! R(L_{7}^{,}\!\!\,^{\mp})=\!\! & \!\!\! G_{7,1}\!-\! G_{7,2}\!+\! G_{8,1}\!-\! G_{8,3}\left(\mp G_{5,4}\right)=\!41\!-\!17\!+\!113\!-\!11(\mp1)=\!126(\mp1).
\end{array}
\]
Let us first examine $\mathcal{L}^{(n')}$ for $n'=1023$. We see
that all $G_{\mu\nu}^{(15)}$ above are present in $\gamma_{\alpha_{\xi}}^{(1023)}\!$:
$\gamma_{\alpha_{1}}^{(1023)}=1$, $\gamma_{\alpha_{3}}^{(1023)}=3$,
$\ldots$, $\gamma_{\alpha_{113}}^{(1023)}=113$. Because LL($G_{\mu\nu}^{(15)}$)
is identical to UR(LL($G_{\mu'\nu'}^{(31)}$)), we may go on checking
whether or not $\gamma_{\alpha_{\xi}}^{(1023)}$ extend to $G_{\mu'\nu'}^{(31)}$
from UL(LL($G_{\mu'\nu'}^{(31)}$)) (\ref{sec:Crotons-as-boundary}).
They do, but include shifts,
\[
\begin{array}{rlcc}
6=\gamma_{\alpha_{6}}^{(1023)} & =G_{9,4} & \!\!\!\!\!\!-\mathit{13}= & \!\!\!\!\!\!\!\!\!19-\mathit{13},\\
29=\gamma_{\alpha_{29}}^{(1023)} & =G_{9,3} & \!\!\!\!\!\!-\mathit{14}= & \!\!\!\!\!\!\!\!\! G_{10,4}-\mathit{14}=43-\mathit{14},\\
100=\gamma_{\alpha_{100}}^{(1023)} & =G_{10,3} & \!\!\!\!\!\!-\mathit{15}= & \!\!\!\!\!\!\!\!\!115-\mathit{15},\\
140=\gamma_{\alpha_{140}}^{(1023)} & =G_{9,2} & \!\!\!\!\!\!-\mathit{15}= & \!\!\!\!\!\!\!\!\! G_{11,4}-\mathit{15}=155-\mathit{15},\\
413=\gamma_{\alpha_{413}}^{(1023)} & =G_{9,1} & \!\!\!\!\!\!-16= & \!\!\!\ldots=G_{12,4}-\mathit{16}=429-\mathit{16},
\end{array}
\]
where the shifts in italics range from $\mathit{13}$ to $\mathit{16}$
and have a presence $\gamma_{\alpha_{\mathit{13}}}^{(1023)}=\mathit{13}$,
$\ldots,$ $\gamma_{\alpha_{\mathit{16}}}^{(1023)}=\mathit{16}$.
The $(3^{\log_{2}(p+1)-3}+1)/2$ members of $G_{\rho}^{(p)}$ that
are less or equal $C_{q}$ ($q=(p-3)/4,p=31$) are those that incur
the shifts (LHS), and they may partner with each other to extend the
scope of $R(\cdot)$. Thus, $100+140=R(L_{8})$, $413-140=R(L_{9}^{+})$
and $413-100+29-6=R(L_{10})$. If, instead of $n'=1023$, $n=1013$
is drawn on, the range of shifts is for three of the above $G_{\mu'\nu'}^{(31)}$
lowered to $\mathit{12}$ down to $\mathit{10}$, 
\[
\begin{array}{rl}
31=\gamma_{\alpha_{31}}^{(1013)}= & G_{9,3}-\mathit{12}\;\,=43-\mathit{12},\\
104=\gamma_{\alpha_{104}}^{(1013)}= & G_{10,3}-\mathit{11}=115-\mathit{11},\\
419=\gamma_{\alpha_{419}}^{(1013)}= & G_{11,3}-\mathit{10}=429-\mathit{10},
\end{array}
\]
the new shifts having the presence $\gamma_{\alpha_{\mathit{12}}}^{(1013)}=\mathit{12}$,
$\ldots,$ $\gamma_{\alpha_{\mathit{10}}}^{(1023)}=\mathit{10}$.
No new $R(\cdot)$ can be formed from the shifted terms (LHS), though.
Note that they are constrained to a (truncated) column, $G_{\mu',3}$.
Another column, namely $G_{\mu',6}$, marks the constraint in the
simple case with no shifts,
\[
\begin{array}{rl}
1= & G_{8,6}=\gamma_{\alpha_{1}}^{(1013)},\\
3= & G_{9,6}=\gamma_{\alpha_{3}}^{(1013)},\\
5= & G_{10,6}=\gamma_{\alpha_{5}}^{(1013)},\\
17= & G_{11,6}=\gamma_{\alpha_{17}}^{(1013)},\\
41= & G_{12,6}=\gamma_{\alpha_{41}}^{(1013)}.
\end{array}
\]
We may, with all due caution, say shifted terms express inertia, here
$-$ while working constructively in UL(LL($G_{\mu'\nu'}^{(31)}$))
for $n'=1023$, they impede expressing terms from UL(LL($G_{\mu'\nu'}^{(31)}$))
for $n=1013$, even though they help preserve $R(L_{6}^{,}\!\!\,^{\mp})$
and $R(L_{7}^{,}\!\!\,^{\mp})$, at least. Whatever the value of $n_{0}$
in the gravitational redshift formulae $n:=\sum_{\nu=1}^{r}p_{\nu}+n_{0}$,
$n':=\ell_{r}+1+n_{0}$, $r+1=n'-n$: equating gravitational and inertial
mass means setting $n_{0}=0$. We may check for terms that satisfy
the $\mathcal{H}$ definition under $\mathcal{L}^{(1013)}\!\!=\!\{G_{\mu\nu}^{(15)}\setminus113\}\cup\{31,104,419\}$.
The complete decomposition of $\mathcal{H}$, the `radiation' felt
by $\mathcal{O}$, into sums of type (\ref{eq:decomp-5/8or9/8}) is
done in Table \ref{tab:Specific-fractions-in-1-1-2}. \newpage{}

\section{\label{sec:Application-to-subatomic}On to subatomic physics}

Speaking of charge and gravity, let us, as a prelude, examine the
ratio of the electrical to the gravitational forces between a proton
and an electron (where a first kind of complementarity comes into
play). Consider the collections formed by $\mu$ Magnus terms $\mathcal{M}_{k}\equiv(2k+1)^{2}(-x)^{k(k+1)/2}\:(0<x<1)$
and the estimated number of protons in the universe, $N=10^{80}$,\vspace{-2mm}
\begin{equation}
\sum^{(\mu)}\mathcal{M}_{k}(N-\mu+1),\label{eq:Magnus}
\end{equation}
from which $x$ is to be determined. The electrical force $F_{e}$
is considered independent of $N$; thus $\mu=1$, \emph{i.e}., only
one Magnus term is there to account for $x_{e}$. Assuming, for the
sake of simplicity, that the boundaries $\Gamma^{(15)}$, $\chi^{(15)}$
are sufficient for the proton-electron system, we make a choice of
the triple $\left(k,2k+1,x_{e}^{-1}\right)$ such that expression
(\ref{eq:Magnus}) forms a least upper bound to the observed ratio
$F_{e}/F_{g}$ under the constraint that only successive Mersenne
numbers are being used. This is fulfilled for $\mathcal{M}_{7}(x_{e})=225(-x_{e})^{28}$
with $k=7$, $2k+1=15$, $x_{e}^{-1}=31$ where (\ref{eq:Magnus})
just assumes the reasonable value $\mathcal{M}_{7}(1/31)\times10^{80}\approx3.92\times10^{40}$.
Thus
\[
x_{e}=\frac{1}{31}.
\]
In contrast, the gravitational force according to Mach's principle
is dependent on all other gravitating bodies in the universe so that,
in this case, $\mu=N$ and expression (\ref{eq:Magnus}) reduces to
$10^{80}$ Magnus sum terms, starting with $k=0$, that are going
to account for $x_{g}$. To good approximation,
\[
x_{g}\approx\Lambda=0.10765\ldots,
\]
the so-called 'one-ninth' constant which is the unique exact solution
of the full Magnus equation $\sum_{k=0}^{\infty}(2k+1)^{2}(-x)^{k(k+1)/2}=0\quad(0<x<1)$.
Next, we come to the croton complementarity mentioned in the introduction,
which plays the part of fine-tuning: 170 croton field values are representable
on $\Gamma^{(15)}$, 40 on $\chi^{(15)}$; but only the difference
in the number of values represented seems to matter (Sect.~\ref{sec:Supplementary-remarks-and}),
leaving $2^{130}$ combinations for the power set of crotons compatible
with charged particles. On the other hand, we have seen that, in \emph{$3D$}
space, the part of interest here of curved spacetime or gravity, there
is just one `wavelength' that is minimal in either \emph{$3D$} space
chunk form, squashed or stretched: $(\lambda_{\textrm{min}})_{3}=19$
(see previous section). The power set allowing just two combinations
in this case, one finds
\begin{equation}
\frac{2^{130}}{x_{e}}{\displaystyle \:\div\:\frac{2}{x_{g}}\approx2.27123\times10^{39},}\label{eq:ratio}
\end{equation}
which coincides with the measured ratio $F_{e}/F_{g}$ to five decimal
places.

\noindent In the squashed $3D$ space chunk case, $(\lambda_{\textrm{min}})_{3}=19$
means there are alignments 
\[
(13+k)L_{3}^{+}+(6-2k)L_{3}+kL_{3}^{-}=241,\quad k=0,1,2,3;
\]

 \noindent and in the stretched $3D$ space chunk case, $(\lambda_{\textrm{min}})_{3}=19$
means 
\[
(16-k)L_{3}^{+}+(1+2k)L_{3}+(2-k)L_{3}^{-}=242,\quad k=0,1,2.
\]
From that we glean that, in general, the only never absent sphere
packing form is $L_{3}^{+}=13$. Crotons $\leq13$ are key to providing
a $3D$ scenery: However long the legs of ever-expanding Mersenne
fluctuations, the feet are always rooted in crotons $\leq13$. For
a Mersenne fluctuation to allow for particulates, a cloud of crotons
$\leq13$ has to keep company with them to administer a background
of 2-$3D$ space chunks. The denser that cloud, the more convincing
the impression of a persistent $3D$ continuum in which a particulate
is embedded.

\noindent We are now prepared to subatomic physics. We have seen
that complementarity of bound\-ary field values affects the power
set of croton combinations admissible in a situation, the theoretical
upper bound for combinations of order 31 being $2^{3^{18}-1}$. It
seems reasonable to associate $\Gamma_{x}^{(31)}\notin\chi^{(31)}$
with nuclear phenomena $-$ whose fundamental laws and constants are
unknown and whose pecularities such as the EMC effect and SRC plateaux
\cite{Higin} have remained puzzling to this day $-$ and reserve
non-complementary croton combinations to quarks and preons. Magnus-type
considerations cannot be expected to apply without qualification.
A safe starting point is to presume that preons carry electric charge,
an assumption that allows to associate crotonic activity to (para)fermionic
forms of particulates, from superordinate levels such as protons and
neutrons down to quarks and quark constituents.

\noindent Oscar Wallace Greenberg envisaged a parafermionic nature
of quarks. But with the advent of QCD, and the experimental findings,
valid to this day, that quarks are pointlike down to $10^{-20}\,$m,
preons, parafermionic or otherwise, have not found much acclaim among
physicists. The consequence of the experimental standoff is that preons,
if they exist, must inhabit extradimensions, do aggregate there and
betray their origin only in short-lived resonances known as quark
flavors. It is known that the up quark carries more momentum than
the down quark, which makes it likely that even the two of them are
not of the same dimensional origin. The following is not meant to
be a worked out general model of hadronic matter $-$ it just contemplates
on the possible mathematical structure of the subatomic onion in the
light of crotonic activity. In what follows we use the notation $f_{n+1}(=2^{n+1}-1)$
to denote (para-)fermionic order and the symbols $\mathrm{c}_{\textrm{up}}^{(f_{n+1})}$
and $\mathrm{c}_{\textrm{down}}^{(f_{n+1})}$ for up-type and down-type
preon charge of that order, respectively.
\begin{conjecture}
\label{con:Preons-of-order}Preons of order $f_{n+1}$ are either
up-type or down-type, $\mathrm{\mathit{preon}}_{\textrm{\emph{up}}}^{(f_{n+1})}$
or $\mathrm{\mathit{preon}}_{\textrm{\emph{down}}}^{(f_{n+1})}$.
The electric charge (in \emph{e}) of up-type items is given by the
expressions $\mathrm{c}_{\mathrm{up}}^{(f_{n+1})}=(f_{n+1}-\sum_{s=0}^{n}f_{s})/\prod_{r=1}^{n+1}f_{r}=(n+1)/\prod_{r=1}^{n+1}f_{r}$,
while down-type items\emph{ }have the charge $\mathrm{c}_{\mathrm{down}}^{(f_{n+1})}=-\sum_{s=0}^{n}f_{s}/\prod_{r=1}^{n+1}f_{r}$
(see Table \ref{tab:Interordinal-preon-model-1-1-1} below). The charge
of up-type items transforms as\emph{ $\mathrm{c}_{\mathrm{up}}^{(f_{n})}=(f_{n+1}-1)\mathrm{c}_{\mathrm{up}}^{(f_{n+1})}+\mathrm{c}_{\mathrm{down}}^{(f_{n+1})}$}
and the charge\emph{ }of down-type items as $\mathrm{c}_{\mathrm{down}}^{(f_{n})}=(f_{n}+1)\mathrm{c}_{\mathrm{down}}^{(f_{n+1})}+f_{n}\mathrm{c}_{\mathrm{up}}^{(f_{n+1})}$\textup{.}

\begin{table}[H]
\caption{\label{tab:Interordinal-preon-model-1-1-1}Mersennian preon charge
model}
\quad{}\quad{}\quad{}%
\begin{tabular}[t]{|c|c|r@{\extracolsep{0pt}.}l|}
\hline 
$f_{n+1}$ & up-type charge $(\textrm{c}_{\textrm{up}})$ & \multicolumn{2}{c|}{down-type charge $(\textrm{c}_{\textrm{down}})$}\tabularnewline
\hline 
\hline 
1 & 1 & \multicolumn{2}{c|}{0}\tabularnewline
\hline 
3 & $\nicefrac{2}{3}$ & \multicolumn{2}{c|}{$-\,\nicefrac{1}{3}$}\tabularnewline
\hline 
7 & $\nicefrac{3}{21}$ & \multicolumn{2}{c|}{$-\,\nicefrac{4}{21}$}\tabularnewline
\hline 
15 & $\nicefrac{4}{315}$ & \multicolumn{2}{c|}{$-\,\nicefrac{11}{315}$}\tabularnewline
\hline 
31 & $\nicefrac{5}{9765}$ & \multicolumn{2}{c|}{$-\,\nicefrac{26}{9765}$}\tabularnewline
\hline 
63 & $\nicefrac{6}{615195}$ & \multicolumn{2}{c|}{$-\,\nicefrac{57}{615195}$}\tabularnewline
${\scriptstyle \cdots}$ & ${\scriptstyle \cdots}$ & \multicolumn{2}{c|}{${\scriptstyle \cdots}$}\tabularnewline
\hline 
\end{tabular}
\end{table}

~
\end{conjecture}
\noindent  The Magnus formalism suggests a connection between $f_{n+1}$
(or $f_{n}$) and $x_{e}^{-1}=31$ $-$ this will not only suffice
for the proton and the neutron (which are assigned the least order
$f_{1}=1$); if the greatest assignment eligible is $f_{n+1}:=x_{e}^{-1}$,
it suffices for three generations of quarks, and if $f_{n}:=x_{e}^{-1}$
is eligible, for a fourth generation as well. Here, the clue to successful
bounds for a representation in terms of kissing numbers comes from
a divisibility postulate for a generation's $L_{\textrm{up}}$: in
addition to being divisible by $f_{n+1}-1$, $L_{\textrm{up}}$ must
contain a genuine prime factor $P_{\mu}>13$ (larger than the $L_3^+$ of $3D$ space) such that 
\begin{equation}
\left|P_{\mu}-L_{\mu_{0}+\mu}\right|=1\quad(\mu=1,2,3).\label{eq:postulate}
\end{equation}
The task is for the $\mu$th generation completed when all its kissing
numbers $L_{\textrm{down}}$ $-$ there are several $-$ which are
divisible by $f_{n}+1$ and have but prime factors less $P_{\mu}$
are identified. Then for each $L_{\textrm{down}}$ found the single
ratio $\frac{L_{\textrm{up}}}{L_{\textrm{down}}}$ could be considered
a bound to the ratio $\frac{m_{\textrm{up}}}{m_{\textrm{down}}}$
in question. However, true to the Magnus ansatz, getting a fine-tuned
result requires taking all contributors into account. Table \ref{tab:Prime-factors-of}
shows how, for each generation, a ratio $\frac{L_{\textrm{up}}}{\Sigma L_{\textrm{down}}}$
can be deduced that bounds the respective measured ratio $\frac{m_{\textrm{up}}}{m_{\textrm{down}}}$
from below. This principle is best understood as a simile to the Magnus
ansatz, where the intra-generational quark-mass ratios $m_{\textrm{u}}/m_{\textrm{d}}$,
$m_{\textrm{c}}/m_{\textrm{s}}$ and $m_{\textrm{t}}/m_{\textrm{b}}$
replace the dimensionless force ratio $F_{e}/F_{g}$.

\noindent  Quark mass is assumed to result from crotonic activity,
and the configurations {\tt c},{\tt s} and {\tt t},{\tt d} make
it clear that this activity has to cover extended spans of orders.
Here, only leading-order crotonic activity is considered in deriving
bounds for intra-generational mass ratios. This implies identifying
where leading-order crotonic activity singles out space chunks that
suit the up-type quark of a generation and other space chunks suiting
the down-type quark. The kissing numbers of the target spaces, $L_{\textrm{up}}$
and $L_{\textrm{down}}$, must in turn show the divisibility properties
demanded in Conjecture \ref{con:To-serve-as}. But that's only a necessary
condition. In the Magnus ansatz, assignment of successive Mersenne
numbers to the triple $(k,2k+1,x_{e}^{-1}$) is essential to getting
a handle on bounds.

\newpage{}
\begin{table}[H]
\caption{\label{tab:Prime-factors-of}Prime factors of $(x_{e})^{-1}$ kissing
numbers; characterstic divisors determine up-type and down-type kissing
numbers that bound measured intra-generational quark mass ratios from
below}
\vspace{0.2cm}
\begin{tabular}{|l|l|c|c|c|c|}
\hline 
\multirow{2}{*}{$m$} & \multirow{2}{*}{$L_{m}$} & \multirow{2}{*}{prime factorization} & \multirow{2}{*}{$\begin{array}{c}
\textrm{divisors}\\
\textrm{{\tt t}}\:[\textrm{{\tt b}}]
\end{array}$} & \multirow{2}{*}{$\begin{array}{c}
\textrm{divisors}\\
\textrm{{\tt c}}\:[\textrm{{\tt s}}]
\end{array}$} & \multirow{2}{*}{$\begin{array}{c}
\textrm{divisors}\\
\textrm{{\tt u}}\:[\textrm{{\tt d}}]
\end{array}$}\tabularnewline
 &  &  &  &  & \tabularnewline
\hline 
\hline 
1 & 2 & 2 &  &  & \tabularnewline
\hline 
2 & 6 & $2\times3$ &  &  & \tabularnewline
\hline 
3 & 12 & $2^{2}\times3$ &  &  & \tabularnewline
\hline 
4 & 24 & $2^{3}\times3$ &  &  & $[2^{2}]$\tabularnewline
\hline 
5 & 40 & $2^{3}\times5$ &  &  & $[2^{2}]$\tabularnewline
\hline 
6 & 72 & $2^{3}\times3^{2}$ &  &  & $[2^{2}]$\tabularnewline
\hline 
7 & 126 & $2\times3^{2}\times7$ &  &  & \tabularnewline
\hline 
8 & 240 & $2^{4}\times3\times5$ & $[2^{4}]$ & $[2^{3}]$ & $[2^{2}]$\tabularnewline
\hline 
9 & 272 & $2^{4}\times17$ & $[2^{4}]$ & $[2^{3}]$ & $[2^{2}]$\tabularnewline
\hline 
10 & 336 & $2^{4}\times3\times7$ & $[2^{4}]$ & $[2^{3}]$ & $[2^{2}]$\tabularnewline
\hline 
11 & 438 & $2\times3\times$%
\framebox{\begin{minipage}[t][0.6\totalheight]{0.03\columnwidth}%
73%
\end{minipage}} &  &  & $2\times3$\tabularnewline
\hline 
12 & 756 & $2^{2}\times3^{3}\times7$ &  &  & \tabularnewline
\hline 
13 & 918 & $2\times3^{3}\times17$ &  &  & \tabularnewline
\hline 
14 & 1422 & $2\times3^{3}\times79$ &  &  & \tabularnewline
\hline 
15 & 2340 & $2^{2}\times3^{2}\times5\times13$ &  &  & \tabularnewline
\hline 
16 & 4320 & $2^{5}\times3^{3}\times5$ & $[2^{4}]$ &  & \tabularnewline
\hline 
17 & 5346 & $2\times3^{5}\times11$ &  &  & \tabularnewline
\hline 
18 & 7398 & $2\times3^{3}\times137$ &  &  & \tabularnewline
\hline 
19 & 10668 & $2^{2}\times3\times7\times$%
\framebox{\begin{minipage}[t][0.6\totalheight]{0.04\columnwidth}%
127%
\end{minipage}} &  & $2\times7$ & \tabularnewline
\hline 
20 & 17400 & $2^{3}\times3\times5^{2}\times29$ &  &  & \tabularnewline
\hline 
21 & 27720 & $2^{3}\times3^{2}\times5\times7\times11$ &  &  & \tabularnewline
\hline 
22 & 49896 & $2^{3}\times3^{4}\times7\times11$ &  &  & \tabularnewline
\hline 
23 & 93150 & $2\times3^{4}\times5^{2}\times23$ &  &  & \tabularnewline
\hline 
24 & 196560 & $2^{4}\times3^{3}\times5\times7\times13$ &  &  & \tabularnewline
\hline 
25 & 197040 & $2^{4}\times3\times5\times821$ &  &  & \tabularnewline
\hline 
26 & 198480 & $2^{4}\times3\times5\times827$ &  &  & \tabularnewline
\hline 
27 & 199912 & $2^{3}\times24989$ &  &  & \tabularnewline
\hline 
28 & 204188 & $2^{2}\times51047$ &  &  & \tabularnewline
\hline 
29 & 207930 & $2\times3\times5\times29\times$%
\framebox{\begin{minipage}[t][0.6\totalheight]{0.04\columnwidth}%
239%
\end{minipage}} & $2\times3\times5$ &  & \tabularnewline
\hline 
30 & 219008 & $2^{7}\times29\times59$ &  &  & \tabularnewline
\hline 
31 & 230872 & $2^{3}\times28859$ &  &  & \tabularnewline
\hline 
\multirow{2}{*}{} & \multirow{2}{*}{$\frac{L_{\textrm{up}}}{\Sigma L_{\textrm{down}}}$} & \multirow{2}{*}{} & \multirow{2}{*}{$\thickapprox40.23$} & \multirow{2}{*}{$\thickapprox12.58$} & \multirow{2}{*}{$\thickapprox0.45$}\tabularnewline
 &  &  &  &  & \tabularnewline
\hline 
\multirow{2}{*}{} & \multirow{2}{*}{$\frac{m_{\textrm{up}}}{m_{\textrm{down}}}$} & \multirow{2}{*}{} & \multirow{2}{*}{$\thickapprox41.86$} & \multirow{2}{*}{$\thickapprox13.58$} & \multirow{2}{*}{$\thickapprox0.48$}\tabularnewline
 &  &  &  &  & \tabularnewline
\hline 
\end{tabular}
\end{table}
 \newpage{}

\noindent Let us clear up the inner workings of Table \ref{tab:Prime-factors-of},
beginning with the third-generation quarks, $\textrm{{\tt t}}$ and
$\textrm{{\tt b}}$. With the top quark, $f_{5}-1=30$ hyperspheres
must fit in a chunk of space such that $f_{5}-1$ and a prime factor
$P_{3}$ satisfying postulate (\ref{eq:postulate}) for some $\mu_{0}$
divide its kissing number without rest. Both is true for the spaces
$23D$ and $29D$, with the candidates $(93150,23)$ and $(207930,239)$
for $(L_{\textrm{up}},P_{3})$, respectively. Only by observing that
the same selection rules must apply to the other generations are we
able to decide that $(207930,239)$ (framed in the table) is the appropriate
pair $-$ a $P_{3}=L_{4}-1$ would not leave place for a $P_{2}$
exceeding 13: $L_{3}\pm1=P_{2}\ngtr13$. The bottom quark is collectively
realized by all subspace chunks in which there are $\,\, f_{4}+1=16\,\,$
hyperspheres such that $f_{4}+1$ divides their kissing numbers without
rest (marked by $[\,]$) and the prime factors involved are less $P_{3}=239$.
For the next lower generation, the pair suiting the up-type quark
is from $19D$, $(10668,127)$, and for the first generation, that
pair \noindent is from $11D$, $(438,73)$. In accordance with postulate
(\ref{eq:postulate}), the $P_{\mu}$ ($\mu=1,2,3$), specify a triple
of successive kissing numbers, $(72,126,240)$. Unsurprisingly 
\begin{equation}
\mu_{0}=\log_{2}\left(x_{e}^{-1}+1\right).
\end{equation}
The kissing numbers of the subspace chunks corresponding to down-type
quarks {\tt s} and {\tt d} too satisfy the required divisibilities
(again marked by $[\,]$).%
\footnote{$\,$ large kissing numbers are an active field of research \cite{Cohn}; those used here are taken from {\small\tt http://www.math.rwth-aachen.de/Gabriele.Nebe/LATTICES/kiss.html}%
}\noindent We now want to hint at the possible existence of a fourth
quark generation. The entries of Table \ref{tab:Quantities-:=00003D-;}
may be used to determine quark family characteristics: 

\begin{table}[H]
\caption{\label{tab:Quantities-:=00003D-;}Quantities $S_{\mu_{0}+\mu-1}$:=
$\sum_{m=1}^{\mu_{0}+\mu-1}(2^{m}-1)$; $L_{\mu_{0}+\mu}$; $P_{\mu}$;
$\Xi(\mu)$ := $6\,\textrm{Prime}(\mu)+(-1)^{\mu}$}
\bigskip{}
\qquad{}\qquad{}\qquad{}\qquad{}\qquad{}%
\begin{tabular}{|c|c|c|c|c|}
\hline 
$\mu$ & $S_{\mu_{0}+\mu-1}$ & $L_{\mu_{0}+\mu}$ & $P_{\mu}$ & $\Xi(\mu)$\tabularnewline
\hline 
\hline 
1 & 57 & 72 & 73 & 11\tabularnewline
\hline 
2 & 120 & 126 & 127 & 19\tabularnewline
\hline 
3 & 247 & 240 & 239 & 29\tabularnewline
\hline 
4 & 502 & 272 & 271 & 43\tabularnewline
\hline 
5 & 1013 & 336 & 337 & 65\tabularnewline
$\vdots$ &  &  &  & \tabularnewline
\hline 
\end{tabular}
\end{table}

\noindent Where $\chi_{\textrm{prime}}(\cdot)$ is the characteristic
function of prime numbers, $L_{\textrm{up}}$ can be said to belong
to the family, and be identified with $L_{\Xi(\mu)}$, if $\chi_{\textrm{prime}}(\Xi(\mu))=1$.
This is obvious for $\mu=1,2,3$. One further notes that the signum
function values, $\textrm{sgn}(P_{\mu}-L_{\mu_{0}+\mu})$ and $\,\textrm{sgn}(S_{\mu_{0}+\mu-1}-P_{\mu})$,
cancel each other out for $\mu=1,2,3$. 

\noindent But said observations hold out for $\mu=4$: The equations

\begin{equation}
\begin{array}{c}
\textrm{sgn}(P_{\mu}-L_{\mu_{0}+\mu})-\textrm{sgn}(P_{\mu}-S_{\mu_{0}+\mu-1})=0,\\
\chi_{\textrm{prime}}(\Xi(\mu))-1=0
\end{array}\quad(\mu_{0}=\log_{2}\left(x_{e}^{-1}+1\right))\label{eq:conservation law}
\end{equation}
\noindent are not violated until at $\mu=5$. What might thus constitute
the quark family's conservation law would predict that $L_{43}$ have
2, 31 and 271 among its prime factors and serve as the $L_{\textrm{up}}$
of a fourth-generation quark {\tt t}$'$.

\noindent It should be noted that the bounds given in Table \ref{tab:Prime-factors-of}
for the first three generations (next-to-last row) are valid for intra-generational
mass ratios only (last row). Ratios of mass for quarks that belong
to different generations are distinctly different. Here, the methods
developed in Sect. \ref{sec:Pregeometric-categories-relevant} take
effect $-$ in particular, the juxtaposition $x$ vs. $x-2$ as applied
to \emph{both} the lower and the upper bound of sum gives a good match
\begin{equation}
\begin{array}{llll}
\textrm{measured} & m_{\textrm{c}}/m_{\textrm{u}}\thickapprox560.84\quad & \textrm{control} & \frac{9}{8}\sum_{l=3}^{8}p_{l}=560.25\\
\:\quad\textrm{ "} & m_{\textrm{t}}/m_{\textrm{c}}\thickapprox135.64\quad & \textrm{\quad\textrm{ "}} & \frac{9}{8}\sum_{l=1}^{6}p_{l}=135
\end{array}\label{eq:intermass}
\end{equation}
\noindent  $\frac{9}{8}\sum_{l=-1}^{4}p_{l}=28.6875$ as a possible
but mathematically ugly weight (in top masses) of the hypothetical
{\tt t}$'$ is an educated guess at best. 

\noindent  How preon charge distributions relate to targets of crotonic
activity is propounded in our second conjecture:
\begin{conjecture}
\label{con:To-serve-as}To qualify as constituents of a superordinate
preon of order $f_{n_{0}+1}$, preon charges must occupy all hyperspheres
of the constituents' packings; the charge multiplets have minimal
order \textup{$f_{\nu>n_{0}+1}$}, or multiples thereof in case there
exists a $\textrm{\emph{LL}}(G_{\alpha\beta}^{(f_{\nu-1})})$ at least
the size of a building block $\binom{a\: b}{c\: d}$.
\end{conjecture}
\noindent We will confine the discussion to the first quark generation
(valence quarks) and focus on their structural configurations as adapted
from Table \ref{tab:Prime-factors-of}. The $\mathtt{u}$ content
is given by $L_{\textrm{up}}=L_{11}=438$, while the $\mathtt{d}$
content would be equal to $\Sigma L_{\textrm{down}}=L_{10}\!+\! L_{9}\!+\! L_{8}\!+\! L_{6}\!+\! L_{5}\!+\! L_{4}=984$
(with a total mixed content 1422). With these assignments, the proton's
total content is $2$$\mathtt{u}$$+$$\mathtt{d}$ $=1860$, and
the neutron's total content $\mathtt{u}$$+$$2$$\mathtt{d}$ $=2406$.
A total mixed content 1422, however, is incompatible with what is
known about the proton: 
\begin{table}[H]
\caption{\label{tab:Structure-of-nucleons}Structure of mixed content of the
valence quarks }
\bigskip{}
\begin{tabular}{|c|c|c|c|}
\hline 
\multicolumn{2}{|c|}{$\begin{array}{c}
\textrm{charge multiplets of order}\, f_{\nu>n_{0}+1}\\
(\textrm{total mixed content}\,1422)
\end{array}$} & \multicolumn{2}{c|}{$\begin{array}{c}
\!\!\!\textrm{charge multiplets of order}\, f_{\nu>n_{0}+1}\\
\!\!\!(\textrm{vetted mixed content}\,1398)
\end{array}\!\!\!$}\tabularnewline
\hline 
\hline 
$=474\times3+0$ & $=474\times3+0$ & $=466\times3+0$ & $=466\times3-0$\tabularnewline
\hline 
$=203\times7+1$ & $=204\times7-6$ & $=199\times7+5$ & $=200\times7-2$\tabularnewline
\hline 
$=94\times15+12$ & $=95\times15-3$ & $=93\times15+3$ & $=94\times15-12$\tabularnewline
\hline 
$=45\times31+27$ & $=46\times31-4$ & $=45\times31+3$ & $=46\times31-12$\tabularnewline
\hline 
$\begin{array}{c}
\Sigma\Delta=40\\
\neq\Sigma p_{l}
\end{array}$ & $\begin{array}{c}
\left|\Sigma\Delta\right|=13\\
\neq\Sigma p_{l}
\end{array}$ & $\begin{array}{c}
\Sigma\Delta=11\\
=\Sigma_{l=1}^{3}p_{l}
\end{array}$ & $\begin{array}{c}
\left|\Sigma\Delta\right|=26\\
=\Sigma_{l=1}^{4}p_{l}
\end{array}$\tabularnewline
\hline 
\end{tabular}
\end{table}

\noindent So we are forced to readjust the $\mathtt{d}$ content
such that the $L_{4}=24$ contribution is eliminated from $\Sigma L_{\textrm{down}}$
$-$ the $\mathtt{d}$ content becomes 960, the total mixed content
1398, the neutron's total content 2358 and the proton's total content
1836. Since our superordinate preon is the nucleon, the superordinate
charge order, according to Table \ref{tab:Interordinal-preon-model-1-1-1},
is $f_{n_{0}+1}=1$, and the constituents, according to Table \ref{tab:Structure-of-nucleons},
have charge multiplets of orders $f_{n_{0}+2}=\!3$, $f_{n_{0}+3}=\!7$,
$f_{n_{0}+4}=\!15$ and $f_{n_{0}+5}=\!31$ (in case of $f_{n_{0}+4},f_{n_{0}+5}$,
multiplet orders may form multiples of 15,31). For $f_{n_{0}+1}=1$,
a $\textrm{LL}(G_{\alpha\beta}^{(1)})$ (read Lower-Left quadrant
of square matrix $(G_{\alpha\beta}^{(1)})$) at least the size of
a building block $\binom{a\: b}{c\: d}$ as demanded by the second
part of Conjecture \ref{con:To-serve-as}, does not exist. So the
$\textrm{c}_{\textrm{up}}$ request for the proton (in units of electric
charge $\textrm{e}$) is simply identified with one triplet ($\nicefrac{2}{3}$,$\nicefrac{2}{3}$,$\nicefrac{-1}{3})\:\rightarrow1$,
and the $\textrm{c}_{\textrm{down}}$ request with one triplet ($\nicefrac{2}{3}$,$\nicefrac{-1}{3}$,$\nicefrac{-1}{3})\:\rightarrow0$;
the remaining $873\!+\!957$ hyperspheres are divided into triplets
of vanishing charge $(\nicefrac{2}{3}$,$\nicefrac{-1}{3}$,$\nicefrac{-1}{3})$
$-$ see first row of Table \ref{tab:Structure-of-nucleons-1-2}: 

\begin{table}[H]
\caption{\label{tab:Structure-of-nucleons-1-2}Structure of protons in multiplet
form 3,7,15$r$ and 31$r$ }
\bigskip{}

\begin{tabular}{|c|c|c|}
\hline 
content & \multicolumn{2}{c|}{charged multiplets}\tabularnewline
\hline 
\hline 
2$\mathtt{u},\mathtt{d}$ $\!\!\!\!\!\!\!\!\!\begin{array}{c}
\\
\\
\\
\end{array}$ & $\begin{array}{c}
\left[291(\nicefrac{2}{3},\nicefrac{-1}{3},\nicefrac{-1}{3})\right.\\
\left.+(\nicefrac{2}{3},\nicefrac{2}{3},\nicefrac{-1}{3})\right]_{2\mathtt{u}}^{876}
\end{array}$ & $\begin{array}{c}
\left[319(\nicefrac{2}{3},\nicefrac{-1}{3},\nicefrac{-1}{3})\right.\\
\left.+(\nicefrac{2}{3},\nicefrac{-1}{3},\nicefrac{-1}{3})\right]_{\mathtt{d}}^{960}
\end{array}$\tabularnewline
\hline 
2$\mathtt{u},\mathtt{d}$ $\!\!\!\!\!\!\!\!\!\begin{array}{c}
\\
\\
\end{array}$ & $\!\!\!\!\begin{array}{c}
\left[123(\underbrace{\nicefrac{3}{21},\nicefrac{\ldots,3}{21}},\underbrace{\nicefrac{-4}{21},\nicefrac{\ldots,-4}{21}})\right.\\
\quad\:{\scriptscriptstyle 4}\qquad\qquad\quad{\scriptscriptstyle 3}\\
\!\!\!\!\!\left.+2(\underbrace{\nicefrac{3}{21},\nicefrac{\ldots,3}{21}},\nicefrac{-4}{21})\right]_{2\mathtt{u}}^{875}\\
\;{\scriptscriptstyle 6}\qquad\qquad\quad
\end{array}$\negthinspace{}\negthinspace{}\negthinspace{} & $\!\!\!\!\!\begin{array}{c}
\left[136(\underbrace{\nicefrac{3}{21},\nicefrac{\ldots,3}{21}},\underbrace{\nicefrac{-4}{21},\nicefrac{\ldots,-4}{21}})\right.\\
\quad\,{\scriptscriptstyle 4}\qquad\qquad\quad\,{\scriptscriptstyle 3}\\
\left.+(\underbrace{\nicefrac{3}{21},\nicefrac{\ldots,3}{21}},\underbrace{\nicefrac{-4}{21},\nicefrac{\ldots,-4}{21}})\right]_{\mathtt{d}}^{959}\\
\!\!\!\!\!\!\!\!\!{\scriptscriptstyle 3}\qquad\qquad\quad{\scriptscriptstyle 4}
\end{array}$\negthinspace{}\negthinspace{}\negthinspace{}\tabularnewline
\hline 
total $\mathtt{p}\!\!\!\!\!\!\!\!\!\begin{array}{c}
\\
\\
\\
\end{array}$ & \multicolumn{2}{c|}{$\!\!\!\!\begin{array}{c}
\left[117(\underbrace{\nicefrac{4}{315},\nicefrac{\ldots,4}{315}},\underbrace{\nicefrac{-11}{315},\nicefrac{\ldots,-11}{315}})+2(\underbrace{\nicefrac{4}{315},\nicefrac{\ldots,4}{315}},\nicefrac{-11}{315})\right.\\
\!\!\!\!\!\!\!\!\!\!{\scriptscriptstyle 11}\qquad\qquad\qquad\:{\scriptscriptstyle 4}\qquad\qquad\qquad\qquad\;{\scriptscriptstyle 29}\\
\left.+(\underbrace{\nicefrac{4}{315},\nicefrac{\ldots,4}{315}},\underbrace{\nicefrac{-11}{315},\nicefrac{\ldots,-11}{315}})\right]_{\mathtt{p}}^{1830}\\
{\scriptscriptstyle 4}\qquad\qquad\qquad\:{\scriptscriptstyle 11}\qquad\quad
\end{array}$}\tabularnewline
\hline 
total  $\mathtt{p}\!\!\!\!\!\!\!\!\!\begin{array}{c}
\\
\\
\\
\end{array}$ & \multicolumn{2}{c|}{$\begin{array}{c}
\textrm{c requests }(\underbrace{\nicefrac{5}{9765},\nicefrac{\ldots,5}{9765}},\underbrace{\nicefrac{-26}{9765},\nicefrac{\ldots,-26}{9765}})\Rightarrow2883>1836\\
\qquad\qquad\quad{\scriptscriptstyle a}\qquad\qquad\qquad\quad{\scriptscriptstyle b}\qquad\qquad\quad{\scriptscriptstyle (a+b=31r)}
\end{array}$}\tabularnewline
\hline 
\end{tabular}
\end{table}

\noindent Nor does a $\textrm{LL}(G_{\alpha\beta}^{(3)})$ the size
of a building block $\binom{a\: b}{c\: d}$ exist, so that the $\textrm{c}_{\textrm{up}}$
and $\textrm{c}_{\textrm{down}}$ requests for the proton would in
a similar manner be fulfilled with simple septets $\rightarrow\nicefrac{2}{3},\nicefrac{-1}{3},0$
out of all formable:
\[
\begin{array}{lcrcc}
{\displaystyle (}\nicefrac{3}{21},\nicefrac{3}{21},\nicefrac{3}{21},\nicefrac{3}{21},\nicefrac{3}{21},\nicefrac{3}{21},\nicefrac{-4}{21}{\displaystyle )} & \rightarrow & \nicefrac{2}{3} &  & (a)\\
(\nicefrac{3}{21},\nicefrac{3}{21},\nicefrac{-4}{21},\nicefrac{-4}{21},\nicefrac{-4}{21},\nicefrac{-4}{21},\nicefrac{-4}{21}) & \rightarrow & \nicefrac{-2}{3} &  & (b)\\
(\nicefrac{3}{21},\nicefrac{3}{21},\nicefrac{3}{21},\nicefrac{-4}{21},\nicefrac{-4}{21},\nicefrac{-4}{21},\nicefrac{-4}{21}) & \rightarrow & \nicefrac{-1}{3} &  & (c)\\
(\nicefrac{3}{21},\nicefrac{3}{21},\nicefrac{3}{21},\nicefrac{3}{21},\nicefrac{3}{21},\nicefrac{-4}{21},\nicefrac{-4}{21}) & \rightarrow & \nicefrac{1}{3} &  & (d)\\
(\nicefrac{3}{21},\nicefrac{3}{21},\nicefrac{3}{21},\nicefrac{3}{21},\nicefrac{-4}{21},\nicefrac{-4}{21},\nicefrac{-4}{21}) & \rightarrow & 0 &  & (e)
\end{array}
\]
One of the septet solutions matching with row 2 of Table \ref{tab:Structure-of-nucleons-1-2}
would be: 
\begin{table}[H]
\caption{\label{tab:Structure-of-nucleons-1-1}Structure of protons in septet
form }
\bigskip{}

\qquad{}\qquad{}\qquad{}%
\begin{tabular}{|c|c|}
\hline 
request & charged hyperspheres divided into septets\tabularnewline
\hline 
\hline 
$2\mathtt{u}$ & $L_{11}+L_{11}^{-}=875$\tabularnewline
\hline 
$\mathtt{d}$ & $L_{10}+L_{9}+L_{8}+L_{6}+L_{5}^{-}=959$\tabularnewline
\hline 
\end{tabular}
\end{table}

\noindent Here, we find centerpiece-free packings as well as two
packings that lack both centerpiece and a peripheral hypersphere.

\noindent A $\textrm{LL}(G_{\alpha\beta}^{(7)})$ at least the size
of $\binom{a\: b}{c\: d}$ does exist (see \ref{sec:Crotons-as-boundary}),
and it is for the first time that it can shape the charge multiplets,
which here are of minimal length $f_{n_{0}+4}=\!15$. It reads $\textrm{LL}(G_{\alpha\beta}^{(7)})=\binom{\boldsymbol{1\!\!\!1}\:1}{\boldsymbol{1}\!\!\!\boldsymbol{1}\:1}$,
and if we denote it by $\binom{G_{3,1}\: G_{3,2}}{G_{4,1}\: G_{4,2}}$,
the $\textrm{c}_{\textrm{up}}$ and $\textrm{c}_{\textrm{down}}$
proton requests take the form
\[
\begin{array}{c}
4(\underbrace{\nicefrac{4}{315},\nicefrac{\ldots,4}{315}},\underbrace{\nicefrac{-11}{315}})\\
{\scriptstyle 2G_{4,1}f_{n+4}-G_{3,1}}\quad\!{\scriptstyle G_{3,1}}
\end{array}
\]

\noindent and
\[
\begin{array}{c}
G_{3,1}(\underbrace{\nicefrac{4}{315},\nicefrac{\ldots,4}{315}},\underbrace{\nicefrac{-11}{315},\nicefrac{\ldots,-11}{315}}),\\
\quad{\scriptstyle 4}\qquad\qquad\qquad{\scriptstyle 11}
\end{array}
\]

\noindent respectively. That some multiplets come in form of proper
multiples of the minimal order $f_{\nu>n_{0}+1}$ implies that only
the total proton content $2\mathtt{u+d}$ matters as regards divisibility
by the minimal order. The appropriately adjusted content must be as
close to 1836 as possible. In case of $f_{n_{0}+4}=\!15$, this is
$1830=122\times15$, so that there remain 117 quindecuplets to fill
with the vanishing charge form 
\[
\begin{array}{c}
(\underbrace{\nicefrac{4}{315},\nicefrac{\ldots,4}{315}},\underbrace{\nicefrac{-11}{315},\nicefrac{\ldots,-11}{315}})\\
\!\!\!\!\!\!\!{\scriptstyle 11}\qquad\qquad\qquad{\scriptstyle 4}
\end{array}
\]
\noindent (see row 3 of Table \ref{tab:Structure-of-nucleons-1-2}).
One of the solutions would read:
\begin{table}[H]
\caption{\label{tab:Structure-of-nucleons-1-1-1}Structure of protons in quindecuplet
form }
\bigskip{}

\ %
\begin{tabular}{|c|c|}
\hline 
request & charged hyperspheres divided into quindecuplets\tabularnewline
\hline 
\hline 
\multirow{2}{*}{- - - $\begin{array}{c}
2\mathtt{u}\\
\mathtt{d}
\end{array}$- - -} & \multirow{2}{*}{$\begin{array}{c}
2L_{11}^{-}=874\\
L_{10}+L_{9}^{-}+L_{8}^{-}+L_{6}^{-}+L_{5}^{-}=956
\end{array}$ }\tabularnewline
 & \tabularnewline
\hline 
\end{tabular}
\end{table}

\noindent For multiplets of minimal order $f_{n+5}=\!31$, the $\textrm{c}_{\textrm{up}}$
and $\textrm{c}_{\textrm{down}}$ proton requests are determined by
the upper left (or, equivalently, lower right) building block of 
\[
\textrm{LL}(G_{\alpha\beta}^{(15)})=\left(\begin{array}{cccc}
G_{5,1} & G_{5,2} & G_{5,3} & G_{5,4}\\
G_{6,1} & G_{6,2} & G_{6,3} & G_{6,4}\\
G_{7,1} & G_{7,2} & G_{7,3} & G_{7,4}\\
G_{8,1} & G_{8,2} & G_{8,3} & G_{8,4}
\end{array}\right)=\left(\begin{array}{cccc}
\boldsymbol{5\!\!\!5} & 3 & 1 & 1\\
\boldsymbol{1\!\!\!1}\boldsymbol{1\!\!\!1} & 5 & 1 & 1\\
41 & 17 & \boldsymbol{5\!\!\!5} & 3\\
113 & 41 & \boldsymbol{1\!\!\!1}\boldsymbol{1\!\!\!1} & 5
\end{array}\right).
\]
\noindent The requests become
\[
\begin{array}{cc}
4(\underbrace{\nicefrac{5}{9765},\nicefrac{\ldots,5}{9765}},\underbrace{\nicefrac{-26}{9765}}) & (4\times682\textrm{ hyperspheres)}\\
{\scriptstyle 2G_{6,1}f_{n+5}-G_{5,1}}\quad\!{\scriptstyle G_{5,1}}
\end{array}
\]

\noindent and
\[
\begin{array}{cc}
G_{5,1}(\underbrace{\nicefrac{5}{9765},\nicefrac{\ldots,5}{9765}},\underbrace{\nicefrac{-26}{9765},\nicefrac{\ldots,-26}{9765}}), & (5\times31\textrm{ hyperspheres)}\\
\quad{\scriptstyle 5}\qquad\qquad\qquad\quad{\scriptstyle 26}
\end{array}
\]
respectively. A total request of 2883 charged hyperspheres exceeds
by far the available total proton content (see Table \ref{tab:Structure-of-nucleons-1-2},
row 4) and directly leads to the below canvassing of nucleon transmutation,
to which we will prefix a view on the $\textrm{c}_{\textrm{up}}$
and $\textrm{c}_{\textrm{down}}$ requests of the neutron $-$ in
the summarized form given in Table \ref{tab:Structure-of-nucleons-1-2-1}:
\begin{table}[H]
\caption{\label{tab:Structure-of-nucleons-1-2-1}Structure of neutrons in multiplet
form 3,7,15 and 31 }
\bigskip{}

\qquad{}\qquad{}%
\begin{tabular}{|c|c|c|}
\hline 
content & \multicolumn{2}{c|}{charged multiplets}\tabularnewline
\hline 
\hline 
$\mathtt{u},\mathtt{\mathrm{2}d}$ $\!\!\!\!\!\!\!\!\!\begin{array}{c}
\\
\\
\\
\end{array}$ & $\begin{array}{c}
[\left.146(\nicefrac{2}{3},\nicefrac{2}{3},\nicefrac{-1}{3})\right]_{\mathtt{u}}^{438}\end{array}$ & $\left[640(\nicefrac{2}{3},\nicefrac{-1}{3},\nicefrac{-1}{3})\right]_{\mathtt{2d}}^{1920}$\tabularnewline
\hline 
$\mathtt{u},\mathtt{\mathrm{2}d}$ $\!\!\!\!\!\!\!\!\!\begin{array}{c}
\\
\\
\end{array}$ & $\begin{array}{r}
\left(438+\epsilon\right)\nmid7\\
{\scriptstyle {\textstyle \epsilon}\:\in\left\{ -1,0,1\right\} }
\end{array}$ & $\begin{array}{r}
\left(1920+\epsilon\right)\nmid7\\
{\scriptstyle {\textstyle \epsilon\:}\in\left\{ -1,0,1\right\} }
\end{array}$\tabularnewline
\hline 
total $\mathtt{n}\!\!\!\!\!\!\!\!\!\begin{array}{c}
\\
\\
\\
\end{array}$ & \multicolumn{2}{c|}{$\begin{array}{c}
\\
\left.157(\underbrace{\nicefrac{4}{315},\nicefrac{\ldots,4}{315}},\underbrace{\nicefrac{-11}{315},\nicefrac{\ldots,-11}{315}})\right]_{\mathtt{n}}^{2355}\\
\;{\scriptscriptstyle 11}\qquad\qquad\qquad\:{\scriptscriptstyle 4}\qquad\quad
\end{array}$}\tabularnewline
\hline 
total $\mathtt{n}\!\!\!\!\!\!\!\!\!\begin{array}{c}
\\
\\
\\
\end{array}$ & \multicolumn{2}{c|}{$\begin{array}{c}
\\
\left.76(\underbrace{\nicefrac{5}{9765},\nicefrac{\ldots,5}{9765}},\underbrace{\nicefrac{-26}{9765},\nicefrac{\ldots,-26}{9765}})\right]_{\mathtt{n}}^{2356}\\
{\scriptscriptstyle 26}\qquad\qquad\qquad\quad{\scriptscriptstyle 5}\qquad\qquad
\end{array}$}\tabularnewline
\hline 
\end{tabular}
\end{table}

\noindent As, in terms of centerpiece-free hypersphere packings,
the total $\mathtt{n}$ content is 2358, we find that for one $L_{\textrm{up/down}}+\epsilon\:({\textstyle \epsilon}\:\in\left\{ -1,0,1\right\} )$
no more than four counteracting packings lacking both the centerpiece
and a peripheral hypersphere are needed in row 3, and no more than
three in row 4. For two $L_{\textrm{up/down}}$'s with a positive
epsilon, row 3 would cease to work. What doesn't work in the first
place is case row 2.\newpage{}

\subsection{\label{sub:Beta}Beta decay }

We are nearing a position to discuss spontaneous $\beta$-decay. In
order that preon levels $f_{\nu>n_{0}+1}$ get activated, some constituents
$L_{m}\!\in\!\left\{ L_{\textrm{up}}\right\} \cup\left\{ L_{\textrm{down}}\right\} $
have to give way to $L_{m}^{\pm}$. We gather from Tables \ref{tab:Structure-of-nucleons}
and \ref{tab:Structure-of-nucleons-1-2-1} that
for row 2), a neutron excited to preon charge level $f_{n_{0}+3}\!=\!7$
loses realizability and undergoes a transmutation. The closest-to-normal,
yet under the shape condition
\[
\mathtt{n}\sim L_{11}\!+\!\epsilon_{11}+2\left(\Sigma_{m=5}^{6}(L_{m}\!+\!\epsilon_{m})+\Sigma_{m=7}^{10}(L_{m}\!+\!\epsilon_{m})\right)\qquad(\epsilon_{m}\in\{-1,0,1\})
\]
unrealizable $\mathtt{u},\!\mathtt{\mathrm{2}d}$ content would be $63\!\times\!7+274\!\times\!7\!=\!337\!\times\!7\!=\!2359$.
Realizability is afforded in terms of a proton under $f_{n_{0}+3}\!=\!7$
with $\mathtt{u}\!+\!\mathtt{u'}\!,\mathtt{d}$ content $\mathrm{125\!\times\!7\!+\!137\!\times\!7\!=\!875\!+\!959=1834}$.
The shape condition satisfied here is 
\[
\mathtt{p}\sim(L_{11}\!+\!\epsilon_{11}+L'_{11}\!+\!\epsilon_{11}')+\Sigma_{m=5}^{6}(L_{m}\!+\!\epsilon_{m})+\Sigma_{m=7}^{10}(L_{m}\!+\!\epsilon_{m})\quad(\epsilon_{m}\in\{-1,0,1\}).
\]
Alongside the transmutation comes a difference $525$, split in familiar
manner into an antineutrino- and electron part:
\[
\begin{array}{cccccc}
2359\rightarrow & 1834 & + & 526 & - & 1.\\
\mathtt{n} & \mathtt{p} &  & \bar{\nu} &  & e^{{\scriptscriptstyle -}}
\end{array}
\]
\noindent For row 4), a proton excited to preon charge level $f_{n_{0}+5}=\!31$
similarly loses realizability and undergoes the reverse transmutation
proton $\rightarrow$ neutron. The enormous, unrealizable total content
amounts to $2883$. Realizability is provided in terms of a neutron
under $f_{n_{0}+5}=\!31$, requiring a moderate total content $2356$;
the transmutation is accompanied by a difference $527$, split into
a neutrino- and positron part: 
\[
\begin{array}{cccccc}
2883\rightarrow & 2356 & + & 526 & + & 1.\\
\mathtt{p} & \mathtt{n} &  & \nu &  & e^{{\scriptscriptstyle +}}
\end{array}
\]

\noindent Once again, we find a juxtaposition $x$ vs. $x-2$: the
second reaction occurs under $f_{n_{0}+5}$, the first under $f_{n_{0}+3}.$
Note also that the lepton companion remains implicit via a $\mp1$
charge correction: the electron imparts a unit decrement on the antineutrino,
the positron a unit increment on the neutrino. We will come back to the juxtaposition shortly,
but first want to reveal a pecularity of transmutations for collections of protons or neutrons.
When a septet regime becomes dominant, with mostly  {\tt d}  content   $959$, admixtures of triple-based {\tt d}  content  $960$ are tolerated,  but only in delicately balanced quota form: Three 959s permit one 960 to join in, the tolereance being
\[ \frac{959+959+959+960}{4\times7}=137.03571... \leq \textrm{\textgreek{a}}{^{-1}}.
\]
And three hundred eighty-three 959s go well together with one hundred twenty-nine 960s:
\[ \frac{383\times959+129\times960}{512\times7}=137.0359933... \leq \textrm{\textgreek{a}}{^{-1}}
\]
or generally,
\[ \frac{(2^{7t+2}-\sum_{r=0}^{t}{2^{7r}})\times959+\sum_{r=0}^{t}{2^{7r}}\times960}{2^{7t+2}\times7} \le \textrm{\textgreek{a}}{^{-1}}.
\]
Swap  $2^{7t+2}$ for  $2^{7t+3}$ and you get an almost identical formula for 959s tolerating 961s:
\[ \frac{(2^{7t+3}-\sum_{r=0}^{t}{2^{7r}})\times959+\sum_{r=0}^{t}{2^{7r}}\times961}{2^{7t+3}\times7} \le \textrm{\textgreek{a}}{^{-1}}.
\]

\noindent Put another way, sort of a  Mach principle for weak processes  would be responsible for keeping deviations from the inverse fine structure constant small.

\noindent It is known that there exist several neutrino flavors $-$
at least two beyond the electron neutrinos.  As we shall see, identifying these additional neutrino flavors may shed light on the very origin of leptons. \newline The successor relations
\[
\begin{array}{c}
\textrm{succ}{}_{-}(x)=2x+3,\\
\textrm{succ}{}_{0}(x)=2x+2,\\
\textrm{succ}{}_{+}(x)=2x+1,
\end{array}
\]
when applied to the start values 285,286,287, yield
\[
\begin{array}{c}
\left(\textrm{succ}{}_{-}\right)^{4}(285)=L_{16}+285,\\

\left(\textrm{succ}{}_{0}\right)^{4}(286)=L_{16}+286,\\

\left(\textrm{succ}{}_{+}\right)^{4}(287)=L_{16}+287.
\end{array}
\]

 \noindent The picture that forms is that the values 285,287 $-$ as well as their split form 286 when the lepton emerges $-$ may be thought
of as  neutrino family members settling relative to level 0 and the values 4605,4606,4607
as their cousins settling relative to level $L_{16}$. A middle one,
$L_{8}=240$, would then be responsible for the members 525,526,527.

\noindent They can, in the first place, be understood as an  `electroweak` phenomenon:
linear combinations of  $\frac{5}{8}\;p_l$ with  $\frac{9}{8}\;p_l$, for instance, yield
\[
\begin{array}{c}
-\frac{5}{8}\;p_2+\frac{9}{8}\;p_8=-1\frac{7}{8}+286\frac{7}{8}=285,\\
\\
-\frac{5}{8}\;p_2+\frac{9}{8}\;p_{12}=-1\frac{7}{8}+4606\frac{7}{8}=L_{16}+285,\\
\\
\frac{5}{8}(p_1+p_2)+\frac{9}{8}(-p_2+p_8)=\frac{5}{8}+1\frac{7}{8}-3\frac{3}{8}+286\frac{7}{8}=286,\\
\\
\frac{5}{8}(p_1+p_2)+\frac{9}{8}(-p_2+p_{12})=\frac{5}{8}+1\frac{7}{8}-3\frac{3}{8}+4606\frac{7}{8}=L_{16}+286,\\
\\
-\frac{5}{8}p_3+\frac{9}{8}(-p_2+p_3+p_8)=-4\frac{3}{8}-3\frac{3}{8}+7\frac{7}{8}+286\frac{7}{8}=287,\\
\\
-\frac{5}{8}p_3+\frac{9}{8}(-p_2+p_3+p_{12})=-4\frac{3}{8}-3\frac{3}{8}+7\frac{7}{8}+4606\frac{7}{8}=L_{16}+287.\\
\end{array}
\]
According to QFD, however, nucleon transmutation is just a weak process, so, as far as 525,526,527 are concerned, the weak part of  quasi-supersymmetric relations should suffice for their representation. While quasi-supersymmetric partners of the electromagnetic sector are separated by a residual `energy' $\frac{5}{2}\hbar$, that residue    becomes $\frac{9}{2}\hbar$ in the weak sector. (The apostrophes indicate that, in the absence of a definition of a time unit, those quantities arent't but preforms of energy packets.) The idea now is to look whether neutrinos or electrons crop up as resonances of that particular  residue. Individually, we observe
\[
\begin{array}{c}\frac{2}{9}\times285=63\frac{3}{9},\qquad \frac{2}{9}\times286=63\frac{5}{9}, \qquad\frac{2}{9}\times287=63\frac{7}{9},\\
\\
\frac{2}{9}\times525=116\frac{6}{9},\qquad\frac{2}{9}\times526=116\frac{8}{9},\qquad\frac{2}{9}\times527= 117\frac{1}{9},\\
\\
\frac{2}{9}\times4605=1023\frac{3}{9},\qquad\frac{2}{9}\times4606=1023\frac{5}{9}, \qquad\frac{2}{9}\times4607=1023\frac{7}{9}.
\end{array}
\]
Yet, together they   smoothly yield the {\tt d} content 960\,:

\noindent \smallskip{}

 \includegraphics[scale=0.66]{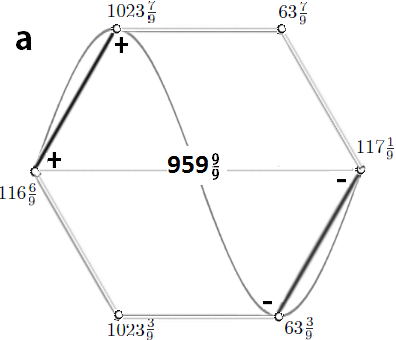}\quad\includegraphics[scale=0.66]{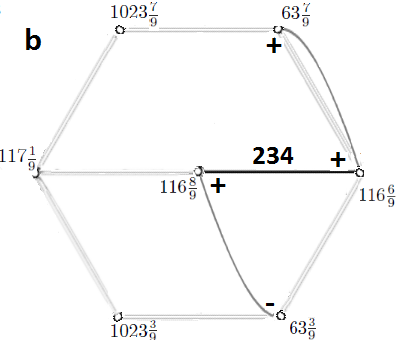}

 \bigskip{}

\noindent The {\tt d}-content transition $960\rightarrow959$ (signaling the transition from  a triple-based  to a septet-based regime, $f_{n_0+2} \rightarrow f_{n_0+3}$, in beta decay) is obviously orchestrated by  {\em all} members of the neutrino family  (diagram a). Neither 960 nor 959, however, are resonances (integers divisible by  residue 9/2 without rest). The picture doesn't get complete until we take note of   the antineutrino involved in the process. In diagram b,
  we can see  what's going on:  the antineutrino, assisted by  285,287, assumes its split form 526 after the electron is  expelled by the nucleon, $525\rightarrow526+e^-$ $-$ all of this expressed in the language of quasi-supersymmetric residual energy where the characteristic resonance     $(116\frac{6}{9}+116\frac{8}{9})+(63\frac{7}{9}-63\frac{3}{9})=234$ is used. 
 
\noindent The {\tt d}-content transition $960\rightarrow961$, which is the result of the transition $f_{n_0+2} \rightarrow f_{n_0+5}$ in fusion, engages all neutrino family members as well (diagram c). Since electron and antielectron are  the same except for charge, the resonance, indicating  the neutrinos's split form
 $527\rightarrow526+e^+$, must be the same. By contrast, though, it's a standalone now  (diagram d) $-$ the contributions of the other family members just cancel each other out: $117\frac{1}{9}+116\frac{8}{9}\pm(1023\frac{7}{9}-1023\frac{3}{9})\mp(63\frac{7}{9}-63\frac{3}{9})=234$. 

 \bigskip{}

 \includegraphics[scale=0.66]{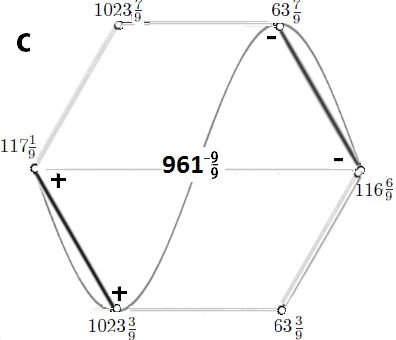} \quad\includegraphics[scale=0.66]{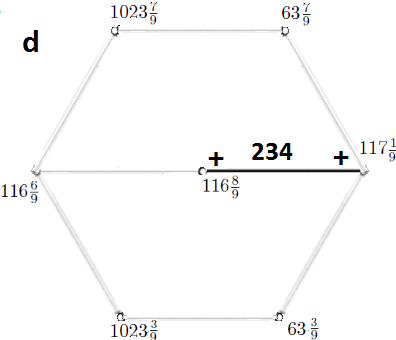}
 \bigskip{}

\noindent We may behold that electrons and antielectrons seem to owe their existence to a resonance of quasi-supersymmetric residual energy.

\noindent Where should we look for a link between neutrinos and Mersenne fluctuations?
 A preliminary inspection of the  $\between^{(n)}_{\alpha}$ incidences of non-split level-0 and level-$L_8$ antineutrinos/neutrinos -- 285,287 and 525,527 --
yields a rather uniform behavior. In order to make them  comparable, we chose a setting    that allows the numbers of peak incidences of 285 and of 525 to coincide $-$ which occurs for a level-0 flavor range of  $7095\!\leq\! n\!\leq13996$ vs. level-$L_8$ flavor range of  $2902\leq n\leq13996$. (Inspection of incidences of the  level-$L_{16}$ flavor, 4605,4607, is beyond the scope of present-day online calculators.) The legend to the table is as follows: an {\tt\bf l} stands  for left-leg incidences, divided into those with failed   and those with successful right-leg captures (underlined), {\tt\bf p} for peak incidences, {\tt\bf r} for right-leg incidences, again divided into those deriving from false left-legs (284 or 288, 524 or 528)  and proper ones (286, 526) and {\tt\bf lr} for left-legs  that make it to the right leg  with their identity preserved. 

\bigskip{}

\qquad{}\qquad{}\qquad{}%
\begin{tabular}{|c|c|c|c|}
\hline 
{\tt\bf l} & {\tt\bf p} & {\tt\bf r} & {\tt\bf lr}  \tabularnewline
\hline\hline
 \multicolumn{4}{|c|}{$I(\,\between_{\alpha}^{(n)}=285)\quad(7095\leq n\leq13996)$} \tabularnewline

$13=8+{\underline5}$  &21  &$8={\underline4}+4$ &23 \tabularnewline
$20\%$ &$32\%$ &$13\%$ &$35\%$ \tabularnewline
\hline \hline 
 \multicolumn{4}{|c|}{$I(\,\between_{\alpha}^{(n)}=287)\quad(7095\leq n\leq13996)$} \tabularnewline

$13=7+{\underline6}$  &25  &$8={\underline8}+0$ &28 \tabularnewline
$18\%$  & $33\%$  & $11\%$ &  $38\%$  \tabularnewline
\hline \hline 
 \multicolumn{4}{|c|}{$I(\,\between_{\alpha}^{(n)}=525)\quad(2902\leq n\leq13996)$} \tabularnewline
$3=2+{\underline1}$  &21  &$4={\underline2}+2$ &14 \tabularnewline
$7\%$  &$50\%$   & $10\%$   &$33\%$    \tabularnewline
\hline \hline 
 \multicolumn{4}{|c|}{$I(\,\between_{\alpha}^{(n)}=527)\quad(2902\leq n\leq13996)$} \tabularnewline
 
$4=3+{\underline1}$  &19  &$3={\underline0}+3$ &12 \tabularnewline
$11\%$     & $50\%$      & $8\%$     & $31\%$     \tabularnewline
\hline 
\end{tabular}

\bigskip

\noindent  With due caution, one may conclude  (see {\tt\bf lr} column) that only about a third of  (anti-)neutrinos survive  the cruise (in expansion time)  from  left leg  to  right leg, which is in line with the empirical evidence of  neutrino oscillations:

 \bigskip{}

 \includegraphics[scale=0.66]{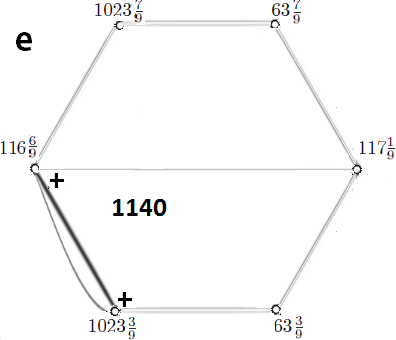} \quad\includegraphics[scale=0.66]{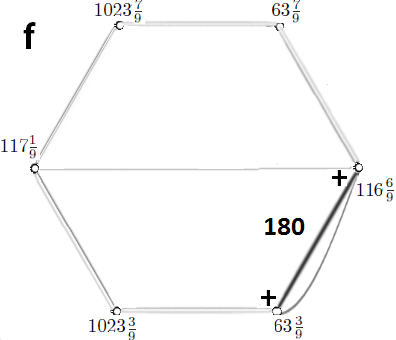}

\smallskip{}
 \includegraphics[scale=0.66]{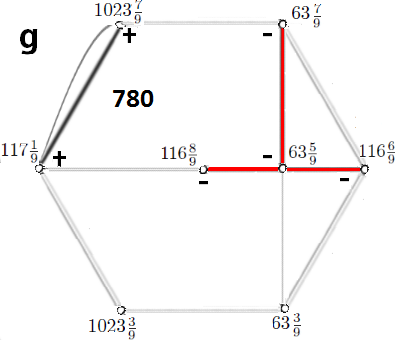} \quad\includegraphics[scale=0.66]{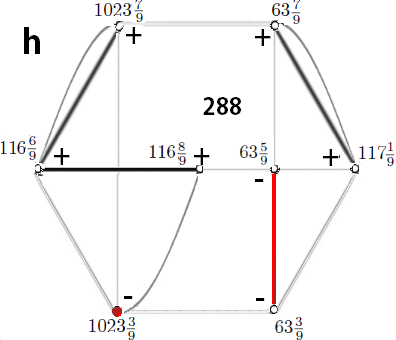}

 \smallskip{}

\noindent The oscillations $4605\leftrightarrow525$ and $4607\leftrightarrow527$  use the  non-resonances 1140 and 780  (diagrams e and g). In contrast, oscillations $285\leftrightarrow525$ and $287\leftrightarrow527$ use the resonances $234\!\mp\!54$ (diagrams f and h).   In  e and g, the diagram-f resonance  $234\!-\!54$   occrs atop or beneath 960 (our {\tt d} content), which underpins the unity of neutrinos. While oscillations in e and f occur independently (irrespective of (anti-)neutrino capture on other levels), (anti-)neutrino capture is a must for oscillations occurring in g and h. Due to that, even though incidences of  4605 and 4607 are left implicit, they  make themselves felt  by rollbacks  (red lines): the  oscillations in question $-$ one  $4607\rightarrow527$, the other $287\rightarrow527$ $-$ entail  $(-116\frac{8}{9}-116\frac{6}{9})+(116\frac{6}{9}+116\frac{8}{9})=0$, as well as rollbacks in the form of one fewer capture  of level-0 neutrino 287 (diagram g) and antineutrino 285 (diagram h) each $-$ a behavior that could be misunderstood  as CP violation if it weren't for  one further rollback in 4605 (red point). 
For instance, in the above table's {\tt\bf p} column, the sum of neutrino  peak incidences, $44=25+19$, could be changed to  $46=23+23$ via oscillation and  capture rollback $-$ as shown in blue in the table below  $-$   to  bring the mean value of the {\tt\bf lr} incidences'  share in  all incidences of   285,\,287,\,525 and 527  closer to one third $-$ 1.378422/4 vs.\,1.381347/4 (albeit on an almost doubled standard deviation, 0.04789 vs.\,0.0269):
$\!\!$\footnote%
{$\;$incidences incorporating oscillations  regularly may in turn  be guided by  spread avoidance}
\bigskip{}

\qquad{}\qquad{}\qquad{}%
\begin{tabular}{|c|c|c|c|}
\hline 
{\tt\bf l} & {\tt\bf p} & {\tt\bf r} & {\tt\bf lr}  \tabularnewline
\hline\hline
 \multicolumn{4}{|c|}{$I(\,\between_{\alpha}^{(n)}=285)\quad(7095\leq n\leq13996)$} \tabularnewline

$11=8+\textcolor{blue}{\underline3}$  &23  &$7={\underline4}+3$ &23 \tabularnewline
$17\%$ &$36\%$ &$11\%$ &$36\%$ \tabularnewline
\hline \hline 
 \multicolumn{4}{|c|}{$I(\,\between_{\alpha}^{(n)}=287)\quad(7095\leq n\leq13996)$} \tabularnewline

$11=7+\textcolor{blue}{\underline4}$  &\textcolor{blue}2\textcolor{blue}3  &$8={\underline8}+0$ &28 \tabularnewline
$16\%$  & $33\%$  & $11\%$ &  $40\%$  \tabularnewline
\hline \hline 
 \multicolumn{4}{|c|}{$I(\,\between_{\alpha}^{(n)}=525)\quad(2902\leq n\leq13996)$} \tabularnewline
$3=2+{\underline1}$  &21  &$4={\underline2}+2$ &14 \tabularnewline
$7\%$  &$50\%$   & $10\%$   &$33\%$    \tabularnewline
\hline \hline 
 \multicolumn{4}{|c|}{$I(\,\between_{\alpha}^{(n)}=527)\quad(2902\leq n\leq13996)$} \tabularnewline
 
$4=3+{\underline1}$  &\textcolor{blue}2\textcolor{blue}3  &$3={\underline0}+3$ &12 \tabularnewline
$10\%$     & $55\%$      & $7\%$     & $28\%$     \tabularnewline
\hline 
\end{tabular}---

\bigskip

\noindent The  link between neutrino capture and -oscillations opens up the possibility of a quasisupersymmetric way to express the `weak' Mach principle  $-$ meaning that in septet-dominated collections of protons or neutrons, triple-based admixtures are tolerated in quota form).  For  this to become an option, an integer multiple of {\small 9/2} divisible by 7 must be at hand that can, in the sequel, be juxtaposed to neutrino capture-/oscillation-related integer multiples $- $ which exhibit divisibility   by 3 throughout. One, the smallest, with divisibiliy 7 lies  in the interior of the hexagon (diagram {\tt i}):  $1204=63\frac{5}{9}+116\frac{8}{9}+1023\frac{5}{9}$; the other, of double  value, forms the hexagon's periphery (diagram j):

\bigskip{}
\includegraphics[scale=0.66]{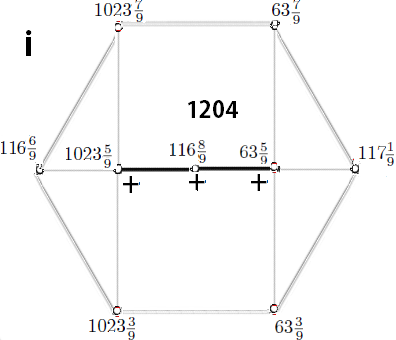}\includegraphics[scale=0.66]{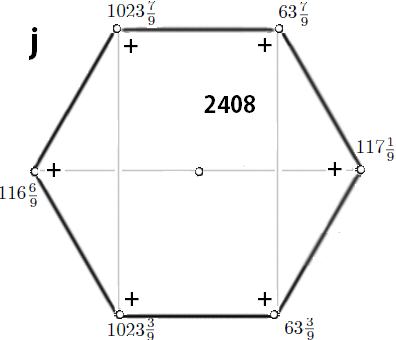}

 \bigskip{}

\noindent One further condition is a duality construction linking nucleon transmutation  to the hypothetisized `weak' Mach principle in quasisupersymmetric form. Parity considerations lead us to separate  integer parts from fraction parts in  the mixed fractions  in diagrams a and c $-$ so that the former represent {\tt d} content and the latter  positron (electron) charge $q\,$:
\[
\phi(a,b):= 960+(a-b)+\epsilon_{ab}\,q.
\]
Thus, $\phi(116,117)=959+\frac{9}{9}$, $\phi(117,116)=961+\frac{\!\!-9}{9}$.
The `weak' Mach principle's  function, on the other hand, can be formulated in form of integer multiples (of {\small 9/2})  with divisibilities 7 first and  3 second:
\[ \psi^{(t)}(a,b):= \frac{\mid\! a-b\!\mid(2^{7t+2}-\sum_{r=0}^{t}{2^{7r}})\, a+\sum_{r=0}^{t}{2^{7r}}\, b}{\mid\! a-b\!\mid2^{7t+2}\times7
}.
\]
Now the duality consists in  decimal places (instead of fraction parts) that are to be  separated from integer parts: 
 \[
\psi^{(t)}(959,960)=137+0.035.., 
\]
\[
 \psi^{(t)}(1204,960)=172-0.035.. , \qquad \psi^{(t)}(2408,960)=344-0.035.. ,
\]
as well as in integer parts realizing the  `weak' Mach principle's  inherent $7t$ vs. $\!(7t+2)$ juxtaposition in a subtler, $t\,$-independent way:
\[
T_3:\sum_{\nu=1}^3 \left\lceil\sigma_\nu\right\rceil\mapsto137, \qquad T_3+\left\lceil\sigma_4\right\rceil=172\,;
\]
\[ 
T_5:\left\lceil\sigma_5\right\rceil\mapsto137, \qquad\sum_{\nu=1}^4 \left\lceil\sigma_\nu\right\rceil+T_5+\left\lceil\sigma_6\right\rceil=344.
\]
These  are conserved quantities:  substitute $b\in\left\{180,234,288,780,1140\right\}$ for 960 in $\psi^{(t)}$ and nothing changes $-$ the important thing is $7\vert a,\;3\vert b\;$ and $b<1204$ or $2408$. In other words, neutrino oscillations obey $\psi^{(t)}\,$-invariance.

\noindent We now show how  transformation of time-like into space-like  refinement informs neutrinos with identity $-$ an  identity that in expansion time figures as width. Checking  Mersenne fluctuations for incidences of  $L_{16}=4320$ provides a first hint. It turns out that, under state-of-the-art computing conditions, there's only one incidence of $\between_{\alpha}^{(n)}=4320$, namely  for  $n=14051$, a prime number. Perhaps, then, for transformation between two kinds of refinement to take place, incidences  have to be in a special category of Mersenne fluctuations defined by
\[
\between^{(\pi_r^{(1)})}_\alpha=L _m\quad\left(\between^{(\pi_r^{(1)}+2)}_{\alpha'}=L_m^-,\; m\in\left\{\textrm{\small8,16}\right\}\right).
\]
Where $\pi_{r}$ means  $r$\hspace{1pt}{\small th} prime, the space-like refinement found for the level-$L_8$ incidence  would be similar in size to the index $r$ indicating  the time-like refinement of the level-$L_{16}$ incidence. 
 We   find $\between_{526}^{(14051)}=4320$, $\between_{536}^{(14053)}=4319$ (see table left, below).  
14051  is equal the prime-indexed prime $\pi_{1657}^{(1)}$. Under the recursive definition $\pi_r^{(1)}:=\pi_r$, $\pi_r^{(k)}=\pi_{\pi_r^{(k-1)}}$, it has a follower, $1657=\pi_{160}^{(2)}$. The new index is resolvable into   factors, $160=2^5\times5$, so the recursion stops at $k=2$. A  $\between$ with prime $n$  occurs  more frequently for level  $L_8=240$. Among those incidences, one  stands out (see table  right): $\between_{1605}^{(3739)}=240$, $\between_{1583}^{(3741)}=239$. Since 3739 is equal $\pi_{522}^{(1)}$, the Mersenne fluctuation belongs in the same category, albeit as one short of recursion:  522 is composite, $2\times3^2\times29$, so  $k=1$.
The relationship
$\pi_{r}^{(k)}\sim r\left(\log r\right)^k$ suggests  that the quantity $k$ might matter. The situation is this: On the one hand, there's 3739's relatively advanced space-like refinement $\beta=1605$, seeking  an appropriately advanced  time-like refinement $\pi_{1657}^{(1)}=14051$; on the other hand, the  lesser time-like refinement $3739=\pi_{522}^{(1)}$ is after an appropriately lesser  space-like refinement, finding it in 14051's $\alpha=526$. In order to make ends meet, the transformation of prime index into space-like refinement and vice versa needs   a joint   $-$ 52 or 4. One that turns out to be controlled by a  juxtaposition in upper bounds of sum: where $k=2$, $k\sum_{l=1}^4p_l=52$; where $k=1$, $k\sum_{l=1}^2 p_l=4$. 
Given the  two $n$'s' compliance with the category, the transformation could not happen without the space-like refinements' active assistance $- $ it's all about their (locally)  extremal values! There's only one prime $n$ in the relevant part of the Mersenne fluctuation on the  right, and   space-like refinement assumes its  maximum right there: $\beta=1605$. By contrast,
the relevant part of the Mersenne fluctuation  on the left has two prime $n$'s: 14051 and 14057, associated with the space-like refinements $\alpha=526$ and $560$, respectively.  Yet, what gets singled out is the  minimum, $\alpha=526$ (figuring, in expansion time, we recall, as the width of an unreferenced electron/antielectron neutrino):
 \bigskip{}

\label{pagepoin}
\quad{}%
\begin{tabular}{|c|c|c|c|}
\hline 
$n$ & $n$ prime & $\alpha$ & $\between_{\alpha}^{(n)}$\tabularnewline
\hline 
\hline 
$\vdots$ &  &  & \tabularnewline

14046 &  & 548 & 134\tabularnewline
\hline 
14047 &  & 554 & 269\tabularnewline
\hline 
14048 &  & 584 & 539\tabularnewline
\hline 
14049 &  & 560 & 1079\tabularnewline
\hline 
14050 &  & 528 & 2159\tabularnewline
\hline 
14051 & \Checkmark{} & {\bf526} & {\bf4320}\tabularnewline
\hline 
14052 &  & 532 & 8640\tabularnewline
\hline 
14053 &  & 536 & 4319\tabularnewline
\hline 
14054 &  & 532 & 2159\tabularnewline
\hline 
14055 &  & 550 & 1079\tabularnewline
\hline 
14056 &  & 534 & 539\tabularnewline
\hline 
14057 & \Checkmark{} & 560 & 269\tabularnewline
\hline 
14058 &  & 538 & 134\tabularnewline

$\vdots$ &  &  & \tabularnewline
\hline 
\end{tabular}
\begin{tabular}{|c|c|c|c|}
\hline 
$n$ & $n$ prime & $\beta$ & $\between_{\beta}^{(n)}$\tabularnewline
\hline 
\hline 
  &  &  & \tabularnewline
 
 $\vdots$ &  &  & \tabularnewline
3735 &  & 1537 & 14 \tabularnewline
\hline
3736 &  & 1559 & 29\tabularnewline
\hline 
3737 &  & 1597 & 60\tabularnewline
\hline 
3738 &  & 1595 & 120\tabularnewline
\hline 
3739 & \Checkmark{} & {\bf1605} & {\bf240}\tabularnewline
\hline 
3740 &  & 1593 & 480\tabularnewline
\hline 
3741 &  & 1583 & 239\tabularnewline
\hline 
3742 &  & 1595 & 119\tabularnewline
\hline 
3743 &  & 1599 & 59\tabularnewline
\hline 
3744 &  & 1531 & 29\tabularnewline
\hline
3745 & & 1533 &  14\tabularnewline
 $\vdots$ &  &  & \tabularnewline

& &  & \tabularnewline
\hline 
\end{tabular}
\bigskip{}

\noindent  The  conditions  for the  correlation  to imbue the unreferenced $\textrm {e}^{\pm}$ neutrino with identity can probably be generalized: 
\[
r_1\sim\beta,\quad r_2\sim \alpha\quad \left(\alpha \in\{\textrm{\small 286,526,4606}\}\right).
\]  
\noindent For the time being, checking till $\alpha=4606$ or for an $n$ $\sim$ $\pi_{4606}=44263$ is hard, or simply not yet computable. We  go ahead checking instead how neutrinos, once  they have come into existence, perform in expansion time. To this end, we demand that the previously defined  successor relations  be
themselves translated into Mersenne fluctuations (Eq. (\ref{eq:left-leg})).
 One can demonstrate
this with generic Mersenne fluctuations 
\[\textrm{mf}\:\left(\nu=n-r,\ldots,n+r;\, r<n,\, r\in\mathbb{N}_{0}\right)
\]
for homogeneous $\epsilon$. For $\epsilon=1$ and $\epsilon=0$,
we respectively get 
\[
\begin{array}{c}
\textrm{succ}{}_{0}\quad\mapsto\quad\textrm{mf leftleg}\quad(\ldots,286,574,1150,2302,\boldsymbol{4606},\ldots),\\

\textrm{succ}{}_{+}\quad\mapsto\quad\textrm{mf leftleg}\quad(\ldots,287,575,1151,2303,\boldsymbol{4607},\ldots).
\end{array}
\]
Interestingly, though, when the first successor relation, $\textrm{succ}{}_{-}$,
is translated into the remaining mf leftleg progression (with homogeneous
$\epsilon=-1$), we get 
\[
\textrm{succ}{}_{_{{\textstyle -}}}\quad\mapsto\quad\textrm{mf leftleg}\quad(\ldots,285,570,1140,2280,\boldsymbol{4560}(=L_{16}\!+\! L_{8}),\ldots).
\]
 This could be interpreted as follows. On the rare occasion a Mersenne
fluctuation homogeneous in $\epsilon=1(0)$ hits 4606(7), we have
(in expansion time) an oscillation of a neutrino settled relative
to 0 to one settled relative to $L_{16}$. On the likewise rare occasion
a Mersenne fluctuation homogeneous in $\epsilon=-1$ hits 4560, the
middle level $L_{8}=240$ is established for 525/526/527 $-$ the
familiar unreferenced electron/antielectron neutrino. 
\noindent This  gives room to  an   auto-correlative category  which we  want to explore first, vital as it is  to a genuine interpretation of  neutrinos' performance in expansion time:
\[
\between^{(2\pi_r)}_\alpha=L _{16}+L_{8},\quad\between^{(2\pi_r+2)}_{\beta}=L_{16}+L_{8}.
\]
How the correlation  works out  can be probed on the case $\between^{(9634)}_{575}=\,\between^{(9636)}_{589}=4560$. A first novelty is that $k$ is no longer of importance $-$ we deal with factors only. 9634 evolves into $2\times4817$ and 9636 into  $2^2\times3\times11\times73$. After accepting (largest) factors instead of  `standalone' primes as  assignments to $r$,   dissimilarities  catch the eye:   $\vert73-4817\vert\nsim\vert589-575\vert$ (where similarity is defined by  the relation  $\left\vert r_b-r_a \right\vert\sim\left\vert \beta-\alpha \right\vert \,$). Obviously an $r_a$ that's    more similar to  $r_b$ than is $\pi_r$   is needed. The  way to a solution leads down the index chain: 4817 is  the 649{\small\hspace{0.8pt}th} prime. $649$  has an  evolution on its own, $11\times59$, so we get a new  largest factor that we can  assign, $r_a:=59$. Astonishingly, the result is not merely nighness $-$ it's a perfect match: $\left\vert r_b-r_a \right\vert=\left\vert \beta-\alpha \right\vert \,$.  Note that the space-like refinement $\alpha=575$ figures as width in expansion time $-$ as  the successor of width 287 according to $\textrm{succ}_+$. So the auto-correlation is one between a  positively-charged-level-0-lepton neutrino 287 that is about to morph into an antielectron neutrino 527 or vice versa. The first part of the transformation  sees  the light: 
\[ -{r_a}+\alpha+\sum_{l=1}^3p_l=-{r_b}+\beta+\sum_{l=1}^3 p_l=527.
\]
 For the  second part, one has to consider the successorship $287\rightarrow575$  causes  two pairs $(r,\gamma)$ to emerge rather than  one: 
 \[(r_a+\alpha)+(-r_b-\beta)+\sum_{l=1}^3 p_l+\sum_{l=1}^5 p_l+\sum_{l=1}^7 p_l=287. 
 \]
 Hardly surprisingly, transforming  level 0 into $L_8$ and   $L_{16}$ into $L_8$, or vice versa, obeys a cognate set of rules: 
 \[(-r_b+\alpha-\sum_{l=1}^8 p_l)+L_8=L_8, \]\[
 L_{16}-(-r_b+r_a+\sum_{l=1}^{11} p_l+\sum_{l=1}^{3} p_l)=L_{16}-(-\beta+\alpha+\sum_{l=1}^{11} p_l+\sum_{l=1}^{3} p_l)=L_8.
 \]
  Even interpolation to  intermediate level becomes  possible:%
\footnote{$\;$ What's more, there are representations of {\tt d} quark constituents and of the {\tt u} quark  tightly knit to the second transformation rule  ``$\!\!$ two pairs $(r,\gamma)$ in case $527\mapsto287$\,'' $-$ with generally more sum terms involved:\\
 $(\beta-r_a-\alpha+r_b)-\sum_{l=1}^2p_l=L_4$, \\ $(-\beta+r_a+\alpha-r_b)+(\sum_{l=1}^5p_l+\sum_{l=1}^3p_l)=L_5$, \\ $(\beta+r_a-\alpha-r_b)+(\sum_{l=1}^5p_l+\sum_{l=1}^3p_l+\sum_{l=1}^2p_l)=L_6$,\\$(\beta+r_a+\alpha-r_b)-(\sum_{l=1}^8p_l+\sum_{l=1}^7p_l+\sum_{l=1}^6p_l+\sum_{l=1}^4p_l+\sum_{l=1}^3p_l+\sum_{l=1}^2p_l)=L_8$, \\$(\beta-r_a+\alpha+r_b)-(\sum_{l=1}^8p_l+\sum_{l=1}^7p_l+\sum_{l=1}^6p_l+\sum_{l=1}^4p_l+\sum_{l=1}^3p_l)=L_9$,\\$(\beta-r_a-\alpha+r_b)+(\sum_{l=1}^7p_l+\sum_{l=1}^5p_l+\sum_{l=1}^2p_l)=L_{10}$;\\$(-\beta+r_a+\alpha-r_b)+\sum_{j=1}^7\left(\sum_{l=1}^{j}p_l\right)=L_{11}$.}
 \[
 L_{16}-(r_a-\beta+\sum_{l=1}^{11} p_l+\sum_{l=1}^{3} p_l)=L_{12}.
 \].
\vspace*{4mm}

\noindent Generally, Mersenne fluctuations  inhomogeneous in $\epsilon$ will prevail, which  changes the interpretation
slightly. For owners of the WolframAlpha app, there's the chance
to check the (perhaps typical) behavior of the  Mersenne fluctuations' ensemble \,$\left(\between_{a,b,c}\right)_{\alpha_{n},\beta_{n},\gamma_{n}}^{(n)}$ between 12461 and 12469:

\medskip{}

\quad{}%
\begin{tabular}{>{\centering}p{4mm}>{\centering}p{3mm}>{\centering}p{3mm}>{\centering}p{5mm}>{\centering}p{5mm}>{\centering}p{6mm}>{\centering}p{5mm}>{\centering}p{5mm}>{\centering}p{4mm}>{\centering}p{4mm}>{\centering}p{5mm}}
 & {\tiny
12461} & {\tiny
12462} & {\tiny
12463} & {\tiny
12464} & {\tiny
12465} & {\tiny
12466} & {\tiny
12467} & {\tiny
12468} & {\tiny
12469} & \tabularnewline
 &  &  &  &  & $4606_{a}$, &  &  &  &  & \tabularnewline
 &  &  &  & $2303_{a}$, &  & $2302_{a}$, &  &  &  & \tabularnewline
 &  &  & $1151_{a}$, &  &  &  & $1150_{a}$, &  &  & \tabularnewline
 &  & $575_{a}$, &  &  &  &  &  & $575_{a}$, &  & \tabularnewline
$(\ldots$, & $287_{a}$, &  & $\,$$\,286_{b}$, &  &  &  & $\;\,285_{b},$  & $\,287_{c}$, & $\,287_{a}$, & $\,\ldots).$\tabularnewline
\end{tabular}
\bigskip{}

\noindent In this scenario, we recognize 1)  a level-0
neutrino $\left(\between_{a}\right)_{\alpha_{12461}}^{(12461)}=287$
that loses its association with a positively charged level-0 lepton at $n=12463$:
\[\left(\between_{b}\right)_{\beta_{12463}}^{(12463)}=286.
\]
The unreferenced neutrino then oscillates to 
\[\left(\between_{a}\right)_{\alpha_{12465}}^{(12465)}=4606(=L_{16}+286).
\]
 2) After an episode of mutual annihilation of the negatively and the positively charged level-0 leptons
 \[\left(\between_{b}\right)_{\beta_{12467}}^{(12467)}=285,
\]
\[\left(\between_{c}\right)_{\gamma_{12468}}^{(12468)}=287,
\]
the unreferenced level-$L_{16}$ neutrino oscillates back to level 0 and resumes association  with the abandoned positively charged level-0 lepton, thus preserving lepton number:
\[\left(\between_{a}\right)_{\alpha_{12469}}^{(12469)}=287.
\]

\noindent This interactive scenario requires some backing up, which  the underlying  correlations furnish aptly.  The unreferenced level-$L_{16}$ neutrino \[
(\between_a)_{182}^{(12465)}=4606
\]
 feels an affinity for what is its ancestor under homogeneous $\textrm{succ}_0$ conditions, the  level-0 neutrino 286; as the table below shows, this tendency finds its expression in the correlation with  
\[(\between_b)_{54}^{(12463)}=286,
\] 
assisted by a second  correlation between  
\[(\between_a)_{x}^{(y)}=287\quad (x,y)\in\left\{(\textrm{\small176,12461}),(\textrm{\small178,12469}\right\}\quad\textrm{and}\quad(\between_c)_{297}^{(12468)}=287.
\] 
That being a correlation with a third party, however, the interaction couldn't proceed without an agent  linking  c to b $-$ a mirroring of  the space-like refinement $\gamma_{12467}=285$  in  width  $(\between_b)_{52}^{(12467)}=285\;\textrm{(in bold)}$. The relevant transformations are:  $\!$\footnote{$\,\;$note the ever-recurring upper-bound-of-sum juxtaposition  $x$ vs. $\!\!\!$ $(x\pm2)$}$^,$\footnote{$\!\;$ the internal transformation a $\mapsto$ a is given by $(r_{12461}+\alpha_{12461})+(r_{12462}+\alpha_{12462})+\dots+(r_{12464}+\alpha_{12464})+(r_{12465}-\alpha_{12465})+(r_{12466}+\alpha_{12466})+\dots+(r_{12469}+\alpha_{12469})+\sum_{j=2}^9\left(\sum_{l=1}^j p_l\right)=4606$}

\noindent  (a $\mapsto$ b) 
\[(-r_{12463}+\beta_{12463})+(-r_{12465}+\alpha_{12465})+(r_{12467}+\gamma_{12467}) \;+\]
\[\sum_{l=1}^3p_l-\sum_{l=1}^2p_l+\sum_{l=1}^1p_l=286,
\]
(a  $\mapsto$ c)
 
 \[  (-r_{12461}+\alpha_{12461})\!+(-r_{12469}\!+\alpha_{12469})\!+(r_{12468}
 +\gamma_{12468})\;-\]\[\sum_{l=1}^5p_l\!-\sum_{l=1}^4p_l\!+\sum_{l=1}^3p_l=287.
 \]

\smallskip{}
\begin{tabular}{|c|r|c|c|c|} 
\hline 
$n$ & factorization  & a & b & c \tabularnewline
\hline 
\hline 
12461 & $17\times733$ & $\between_{176}^{(n)}=\fbox{\fbox{287}}$& & \\
  & {\small733}$=\pi_{130}$ & & & \tabularnewline
  \hline
  12462 & $2\times 3\times31\times67$ &$\between_{180}^{(n)}{\!\,\!}=575$ & & \tabularnewline
  \hline
  12463 & $11^2\times103$ & $\between_{190}^{(n)}{\!\,\!}=1151$ & $\between_{54}^{(n)}=\fbox{286}$  & \tabularnewline
  \hline
  12464 & $2^4\times19\times41$ & $\between_{174}^{(n)}{\!\,\!}=2303$ & $\between_{50}^{(n)}{\!\,\!}=572$ & \tabularnewline
  \hline
  12465 & $3^2\times5\times277$ & $\between_{182}^{(n)}=\fbox{4606}$ & $\between_{58}^{(n)}{\!\,\!}=1144$ &
\tabularnewline
 \hline
 12466 & $2\times23\times271$ & $\between_{200}^{(n)}{\!\,\!}=2302$ & $\between_{64}^{(n)}{\!\,\!}=571$ & \tabularnewline
 \hline
 12467 & $7\times13\times137$ & $\between_{180}^{(n)}{\!\,\!}=1150$ & $\between_{52}^{(n)}{\!\,\!}=\textrm{\bf285}$ & $\between_{\textrm{\bf285}}^{(n)}{\!\,\!}=143$ \tabularnewline
 \hline
 12468 & $2^2\times3\times1039$ &$\between_{182}^{(n)}{\!\,\!}=575$ &  & $\between_{297}^{(n)}=\fbox{\fbox{287}}$ \\
  & {\small1039}$=\pi_{175}$ & & & \tabularnewline
  \hline
  12469 & $37\times337$ & $\between_{178}^{(n)}=\fbox{\fbox{287}}$ & & $\between_{293}^{(n)}{\!\,\!}{\!\,\!}=143$ 
\tabularnewline
\hline
  & \multicolumn{4}{|l|}{\footnotesize \,\, }  \tabularnewline 
  & \multicolumn{4}{|l|}{\footnotesize if the largest factor $\gg L_{11}$, then its prime index rather than  the  factor}  \tabularnewline   

  & \multicolumn{4}{|l|}{\footnotesize itself is assigned to $r$}    \tabularnewline 
   \hline

\end{tabular}

\bigskip{}

\noindent The next 
Mersenne fluctuations' ensemble, $\left(\between_{a,b,c}\right)_{\alpha_{n},\beta_{n},\gamma_{n}}^{(n)}$,
 $n\in\left\{\textrm{\small 9630},{\scriptstyle \ldots},\textrm{\small9640}\right\}$, was chosen on purpose as in their midst reigns the `auto-correlative' category already discussed:

\bigskip{}
\hspace{-10mm}%
\begin{tabular}{c>{\centering}p{12mm}>{\centering}p{12mm}>{\centering}p{5mm}>{\centering}p{5mm}>{\centering}p{5mm}>{\centering}p{5mm}>{\centering}p{5mm}>{\centering}p{5mm}>{\centering}p{5mm}>{\centering}p{5mm}>{\centering}p{6mm}}
 & $\,\,${ \tiny
9630} & { \tiny
9631} & { \tiny
9632} & { \tiny
9633} & { \tiny
9634} & { \tiny
9635} & { \tiny
9636} & { \tiny
9637} & { \tiny
9638} & { \tiny
9639} & { \tiny
9640}\tabularnewline
 &  &  &  &  &  & $9121_{a}$, &  &  &  &  & \tabularnewline
 &  &  &  &  & $4560_{a}$, &  & $4560_{a}$, &  &  &  & \tabularnewline
 &  &  &  & $2279_{a}$, &  &  &  & $2279_{a}$, &  &  & \tabularnewline
 &  &  & $1139_{a}$, &  &  &  &  &  & $1139_{a}$, &  & \tabularnewline
 &  & $\underline{526}_{b}$,$569_{a}$ &  &  &  &  &  &  &  & $569_{a}$, & $\underline{524}_{c}$\tabularnewline
 & $({\scriptstyle \ldots},\underline{284}_{a}\quad$ &  &  &  &  &  &  &  &  &  & $\underline{284}_{a},{\scriptstyle \ldots})$\tabularnewline
\end{tabular}

\bigskip{}

\noindent In this ensemble, two levels, one 0, the other
240, are confronted with each other, the middle one (level
240) being established (in expansion time) at $n=9634$ and $n=9636$.
Particle pairs created in this interval can be 1) identical leptons,
one moving forward in time, the other backward; 2) antiparticle neutrinos
with their lepton partners, one pair moving forward in time, the other
backward. We will discuss case one first. There are two level-240
leptons (electrons) created at $n=9635$, one moving forward and transmuting
into a level-0 lepton at $n=9636$, the other moving backward and
incurring transmutation at $n=9634$. They continue their respective
paths symmetrically till they reach, after 10 steps in all, a scalariform
bottom layer that runs from $n=9630$ to $n=9640$. Where the bottom
was hit, a correction term $-1$ is deposited. Now comes case two.
The pair in question may be a level-0 lepton plus its antiparticle
neutrino. They hit the bottom, at $n=9630$ in backward mode and $n=9640$
in forward mode, depositing $-1$ and $286$ each time. Each partner
performs 5 steps till it reaches the bottom, that is 10 steps per
pair. Things get more complicated if, while ending up invariably as
an electron, which kind of lepton is created depends on the mode.
In backward mode, it is created rightmost as an electron at $n=9635$,
transmutes into a level-0 lepton at $n=9634$, bounces back to reacquire
electron status at $n=9635$ and continues its path toward the bottom,
hitting it after 7 steps at $n=9640$ and depositing there $-1$;
the antielectron neutrino is created leftmost at $n=9634$; it hits
the bottom after $10-7$ steps, that is at $n=9631$, and deposits
526 there. In forward mode, the lepton is created leftmost as a level-0
lepton at $n=9634$; it then transmutes into an electron at $n=9635$,
continuing its path till it reaches the bottom, which it hits after
6 steps at $n=9640$ and deposits $-1$ there; the antielectron neutrino
is created rightmost at $n=9636$; it hits the bottom after $10-6$
steps at $n=9640$ and deposits 526 there. We behold from the above:
a) there is an invariant number (10 in this case) which characterizes
both the number of steps and how many times the bottom is hit, and
b) what's shown at the bottom  by $\left(\between_{a,b,c}\right)_{\alpha_{n},\beta_{n},\gamma_{n}}^{(n)}$
is a superposition (underlined in the graph) of the possible outcomes:
$-1+(-1+286)$ at $n=9630,9640$; 526 at $n=9631$ and $-1+(-1+526)$
at $n=9640$.

\subsubsection{\label{sub:timecrys}Time crystals}

\noindent Historically, the model of three separate flavors of neutrinos had to be revised when it was discovered that they morph into one another (see prior discussion).  Neutrinos, consequently, weren't considered massless entities moving at light speed anymore $-$ their oscillations inextricably implied they possess mass, although in very tiny amounts. Had they been truly massless, time $-$ according to special relativity  $-$ would not flow and space not have any extension to them, nor could there be any time-like or space-like refinement. Alternatively, endowed with a tiny mass, they would fill a limited volume of time and space or more to the point, experience a finite number of time-like ({\em tr}) and space-like ({\em sr}) refinements like in a time crystal. If so, not only must we look for neutrino incidences  in  particulate form $\between^{()}$ but also take  account of how they originate in  $\gamma^{()}$. The relationship to bring these aspects together would be  a quadratic one in expansion time:  $\gamma^{()}$ at $n_0$, and $\between^{()}$ at   $n_0^2\left.{-\nu_1}\atop{+\nu_2}\right.$.

\noindent As the first table below shows, the mystery can  at least in part be solved for the electron antineutrino (a right-handed entity, according to current thinking): The  Mersenne fluctuations mf I and mf II exhibit a pecularity in that  their  otherwise pulsating space-like refinements  get stalled at   `islands  of  repose'  $-$  space-like correspondences of the `now' of $n_0$  whose locations are strictly choreographed and each of which spans two units: a `hillside,' and a `bayside' with upstream `anchorage.'   Of same height $h=11$, mf I and mf II are aligned  Mersenne fluctuations. Yet,  being of same height, or being aligned, is not a  necessary condition for Mersenne fluctuations to possess islands. As the second table   reveals, out of more spread-out  (but still quadratically related) fluctuations mf III and mf IV, one holds an  island as well. As for the aligned case,   the {\em sr}\,s of mf\,I  group themselves near $\tilde\alpha_{\textrm{\,I}}=1445$ and those of mf\,II near $\tilde\alpha_{\textrm{\,II}}=17$, the reason being  they follow the two factors of $n_0$, $f_1=5$ and $f_2=17$, such that
\[\textrm{}	\;\tilde\alpha_{\textrm{\,I}}= n_0^2/f_1,\quad \tilde\alpha_{\textrm{\,II}} =\tilde\alpha_{\textrm{\,I}}/n_0=f_2.
\]
Between $f_1$ and $f_2$, there lie the primes  7, 11 and 13. The unaligned case (second table) takes  account of the  remaining heights $-$ the aligned fluctuations mf I and mf II had $h=11$;   mf III and mf IV, on the other hand, come with $h'=7$ and  $h''=13$, respectively. The {\em sr}\,s of mf\,III  group themselves near 756 and the {\em sr}\,s of mf\,IV near 247:
\[\textrm{}	\;\tilde\alpha_{\textrm{\,III}}=756 =2^2\times3^3\times7,\quad \tilde\alpha_{\textrm{\,IV}} =247 = 13\times19.
\]For every Mersenne fluctuation, $\alpha$ has a  parity, and the number of islands for  mf\,I,...,\hspace{1pt}mf\,IV depends on whether $\tilde\alpha$ respects or violates that parity:
 \[(\tilde\alpha+\alpha)\; \textrm{odd} \Rightarrow \textrm{odd}\,\#\;\textrm{of islands},\quad(\tilde\alpha+\alpha)\; \textrm{even} \Rightarrow  \textrm{even}\,\#\;\textrm{of islands}.
 \]
 Checking  the tables for that we get one island each for mf I and mf IV   and two islands for mf II and none for mf III. A second observation concerns  the average behavior of  $sr\textrm{s}$: In mf II, $\alpha$\hspace{0.5pt}s  on average  are \hspace{-1pt}4 {\em sr}\hspace{1pt}s lower  than $\tilde\alpha_{\textrm{\hspace{-1pt}\,II}}$, and the $\alpha$\hspace{0.5pt}s of mf \hspace{-1pt}I even \hspace{-1pt}39 $sr\textrm{s}$ lower  than $\tilde\alpha_{\textrm{\,I}}$. And for  the second pair: mf IV's $\alpha$\hspace{0.5pt}s   1  {\em sr} lower than $\tilde\alpha_{\textrm{\,IV}}$, mf III's $\alpha$\hspace{0.5pt}s   8  {\em sr}\hspace{1pt}s lower than $\tilde\alpha_{\textrm{\,III}}$. More tellingly,  the
 ratio of deviations   $(\bar\alpha_X-\tilde\alpha_X)/(\bar\alpha_{X\!+1}-\tilde\alpha_{X\!+1})$  leads to two characteristic  numbers: $X=1$ to an approximation of  the relevant fluctuations' height minus one; $X=3$ to the number $2^3$ that, as will be elucidated further below, is key to integrating the results  of the second pair of fluctuations in the environment of mf I:
\[\left\lceil\frac{\bar\alpha_{\textrm{\,I}}-\tilde\alpha_{\textrm{\,I}}}{\bar\alpha_{\textrm{\,II}}-\tilde\alpha_{\textrm{\,II}}}\right\rceil=\left\lceil\frac{-39}{-4}\right\rceil=h-1=10
\]
and
\[\frac{\bar\alpha_{\textrm{\,III}}-\tilde\alpha_{\textrm{\,III}}}{\bar\alpha_{\textrm{\,IV}}-\tilde\alpha_{\textrm{\,IV}}}=\frac{-8}{-1}=8.
\] 
 We further note  that successfully formed islands possess indexes (that contain factors) forming a sequence of consecutive primes:
\[\textrm{mf IV:}\quad\textcolor{blue}{242}=2\times\mathbf{\textcolor{blue}{11}}^2;\quad\textrm{mf II:}\quad \mathbf{\textcolor{blue}{13}},\quad\mathbf{\textcolor{green}{17}};\quad\textrm{mf I:}\quad\textcolor{blue}{1406}=2\times\mathbf{\textcolor{blue}{19}}\times37.
\]
Contrary to the last three island indexes, which coicide  with   either $\tilde\alpha$ or $\bar\alpha$, 242 does not seem to have  such a connection. We may however use the latter triple of indexes as a starting point for  a consolidated consideration of a pair of Mersenne fluctuations of same height and hypothetically perfect bilateral symmetry $-$ as already borne out  regarding mf II, but still wanting in mf I. Summands $(h-1=10)$ neglected, a perfect template for the latter would be the fluctuation

 \noindent{ itf} :=

\noindent$\;1,3(+\!10),7,15(+\!10),$ 

 $31, 64(+\!10), 130,262(+\!10), 525, 1051(+\!10),2103, 1051(+\!10),\dots, 64(+\!10), 31,$

\noindent$\;15(+\!10),7,3(+\!10),1$

\smallskip{}
\noindent It uses the premise that the electron antineutrino breaks $\epsilon$-symmetry: Terms left of 31 grow with $\epsilon=0$, their followers left of 525  with $\epsilon=1$, whereupon terms grow again with $\epsilon=0$, peaking at 2103 and then  descending like they were mirrored.   We know from inspecting the relevant columns of the first table which island is involved with mf I and which with mf II. However, an input from $\gamma^{()}$, with no  information other than  $h-1$, strives to bring forth the best  of both  worlds by  alternatingly adding  $h-1$  and  treat itf  as {\em one} entity (ideal typical fluctuation). Under these conditions, there are correlations that allow  identifying a factor of a 2-sum of itf terms that would coincide with (a factor of) an island index. For an island located left of the peak, both sum terms must lie between $31$ and $525$ ($\epsilon$ unbroken), while an island located at the peak or right of it means one  term  $<525$, the other  $>525$, or broken $\epsilon$. Thus
\[\hspace*{-2.2cm}(130+74)=2^2\times3\times\textcolor{green}{17}\quad\widehat{=}\quad\textrm{island index}\;\;\textcolor{green}{17};
\]
\[(2103+272)=5^3\times\textcolor{blue}{19}\quad\widehat{=}\quad\textrm{island index}
 \;\,{1406=2\times\textcolor{blue}{19}\times37};
 \]
\[\hspace*{-1.4cm}(1061+31)=2^2\times3\times7\times\textcolor{blue}{13}\quad\widehat{=}\quad\textrm{island index }
\;\textcolor{blue}{13}.
\]  However convincingly this recipe reproduces  the former results, we will  later (along the lines of Conjecture \ref{con:To-serve-as} regarding the proton's multiplet construction) learn about a more encompassing way to derive them.

\noindent Now that we know a little about island indexes, we may address the question of how information about `particle anchorages' and `baysides'   is gained (the former in bold in the table).  The intersection of  dom\,itf  with the union of the domains of mf I, mf II and $\gamma^{(n_0)}$  yields $\{130,272\}$.  130 has an mf I incidence $\between^{(n_0^2-2)}_{1418}=\underline{130}$. 272, on the other hand, has two  incidences in mf II, $\between^{(n_0^2-1)}_{19}=\;\between^{(n_0^2+5)}_{13}=\underline{272}$, both of which are legit in that their mf I counterparts, $\between^{(n_0^2-1)}_{1402}=\;\between^{(n_0^2+5)}_{1452}=262$, don't interfere with an island and $-$ from the point of view of  summand  $10$ $-$ have the right  difference $272-262$. We are now in a position to re-interpret the above incidences as follows. While $\gamma^{(n_0)}_{1061}$ maps to $\between^{(n_0^2-2)}_{1418}$ (mf I), $\gamma^{(n_0)}_{681}$ maps to $\between'^{(n_0^2-1)}_\mathbf{19}$  as well  as to $\between'^{(n_0^2+5)}_\mathbf{13}$ (mf II)  $-$  a multivalued, or multifunction, in essence.  
 We learn that the double  image coincides with two anchorages at two  characteristic  {\em sr}\,s, whilst the single  image lacks this kind of connection. A possible physical argument  might run thus:  Both  272 and  130 have a radiation  origin in the first place. Once multifunctionally mapped to $\between^{(n)}$, however, they become time crystals:  $272(=\!L_9)$  a {\tt d}-quark constituent, a particle confined by two islands that form because of mf II'\hspace{0.5pt}s  perfect bilateral symmetry;  130 a preform of the electron antineutrino which, because of the lesser symmetry of  mf I, is prevented  from becoming a real particle 525,
hindered by the one island present.
  
  \noindent A finite, interval-nesting recursion limits the {\em sr} ranks defining $\gamma^{(n_0)}$'s domain:
\[  b_1=13, 
\quad
b_{k} = \left\lfloor b_{k-1}/2\right\rfloor\quad (1<k\leq3).
\] 

\noindent (The recursion $-$ as a segment of a longer recursion $-$  is also present  in mf I (as for the explicit form, see red terms.) And   $b_k$ squared  apply to
\[ \beta^{(k)} :=b_k^2+\delta_{b_k,\textrm{even}}+2^{8+k}\quad (1\leq k\leq3),
\]

\smallskip{}\hspace{-6mm}
\begin{tabular}{|c|l|l|} 
\hline
  \multicolumn {3} {|l|}{\footnotesize \,\, }  \tabularnewline
  \multicolumn{3}{|l|}{{\footnotesize \,\, }  {\small  \label{cohere}Anchorage within range}  $n=n_0^2$ ${-\nu_1}\atop{+\nu_2}$$=7225$${-\nu_1}\atop{+\nu_2}$}  \tabularnewline 
  $\Big.n$&\small $(h=11)$ mf I &\small $(h=11)$ mf II \tabularnewline

\hline
\hline
$n_0^2-8$ &  $\between_{1406}^{(n)}=1$ & $\between^{(n)}_{11}=1$    \tabularnewline
  \hline
 $n_0^2-7$ &  $\between_{1388}^{(n)}=3$ & $\between^{(n)}_{17}=3$    \tabularnewline
  \hline
  $n_0^2-6$ &  $\between_{1416}^{(n)}=7$ & $\between^{(n)}_{15}=7$    \tabularnewline
  \hline
$n_0^2-5$ &  $\between_{1392}^{(n)}=15$ & $\between^{(n)}_{13}=16$    \tabularnewline
  \hline 
$n_0^2-4$ &  $\between_{1390}^{(n)}=32$ & $\between^{(n)}_{15}=33$    \tabularnewline
  \hline
 $n_0^2-3$ &  $\between_{1440}^{(n)}=64$ & ${\between}^{(n)}_{\textcolor{green}{17}}=67$    \tabularnewline
  \hline
$n_0^2-2$ &  $\between_{1418}^{(n)}=\underline{130}(\neq2^7+\textcolor{red}{3})$ & ${\between}^{(n)}_{\textcolor{green}{17}}=136$    \tabularnewline
  \hline
$n_0^2-1$ &  $\between_{1402}^{(n)}=262(=2^8+\textcolor{red}{6})$  & ${\between}^{(n)}_{\mathbf{19}}=\underline{\textrm{\bf272}}$    \tabularnewline
  \hline
$n_0^2\;\;$ &  $\between_{1372}^{(n)}=525(=2^9+\textcolor{red}{13})$  & ${\between}^{(n)}_{15}=546$    \tabularnewline
 \hline
$n_0^2+1$ &  $\between_{1380}^{(n)}=1051$ & ${\between}^{(n)}_{11}=1093$    \tabularnewline
 \hline
$n_0^2+2$ &  $\between_{\textcolor{blue}{1406}}^{(n)}=2103$ & ${\between}^{(n)}_{9}=2187$    \tabularnewline
 \hline
$n_0^2+3$ &  $\between_{\textcolor{blue}{1406}}^{(n)}=1051$ & ${\between}^{(n)}_{13}=1093$    \tabularnewline
  \hline
$n_0^2+4$ &  $\between_\textrm{\bf1410}^{(n)}=\textrm{\bf525}(=2^9\!+\!\textcolor{red}{13})$ & ${\between}^{(n)}_{11}=546$    \tabularnewline
   \hline
$n_0^2+5$ &  $\between_{1452}^{(n)}=262(=2^8+\textcolor{red}{6})$ & ${\between}^{(n)}_{\mathbf{13}}=\underline{\textrm{\bf272}}$    \tabularnewline
   \hline
$n_0^2+6$ &  $\between_{1418}^{(n)}=131(=2^7+\textcolor{red}{3})$  & ${\between}^{(n)}_{\textcolor{blue}{13}}=136$    \tabularnewline
    \hline
    $n_0^2+7$ &  $\between_{1394}^{(n)}=65$  & ${\between}^{(n)}_{\textcolor{blue}{13}}=67$    \tabularnewline
       \hline
       $n_0^2+8$ &  $\between_{1424}^{(n)}=32$ & ${\between}^{(n)}_{15}=33$    \tabularnewline
         \hline
       $n_0^2+9$ &  $\between_{1400}^{(n)}=16$ & ${\between}^{(n)}_{11}=16$    \tabularnewline
         \hline   
       $n_0^2+10$ &  $\between_{1416}^{(n)}=7$ & ${\between}^{(n)}_{9}=7$    \tabularnewline
         \hline      
       $n_0^2+11$ &  $\between_{1366}^{(n)}=3$ & ${\between}^{(n)}_{13}=3$    \tabularnewline
         \hline
       $n_0^2+12$ &  $\between_{1430}^{(n)}=1$ & ${\between}^{(n)}_{7}=1$    \tabularnewline
       
                  \hline\end{tabular}
\begin{tabular}{|l|}\hline 
  ${\big.}\atop{\textrm{\,\small  where}\;{\textstyle n_0=85}}$\tabularnewline$\Bigg.\qquad\gamma^{(n_0)}_{\beta^{(k)}}$
  \tabularnewline 
   \hline\hline
   $\Big.$\tabularnewline     $\Big.$\tabularnewline   $\Big.$\tabularnewline   $\Big.$\tabularnewline   $\big.$\tabularnewline  
      $\!\gamma^{(n_0)}_{2057}=74(=\!64\!+\!h-\!1)$ 
   \tabularnewline  
     $\!\gamma^{(n_0)}_{1061}=130$ \tabularnewline  
     $\!\Big.\gamma^{(n_0)}_{681}=272(=\!262\!+\!h\!-\!1)$
  \tabularnewline   \tabularnewline \tabularnewline 
  
     \tabularnewline     \tabularnewline 
    $$\big.\tabularnewline   
     
      $\;\textrm{\small681=}2^9\textrm{\small+\textcolor{red}{13}}^2$\tabularnewline 
  $\Big.\textrm{\small1061=}2^{10}\textrm{\small+\textcolor{red}{6}}^2\!\!+\!1\,$\tabularnewline 
  $\Big.\textrm{\small2057=}2^{11}\textrm{\small+\textcolor{red}{3}}^2$\tabularnewline $\Big.$   \tabularnewline  
   $\Big.$\tabularnewline   $\Big.$\tabularnewline   $\Big.$\tabularnewline  $\big.$\tabularnewline 
  \tabularnewline 
    \hline
\end{tabular}
\bigskip{}

\noindent leading to           
\[ \gamma^{(n_0)}_{\beta^{(1)}}:=\gamma^{(n_0)}_{681}=272,\quad
                   \gamma^{(n_0)}_{\beta^{(2)}}:=\gamma^{(n_0)}_{1061}=130
                   \quad
                   (\gamma^{(n_0)}_{\beta^{(3)}}:=\gamma^{(n_0)}_{2057}=74).
 \]
For the  index triple  $ \beta^{(k)}$ generated in this way to be used in $ \beta^{(k)}\mapsto\,\between^{()}$, it has to undergo   a permutation  $\varsigma$ in  $k=1,2,3$, and the permuted triple a {\em tr}\hspace{1pt}-shift $\tau$. Assigning  $\textrm{sgn}()=-1$ to  such an operation if it is odd and  $\textrm{sgn}()=+1$ if it is even, one verifies  that there's  positive product parity: $\textrm{sgn}(\varsigma)\cdot\textrm{sgn}(\tau)=1$.  The above formulae also suggest that  chirality  is genuinely brought about by the direction the  information is processed:  The course taken by $2^n$ in the first column  is  retrograde relative to the one taken in the third column. Opposite chirality (as associated with electron neutrino)  just  could mean  `information  processed in `same-sense' manner,
\[ \textrm{527}=2^9\!+\!\textcolor{red}{15},\;
\textrm{263}=2^8\!+\!\textcolor{red}{7},\quad \textrm{alongside}\quad\beta^{(1)}=2^9+\textcolor{red}{15}^2,\;
 \beta^{(2)}=2^{8}\textrm{\small+\textcolor{red}{7}}^2.
 \]
Now we are in a position to continue the discussion of  island indexes for mf I, mf II in a comprehensive manner that includes  the electron neutrino and elucidates how $\tilde{\alpha}_\textrm{III},\tilde{\alpha}_\textrm{IV}$ are integrated in the environment of mf I. Numbers that are  key to the subject of  Sect.\ref{GJnum}  become essential here as well:
\[{\cal G}=\{3,5,11,17,19,41,43,113,115,155,429\},
\quad
{\cal J}=\{3,13,15,117,143,149\}.
\]
They are derived from the same mathematical machinery that governs  the structure of protons in multiplet form (see Conjecture \ref{con:To-serve-as} and AppendixA). Indeed, for the electron antineutrino, the island indexes (or their decisive factor) follow from the recursion
\[2^3+5=13,\quad 2^2+13=17,\quad 2+17=19,
\]
while the corresponding (reversely oriented) recursion for the electron neutrino reads
\[-2+17=15,\quad -2^2+15=11,\quad -2^3+11=3.
\]
The indexes are unique in that they share  a ${\cal GJ}$ signature:
\[{{\cal G}\atop5}\;{{\cal J}\atop13}\;
{{\cal G}\atop17}\;{{\cal G}\atop19}\qquad\qquad{{\cal G}\atop17}\;{{\cal J}\atop15}\;{{\cal G}\atop11}\;{{\cal G}\atop3}
\]
Another signature is moved center stage when considering `baysides.' There's a surprising connection between them and  $f_1$ and $f_2$:
The first four Fermat primes being 3, 5, 17 and 257, the first four Mersenne primes  3,\,7,\,31 and 127, and  $n_0=f_1\times f_2=5\times17$, we get the two mf II bayside {\em tr}\hspace{1pt}s, $n_0^2+6=7231=7\times1033$ and $n_0^2-2=7223=31\times233$, both of which  adopt  a consistent signature
\[{{0}\atop3}\;{{1}\atop5}\;
{{1}\atop17}\;{{0}\atop257}\qquad\qquad{{0}\atop3}\;{{1}\atop7}\;{{1}\atop31}\;{{0}\atop127}
\]
 and provide a profound link to the $\cal GJ$ signature  in that $(n_0^2+6)-(n_0^2-2)$ and $({n^*_0}^2-{\nu^*_1})-({n^*_0}^2+{\nu^*_2})$ respectively mark the starting point  $2^3$  and end point  $-2^3$ of the $\cal GJ$ recursions, underpinning for the first time the  integrative power of
\[\frac{\bar\alpha_{\textrm{\,III}}-\tilde\alpha_{\textrm{\,III}}}{\bar\alpha_{\textrm{\,IV}}-\tilde\alpha_{\textrm{\,IV}}}=\frac{-8}{-1}=2^3.
\]  
It  also brings with it the motivation to  follow further consolidation options. As a look-up of the second table makes it clear, mf I, mf III and mf IV can be consolidated in a parity preserving manner:

\[\quad\alpha_\textrm{I}\!+\alpha_\textrm{III}\!+\alpha_\textrm{IV} \quad\textrm{even}\qquad\alpha_\textrm{I}\!+\tilde\alpha_\textrm{I}\!+\alpha_\textrm{III}\!+\tilde\alpha_\textrm{III}\!+\alpha_\textrm{IV}\!+\tilde\alpha_\textrm{IV} \quad\textrm{even}
\]
because    $\tilde\alpha_\textrm{I}\!+\tilde\alpha_\textrm{III}\!+\tilde\alpha_\textrm{IV}=1445+247+756$. This means that mf I, when consolidated with mf III and mf IV, gains a new asset in form of the  originally mf IV-owned island with index $\tilde{\alpha}_\textrm{IV}=242$  $-$ a kind of `dry dock' for the  unreferenced level-0 neutrino 286:

 \smallskip {}
\begin{tabular}{|c|l|l|} 
\hline 
  \multicolumn {3} {|l|}{\footnotesize \,\, }  \tabularnewline
  \multicolumn{3}{|l|}{{\footnotesize \,\, }  {\small Mersenne fluctuations $\between^{(n)}$} \,for  $n=n_0^2$ ${-\nu'_1}\atop{+\nu'_2}$ } \tabularnewline 
  $\Big.n$&\small $(h'=7)$ mf III &\small $(h''=13)$  mf IV \tabularnewline

\hline
\hline $n_0^2-26$ & & $\between_{250}^{(n)}=1$      \tabularnewline
  \hline 
  $n_0^2-25$ & &  $\between_{238}^{(n)}=2$      \tabularnewline
  \hline
 $n_0^2-24$ & &  $\between_{248}^{(n)}=4$     \tabularnewline
  \hline
 $n_0^2-23$ & &  $\between_{260}^{(n)}=8$     \tabularnewline
  \hline
 $n_0^2-22$ & &  $\between_{254}^{(n)}=17$      \tabularnewline
  \hline
 $n_0^2-21$ & &  $\between_{252}^{(n)}=35$      \tabularnewline
  \hline
 $n_0^2-20$ & &  $\between_{238}^{(n)}=71$      \tabularnewline
  \hline
 $n_0^2-19$ & &  $\between_{250}^{(n)}=143$      \tabularnewline
  \hline
 $n_0^2-18$ & &  $\between_{\textcolor{blue}{242}}^{(n)}=\mathbf{286}$       \tabularnewline
  \hline
 $n_0^2-17$ & &  $\between_{\textcolor{blue}{242}}^{(n)}=573$      \tabularnewline
  \hline
 $n_0^2-16$ & &  $\between_{226}^{(n)}=1148$      \tabularnewline
  \hline
 $n_0^2-15$ & &  $\between_{242}^{(n)}=2297$      \tabularnewline
  \hline
   $n_0^2-14$ & &  $\between_{238}^{(n)}=4595$     \tabularnewline
  \hline
$n_0^2-13$ &  & $\between_{256}^{(n)}=2297$   \tabularnewline
  \hline
 $n_0^2-12$ &   &  $\between_{254}^{(n)}=1148$    \tabularnewline
  \hline
  $n_0^2-11$ &  $\between_{754}^{(n)}=1$ &   $\between_{244}^{(n)}=573$  \tabularnewline
  \hline
$n_0^2-10$ &  $\between_{748}^{(n)}=2$     &  $\between_{252}^{(n)}=\mathbf{286}(=2^8+\textcolor{red}{30})$\tabularnewline
  \hline 
$n_0^2-9$ & $\between_{750}^{(n)}=6$  & $\between_{234}^{(n)}=\underline{142}(=2^7+\textcolor{red}{14})$    \tabularnewline
  \hline
 $n_0^2-8$ & $\between_{738}^{(n)}=14$  &  $\between_{238}^{(n)}=71(=2^6+\textcolor{red}{7})$   \tabularnewline
  \hline
$n_0^2-7$ &$\between_{750}^{(n)}=29$   &  $\between_{240}^{(n)}=\underline{35}(=2^5+\textcolor{red}{3})$   \tabularnewline
  \hline
$n_0^2-6$ &  $\between_{770}^{(n)}=59$   &  $\between_{252}^{(n)}=17$ \tabularnewline
  \hline
$n_0^2-5$ & $\between_{728}^{(n)}=119$   &   $\between_{242}^{(n)}=8$ \tabularnewline
 \hline
$n_0^2-4$ &$\between_{756}^{(n)}=\underline{\mathbf{239}}$ &   $\between_{254}^{(n)}=4$     \tabularnewline
 \hline
$n_0^2-3$ & $\between_{770}^{(n)}=119$ &   $\between_{250}^{(n)}=2$    \tabularnewline
 \hline
$n_0^2-2$ &   $\between_{742}^{(n)}=59$  & $\between_{254}^{(n)}=1$   \tabularnewline
  \hline
$n_0^2-1$ &    $\between_{744}^{(n)}=29$ &  \tabularnewline
   \hline
$n_0^2\;\;$ &    $\between_{730}^{(n)}=14$ &\tabularnewline 
  \hline
$n_0^2+1$ &      $\between_{754}^{(n)}=6$  &  \tabularnewline
    \hline
    $n_0^2+2$ &     $\between_{744}^{(n)}=2$  & \tabularnewline
       \hline
       $n_0^2+3$ &    $\between_{742}^{(n)}=1$  & \tabularnewline
         \hline

                  \hline\end{tabular}\,\nolinebreak
\begin{tabular}{|l|}\hline 
  ${\Big.}\atop{\textrm{\,\small  where}\;{\textstyle n_0=85}}$\tabularnewline$\Bigg.\qquad\gamma^{(n_0)}_{\beta^{(k)}}$
  \tabularnewline 
   \hline\hline
  $\Big.$\tabularnewline    $\Big.$
 \tabularnewline  $\Big.$ \tabularnewline   $\Big.$  \tabularnewline   \tabularnewline   \tabularnewline   
\tabularnewline   \tabularnewline   \tabularnewline   \tabularnewline   \tabularnewline     $\Big.$  \tabularnewline   \tabularnewline   \tabularnewline    $\Big.$\tabularnewline           $\Big.$\tabularnewline     $\Big.$\tabularnewline  $\!(\gamma^{(n_0)}_{514}=298=286\!+\!h''-\!1)$ $\Big.$ \tabularnewline  $(\gamma^{(n_0)}_{764}=142)$$\Big.$\tabularnewline 
    $(\gamma^{(n_0)}_{486}=83=71\!+\!h''-\!1)$   \tabularnewline  $ (\gamma^{(n_0)}_{286}=35)$  $\Big.$\tabularnewline  $\Big.$\tabularnewline  $\big.$  \tabularnewline$\;\gamma^{(n_0)}_{1364}=239$

     $\Big.$ \tabularnewline  
 $\Big.$ \tabularnewline   \tabularnewline   \tabularnewline   \tabularnewline   \tabularnewline    
  \tabularnewline       \tabularnewline     $\Big.$      \tabularnewline  
    \hline
\end{tabular}

\newpage

 \noindent  Yet the new island must  be compatible with   an electron antineutrino 525 that has at least   reached its virtual phase. Writing $\alpha(n,\between^{(n)})$ and  $\beta(n,\gamma^{(n)})$ for the $\between$- and $\gamma$-related index functions, respectively, we note this is indeed the case for
  \[(n^2_0-10)-(n^2_0-11)>0\Rightarrow
  \alpha(n^2_0\!-10,286)-\alpha(n^2_0\!-11,573)=252-244>0
 \]
 as opposed to 
 \[(n^2_0-18)-(n^2_0-19)>0\;\Rightarrow \alpha(n^2_0\!-18,286)-\alpha(n^2_0\!-19,143)=242-250\ngtr0.
 \]
In a similar vein, after the consolidation $525 =(\between^{(n^2_0-10)}_{252}\!\!+\gamma^{(n_0)}_{1364})=(286+239)$,
  \[n^2_0-n_0>0\Rightarrow\alpha(n^2_0,525)-\beta((n_0,239)=1372-1364=2^3
 \] 
  as opposed to a potentially  unhindered electron-antineutrino 
  \[(n^2_0+4)-n_0>0\Rightarrow\alpha(n^2_0+4,525)-\beta(n_0,239)=1410-1364\neq2^3,
 \]  
 which proves that $2^3$  does not only  help tell apart electron antineutrinos from electron neutrinos but also  shed light on the circumstances of their provenance.
 
  \noindent As generalizing conditions could be named: (a) a `one island only' Mersenne fluctuation harboring an (anti-)neutrino flavor  in preform; (b) two islands in its twin that  confine a down-type quark constituent; (c) both fluctuations spanning a range $n_0^2\left.{-\nu_1}\atop{+\nu_2}\right.$ relative to the  mapping $\gamma^{(n_0)}$'s  $n_0$.
\medskip {}

  \vspace*{-0.6mm}\noindent If Table \ref{tab:Prime-factors-of} holds true, we can  elucidate this further by analyzing the complexities  quarks (or constituents in case they are of down type)  incur, and relate them  to one another.
To this end, we take the secondary diagonal of $\mathrm{LL}\,(G_{\mu\nu}^{(31)})$ (cf. Appendix A) and form\footnote{$\;$edge entries 2430289 and 1  trimmed off} the 6-tuple $\textrm{SD}=(227089,15297,1633,113,17,1)$. We may, in the same vein, form from  the main-diagonal representatives of $\mathrm{LL}\,(G_{\mu\nu}^{(31)})$, $\mathrm{LL}\,(G_{\mu\nu}^{(15)})$ and  $\mathrm{LL}\,(G_{\mu\nu}^{(7)})$  a 3-tuple $\textrm{MD}=(429,5,1)$. 

\noindent 

\[
\mathrm{LL}\,(G_{\mu\nu}^{(31)})=\qquad\qquad\qquad\qquad
\]
$\quad\left(\begin{array}{c}
\\
\\
\\
\\
\\
\\
\\
\\
\end{array}\right.$\hspace{-0.3cm}%
\begin{tabular}{cccc|cccc}
$\underline{\mathbf{429}}$ & \textcolor{lightgray}{155} & \textcolor{lightgray}{43} & \textcolor{lightgray}{19} & $\underline{\mathbf{5}}$ & \multicolumn{1}{c|}{\textcolor{lightgray}{3}} & $\underline{\textrm{\bf1}}$ & \textcolor{lightgray}{1}\tabularnewline
\textcolor{lightgray}{1275} & $\underline{\mathbf{429}}$ & \textcolor{lightgray}{115} & \textcolor{lightgray}{43} & \textcolor{lightgray}{11} & \multicolumn{1}{c|}{$\underline{\mathbf{5}}$} & {\bf1} & $\underline{\textrm{\bf1}}$\tabularnewline
\cline{7-8} 
\textcolor{lightgray}{4819} & \textcolor{lightgray}{1595} & $\underline{\mathbf{429}}$ & \textcolor{lightgray}{155} & \textcolor{lightgray}{41} & \bf{17} & $\underline{\mathbf{5}}$ & \textcolor{lightgray}{3}\tabularnewline
\textcolor{lightgray}{15067} & \textcolor{lightgray}{4819} & \textcolor{lightgray}{1275} & $\underline{\mathbf{429}}$ & {\bf113} & \textcolor{lightgray}{41} & \textcolor{lightgray}{11} & $\underline{\mathbf{5}}$\tabularnewline
\cline{5-8} 
\textcolor{lightgray}{58781} & \textcolor{lightgray}{18627} & \textcolor{lightgray}{4905} & \multicolumn{1}{c}{\bf1633} & $\underline{\mathbf{429}}$ & \textcolor{lightgray}{155} & \textcolor{lightgray}{43} & \textcolor{lightgray}{19}\tabularnewline
\textcolor{lightgray}{189371} &\textcolor{lightgray}{58781}& {\bf15297} & \multicolumn{1}{c}{\textcolor{lightgray}{4905}} & \textcolor{lightgray}{1275} & $\underline{\mathbf{429}}$ &\textcolor{lightgray}{115} & \textcolor{lightgray}{43}\tabularnewline
\textcolor{lightgray}{737953} & {\bf227089} & \textcolor{lightgray}{58781} & \multicolumn{1}{c}{\textcolor{lightgray}{18627}} & \textcolor{lightgray}{4819} & \textcolor{lightgray}{1595} & $\underline{\mathbf{429}}$ &\textcolor{lightgray}{155}\tabularnewline
\textcolor{lightgray}{2430289} & \textcolor{lightgray}{737953} & \textcolor{lightgray}{189371} & \multicolumn{1}{c}{\textcolor{lightgray}{58781}} & \textcolor{lightgray}{15067} &\textcolor{lightgray}{4819} &\textcolor{lightgray}{1275} & $\underline{\mathbf{429}}$\tabularnewline
\end{tabular}\hspace{-0.2cm}$\left.\begin{array}{c}
\\
\\
\\
\\
\\
\\
\\
\\
\end{array}\right)$
\newpage

 \noindent In a second step,  relevant 6- and 3-cube-complex terms  are summed up so that, plus/minus a further term ($o_\nu\in M_{5/8}$ for a quark and  $\sigma_\nu\in M_{9/8}$ for a neutrino), complies with the target value:

\medskip{}
\noindent
\hspace*{-1.5mm}\begin{tabular}{lll}
 $L_{29}\cong\textrm{\tt t}$: & $\!\!(1,-1,-1,-1,1,1)\cdot \textrm{SD}^t +(1,-1,1)\cdot\textrm{MD}^t\;-\;o_{10}$ &$(=207930)$\tabularnewline
  $L_{19}\cong\textrm{\tt c}$: & $\!\!(0,0,0,0,0,0)\cdot \textrm{SD}^t +(1,0,0)\cdot\textrm{MD}^t\;+\;o_{12}$ &$(=10668)$\tabularnewline
  $L_{11}\cong\textrm{\tt u}$: & $\!\!(0,0,0,0,0,0)\cdot \textrm{SD}^t +(1,0,0)\cdot\textrm{MD}^t\;+\;o_{2}$ &$(=438)$\tabularnewline
    $L_{16}\in\textrm{\tt b}$: & $\!\!(0,0,1,1,1,-1)\cdot \textrm{SD}^t +(0,0,-1)\cdot\textrm{MD}^t\;+\;o_{10}$ &$(=4320)$\tabularnewline
      $L_{10}\in\textrm{\tt b,s,d}$: & $\!\!(0,0,0,0,1,0)\cdot \textrm{SD}^t +(0,0,0)\cdot\textrm{MD}^t\;+\;o_{7}$ &$(=336)$\tabularnewline
        $L_{9}\in\textrm{\tt b,s,d}$: & $\!\!(0,0,0,1,0,0)\cdot \textrm{SD}^t +(0,0,0)\cdot\textrm{MD}^t\;+\;o_{6}$ &$(=272)$\tabularnewline
          $L_{8}\in\textrm{\tt b,s,d}$: & $\!\!(0,0,0,1,1,0)\cdot \textrm{SD}^t +(1,0,0)\cdot\textrm{MD}^t\;-\;o_{7}$ &$(=240)$\tabularnewline
                    $L_{6}\in\textrm{\tt d}$: & $\!\!(0,0,0,1,-1,0)\cdot \textrm{SD}^t +(0,-1,0)\cdot\textrm{MD}^t\;-\;o_{3}$ &$(=72)$\tabularnewline
                  $L_{5}\in\textrm{\tt d}$: & $\!\!(0,0,0,0,0,0)\cdot \textrm{SD}^t +(0,0,1)\cdot\textrm{MD}^t\;+\;o_{4}$ &$(=40)$\tabularnewline
                     $L_{4}\in\textrm{\tt d}$: & $\!\!(0,0,0,0,0,0)\cdot \textrm{SD}^t +(0,1,0)\cdot\textrm{MD}^t\;+\;o_{3}$ &$(=24)$\tabularnewline
 \end{tabular}

\noindent\hspace*{-1.5mm}\begin{tabular}{lll}\tabularnewline
level-$L_{16}\;\textrm{neutrino}$: & $\!\!(0,0,0,0,0,0)\cdot \textrm{SD}^t +(0,0,-1)\cdot\textrm{MD}^t\;+\;\sigma_{11}$ &$(=4606)$\tabularnewline
level-$L_{8}\;\textrm{neutrino}$: & $\!\!(0,0,0,1,0,0)\cdot \textrm{SD}^t +(1,0,1)\cdot\textrm{MD}^t\;-\;\sigma_{3}$ &$(=526)$\tabularnewline
level-$0\;\textrm{neutrino}$: & $\!\!(0,0,0,0,0,0)\cdot \textrm{SD}^t +(0,0,-1)\cdot\textrm{MD}^t\;+\;\sigma_{7}$ &$(=286)$.\tabularnewline
 \end{tabular}
 
 \vspace*{2mm}
 \noindent One notices that  the bottom quark ({\tt b}), strange quark ({\tt s}) and down quark ({\tt d}) share ($L_8,L_9,L_{10}$) as constituents. There are additions: the bottom quark includes $L_{16}$, too, and the down quark includes $L_4$, $L_5$ and $L_6$ (see Table \ref{tab:Prime-factors-of}). Representational differences notwithstanding (level-0 and level-$L_{16}$  neutrinos  don't refer to SD, while  level-$L_8$ or electron neutrinos do), these  flavors could, at least in principle, occur in conjunction with {\tt b,c,d} constituents in time-crystal contexts. There's one caveat, though.  What we previously called   {\tt u}+{\tt d} content, $1422$ $=L_{11}+(L_4+L_5+L_6+L_8+L_9+L_{10})$, leads to two additional representations: $1422=o_9+\sigma_6$ and $1422=L_{14}$. When we take  orders (subscripts)  as  equivalents of Euclidean dimension, the first is a bicube complex in 9+6 dimensions, the second, and created holographically  in one dimension less, is a densest sphere packing in 14D space without centerpiece (`mini black hole').\footnote{\,A more technical description is given on page \pageref{techpoin}} Neither form of representation works with  {\tt s+c}- or {\tt t+b} content.  So it's entirely possible that morphing into level-$L_8$ is the only way for  level-0 and level-$L_{16}$ flavors to be able to couple to quark matter, which  means that the only pairings one would expect to occur  in  time-crystal contexts are those of 525, 526 or 527 with $L_8,L_9 ,L_{10}$, such as the combination shown above: $525(\cong286+239)$ with $L_9(=272)$.

\subsubsection{\label{sub:rep-source-sink}Representations and `source-' and `sink' dimensions}  
 The  additional representations  are closely related to what in Sect.\,\ref{sub:Duality-controls:-the-1}, Table \ref{tab:Key-particle-creation-related-2} we will call `source' and `sink' dimensions. Their  inclusion in turn sets the stage for  a second  bicube complex (of 47+10  rather than  9+6 dimensions). Moreover $-$ regarding    quark/neutrino calculations $-$ they are found to stipulate the contributions' dimensions, contributions that do or do not cancel each other out.  (For corresponding calculations that lie beyond quarks or neutrinos, cf.$\,$\ref{sec:A-special-cube}.) Only an excerpt  of Table \ref{tab:Key-particle-creation-related-2} is reproduced here:

\hspace*{1.5mm}\begin{tabular}{cc}
 \multicolumn {2}{l}{$p$-based `source' $(\nu) $ and `sink' dimensions $(\nu')$:}\tabularnewline
 $\qquad\nu\qquad\qquad\qquad\qquad$ &$\nu'$\tabularnewline \tabularnewline
 $\qquad4\qquad\qquad\qquad\qquad$ &$3$\tabularnewline  
 $\qquad10\qquad\qquad\qquad\qquad$ &$7$\tabularnewline
 $\qquad\vdots\qquad\qquad\qquad\qquad$ &$\vdots$\tabularnewline  \tabularnewline
 \multicolumn {2}{l}{$o$-based `source' $(\nu) $ and `sink' dimensions $(\nu')$:}\tabularnewline
 $\qquad\nu\qquad\qquad\qquad\qquad$ &$\nu'$\tabularnewline \tabularnewline
 $\qquad12\qquad\qquad\qquad\qquad$ &$18$\tabularnewline  
 $\qquad21\qquad\qquad\qquad\qquad$ &$29$\tabularnewline
 $\qquad\vdots\qquad\qquad\qquad\qquad$ &$\vdots$\tabularnewline  \tabularnewline

 \end{tabular}
 
\noindent From the $o$-based `source' dimensions arises $\Delta_o=21-12=9$, and from the $p$-based `source' dimensions $\Delta_p=10-4=6$, which together form the `outer' dimensions of the 9+6 bicube complex.  Its `inner' dimensions arise analogously, namely from  the $o$-based and $p$-based `sink' dimensions $\Delta'_o=29-18=11$ and $\Delta'_p=7-3=4$, respectively. The dimension 11 however is peculiar in that  it assumes  parforce the   form $p_3+p_2+p_1\rightarrow7+3+1$. The ensuing bicube complex administers  in a physically plausible way  the vital contributions of both the non-({\tt u\,+\,d}) and the {\tt d}-minimal sphere packings plus the way $\textrm{e}^{\scriptscriptstyle+}$ and $\textrm{e}^{\scriptscriptstyle-}$  come into play:
 \[\Big[\textrm{e}^{\scriptscriptstyle+}+\big(1,1,1,1,1,1,1,1,1\big)\cdot\big(L_7+L_3^{\scriptscriptstyle+}+L_1^{\scriptscriptstyle+},L_7+L_3^{\scriptscriptstyle+}+L_1^{\scriptscriptstyle+},\dots,L_7+L_3^{\scriptscriptstyle+}+L_1^{\scriptscriptstyle+}\big)^t\Big]+\]\[\Big[\big(1,1,1,1,1,1\big)\cdot\big(L_4,L_4,\dots,L_4\big)^t+\textrm{e}^{\scriptscriptstyle-}\Big]=o_{\Delta_o}+\sigma_{\Delta_p}=o_9+\sigma_6=1422.\]
 The pattern substitutions $7+3+1\rightarrow L_7+ L^{\scriptscriptstyle+}_3+L^{\scriptscriptstyle+}_1$, $4\rightarrow L_4$ are a key ingredient here. They apply again when it comes to constructing the  47+10, dual  bicube complex:  From the $o$-based `sink' dimensions we get $\Sigma'_o=18+29=47$, and from the $p$-based `sink' dimensions we get $\Sigma'_p=3+7=10$, which become the  `outer' dimensions of the bicube complex.  Just as 3 and 7 have the status of  $p$ numbers, 47 gets its status by $o$-numbers: $9+38=o_2+2o_3$. The new `inner' dimensions arise analogously, namely from  the `source' dimensions $\Sigma_o=12+21=33$ and $\Sigma_p=4+10=14$, respectively. 33 being a multiple of 11, it's a case for parforce treatment: $33\rightarrow3(7+3+1)$; using  the shorthands  ${\cal L}^{\scriptscriptstyle\_\_\,+}_{7\;3\;1}$ and  ${\cal L}^{\scriptscriptstyle\_+\,+}_{7\;3\;1}$ for the  substitutions $7+3+1\rightarrow L_7+L_3+L^{\scriptscriptstyle+}_1$ and $7+3+1\rightarrow L_7+L^{\scriptscriptstyle+}_3+L^{\scriptscriptstyle+}_1$, respectively, we obtain
 \[\!{\Big[\!\atop\;}{\big(\atop\;}{\underbrace{1,1,{\scriptstyle\dots},\!1,}\atop{\scriptscriptstyle9}}{\underbrace{1,1,{\scriptstyle\dots},\!1}\atop{\scriptscriptstyle38}}{\big)\atop\;}
 {\!\cdot\big(3{\cal L}^{\scriptscriptstyle\_+\,+}_{7\;3\;1},\atop\;}{3{\cal L}^{\scriptscriptstyle\_+\,+}_{7\;3\;1},\atop\;}{{\scriptstyle\dots},\atop\;} {3{\cal L}^{\scriptscriptstyle\_+\,+}_{7\;3\;1},3{\cal L}^{\scriptscriptstyle\_\_\,+}_{7\;3\;1},\atop\;}{3{\cal L}^{\scriptscriptstyle\_\_\,+}_{7\;3\;1},{\scriptstyle\dots},3{\cal L}^{\scriptscriptstyle\_\_\,+}_{7\;3\;1}\big)^t\Big]+\atop\;}\]\[
 {\!\!\Big[(1,1,1,1,1,1,1,1,1,1)\!\cdot\!(L_{14},L_{14},{\scriptscriptstyle\dots}, L_{14})^t\Big]=\atop\;}\qquad\qquad\qquad\qquad\quad\quad\quad\]\[\,\quad34128=L_4(o_9+\sigma_6)=24({\tt u+d})=8({\tt p}+{\tt n}).\qquad\qquad\qquad\qquad\qquad\qquad\qquad\]
\noindent As  terms $L_7$  in the 9+6 bicube complex  get incrementally  replaced by  $L_7^{\scriptscriptstyle+}$, an electron shell may build up and mark the transition  subatomic-to-atomic. For: every time this occurs, an additional $\textrm{e}^{\scriptscriptstyle-}$  pops up, infused with electric charge -2/3 to compensate for the  enhancement $L_7\!\rightarrow \!L^+_7$, of charge increase +2/3; because of an overall factor 3/2 attached to a nucleon (3{\tt u}+3{\tt d}, for instance, matches the deuteron's proton-neutron pair) the compensation is  effectively of unit charge. This holds true  as well  for the 47+10  bicube complex where the replacements are $L_{14}\rightarrow L_{14}^{\scriptscriptstyle+}$ (showing that the causes for an $\textrm {e}^{\scriptscriptstyle-}$ instance might be of  a more speculative than just hadronic nature).

\noindent What should be noted as well: the $o$-type contributions in the above SD- and MD-based quark calculations fall in two categories: those that cancel each other out and those that do not.
The former category feeds from `sink' dimensions, $3\rightarrow\pm o_3$, $7\rightarrow\pm o_7$ and $3+7\rightarrow\pm o_{10}$, while $o$-type contributions  remaining strictly positive  always feed from `source' dimensions: $12\rightarrow o_{12}$,   $21-19\rightarrow o_{2}$ (a hybrid since, according to Table \ref{tab:Key-particle-creation-related-2}, 21 is  $o\,$-based whereas 19 is $p$-based), $10-4\rightarrow o_6$ and $4\rightarrow o_4$. In contrast, the $\sigma$-type contributions in SD- and MD-based neutrino calculations belong exclusively, and hence chirally, in the `sink' category:  $29-18\rightarrow+\sigma_{11}$ (this time spared with parforce treatment!), $3\rightarrow-o_3$ and $7\rightarrow+o_7$.

\noindent Summarizing, the bicube complex $34128=3\times8\times1422$ is spanned by a whopping $47+10=57$ dimensions.\footnote{\,As for strategies for detecting (or inferring) incidences of $\between^{(n)}[\,]=34128$, one should look for amplitudes which MF-wise are comparable in size with $L_m$s and keep them company. In other words: for incidences of $L_m$s with low powers of two as coefficients. For instance, with  as few companions as possible, $2L_{10}+34128=2L_{20}$ (two companions), $2L_{19}+34128=2L_{21}+L_4\,(=55464)$ or
$4L_5+L_9+34128=8L_{16}\,(=34560)$ (three companions each).\, If the refinement for $\between^{(n)}[\alpha]=34128$ is too hard to compute, maybe one of the companions' is not.} Together with the  $1422=(o_9+\sigma_6)$ bicube complex with   parsimonious $9+6=15$ dimensions, the two form the basis of a rich atomic tableau.  Effected by means of `source' and `sink' dimensions, the elemental view of {\tt u+d} content, $L_{14}=1422=L_{11}+(L_4+L_5+L_6+L_8+L_9+L_{10})$, experiences a reinterpretation
as $(o_9+\sigma_6)$
which is  marked by the emergence of  electrons and neutrinos. 
 \subsubsection{\label{sub:tetrapetalia}On to the tetrapetalia and transactional spacetime}
 In principle, the  reinterpretation could take the form of  a phase transition.  We therefore propose a spacetime construct \`a la Minkowski,  a 4D setting with space components $x,y,z$ and  time component  $t_\nu$. As a measure of 3D volume  ($V$)  we use  tetrahedral numbers  ({\em Te\,}$_v$) whose (discrete) values are determined by $p_n$ and $\sigma_\tau$:
 \[\begin{array}{rl}V:=\textrm{\em Te\,}_v=-p_n +\sigma_\tau & (v=6\cdot2^n-1,\tau=4+3n; \quad n=0,1,2, ...).\end{array}\]
 Led by the perspective that the (9+6)-dimensional bicube complex can,  holographically,  be pictured by a 14D Euclidean quantity like $L_{14}$, we choose to drop  one space component and  introduce the transformation equation  \[t_\nu=(x^2+y^2)^2+\frac{\sqrt{x^2+y^2-1}}{ x^2+y^2}\] which has what we call a {\em tetrapetalia} as  Cartesian plot $-$ one with solutions that unfold as four petals for their real part  and pistils for their imaginary part. If this could be named the primitive transformation equation, its companion, the conjugated equation  
\[\bar t_\nu=(x^2+y^2)^2-\frac{\sqrt{x^2+y^2-1}}{x^2+y^2},\] would form a similar tetrapetalia, again with four petals as real part, but a  stem  as imaginary part:

\hspace*{4cm} {\bf\tt Tetrapetalia}
 
 \smallskip {}
 
 \smallskip{}

\hspace*{1cm}{\bf\tt (primitive)}\hspace*{6cm}{\bf\tt (conjugated)}\newline
\includegraphics[scale=0.33]{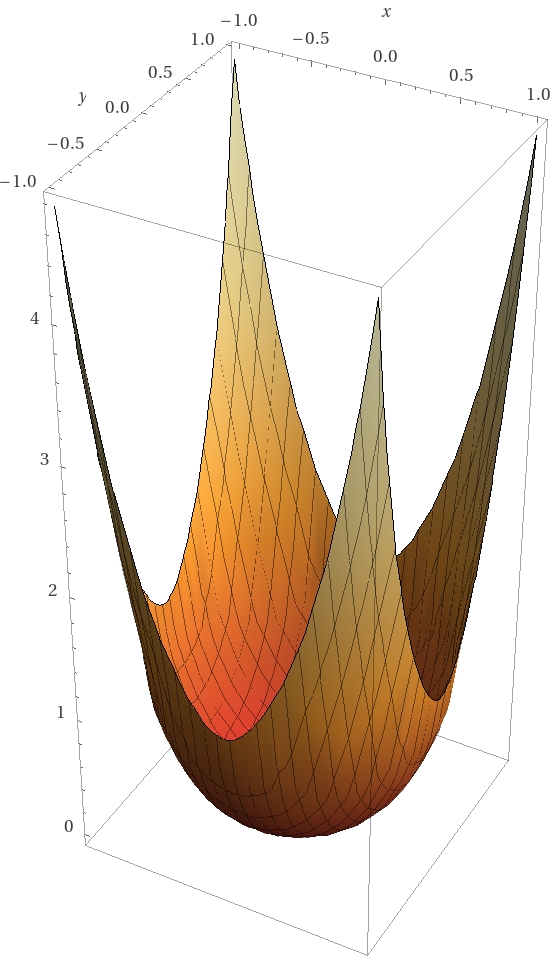}\includegraphics[scale=0.33]{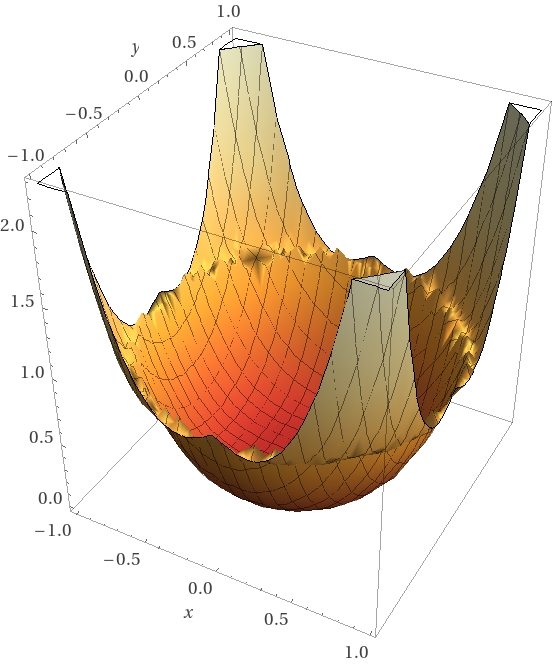}

\includegraphics[scale=0.33]{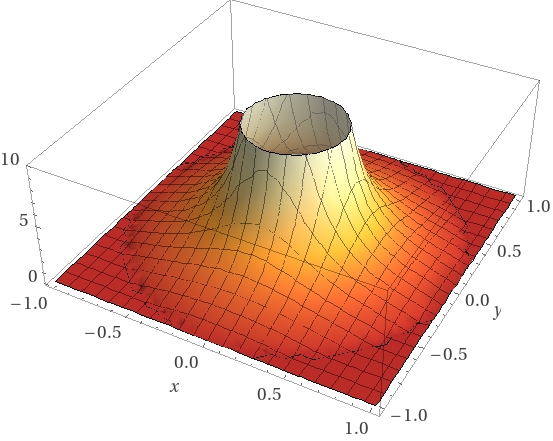}\includegraphics[scale=0.33]{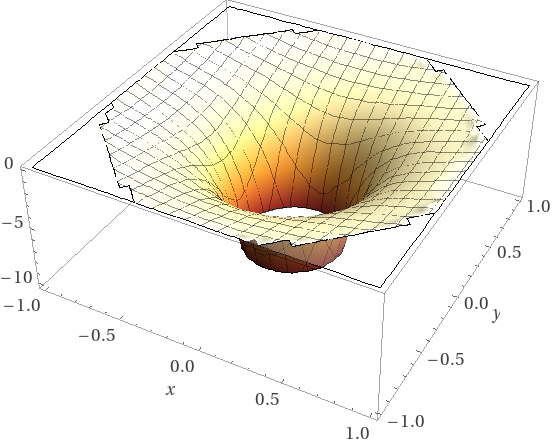}

 \noindent One notes that for the real part, solutions tend to climb  the $x=\pm y$ corners, so it's worthwhile to simplify further and check solutions along $y=x$:

 
 \noindent$t_\nu=(x^2+y^2)^2\!+\!\frac{\sqrt{x^2+y^2-1}}{x^2+y^2}\;$  $(y=x)$\quad$\bar t_\nu=(x^2+y^2)^2\!-\!\frac{\sqrt{x^2+y^2-1}}{x^2+y^2}\quad$  $(y=x)$

 \smallskip{}

\noindent\includegraphics[scale=0.8]{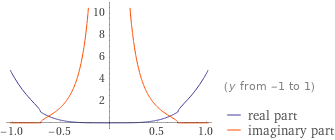}\quad\includegraphics[scale=0.8]{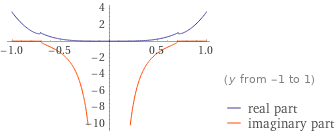}

\smallskip{}

\noindent\includegraphics[scale=0.8]{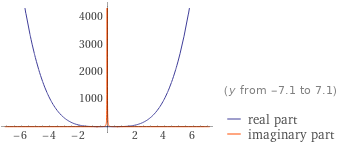}\quad\includegraphics[scale=0.8]{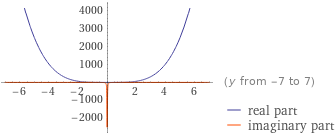}

 \bigskip{}

\noindent In Minkowskian terms ($t_\nu\sim t_{\textrm{M}},\bar t_\nu\sim t_{\textrm{M}}$), both the primitive and the conjugated equations'  large-scale behaviors  obey causality while the conjugate's small-scale behavior violates it  $-$ with lift-off from the   future light cone at $t_\textrm{M}=1$ and reentry at ${t'}_{\textrm{M}}=1.104...$.\label{causalpoin} This may or may not be an intriguing feature in itself; what we are more interested in, however, is what we call the tetrapetalia's  $1/\sqrt2$-order -- the same order that makes  the future-light-cone line  $t_{\textrm{M}}=\sqrt2x$.
Equating $y=x$, then writing \[\theta:X\rightarrow T,\;\bar\theta:X\rightarrow T,\;\theta^{-1}:T\rightarrow X,\;{\bar\theta}^{-1}:T\rightarrow X\quad (x\in X,t_\nu\in T) \] for the respective primitive and conjugate tetrapetalic transformations, we would normally expect a natural  input  $x$ to yield an irrational output $\theta(x)$ or $\bar\theta(x)$.  But exceptions  there are. With natural solutions to the roots $w=\sqrt{2x^2-1}$,  $x=1,w=1$, $x=5,w=7$, $x=29,w=41$,  $x=169,w=239\dots$,   say, we not only get rational transformation results but also sequences with a ratio $x_n/w_n\rightarrow1/\sqrt2$  that exert  their influence  on the whole irrational output. That $1/\sqrt2$  order we'll associate  with the interconvertent behavior of neutrinos in transit. \newline Let us preface the analysis with a discussion of the decay patterns of neutrinos under transformations of the tetrapetalia. We denote
the  continued fraction expansions of   $X\rightarrow T$-transformations  $\theta(x)[\beta]$ and  ${\bar\theta}(x)[\beta]$ and those of   $T\rightarrow X$-transformations  $\theta^{-1}(t_\nu)[\alpha]$ and  ${\bar\theta}^{-1}(t_\nu)[\alpha]$,  respectively, and, for the sake of completeness,   the associated refinements  by $\beta(x)$, $\bar\beta(x)$, $\alpha(t_\nu)$, $\bar\alpha(t_\nu)$. As with every transformation, we  have to look for conserved quantities,  then. 
Concentrating, for a start, on   $L_{14}=1422$,  we find
\[\begin{array}{ll}
\theta^{-1}(1422)=4.34207...,\quad {\bar\theta}^{-1}(1422)=4.34232...\end{array}\] 
as positive solution to the $T\rightarrow X$ transformation, 
and
\[\begin{array}{ll}
\theta(1422)\!=\!1.635529481222400...\!\times\!10^{13},&  \bar \theta(1422)\!=\!1.635529481222399...\!\times\!10^{13}
\end{array}\]
as solutions to the  $X\rightarrow T$ transformation. We anticipate  that the pecularity that solutions diverge after the third decimal place in one case and not until  the 13th  in the other, will entail a crucial divide.
The ensuing neutrino 526-related information reads\label{0th-approx}\[\begin{array}{rr}
{\bar\theta}(1422)[2168]=526,\\ 
\bar\theta^{-1}(1422)[42]=526,\end{array}\]
which means there is one conserved quantity {\em per se}:  the unreferenced, electron neutrino-related information 526. 
Finding the counterpart for the primitive transformation $\theta(1422)$  is easy  because the twelve matching decimals are simply the consequence of the primitive refinements advancing in step with the conjugated ones: $\Delta[\bar\beta_{526}(1422),\beta_{526}(1422)]=1$.\footnote 
{\,\,the purported relationship can be demonstrated for the simplest nontrivial case $x=3$: $\theta(3)=324+\frac{\sqrt{17}}{18}$ $=324+\frac{0.229...}{1}$ $=324+\frac{1}{4+\frac{0.366...}{1}}$ $=\cdots,$ whereas $\bar\theta(3)=324-\frac{\sqrt{17}}{18}$
$=323+\frac{18-\sqrt{17}}{18}$ $=323+\frac{0.7709...}{1}$ $=323+\frac{1}{1+\frac{0.2971...}{1}}$  $=323+\frac{1}{1+ \frac{1}{3+ \frac{0.3656...}{1}}}$. As differences like $0.366...- 0.3656...$ are getting ever smaller in the sequel, the associated cf terms stay in step.} The meager three  matching decimal places in the inverse case $\theta^{-1}(1422)$, on the other hand, adumbrate  high refinemental dispersion $-$ where 
$\alpha_{526}(1422),\bar\alpha_{526}(1422)$ signify those refinements:
$\alpha_{287}(1278)\gg\bar\alpha_{287}(1278)$,
  $\alpha_{526}(1422)\gg\bar\alpha_{526}(1422)$.\footnote{\,Here, $\bar\alpha_\textrm{neutrino}<2500$ $\rightarrow\alpha_\textrm{neutrino}>2500$, or  $\alpha_\textrm{neutrino}<2500$ $\rightarrow\bar\alpha_\textrm{neutrino}>2500$.}
This more obvious characteristics of the tetrapetalia aside, nothing  points to the existence of further conserved quantities. And yet there  are them $-$ revealed, as we shall see,  by particularly combined  input pairs ($x,t_\nu$) and their associated refinements (see \ref{sec:Minkowski}).
 For instance,  recall that  the (9+6) bicube complex has the characteristic partition $1422=1278+144$. Now form from the first term  the pseudo-Mersenne number $2\cdot1278-1=2555$  and you find to your surprise   that this value  serves as an identity in a combined (286,526)  system:
\[\begin{array}{ccc}\;\bar\theta(1422)[387]=286,& \bar\theta(1422)[2168]=\!526,&(387+2168=2555)\\
\end{array}\] 
We can deal with this system in the following way:
 \[\begin{array}{rrl}
 \textrm{for primitive neutrino 526 :}& 2167 &\textrm{divides}\;14^{20}-1 \\  
\textrm{for primitive neutrino 286 :}& 386 &\textrm{divides}\;81^{4}-1 \\
\textrm{for conjugated neutrino 526 :}&2168  & \textrm{divides}\;29^{6}-1 \\ 
 \textrm{for conjugated neutrino 286 :}&387  & \textrm{divides}\;44^{2}-1 \\
 \end{array} \] 
 That is,  neutrino-related    refinements would conspire with  dimensions $d$ through a pseudo-Mersennian form $\textrm{B}=d^n-1$, dimensions  that assume a special meaning pairwise: \[\begin{array}{rrr} 14+29=43, & 14+44=58, & 29+44=73, \\
 14+81=95, & 29+81=110,& 44+81=125.\end{array}\]In each row,  sum values are spaced 15 units apart, a reminder of  the (9+6) structure of the bicube complex. The sum values in the middle, as if attracted,  approach the Sophie Germain prime 53 as well as the associated safe prime 107. The deviations 5 and 3 from these primes will concern us in a moment, but here we want to highlight the role of the first and the middle sums for the combined (286,526) system under $X\rightarrow T$ transformations:  \[\begin{array}{r}14+29\equiv43
\;(\textrm{mod}\,\,107)
 \\
\beta_{286}(1422)+\beta_{526}(1422)+14+44\equiv43
\;(\textrm{mod}\,\,107) \\
\\14+81\equiv95
\;(\textrm{mod}\,\,107) \\
\bar\beta_{286}(1422)+\bar\beta_{526}(1422)+29+81\equiv97
\;(\textrm{mod}\,\,107) \\
\end{array}\] In both the primitive and the conjugated expressions,  the remainders from the second row can be reinterpreted as the result of subtracting $L_4$  107 times $-$  proof that  the full bicube complex (the 144-part encompassing basis elements, $ L_4$ included) is involved in each step. The value 43 is a conserved quantity in the primitive case, and while rest mass is a conserved quantity in spacetime,   43  constitutes a preserved `sink' dimension. By contrast, the conjugated case (the one  tied to the identity 2555) comes up with a mismatch: $-$ 95 vs. 97. 
Applying the second modulus to the remainders, we get\[\begin{array}{c}43\equiv43
\;(\textrm{mod}\,\,53)
 \\
95\equiv42
\;(\textrm{mod}\,\,53) 
\\
97\equiv44
\;(\textrm{mod}\,\,53) 
\end{array}\] and find these refinements mirrored in \[\begin{array}{ccc}{\bar\theta}^{-1}(1422)[42]\!=526,&{\bar\theta}^{-1}(1278)[309]\!=287&\big(
309\equiv\!44\,(\textrm{mod}\,53)\big).\end{array}\]  
Even if we were to conclude that separately conserved quantities like 42 and 44 characterize `decay' in some form, the atomic context signals that the story could not end without including $\beta+$ decay $-$ which 287 certainly does not stand for. In fact, it turns out that a more complete version recovers 43 as a conserved quantity.  The key  to the solution is  the $m$ vs. $\!\!m\!-\!2$ juxtaposed difference  between  (the1278-part of) $L_{14}$ and $L_{12}$ $(=1278-756)$, which lies close to $\bar p_9=512$ so that we may check the vicinity of $t_\nu=512$. First, in a prepared-for-interaction manner, the (9+6) bicube complex forms as\newline
(1278-part) \[\begin{array}{lcl}
\big(1,1,1,1,1,1,1,1,1\big)&\!\!\!\!\!\!\!\!\cdot
&\!\!\!\!\!\!\!\! \big({\cal L}^{\scriptscriptstyle+\,+\,+}_{7\;3\;1},{\cal L}^{\scriptscriptstyle+\,+\,+}_{7\;3\;1},{\cal L}^{\scriptscriptstyle+\,+\,+}_{7\;3\;1},{\cal L}^{\scriptscriptstyle+\,+\,+}_{7\;3\;1},\\&  &\!\!\!\!{\cal L}^{\scriptscriptstyle\_\,+\,+}_{7\;3\;1},\\
& &\!\!\!\!{\cal L}^{\scriptscriptstyle\_\_\;+}_{7\;3\;1},{\cal L}^{\scriptscriptstyle\_\_\;+}_{7\;3\;1},{\cal L}^{\scriptscriptstyle\_\_\;+}_{7\;3\;1},{\cal L}^{\scriptscriptstyle\_\_\;+}_{7\;3\;1}\big)^t=1278;\end{array}\]
(144-part) \[\begin{array}{lcl} 
(1,1,1,1,1,1)&\!\!\!\!\!\cdot
&\!\!\!\!\! (L^+_{4},L^+_{4},L^+_{4},L^-_{4},L^-_{4},L^-_{4})^t=144.
\end{array}\] 
In the vicinity of $t_\nu=512$, we find ${\theta}^{-1}(544)[427]\!=286$ which has  pendants in the projections
\[\begin{array}{rcl}
\big(0,0,0,0,0,0,1,1,1\big)&\!\!\!\!\!\!\!\!\cdot
&\!\!\!\!\!\!\!\! \big({\cal L}^{\scriptscriptstyle+\,+\,+}_{7\;3\;1},{\cal L}^{\scriptscriptstyle+\,+\,+}_{7\;3\;1},{\cal L}^{\scriptscriptstyle+\,+\,+}_{7\;3\;1},{\cal L}^{\scriptscriptstyle+\,+\,+}_{7\;3\;1},\\&  &\!\!\!\!{\cal L}^{\scriptscriptstyle\_\,+\,+}_{7\;3\;1},\\
& &\!\!\!\!{\cal L}^{\scriptscriptstyle\_\_\;+}_{7\;3\;1},{\cal L}^{\scriptscriptstyle\_\_\;+}_{7\;3\;1},{\cal L}^{\scriptscriptstyle\_\_\;+}_{7\;3\;1},{\cal L}^{\scriptscriptstyle\_\_\;+}_{7\;3\;1}\big)^t+\\
(1,1,1,1,1,0)&\!\!\!\!\!\cdot
&\!\!\!\!\! (L^+_{4},L^+_{4},L^+_{4},L^-_{4},L^-_{4},L^-_{4})^t=423+121=544
\end{array}\] 
and (within the 1278-part)
\[\begin{array}{rcl}
\big(1,1,0,0,0,0,0,0,0\big)&\!\!\!\!\!\!\!\!\cdot
&\!\!\!\!\!\!\!\! \big({\cal L}^{\scriptscriptstyle+\,+\,+}_{7\;3\;1},{\cal L}^{\scriptscriptstyle+\,+\,+}_{7\;3\;1},{\cal L}^{\scriptscriptstyle+\,+\,+}_{7\;3\;1},{\cal L}^{\scriptscriptstyle+\,+\,+}_{7\;3\;1},\\&  &\!\!\!\!{\cal L}^{\scriptscriptstyle\_\,+\,+}_{7\;3\;1},\\
& &\!\!\!\!{\cal L}^{\scriptscriptstyle\_\_\;+}_{7\;3\;1},{\cal L}^{\scriptscriptstyle\_\_\;+}_{7\;3\;1},{\cal L}^{\scriptscriptstyle\_\_\;+}_{7\;3\;1},{\cal L}^{\scriptscriptstyle\_\_\;+}_{7\;3\;1}\big)^t=286.
\end{array}\] In fully interactive form, the bicube complex is enriched with 11 electrons or positrons. Then a further case in point, namely  \[{\theta}^{-1}(532)[446]\!=527,\] is found to conform to the projections
\[\begin{array}{rcl}
\big(0,0,0,0,0,0,1,1,1\big)&\!\!\!\!\!\!\!\!\cdot
&\!\!\!\!\!\!\!\! \big({\cal L}^{\scriptscriptstyle+\,+\,+}_{7\;3\;1},{\cal L}^{\scriptscriptstyle+\,+\,+}_{7\;3\;1},{\cal L}^{\scriptscriptstyle+\,+\,+}_{7\;3\;1},{\cal L}^{\scriptscriptstyle+\,+\,+}_{7\;3\;1},\\&  &\!\!\!\!{\cal L}^{\scriptscriptstyle\_\,+\,+}_{7\;3\;1},\\
& &\!\!\!\!{\cal L}^{\scriptscriptstyle\_\_\;+}_{7\;3\;1},{\cal L}^{\scriptscriptstyle\_\_\;+}_{7\;3\;1},{\cal L}^{\scriptscriptstyle\_\_\;+}_{7\;3\;1},{\cal L}^{\scriptscriptstyle\_\_\;+}_{7\;3\;1}\big)^t+\\
({1,1,1,1},0,0)&\!\!\!\!\!\cdot
&\!\!\!\!\! (L^+_{4},L^+_{4},L^+_{4},L^-_{4},L^-_{4},L^-_{4})^t+11\,\textrm{e}^+=\\
& &\quad423+98+11=532;
\end{array}\] 
\[\begin{array}{rcl}
\big(0,0,1,1,1,0,0,0,0,0\big)&\!\!\!\!\!\!\!\!\cdot
&\!\!\!\!\!\!\!\! \big({\cal L}^{\scriptscriptstyle+\,+\,+}_{7\;3\;1},{\cal L}^{\scriptscriptstyle+\,+\,+}_{7\;3\;1},{\cal L}^{\scriptscriptstyle+\,+\,+}_{7\;3\;1},{\cal L}^{\scriptscriptstyle+\,+\,+}_{7\;3\;1},\\&  &\!\!\!\!{\cal L}^{\scriptscriptstyle\_\,+\,+}_{7\;3\;1},\\
& &\!\!\!\!{\cal L}^{\scriptscriptstyle\_\_\;+}_{7\;3\;1},{\cal L}^{\scriptscriptstyle\_\_\;+}_{7\;3\;1},{\cal L}^{\scriptscriptstyle\_\_\;+}_{7\;3\;1},{\cal L}^{\scriptscriptstyle\_\_\;+}_{7\;3\;1}\big)^t+\\
({1,1,1,1},0,0)&\!\!\!\!\!\cdot
&\!\!\!\!\! (L^+_{4},L^+_{4},L^+_{4},L^-_{4},L^-_{4},L^-_{4})^t+ \textrm{e}^+=\\&&428+98+1=527.
\end{array}\] 
As a matter of fact, both instances 286 and 527 would originate in  primitive transformations $\theta^{-1}(t_\nu)$. Our previous find, $\bar\theta^{-1}(1278)[309]=287$, however, is a reminder that conjugated transformations matter no less. Now the juxtaposition that led  from 1278 to $L_{12}\!=\!756$ offers  the opportunity to search for one. Surprisingly, noting  $\bar\theta^{-1}(756)[37]=766$ (a value linked to $1278-\bar p_9=766$), there emerges a new refinement $-$ $\bar\alpha_{766}=37$  $-$
 that acts as the right corrective enabling conservation of 43 in $\beta+$ decay:
\[\begin{array}{ll}\alpha_{286}=427\;\; \textrm{divides}\;\;62^2-1, 
&\alpha_{527}=446\;\; \textrm{divides}\;\;39^3-1.\end{array}\]
Together with $\alpha_{526}=42$, it allows writing
\[\begin{array}{rl}\alpha_{286}+\alpha_{526}+39-\textcolor{red}{37}&\equiv43\;(\textrm{mod}\,107) \\
\alpha_{527}+62-\textcolor{red}{37}&\equiv43\;(\textrm{mod}\,107). \end{array}
\]
What remains is clarifying why 11 positrons showed up where only one would be needed to define the decay product 527. First of all, conservation  works   under the composite modulus operation ``\,first  107, then 53.'' When applied to the composite number $150$ $-$  the sum of $\beta_{286}(1422)+\beta_{526}(1422)$ modulo 107, plus the middle sum mentioned above $-$ $14+44$  $-$, we get the  twin primes 149 and 151 as 150's neighbors. Applying the modulus, 150 becomes equivalent to 43,   a prime  in its turn   sided  by the composite numbers 42 and 44. This looks like a dual $-$ meaning bidirectional $-$ situation:\,\[\begin{array}{cc}\textrm{149 divides}\;42^{9}-1,&\textrm{151 divides}\;44^{25}-1, \\ 
\textrm{42 divides}\;149^{6}-1, &\textrm{44 divides}\;151^{10}-1.\end{array}\]
The exponents encode the primitive Pythagorean triple $5^2=3^2+4^2$: One  exponent  is equal the hypotenuse squared, $25=5^2,$  while the others are $9=3^2$ and, in summed-up form: $4^2=6+10$. One recognizes that the 9 dimensions of the 1278-part of the (9+6) bicube complex are not a chance product; the next higher representation of the atom, 34128, would likely come  with another  number partaking in a Pythagorean triple. On closer inspection one finds that, for $x=34128$  {\em per se},  neutrino-related information  has unreasonably high refinement; rather, we must look for more accessible candidates in the vicinity. While those showing up mostly turn out  to be singles, two instances, both located within the comfortable zone $\beta<2500$, originate from  $x=34143$\label{34143poin} and, like with 286 and 526, are bound up with net charge zero:  $\theta(34143)[1655]=287$, $\theta(34143)[1883]=525$. Again, a primitive Pythagorean triple ($5^2=3^2+4^2$) is involved, as is a multiple of 11: 
\[25\times1278+15\times144+33\,\textrm{e}^+=34143.
\]
 And the sequence goes on. With a 19908-part assuming the task of 1422's 1278-part, and 49 partaking in the primitive Pythagorean triple $7^2+24^2=25^2$:
\[49\times19908+57\times144+77\,\textrm{e}^+=983777.
\]
(With refinement $\beta_{286}=831$ for $\theta(983777)[831]=286$ being more than twice that for 
$\theta(1422)[386]=286$,  $\beta_{526}$  (if there exists a $\theta(983777)[\beta_{526}]=526$ at all), would be out of reach: $\beta_{526}>2\times2167>2500$.) The sequence's organizing principle seems to be $(2^n-1)\!\times\!11\textrm{e}^+$. A further bump in neutrino information to be found, as successor of $L_{14}=1422$, at $L_{21}=27720$ may be due to a similar principle,  $L_{n\times7}$:
\[\begin{array}{rc}\theta(27720)[149]=526,&\;\;\,\theta(27720)[573]=285, \\
\theta(27719)[1267]=287,&\theta(27721)[1182]=286
\end{array}\]  (with an ion character rather than an atom one).  The principles may be interlinked via tetrahedral numbers: $Te_{\,54}=L_{3\times7}=27720$; and $Te_{\,11}=286$ which may  entail `elevenness' of positron (or electron) occurrences. Otherwise,  $Te_{\,5}=35$, $\dots$, $Te_{\,23}=2300$, $Te_{\,47}=18424$, $Te_{\,95}=147440,\dots$ (our tetrapelia volume marks  $Te_{v}=-p_n+\sigma_\tau$),   remain neutrino-wise nondescript.


\noindent Resuming our above line of reasoning, we   now deal with the tetrapetalia's $1/\sqrt2$ order for $x\in\mathbb{N}$ which, in case of irrational output, entails the relationship

\[\theta(x)[n]=\left\{\begin{array}{ll}\bar\theta(x)[n+1]+1\approxeq\lfloor\sqrt2\,x\rfloor\;\textrm{or}\;\lceil\sqrt2\,x\rceil, &(n=1)\\
\bar\theta(x)[n+1] &(n>1).\end{array}\right.\]
This  lets us conveniently find places holding  neutrino-related information:
\[ \theta(372)[1]=\bar\theta(372)[2]+1=526,\qquad\theta(3257)[1]=\bar\theta(3257)[2]+1=4606.
\]
It's easy to see that the relationship makes  neutrino 286 a special case:   \[\nexists \,x\!:\theta(x)[1]=286.\] Instead,
\[  \theta(202)[1]=\bar\theta(202)[2]+1=285,\qquad\theta(203)[1]=\bar\theta(203)[2]+1=287.
\]
Thus,  while interconversions 
to and fro 287, 527 and 4607 go  unimpeded so that they may stand for oscillations in transit
\[287=\theta(203)[1]\; {\rightleftarrows}\; \theta(373)[1]=527\;{\rightleftarrows}\;\theta(3258)[1]=4607,
\]free (in-transit) oscillations  of decremented type remain limited to types 525 and 4605 $-$ with conjugate transport
\[
\bar\theta(372)[2]=525\rightleftarrows\bar\theta(3275)[2]=4605.\]Because $\nexists\,x:\bar\theta(x)[2]=285$, interconversions to and fro this neutrino type are a phenomenon strictly limited to bound systems. This can be demonstrated using the dual  $\nexists x:\theta(x)[1]=286$ as well as 526's  definition reading $526=286+\textcolor{blue}{240}$, where the rhs emerges  in  compact form in
\[286=\bar\theta^{-1}(\textcolor{blue}{240})[49].\]
If this  inverse transformation were somehow akin to  the conjugate transformation of $1/\sqrt2$ tetrapetalia order,  $\theta(170)=[3340840000;\textcolor{blue}{240},2,...]$,   the latter could be converted into  $\bar\theta(169)=[3...;1,\textcolor{blue}{240},...].$ Instead \[\bar\theta(169)=[3262922884;238,239]\] $-$ a signal that irrationality broke down, inhibiting interconversion. 

\noindent Meaningful combinations of  $1/\sqrt2$-  and  normal order are possible, yet.
This  foremost concerns conservation in mixed cases $\big(\{T\rightarrow X\}_\textrm{c},\{X\rightarrow T\}_\textrm{p}\big)$ where p and c stand for primitive and conjugate, respectively. We can directly turn to the combined (287,526)  system  where the sought-after partners would be the previously presented $1/\sqrt2$-order terms ($\bar\theta(203)[1]=287$ and  $\theta(372)[1]=526$):
 \[\begin{array}{rlc}\bar\alpha_{526}(1422)+\beta_{287}(203)=42+1&\equiv43
\;(\textrm{mod}\,\,107)
&43\equiv43
\;(\textrm{mod}\,\,53) 
\\
\bar\alpha_{287}(1278)+\beta_{526}(372)=309+1&\equiv96
\;(\textrm{mod}\,\,107) 
&\;96\equiv43
\;(\textrm{mod}\,\,53).

\end{array}\]
More generally,   the tetrapetalia's  $1/\sqrt2$ order has a perfusing quality:  
\[x_0,\quad x_1=\lfloor\sqrt2\times x_0\rfloor\;\textrm{or}\;\lceil\sqrt2\times x_0\rceil, \quad x_2=\lfloor\sqrt2\times x_1\rfloor\;\textrm{or}\;\lceil\sqrt2\times x_1\rceil,\;\dots,\textrm{\,say,}
\]
 where all depends on the  choice of $x_0$ and the type of rounding, and where we stumble across  cases with primitive Pythagorean triples   interlinked with transactional spactime \[\Big(\{T\rightarrow X\}_\textrm{cc},\{X\rightarrow T\}_\textrm{pc},\{X\rightarrow T\}_\textrm{pc},\{X\rightarrow T\}_\textrm{pc}\Big)\] $-$ which allows more encompassing combined systems to be characterized.   We already noted that many tetrapetalia coordinates  are neutrino-wise nondescript; the same can be said  of $\sigma_\tau$. There are a few notable exceptions, however. A case in point is $2\sigma_{10}=4006$ which, together with $\sigma_{11}=4607$, 
entails an interesting mix of  normal and   `smeared-out' neutrino systems: $(525,527,\widetilde{4605},\widetilde{4607})$. That is, $\sigma_{10}=2303$, first  nondescript when assigned to $(t_\nu)_0$, turns meaningful by way of assigning $(t_\nu)_1:=\lceil\sqrt2\times\sigma_{10}\rceil=3257$. And  the assignments  $x_0:=\sigma_{10}$, $x_1:=\lfloor2\sqrt2\times x_0\rfloor=6513$, $x'_1:=\lceil2\sqrt2\times x_0\rceil=6514$  are no less rewarding, the two being  the hypotenuses of the primitive Pythagorean triples $6513^2=2505^2+6012^2$  and $6514^2=2464^2+6030^2$, respectively. The mixed system that  develops is
\[ \begin{array}{c}\qquad\theta(6514)[624]=525,\quad\bar\theta(4607)[1992]=527,\\
\theta(6513)[224]=4471,\;\bar\theta(6513)[666]=134\;\quad(4471+134=\widetilde{4605}),\\
\;\theta(6513)[1506]=4125,\,\bar\theta(6513)[1647]=482\,(4125+482=\widetilde{4607}),\\
\!\!\!\!\!\theta(6513)[1613]=\bar\theta(6513)[1614]=525,\\
(624\!+\!1992)+(224\!+\!666)+(1506\!+\!1647)+(1613\!+\!1614)\equiv42
\;(\textrm{mod}\,\,107); \\ 
\begin{array}{l}\bar\theta^{-1}(3257)[356]=5351\\ \bar\theta^{-1}(3784)[425]=744\end{array}\left(\begin{array}{c}\lceil(1/\sqrt2)\times5351\rceil=3784,\\5351-\textcolor{red}{744}=\widetilde{4607}\end{array}\right),
\\
\!\quad\beta_{482}-\bar\alpha_{744}=1646-\textcolor{red}{425}\equiv44\;(\textrm{mod}\,\,107), \\
\bar\beta_{4125}+\bar\alpha_{5351}=1507+356\equiv44 \;(\textrm{mod}\,\,107);\\
\;\quad42\equiv42 \;(\textrm{mod}\,\,53),\quad44\equiv44 \;(\textrm{mod}\,\,53).
\end{array} \]
The dichotomy 42 vs.$\,$44 is exactly the one  found  previously  in the atomic context.  More combined systems show up. When assigned to $x_1$, 3257 becomes the hypotenuse of the primitive Pytha\-gorean triple $3257^2=1232^2+3015^2$ and linked to the  combined system  $({286},\!{526},\!{4606})$. And $x_1:=373=\lceil(1/\sqrt2)\times527\rceil $ is the hypotenuse of $373^2=252^2+275^2$ 
and  is linked to the  combined system $({285},{286},{287},{525},{526},{527})$. To say that  neutrinos are nature's way of making transactional  spacetime work demands  the evidence of more scenarios meeting requirements. For details, see body text at the close of \ref{sec:Minkowski} (p.\,\pageref{transpoin}).

\subsubsection{\label{charge-trans} Phase transition and ordered `First 107, then 53'}
 
The Sophie Germain prime 53, together with the safe prime 107,  turns out to be  intertwined with neutrino-based transactional spacetime.  In expansion time, no more than the prime 53 is needed to achieve  neutrino contextualization:\[\begin{array}{lr}
                      \between^{(16)}[149]=286
                                & \qquad\quad\,149 \equiv 43\, (\textrm{mod}\,53),
                           \end{array}\]     
or in slightly varied form:

      \[\begin{array}{lr}
                  \between^{(8)}[159] =\;\between^{(16)}[149] = 286
                                &\quad
                                159+149 \equiv 43\, (\textrm{mod}\,53).
\end{array}\]                 
A more advanced variant consists in a  random integration of  the safe prime 107, accomplished by exploiting phase (see Sect. \ref{sec:croton-phase}):
Instead of $\between^{(8)}[126]=287$,
                  $\between^{(8)}[159]=286,
                  \between^{(8)}[207]=526$,
the phase picture is taken into focus:
\[\begin{array}{rrl}
                 \between_{+-}^{(8)}[119]=287\\
       
                \between_{+-}^{(8)}[149]=287\\
                 \between_{+-}^{(8)}[201]=527\\
                               &\qquad119+149+201 \equiv &41 \,(\textrm{mod}\,107)\\                                                    &\qquad\equiv &45 \,(\textrm{mod}\,53)\;\;\,\\
                              &    & 43 \,(\textrm{on average}).\end{array}\]
Prime 53 seems to stay dominant even if we imbue  the neutrino picture with charge. As we shall learn on page \pageref{weakpoin}, the weak force is closely linked to the kissing number $L_{24}=196560$. The connection becomes visible in the combination of a transcendental entity like the log function and a non-trivial, yet widely considered  normal number such as the Euler constant $\gamma=0.57721..\,$. We find
\[\frac{\gamma \log(2)}{\log(196560)}[1748]=4320,\] which is complemented by \[\frac{2\gamma \log^3(2)}{\log(4320)}[1114]=240.\] Now \[\begin{array}{c}1748\equiv53-1\,(\textrm{mod}\,53)\\\,1114\equiv53+1\,(\textrm{mod}\,53),\end{array}\]
or 2 charges compensating each other nilcongruently $1748+1114\equiv0\,(\textrm{mod}\,53)$. While in that relationship the nature of the (integer-valued) charge remains vague, it definitely takes the form of  {\em electric} charge in case of the bicube complex (\textcolor{blue}{9}+6)(=\textcolor{blue}{15}). By means of the tetrahedral number $T_{12}$ and the triangular number $D_8$,  we get
\[\begin{array}{cc}196560/91=2160\\196560/182=1080\\ 196560/364=540&(T_{12})\\
196560/728=270\\196560/1456=135\\ \\
4320/\textcolor{blue}{9}=480\\ 4320/18=240 \\ 4320/36=120&(D_8)\\4320/72=60\\ 4320/144=30\\ 4320/288=\textcolor{blue}{15}\end{array}\]
(divisor halvings and doublings start at the parenthesized markings and stop if either  a divisor or a quotient    gets  into a $2\mathbb{N}+1$ boundary).
Not only does the reinterpretation
      \[  L_{11} +(L_4+L_5+L_6+L_8+L_9+L_{10}) 
\;\longrightarrow\;
       9(L_7+L_3^++L_1^+)+ 6L_4\]  (in spacetime, phase transition) allow for the representation of {electric} charge, $o_9+\sigma_6\equiv1278+\textrm {e} ^++144+\textrm {e} ^-$, it gets the  amplitudes 
      $\between^{(n)}[\,] $
morphed into amplitudes of transactional spacetime, 
be it        
      $ \theta(\,)[\,]$, $\bar\theta(\,)[\,]$ or $\bar\theta^{-1}(\,)[\,]$,
for which there holds {\em ordered} \label{techpoin}
      `First 107, then 53.'
When the transition is imminent, it comes  under the supervision 
of $(G_{mn})^{(15)}$. 
In a first step,  an optimal container size is determined, according to divisibility properties with respect to pairings or higher-order federations of $L_m$'s. 
The supervision, in a sequence of descending nesting levels, starts with the 1-element side diagonal 
       $\textrm{LL(LL(LL\,}G_{mn}^{(15)})) = 113$:\[
      (L_8+L_{11})/113 = 6.\]
 Which is followed by a test  with the diagonal of 
$\textrm{LL(LL\,}G_{mn}^{(15)})$, $\textrm{ diag}(41,41)$:\label{hfpoin}
                                                      
      \[(L_4+L_5+L_6+L_{11})/41 = 14,
     \quad(L_4+L_5+L_6+L_8+L_9+L_{10})/41 = 24,\]  where the  maximum is 24.
The third in the league  is the diagonal of
     $\textrm{LL\,}G_{mn}^{(15)}$, $ \textrm{diag}(5,5,5,5)$:
             
      \[ \begin{array}{cc}
      (L_5+L_8)/5 = 56,
        & (L_4+L_{10})/5 = 72, \\
       (L_6+L_{11})/5 = 102, 
      &(L_9+L_{11})/5 = 142,
      \end{array}\]
where the winner is 142.
As a result, we get two leading container sizes -- 142 and 24.
With a   payload  of  nine $L_7,L_3^+,L_1^+$ and six  $L_4$, nine containers of size 142  can be filled with  $ L_7+L_3^++L_1^+$ and six containers of size 24 with $L_4$  to achieve  successful transition  to what we  call  1278+144. In a second step,   pairings   are sorted according to their divisibility by 142 or 24. The  incidence matrix below has the  entries 142,\,24 or  \o{}, meaning compatibility with container size 142,\,24  or incompatibility with either, respectively.

\bigskip{}
\hspace*{2cm}\begin{tabular}{|c|ccccccc|}
\hline
          &$L_4$&$L_5$&$L_6$&$L_8$&$L_9$&$L_{10}$&$L_{11}$ \tabularnewline\hline
 $L_4$&24 &\o{} &24 &24 &\o{} &24 &\o{}  \tabularnewline\hline 
 $L_5$&\o{} &\o{} &\o{} &\o{} &24 &\o{} &\o    \tabularnewline\hline
 $L_6$&24 &\o{} &24 &24 &\o{} &24 &\o    \tabularnewline\hline
 $L_8$&24 &\o{} &24 &24 &\o{} &24 &\o    \tabularnewline\hline
 $L_9$&\o{} &24 &\o{} &\o{} &\o{} &\o{} &142    \tabularnewline\hline
  $L_{10}$&24 &\o{} &24 &24 &\o{} &24 &\o{}    \tabularnewline\hline 
  $L_{11}$&\o{} &\o{} &\o{} &\o{} &142 &\o{} &\o{}    \tabularnewline\hline
  \end{tabular}
 
 \bigskip{}
 
 \noindent  $L_9$  is unique -- it has all options at its disposal, which makes it an ideal double agent. That exceptional position is underpinned by the way kissing number contractions (KNC) work. On page \pageref{techpoin}, we noted a  holographic equivalence of the  (9+6) bicube  complex  with $L_{14}$,  expressible by the (pre-transition) equation   $ L_{11} +(L_4+L_5+L_6+L_8+L_9+L_{10}) =L_{14}$, which is nothing but the second in a sequence of KNC's, the first KNC being $L_4+L_6+L_8=L_4+L_5+L_9=L_{10}$. The involved dimensions 
behave as follows {\small
 \[\begin{array}{cccc}\textrm{\tt I}&\textrm{\tt II}&\textrm{\tt III}&\textrm{\tt IV}\\
 {4+6+8}\atop{4+5+9}&10&{5+9+11}\atop{6+8+11}&14\end{array}\]
\noindent }
$\!\!${\tt I+III} sum up to the prominent residue 43; together with {\tt II}, it yields the Sophie Germain prime 53. The moduli to make  \textrm{\tt I} and \textrm{\tt III} congruent to \textrm{\tt II} and \textrm{\tt IV} are 8 and 11, respectively. If, in analogy to a Sophie Germain situation,  a `First 11, then 8' rule is used, the only kissing number to turn out nilcongruent is $L_9$, no matter whether the context is \textrm{\tt I} or  \textrm{\tt III}:
\[\begin {array} {rcc}272&\equiv8&(\textrm{mod }11),\\
                               8&\equiv0&(\textrm{mod }8).
                               \end{array}\]
Distilled from \[\begin{array}{ccr}(L_4+\textcolor{black}{L_5}+L_6+\textcolor{black}{L_{11}})/41 &= &14,\\
   (L_4+\textcolor{black}{L_5}+L_6+L_8+\textcolor{black}{L_9}+L_{10})/41 &= &24,\\
  (\textcolor{black}{L_9}+\textcolor{black}{L_{11}})/5 &= &142 \end{array}\quad\bigg
  \}\;180\] 
and the KNC  $L_4+L_6+L_8=L_{10}$, the double agent can be attributed a   volume  $\frac{180}{4}\times L_9$. It has four equivalents,   $(45\times L_9=)\,510\times L_4=
306\times L_5=170\times L_6=51\times L_8$.  
Proceeding with
\[\begin{array}{c}(L_5+\textcolor{black}{L_{8}})/5 =56,
\end{array}
\]
we get a second  volume, $56\times L_8$, possessing the equivalents $560\times L_4=336\times L_{5}=40\times L_{10}$. (45 is the  triangular number $D_9$,  56 is the  tetrahedral  number $T_6$, which together yield the nice correspondence $o_9\sim D_9$, $\sigma_6\sim T_6$.)
A third  volume originates with
\[\begin{array}{ccr}(L_5+\textcolor{black}{L_{8}})/5 &=&56,
\\
   (L_4+L_{10})/5 &= &72,\\
  ({L_6}+\textcolor{black}{L_{11}})/5 &= &102, 
\end{array}
\quad\bigg  \}\;180
\] and the KNC $L_4+L_5+L_9=L_{10}$, and is attributed to $L_{11}$, $136\times L_{11}(=T_{16}\times73$). It has the equivalents  $2482\times L_4=219\times L_9$. We see that the double agent's volume is standard since it has the most equivalents.

\smallskip{}
\noindent For  large enough $\theta()[\,]$, $\bar\theta()[\,]$, 
$\bar\theta^{-1}()[\,],$ partitioning into configurations 1278+144 and a  remainder $R$ becomes possible, but requires dragging the outlined provision of standard volume   alongside in the process:

       \[ A_\textrm{\small large} = c(9\times142 + 6\times24 + 45\times272)  + R\qquad
                                            (c \in \mathbb{N}_0).\]
For instance, 

    \[   \bar\theta^{-1} (5200)[236] =
       3\,860\,638 = 282(1422+45\times272) + 54\times142 + 286,\]
where the remainder assumes the interesting form 
                   $    6\times1278 + \textrm{ neutrino } 286$.
\bigskip{}

\noindent A modification is warranted in case of odd amplitudes. For bookkeeping reasons, 639--half of 1278--arises in one of two forms,  $5L_7+3L_1^+$ or  $4L_7+9L_3^++6L_1^+.$  Even distribution among those forms means smooth transition -- every twin pair 639+639 behaves like 1278. When the distribution is uneven and the bookkeeping gets out of balance,  the $L_m$s' inherent tendency to try transitioning individually is furthered -- with dramatic consequences, as we shall see.

\noindent If  639 signals bipartition, 426  stands for tripartition, or black boxes that consist of three 142's each. Whichever variant prevails, in order to count as successful it must be supplemented by  144.
This characteristics is encoded in the 639th prime number 4733 which yields
     \[  \bar\theta^{-1} (4733)[426]=144.\]
With respect to the smallest summands present in 639, $3L_1^+ (=9)$ or $6L_1^+ (=18)$, a case of $L_4$ wandering from even- to odd-assignment partners is observed. The  federations
$(L_6+L_8+L_{10})(=72\times9)$,$(L_5+L_9+L_{11} + L_4)(=43\times18)$
come up with container sizes whose  sum  72+43   is equivalent to the 1-element side diagonal 
 $\textrm{LL(UR(UL(LL\,}G_{mn}^{(31)})))$=115, and 43  to  the diagonal of $\textrm{UR(UL(LL\,}G_{mn}^{(31)}))$   ($\textrm{diag }(43,43)$, see $\,$\ref{sec:Crotons-as-boundary}). Such  passing of  control  from normal  to { interordinal}  $(G_{mn}) $ is accompanied by a switch of the partitioning mode -- from the even-$A_\textrm{\small large}$ mode  to the odd-$A_\textrm{\small large}$:
  \[ A_\textrm{\small large}= c(9\times142 + 6\times24 + 45\times272) + d\times639 + R\quad
                                            (d \textrm{ odd};\;c,R \in2\mathbb{N}_0).\]
For instance, 
  \[  \begin{array}{l}
      
  \bar\theta^{-1} (9221)[35] =
    1\,078\,923 \\
    \qquad= 72(1422+45\times272) + 107\times639 +
                           4\times1286 + 126\times136 + 4606\\
                   \qquad = 72(1422+45\times272) + 107\times639 +
                           4\times1286 + 156\times136 + 526\\
                    \qquad = 72(1422+45\times272) + 87\times639 +
                          10\times1286 + 195\times136 + 286.\end{array}\]
The amplitude's odd value is carried on to an odd multiple of 639, leaving one 639 unbalanced.  A tendency towards an individualistic type of transition  then causes  $L_4$, depending on which variant of 639 it encounters, to pursue a distinct ill-opted partnership, rather than partaking in the morphing  process as a whole. With a limit of  in total nine $ L_1^+, L_3^+, L_7$  and six $L_4$, the result includes an exotic form of matter,
\[\begin{array}{lll}L_4&\longrightarrow&L_4\\
L_5&\longrightarrow&L^+_1+L^+_3+L_4\\
L_6&\longrightarrow&3L_4\\
L_8+L_{11}&\longrightarrow&3 L^+_1+3L^+_3 +5 L_7\\
L_9+L_{10}&\longrightarrow&5 L^+_1+5L^+_3+L_4+4 L_7,\end{array}\]
1286+136, for short, whose creation -- as the above amplitude bears out -- goes hand in hand with that of ordinary matter 1278+144, while at the same time managing neutrino oscillation with ease.\footnote {$R$ for remainder might actually be  a misnomer.} Now, in case of even $A_\textrm{\small large}$  a slightly modified partitioning suffices to bring neutrino 536 to the fore too:
\[3\,860\,638 = 282(1422+45\times272) + 
                     39\times142 + 22\times144 + 526;\]
the inclusion of neutrino 4606, however, implies  less ordinary matter:$\,$\footnote {$\,39\times142+22\times144$ and $63\times142+56\times144$ may be named metastable ordinary matter}
\[3\,860\,638 = 281(1422+45\times272) +
                     63\times142 + 56\times144 + 4606.\]
We have scrutinized $13\,000$ CFRs of 
       \[\begin{array}{cc}
           
       \bar\theta^{-1} (t_\nu) &\qquad ( \bar\theta(x)=t_\nu; \;\;1\le t_\nu\le 13\,000), 
      \end{array}\]
$\bar\theta^{-1}()[\,] $, for short. What we found is that, atop eight instances of $\bar\theta^{-1}()[\,]=1422,$ there are three outcomes associated with ordinary matter $(\bar\theta^{-1}()[\,]=1278)$, and five associated with exotic matter $(\bar\theta^{-1}()[\,]=1286).$ While the uncovered copies of ordinary matter yield the expected result under `First 107, then 53,' \[\begin{array}{rr}
\bar\theta^{-1}(4906)[454]=1422 &\bar\theta^{-1}(3766)[33]=1278 \\
\bar\theta^{-1}(9736)[240]=1422 &\bar\theta^{-1}(6198)[478]=1278 \\
\bar\theta^{-1}(10807)[434]=1422&\bar\theta^{-1}(9312)[330]=1278\\ 
\multicolumn{2}{c}{(454+240+434)+(33+478+330)=1969\equiv 43\, (\textrm{mod }107),}
\end{array}\]those of  exotic matter deviate, 
\[\begin{array}{rr}
\bar\theta^{-1}(8838)[228]=1422 &\bar\theta^{-1}(9419)[494]=1286 \\
\bar\theta^{-1}(8346)[361]=1422 &\bar\theta^{-1}(7673)[413]=1286 \\
\bar\theta^{-1}(4912)[302]=1422 &\bar\theta^{-1}(7632)[271]=1286 \\
\bar\theta^{-1}(3587)[285]=1422 &\bar\theta^{-1}(6704)[292]=1286 \\
\bar\theta^{-1}(2665)[311]=1422 &\bar\theta^{-1}(6501)[334]=1286 \\
\multicolumn{2}{c}{(228+361+302+285+311)+(494+413+271+292+334)} \\
\multicolumn{2}{c}{=3291\equiv 81\, (\textrm{mod }107)} \\
\multicolumn{2}{c}{\qquad81\equiv 28\, (\textrm{mod }53),} 
\end{array}\] 
but the failure is a systematic one: The value  $45=180/4$ proved to be constitutive for the establishment of a standard volume.
With five copies of 1286 and three copies of 1278 (in eight  of 1422),  45 now turns a  sum of   multiples  (red) which   fail to make, and (blue) which succeed to make 3291 congruent to $1969$ under the Sophie Germain relationship\label{SGrelpoin}\[\begin{array}{c}3291-\textcolor{red}{30}\times107\,-\,\, \,1\times53\equiv28, \\
3291-\textcolor{blue}{15}\times107-{31}\times53\equiv43.
\end{array}
\]
\noindent Before going any further with that partitioning, let us point out that  exotic matter doesn't have to remain exotic.
In neutrino-based transactional spacetime, a rollback of  ill-formed pairs becomes possible through a particular juxtaposition $d\textrm{ vs. }\!d-2\,$: a decrease of the nesting depth  by two levels -- for the ill-formed pair $L_8+L_{11}\,(=6\times113)$:  from LL(LL(LL$G^{(15)}_{mn}))=113$ to LL$G^{(15)}_{mn}$;  for the 
ill-formed pair  $L_9+L_{10}\,(=32\times19$): from  $\textrm{UR(UR(UL(LL}G_{mn}^{(31)})))\!=\!19$ to LL(LL$G^{(15)}_{mn}$).
Inspired by the tripartition of 1278, that can be demonstrated by way of the template
 \[\begin{array}{cccc}
 \textrm{\small426}&\textrm{\small426}&\textrm{\small426}&\textrm{\small144}\\
 3(L^+_1+L^+_3+L_7),&3(L^+_1+L^+_3+L_7),&3(L^+_1+L^+_3+L_7),&6L_4,
\end{array}\]
which comes alive with two mappings $\{X\rightarrow T\}$ and two  mappings $\{T\rightarrow X\}$ that are linked  to one another by a pointer chain:

 \[\scriptsize\begin{array}{cccc}
\textcolor{lightgray}{426}&\textcolor{lightgray}{426}&\textcolor{red}{426} &\textcolor{lightgray}{144}\\
&&{\theta(426)[1670]\!=\!144}\\ \\
\textcolor{lightgray}{426}&\textcolor{lightgray}{426}&\textcolor{lightgray}{426}
&\textcolor{red}{144}\\
&&&\theta(144)[1663]\!=\!240 \\ \\
\textcolor{red}{426}&\textcolor{lightgray}{426}&\textcolor{lightgray}{426}&\textcolor{lightgray}{144}\\
{\bar\theta^{-1}(240)[239]\!=\!139}\\ \\
\textcolor{lightgray}{426}&\textcolor{red}{426}
&\textcolor{lightgray}{426}&\textcolor{lightgray}{144}\\
&\bar\theta^{-1}(139)[106]\!=\!426\\
&
\end{array}\]
The  chain has four pointers, each  of which pointing to a different   option of how to break free (in blue). The most obvious endeavors are with 426 and 144, ordinary-matter containers  both.

426 (blue):
\[\begin{array}{lll}L_4&\longrightarrow&L_4\\
L_5&\longrightarrow&L^+_1+L^+_3+L_4\\
L_6&\longrightarrow&3L_4\\
L_8+L_{11}&\longrightarrow&
\textcolor{blue}{3 L^+_1+3L^+_3 +3 L_7}+2L_7\\
L_9+L_{10}&\longrightarrow&5 L^+_1+5L^+_3+L_4+4 L_7,\end{array}\]

 144 (blue):\[\begin{array}{lll}L_4&\longrightarrow&L_4\\
L_5&\longrightarrow&\textcolor{blue}{L^+_1+L^+_3}+L_4\\
L_6&\longrightarrow&3L_4\\
L_8+L_{11}&\longrightarrow&
\textcolor{blue}{3 L^+_1+3L^+_3} +5 L_7\\
L_9+L_{10}&\longrightarrow&
\textcolor{blue}{5 L^+_1+5L^+_3}+L_4+\textcolor{black}{4 L_7}.\end{array}\]


\noindent The other two pointers signify attempts at secondary transformations: For 
$L_8+L_{11}$, the objective of a rollback could be   exchanging one $L_8$ for five $L_5$:

\smallskip{}

240 (red and blue):
 \[\quad5(1286+L_5)-\textcolor{red}{L_8}=5\times1278,\textrm{ and } 5(136-L_5)+\textcolor{blue}{L_8}=5\times144;\]

\noindent unfortunately,   only one $L_5$ is available, so this remains an attempt. What is sauce for the goose is sauce for the gander -- a rollback of $L_8+L_{11}$ would induce a rollback of $L_9+L_{10}$, followed by $L_9$ trying to recruit another $L_{11}=438$ in order to transform {\tt ud} to {\tt uu}:

\bigskip{}

139 (blue):
\[L_9+4L^+_1+\textcolor{blue}{5L^+_1+L^+_3+L_7}=L_{11};\] 

\noindent in transactional spacetime, however, neither one nor two $L_{11}$'s can exist by themselves. That attempt is fruitless as well.

\

\noindent The pointer quartet possesses  invariants. After decrease of nesting depth  by two levels,  differences of the quotients that were competing for maximal container size  (page \pageref{hfpoin}) can  be formed,  from pair-associated quotients ($\Delta_\textrm{\tiny p}(q_a,q_b)$) as well as  from the  outcomes of higher-order federations ($\Delta_\textrm{\tiny hf}(q_a,q_b)$):
\[\begin{array}{lc}
 \textrm{for diag(5,5,5,5): }&\Delta_\textrm{\tiny p}(72,56)=16,\quad
 \Delta_\textrm{\tiny p}(142,102)=40,\\
 
 \textrm{for diag(41,41): }&\Delta_\textrm{\tiny hf}(24,14)=10.
\end{array}\]
The differences become invariants for the pointer chain under `First 107, then 53:'
\[\begin{array}{rrl}
    
 426+144+240+139&\equiv &93\, (\textrm{mod}\,107)\\                                                    93&\equiv &40 \,(\textrm{mod}\,53)\;\;\,\\
 1670+1663+239+106&\equiv &40\, (\textrm{mod}\,107)\\
                              & \equiv   &
 \Delta_\textrm{\tiny p}(142,102)\\
                          \end{array}\]

\noindent This could go on forever -- with $144+240+139+426,\,1663+239+106+1670$, and further cyclic permutations. Yet, the chain may steer a different course after the first two steps  -- when the unreferenced neutrino 286 is followed and simple linear combinations of $\Delta_\textrm{\tiny p}(\,),\Delta_\textrm{\tiny hm}(\,)$ are exploited. The best way to explain the ongoings is to track the chain's  evolution in slow motion:

\[\scriptsize\begin{array}{cccc}
\textcolor{lightgray}{426}&\textcolor{lightgray}{426}&\textcolor{red}{426}&\textcolor{lightgray}{144}\\
&&{\theta(426)[1670]\!=\!144}\\
\textcolor{lightgray}{426}&\textcolor{lightgray}{426}&\textcolor{lightgray}{426}&\textcolor{red}{144}\\
&&&\theta(144)[1663]\!=\!240 \\ \\
\textcolor{red}{426}&\textcolor{lightgray}{426}&{426}&\textcolor{lightgray}{144}\\
{\bar\theta^{-1}(240)[49]\!=\!286}\\ \\
\textcolor{lightgray}{426}&\textcolor{red}{426}&\textcolor{lightgray}{426}&\textcolor{lightgray}{144}\\
&{\bar\theta^{-1}(286)[101]\!=\!6471}\\
\end{array}\]
with  yet-to-establish invariants
\[\begin{array}{rrl}
    
 426+144+240+286&\equiv &26\, (\textrm{mod}\,107)\\                                                    
 26&\equiv & \Delta_\textrm{\tiny p}(72,56)+\Delta_\textrm{\tiny hm}(24,14)\\
 1670+1673+49+101&\equiv &59\, (\textrm{mod}\,107)\\              59& \equiv   &
\;\,6\,(\textrm{mod}\,53)\\
                            6  & \equiv   &\Delta_\textrm{\tiny p}(72,56)-\Delta_\textrm{\tiny hm}(24,14).               \end{array}\] 
                             
                              \bigskip{}
                              
                              \noindent Go one step further

\[\scriptsize\begin{array}{cccc}
\textcolor{lightgray}{426}&\textcolor{lightgray}{426}&\textcolor{lightgray}{426}&\textcolor{red}{144}\\
&&&{\theta(144)[1663]\!=\!240}\\ \\
\textcolor{red}{426}&\textcolor{lightgray}{426}&\textcolor{lightgray}{426}&\textcolor{lightgray}{144}\\
\bar\theta{-1}(240)[49]\!=\!286 \\ 
\textcolor{lightgray}{426}&\textcolor{red}{426}
&\textcolor{lightgray}{426}&\textcolor{lightgray}{144}\\
&{\bar\theta^{-1}(286)[101]\!=\!6471}\\ \\
\textcolor{lightgray}{426}&\textcolor{lightgray}{426}&\textcolor{red}{426}&\textcolor{lightgray}{144}\\
&&{\theta(6471)\left\{\begin{array}{c}
    \;\textrm{[}1602\textrm{]}=\;144\\                 \quad\textrm{[}816\textrm{]}=1278
\end{array}\right.}
\end{array}\]
and you get the surprising result: 
\[\begin{array}{rrl}
    
 144+240+286+6471&\equiv &26\, (\textrm{mod}\,107)\\                                                    
 26&\equiv &\Delta_\textrm{\tiny p}(72,56)+\Delta_\textrm{\tiny hm}(24,14)\\ \\
\textrm{(i) } 1673+49+101+1602&\equiv &98\, (\textrm{mod}\,107)\\              98& \equiv   &
45\,(\textrm{mod}\,53)\\
\textrm{(ii) }\;1673+49+101+816&\equiv &61\, (\textrm{mod}\,107)\\              61& \equiv   &
\,8\,(\textrm{mod}\,53)\\ \\
\frac{(1673+49+101+1602)+(1673+49+101+816)}{2}&\equiv &26\, (\textrm{mod}\,107)\\                                                            
26&\equiv & \Delta_\textrm{\tiny p}(72,56)+\Delta_\textrm{\tiny hm}(24,14).
 %
   \end{array}\]  
What started out with a simply covered  chain  ends up in double cover ({\tt i+ii}),  if this case of successful secondary transformation $1286+136\rightarrow 1278+144$ is to obey   invariant 26.

\noindent One further  side note on neutrinos is appropriate. We learned that $ L_{16}(=4320) $ is directly related to neutrino 526:
        $  \between^{(14051)}[526]=4320$ --
see page \pageref{pagepoin}, where we also noted that 14051 is a prime number. A similar cornerstone exists in case of neutrino-based transactional spacetime: 
The 1278th prime number 10453 encodes the relationship 
       \[\bar\theta^{-1}(10453)[53]=4320.
       \]
With an MF-wise exposition $\,\bar{B}_{t_\nu}^{(z)}:=2^{-z}\bar\theta^{-1}(10453)\;(z\in\mathbb{Z})$, that relationship extends to
        \[ \begin{array}{c}\!\qquad \bar{B}_{t_\nu}^{(-12)}[295] = 286,\quad\bar{B}_{t_\nu}^{(8)}[26] = 526\\ \\
\!\!\qquad\;53\equiv0\,(\textrm{mod }\;\,53)\\
\!\qquad\textrm{\hspace{2mm}}295+26\equiv0\,(\textrm{mod }107).\textrm{\hspace{1cm}}
\end{array}\]
The new code not only confirms what we know about neutrinos and  the role of the ordered Sophie Germain  condition `First 107, then 53'
as a constitutional element of transactional spacetime -- it has  implications that go beyond Sophie Germain. Provided  that $\theta() $, $\bar\theta()$ adhere to the $1/\sqrt2$ order and 
defining $B_x^{(z)}=2^{-z}\theta(x)$, $\bar B_x^{(z)}=2^{-z}\bar\theta(x)$, 
one makes the observation: after they reach the peak, MF terms $B_x^{(z)}[\,]$, $\bar B_x^{(z)}[\,]$  break out of their stationary refinement position [1] or violate stationarity in  [1] or [2] in some other way. There's more to it:  resuming the argument of how  the multiples \textcolor{red}{30} and \textcolor{blue}{15} gain criticality with respect to congruence  (page \pageref{SGrelpoin}), the stationarity numbers (\# of MF terms that are faithful to stationarity) turn out to be dual in that respect. 
Earlier, we noted that one {\tt u+d} alone is inappropriate for representing  a nucleon. Instead, three {\tt u+d}'s are needed to match the deuterium's proton-neutron pair, which amounts to 45 containers 142. This relates to and at the same details what we -- once combined -- register at the  innerts of the MF-like expositions of $\theta(\lfloor(4320/\sqrt2\rfloor) $ and $\bar\theta(\lceil(4320/\sqrt2\rceil)$:   $\textcolor{blue}{15}$ $B_x^{(z)}[\,]$ stationary on position [1],  $\textcolor{red}{30}$ $\bar B_x^{(z)}[\,]$ stationary on position [2]
 while maintaining unity on position [1] -- both until they peak, which suggests the deuterium's nucleons form in an asymmetrical manner dictated by this type of MF-wise confined $\theta$-$\bar\theta$ environ. 
While the  previous analysis  tells us that the same numbers appear in the production ratio of $\bar\theta^{-1}(\,)[\,]=1278$ vs. $\bar\theta^{-1}(\,)[\,]=1286$ with and without congruence to one another:



\hspace*{2.0cm}
\begin{tabular}{|c|c|c|} 
\hline
\multicolumn{3}{|c|}{$\,B_x^{(z)}:=2^{-z}\theta(3055)$}
\tabularnewline
\multicolumn{3}{|c|}{\small{$\bar{B}_{x}^{(z)}:=2^{-z}\bar\theta(3056)$}}
\tabularnewline
\multicolumn{3}{|c|}{
}
\tabularnewline
$\,z\,$ &$B_x^{(z)}[1]$&$\bar{B}_x^{(z)}[1]$,$\bar{B}_x^{(z)}[2]$
\tabularnewline\hline\hline
-12&     $\quad1$& $-$
\tabularnewline\hline\hline
-11&     $\quad2$&1, $\;1$
\tabularnewline\hline
-10&     $\quad4$&1,$\;3$
\tabularnewline\hline
-9&      $\quad8$&1,$\;7$
\tabularnewline\hline
-8&     16&1,$\;15$
\tabularnewline\hline
-7&     33&1,$\;32$
\tabularnewline\hline
-6&     67&1,$\;66$
\tabularnewline\hline
-5&     135&1,$\;134$
\tabularnewline\hline
-4&     270&1,$\;269$
\tabularnewline\hline
-3&     540&1,$\;539$
\tabularnewline\hline
-2&    1080&1,$\;1079$
\tabularnewline\hline
-1&   2160&1,$\;2159$
\tabularnewline\hline
0&    4320&1,$\;4320$
\tabularnewline\hline
$1$&    8640&1,$\;8642$
\tabularnewline\hline
$2$&  $\quad\qquad17281\; (\textcolor{blue}{\tt peak})$& $\;1,17286$
\tabularnewline\hline

\multicolumn{3}{|c|}{}
\tabularnewline
\multicolumn{3}{|c|}{\small continued on next page}
\tabularnewline
\multicolumn{3}{|c|}{}
\tabularnewline\hline
 \end{tabular}

\begin{tabular}{|c|c|c|} 
\hline
\multicolumn{3}{|c|}{}
\tabularnewline
\multicolumn{3}{|c|}{$\,B_x^{(z)}:=2^{-z}\theta(3055)$}
\tabularnewline
\multicolumn{3}{|c|}{$\bar{B}_{x}^{(z)}:=2^{-z}\bar\theta(3056)$}
\tabularnewline
\multicolumn{3}{|c|}{}
\tabularnewline
$\,z\,$ &$B_x^{(z)}[\,]$&$\bar{B}_x^{(z)}[1]$,$\bar{B}_x^{(z)}[2],\dots$
\tabularnewline\hline\hline
$3$&  8640\,@\,3&1,$\;34573$
\tabularnewline\hline
$4$&  4319\,@\,3&$\cdots$
\tabularnewline\hline
$5$&  2159\,@\,3&
\tabularnewline\hline
$6$&  1079\,@\,3&
\tabularnewline\hline
$7$&  539\,@\,3&
\tabularnewline\hline
$8$&  269\,@\,3&
\tabularnewline\hline
$9$&  134\,@\,5&
\tabularnewline\hline
$10$&  67\,@\,5&
\tabularnewline\hline
$11$&  33\,@\,5&
\tabularnewline\hline
$12$&  16\,@\,11&
\tabularnewline\hline
$13$&   7\,@\,13&
\tabularnewline\hline
$14$&   3\,@\,9&$\cdots$
\tabularnewline\hline
$15$&   1\,@\,7&1,$\;141617946$
\tabularnewline\hline
$16$&{\tiny short\,right\,leg}          &1,$\;283235893$
\tabularnewline\hline
$17$&          &1,$\;566471787$
\tabularnewline\hline
$18$&          &$\!\!\qquad\qquad1$,$\;1132943575\; (\textcolor{red}{\tt peak})$
\tabularnewline\hline
$19$&          &2,$\;566471787$
\tabularnewline\hline    
$\cdots$&       &$\cdots$     
\tabularnewline\hline
$26$&.         &256, $\;4425560$  
\tabularnewline\hline
$27$&     & $1, 1, 127, 1, 1, 2212779$ @6
\tabularnewline\hline
$28$&     & $ 3, 1, 63, 3, 1, 1106389$ @6
\tabularnewline\hline
$29$&     & $ 1, 1, 1, 2, 31, 1, 1, 1, 2, 553194$ @10
\tabularnewline\hline
$30$&     & $ 1, 4, 2, 1, 15, 3, 5, 276597$ @8
\tabularnewline\hline
$31$&     & $ 1, 9, 1, 2, 7, 1, 1, 1, 10, 138298$ @10
\tabularnewline\hline
$32$&     & $ 2, 4, 1, 5, 3, 1, 4, 1, 4, 1, 1, 69148$ @12
\tabularnewline\hline
$33$&     & $ 1, 2, 1, 1, 1, 11, 1, 1, 10, 1, 1, 1, 3, 34574$ @14
\tabularnewline\hline
$34$&     & $  2, 1, 3, 24, 22, 3, 1,1, 1, 17286$ @10
\tabularnewline\hline
$35$&     & $  1, 2, 7, 12, 44, 1, 1, 4, 1, 8642$ @10
\tabularnewline\hline
$36$&     &  $ 1, 5, 3, 1, 1, 5, 1, 1, 21, 1, 3, 2, 2, 4321$ @14
\tabularnewline\hline
$\cdots$&  &$\cdots$
\tabularnewline\hline
$41$&     &   $ 1, 1, 1, 6, 1, 1, 4, 3, 1, 21, 1, 1, 5, 1, 2, 11, 1, 134$ @18
\tabularnewline\hline
$42$&     &   $ 3, 14, 9, 1, 1, 44, 12, 2, 1, 5, 2, 67$ @12
\tabularnewline\hline
$43$&     &   $ 1, 1, 1, 28, 4, 1, 3, 22, 24, 1, 2, 2, 1, 2, 1, 33$ @16
\tabularnewline\hline
$44$&     &   $  1, 4, 1, 13, 1, 1, 1, 1, 7, 11, 49, 2, 2, 6, 1, 16$ @16
\tabularnewline\hline
$45$&     &   $   1, 10, 1, 6, 3, 3, 3, 1, 1, 5, 98, 1, 4, 3, 2, 8$ @16
\tabularnewline\hline
$46$&     &   $  1, 22, 1, 2, 1, 1, 1, 6, 1, 1, 3, 2, 1, 1, 48, 1, 9, 1, 1, 2, 1, 3$ @22
\tabularnewline\hline
$47$&     &   $   2, 11, 2, 1, 3, 14, 7, 1, 3, 1, 23, 1, 20, 6, 1, 1$ @16
\tabularnewline
& &
\tabularnewline\hline
\end{tabular}

\newpage\noindent 
{Regarding controls, we may add that $\theta$'s MF terms would suffice to represent all engagements of  LL$G_{mn}^{(15)}$. In case of more involved nucleonic formation, $\bar\theta$'s MF terms would be drawn in instead.  For instance, to conform to an engagement of 
$\textrm{LL(LL(LL(LL}G_{mn}^{(31)})))= 2\,430\,289,$ nine different $\bar B_x^{(z)}[2]$ would be drawn in: 
 
\hspace*{0.2cm}{\small$2\,430\,289 =
 2\,212\,779+138\,298+69148+8642+1079+269+66+7+1.$}}

\noindent A treatise on Mersenne fluctuations cannot be complete without a hypothesis about the nature of life, however tentative.
The relationship between  MF-wise exposition $\bar  B_{t_\nu}^{(z)}:=2^{-z}\bar\theta^{-1}(t_\nu)$ and the conjugate tetrapetalia   is  peculiar. As  $z$  progresses from $\mathbb{Z}^-$ to $\mathbb{Z}^+$, it  provides   $\{X\rightarrow T\}_\textrm{pc}$ processes retrogradely with a rich canvas of  stimuli:\footnote
{\,from the eight copies of 1422 found among the 13\,000 CFRs of $\bar\theta^{-1}(t_\nu)$, six peak in their respective MF expositions; the  two that don't, have MF expositions which for $z<0$ realize an asymmetric subdivision of a pointer to the important transformation $\theta(34143)[\,]$  dealt with intensively on page \pageref{34143poin}: where $\bar  B_{t_\nu}^{(z)}:=2^{-z}\bar\theta^{-1}(2665),\bar  B_{t_\nu'}^{(z)}:=2^{-z}\bar\theta^{-1}(3587)$, $\bar  B_{t_\nu}^{(-3)}[\,]=\textcolor{blue}{11381},\bar B_{t_\nu'}^{(-4)}[\,]=\textcolor{red}{22762}.$}
Far-flung places on the conjugate tetrapetalia's $x=\pm y$  corners   get tightly packed together, from the distant future  to  the causality gap 
 to the complex (largely
imaginary) past. 
A rudimentary definition of life then goes as follows:

\noindent Proposition:

\vspace*{-0mm}\hspace*{0.1cm}\begin{minipage}{9.9cm} Life consists of a standby process 
activated by some amplitude $\bar\theta^{-1}(t_\textrm{act})[\,]$, provided a preselected MF emulation's resemblance to the  MF-wise  exposition of $\bar\theta^{-1}(t_\textrm{act})[\,]$ surpasses a  specifiable threshold.\end{minipage}

\noindent There exist restrictions:

\vspace*{-0mm}
\noindent(i)\,\, \begin{minipage}{10cm}{The MF emulation distills its 
 terms exclusively  from $ \theta() $ and $\bar\theta() $  adhering to the $1/\sqrt2$ order;
the emulation has no feet: $\theta(1)[1]$, $\bar\theta(1)[2]$  do not exist, only $\theta(2)[1]=3$ and $\bar\theta(2)[2]=2$.}\end{minipage}

\noindent(ii)\, \begin{minipage}{10cm}
Whether or not the MF-wise exposition of $\bar\theta^{-1}(t_\textrm{act})$ has symmetric legs,   only a fragment of the MF emulation will match -- ideally, the lower cut-off will be at  $\bar\theta(2)[2]=2$ and $\theta(2)[1]=3$. \end{minipage}

\noindent(iii)\,\begin{minipage}{10cm}
An MF emulation may pursue different  paths. Promising or not, they must meet the premise  that left-leg node and right-leg node assignments  are in that or reverse order: 

\smallskip{}
\noindent {\tt first}
$\qquad\qquad\begin{array}{c}\bar\theta(x_0)[2] \textrm{ \,and \,}\theta(x_0)[1],\end{array}$ 

\smallskip{}
\noindent {\tt then}
$\qquad\begin{array}{c}\bar\theta(\textrm{succ}(x))[2] \textrm{ \,and \,} \theta(\textrm{succ}(x))[1];\end{array}$

\smallskip{}
\noindent where 
$\quad x_0=\big\lfloor\frac{1}{\sqrt2}\bar\theta^{-1}(t_\textrm{act})[\alpha]\big\rfloor$ or 
$x_0=\big\lceil\frac{1}{\sqrt2}\bar\theta^{-1}(t_\textrm{act})[\alpha]\big\rceil$
\smallskip{}

\noindent and, as adapted from Eqs.\,\ref{eq:left-leg} and \ref{eq:right-leg}, 

\[\textrm{succ}(x)=\left\{
\begin{array}{lll}\lfloor x/2\rfloor-\delta&\textrm{\small towards feet}&\quad\!(\delta\in\{0,1\},\;3<x)\\
2x+1+\epsilon&\textrm{\small else}&\quad\;(\epsilon\in\{-1,0,1\}).
\end{array}\right.\]


\end{minipage}
\smallskip{}

\noindent With the MF-wise exposition of $\bar\theta^{-1}(t_\textrm{act})[53]=4320$, where $t_\textrm{act}=10453$, 
\smallskip{}

\renewcommand{\arraystretch}{1.2}
\noindent a{\scriptsize)}\begin{tabular}{c}
$\bar B_{t_\textrm{act}}^{(4)}[\,]=69146$\tabularnewline
$\bar B_{t_\textrm{act}}^{(3)}[\,]=34573$\quad
$\bar B_{t_\textrm{act}}^{(5)}[\,]=34572$\tabularnewline
$\bar B_{t_\textrm{act}}^{(2)}[\,]=17286$\quad\quad\quad
$\bar B_{t_\textrm{act}}^{(6)}[\,]=17286$\tabularnewline
$\bar B_{t_\textrm{act}}^{(1)}[\,]=8642$\quad\quad\quad\quad\quad
$\bar B_{t_\textrm{act}}^{(7)}[\,]=8643$\tabularnewline
$\bar B_{t_\textrm{act}}^{(0)}[\,]=4320$\quad\quad\quad\quad\quad\quad\quad
$\bar B_{t_\textrm{act}}^{(8)}[\,]=4321$\tabularnewline
$\bar B_{t_\textrm{act}}^{(-1)}[\,]=2160$\quad\quad\quad\quad\quad\quad\quad\quad\quad
$\bar B_{t_\textrm{act}}^{(9)}[\,]=2160$\tabularnewline
$\bar B_{t_\textrm{act}}^{(-2)}[\,]=1079$\quad\quad\quad\quad\quad\quad\quad\quad\quad\quad\quad
$\bar B_{t_\textrm{act}}^{(10)}[\,]=1080$\tabularnewline
$\bar B_{t_\textrm{act}}^{(-3)}[\,]=539$\quad\quad\quad\quad\quad\quad\quad\quad\quad\quad\quad\quad\quad
$\bar B_{t_\textrm{act}}^{(11)}[\,]=539$\tabularnewline
$\bar B_{t_\textrm{act}}^{(-4)}[\,]=269$\quad\quad\quad\quad\quad\quad\quad\quad\quad\quad\quad\quad\quad\quad\quad
$\bar B_{t_\textrm{act}}^{(12)}[\,]=269$\tabularnewline
$\bar B_{t_\textrm{act}}^{(-5)}[\,]=134$\quad\quad\quad\quad\quad\quad\quad\quad\quad\quad\quad\quad\quad\quad\quad\quad\quad
$\bar B_{t_\textrm{act}}^{(13)}[\,]=134$\tabularnewline
$\!\!\bar B_{t_\textrm{act}}^{(-6)}[\,]=66$\quad\quad\quad\quad\quad\quad\quad\quad\quad\quad\quad\quad\quad\quad\quad\quad\quad\quad\quad
$\bar B_{t_\textrm{act}}^{(14)}[\,]=66$\tabularnewline
$\!\!\bar B_{t_\textrm{act}}^{(-7)}[\,]=32$\quad\quad\quad\quad\quad\quad\quad\quad\quad\quad\quad\quad\quad\quad\quad\quad\quad\quad\quad\quad\quad
$\bar B_{t_\textrm{act}}^{(14)}[\,]=33$\tabularnewline
$\!\!\dots$ \qquad \qquad \qquad\qquad  \qquad \qquad\qquad  \qquad \qquad \qquad \qquad \qquad\qquad$\dots$
\end{tabular}\newline


\smallskip {}

\noindent and  one of the preselected MF emulations meeting (i)-(iii) equaling 

\noindent \newline
b{\scriptsize)}\begin{tabular}{c}
$\qquad\textcolor{brown}{\bar\theta({48896})[2 ]={69148}\,(+2)}$\tabularnewline
$\quad\textcolor{brown}{\bar\theta(24448)[2 ]={34573}(+0)\; \theta({24448})[1 ]={34574}(+2)}$\tabularnewline
$\quad\bar\theta(12224)[2 ]=\textcolor{black}{17286}\,(+0)\quad\, \theta(\textcolor{black}{12224})[1 ]=\textcolor{black}{17287}\,(+1)$\tabularnewline
$\bar\theta(\textcolor{black}{6112})[ 2]=8642\,(+0)$\quad\quad\quad\,
$\theta(\textcolor{black}{6112})[1 ]=\textcolor{black}{8643}\;(+0)$\tabularnewline
$\bar\theta(3056)[2]=4320\;(+0)$\quad\quad\quad\quad\quad
$\theta(3056)[1 ]=\textcolor{black}{4321\;(+0)}$\tabularnewline
$\bar\theta(1528)[2 ]=2159\;(-1)$\quad\quad\quad\quad\quad\quad
$\theta(1528)[1 ]=\textcolor{black}{2160\;(+0)}$\tabularnewline
$\bar\theta(764)[2 ]=1079\;(+0)$\quad\quad\quad\quad\quad\quad\quad\quad\,
$\theta(764)[1 ]=\textcolor{black}{1080\;(+0)}$\tabularnewline
$\bar\theta(382)[2]=539\;(+0)$\quad\quad\quad\quad\quad\quad\quad\quad\quad\quad\,
$\theta(382)[1]=\textcolor{black}{540\;(+1)}$\tabularnewline
$\bar\theta(191)[2]=269\;(+0)$\quad\quad\quad\quad\quad\quad\quad\quad\quad\quad\quad\quad
$\theta(191)[1]=\textcolor{black}{270\;(+1)}$\tabularnewline
$\bar\theta(95)[2]=133\;(-1)$\quad\quad\quad\quad\quad\quad\quad\quad\quad\quad\quad\quad\quad\quad
$\theta(95)[1]=\textcolor{black}{134\;(+0)}$\tabularnewline
$\!\!\!\!\!\textcolor{black}{\bar\theta(47)[2]=65\;(-1)}$\qquad\quad\quad\quad\quad\quad\quad\quad\quad\quad\quad\quad\quad\quad\quad
$\textcolor{black}{\theta(47)[1]=66\;(+0)}$\tabularnewline
$\!\!\!\!\!\!\textcolor{brown}{\bar\theta(23)[2]=31\;(-1)}$\qquad\quad\quad\quad\quad\quad\quad\quad\quad\quad\quad\quad\quad\quad\qquad
$\textcolor{brown}{\theta(23)[1]=32\;(-1)}$\tabularnewline
$\textcolor{black}{\dots}$ \qquad\qquad \qquad \qquad \qquad \qquad \qquad \qquad \qquad \qquad \qquad$\textcolor{black}{\dots}$
\end{tabular}

\noindent  we get to see both agreement and deviations. If  we allow for deviations $\pm1$ across a level,      the upper cut-off lies at $x=24448$  and the lower cut-off at $x=23$. Whilst the emulation is  seven storeys less high by that measure, there remain nine  which show the rudimentary agreement necessary for  life to exist in the proposed form. For a reverse-order example, see the close of \ref{sec:Minkowski}, p.\,\pageref{actpoin}.

\renewcommand{\arraystretch}{1}

\subsection{\label{sub:The-plateau-effect}The plateau effect}

The incidences $I\!\left(\between_{\,\alpha}^{(n)}=x_{m}\mid x_{m}\!\in\!\left\{ L_{m}^{+},L_{m},L_{m}^{-}\right\} \!;n<\!6262\right)$
(Table \ref{tab:Specific-fractions-in-1-2-1}) follow a hyperbolic
law $1/x_{m}$, modulated by a weight function $w(m)$ so that $I=w(m)/x_{m}$.
The $I$'s to be examined are based on a non-automated read-out of
online CFR calculations, so counts were restricted to a manageable
scope $m\geq8$. For the discussion in this section, we impose the
further restriction $8\leq m\leq11$, so that only the incidences
of the higher down quark constituents%
\footnote{$\;$as can be seen from Table \ref{tab:Prime-factors-of}, there
are more down quark constituents (at $m=4,5,6$), which require automated
read-out because of their high incidence rate %
} and the up quark variants will be considered: \bigskip{}

\quad{}\qquad{}\qquad{}%
\begin{tabular}{|c|c|c|c|}
\hline 
$m$  & $I(\between_{\alpha}^{(n)}=L_{m}^{+})$  & $I(\between_{\alpha}^{(n)}=L_{m})$  & $I(\between_{\alpha}^{(n)}=L_{m}^{-})$\tabularnewline
\hline 
\hline 
$\vdots$  &  &  & \tabularnewline
\hline 
8  & 249  & 184  & 331\tabularnewline
\hline 
9  & 278  & 161  & 219\tabularnewline
\hline 
10  & 164  & 101  & 139\tabularnewline
\hline 
11  & 110  & 65  & 98\tabularnewline
\hline 
\end{tabular}

\bigskip{}

\noindent The most natural approach to interpret these numbers is to treat
their respective weight function as an identity: 
\[
w(m)\equiv I\, x_{m},
\]
implying there's a $3^{\textrm{rd}}$ degree polynomial 
\[
a_{3}\mu^{3}+a_{2}\mu^{2}+a_{1}\mu+a_{0}
\]
which yields a perfect $R^{2}=1$ fit of the data. For $m_{0}=8$,
for instance, the cubic fit for $I\, L_{m_{0}+\textrm{offset}}$,
\[
((m_{0},184L_{m_{0}}),(m_{0}\!+\!1,161L_{m_{0}+1}),(m_{0}\!+\!2,101L_{m_{0}\!+\!2}),(m_{0}\!+\!3,\,65L_{m_{0}+3})),
\]
is $(R^{2}=1)$ perfect with the polynomial 
\[
2\,313\,\mu^{3}-67\,195\,\mu^{2}+640\,026\,\mu-1\,959\,824\,.
\]
Perfect fits with $3^{\textrm{rd}}$ degree polynomials are also achieved
for $I\, L_{m_{0}+\textrm{offset}}^{\pm}$. These polynomials are
robust in the sense that their `fitness' persists even if $n$ is
increased beyond the $n\leq6262$ horizon of Table \ref{tab:Specific-fractions-in-1-2-1}
(then, of course, with a slightly different set of coefficients $a_{\nu}$).
In the above example, $\mu$ was chosen as ${m_{0}+\textrm{offset}}$.
But the polynomials are endowed with a `translational' symmetry to
the effect that a linear substitution $\mu\rightarrow\mu'$ leaves
the respective leading coefficient $a_{3}$ and the shape of the plot
invariant. The unit shift applied $s$ times yields a plateau index
$s\,$; when applied eight times, the plateau index has reached the
level $m_{0}$ where, in terms of $m_{0}-s$, the fit becomes 
\[
\textrm{cubic fit }((0,184\, L_{m_{0}}),(1,161\, L_{m_{0}+1}),(2,101\, L_{m_{0}+2}),(3,\,65\, L_{m_{0}+3}))
\]
with a polynomial 
\[
2\,313(\mu')^{3}-11\,683(\mu')^{2}+9\,002\mu'-44\,160\,
\]
as plotted in the figure:

\noindent \bigskip{}
 \quad{}\quad{}\qquad{}\includegraphics{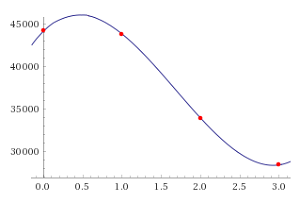}\bigskip{}

\noindent We note however that the next-to-leading coefficient $a_{2}$
varies $-$ as do the remaining $a_{\nu}$. The table below shows
the variations of $a_{2}$ in dependence of $m_{0}-s$:\bigskip{}

\quad{}\quad{}\qquad{}\global\long\def\arraystretch{1.25}
\begin{tabular}{|c|c|c|c|}
\hline 
$m_{0}-s$  & $a_{2}\textrm{ for }I(L_{m}^{+})$  & $a_{2}\textrm{ for }I(L_{m})$  & $a_{2}\textrm{ for }I(L_{m}^{-})$\tabularnewline
\hline 
\hline 
8  & $-243\,971$  & $-67\,195$  & $-\frac{11\,645}{2}$\tabularnewline
\hline 
7  & $-\frac{437\,783}{2}$  & $-60\,256$  & $-4788$\tabularnewline
\hline 
6  & $-193\,812$  & $-53\,317$  & $-\frac{7507}{2}$\tabularnewline
\hline 
5  & $-\frac{337\,465}{2}$  & $-46\,378$  & $-2719$\tabularnewline
\hline 
4  & $-143\,653$  & $-39\,439$  & $-\frac{3369}{2}$\tabularnewline
\hline 
3  & $-\frac{237\,147}{2}$  & $-32\,500$  & $-650$\tabularnewline
\hline 
2  & $-93\,494$  & $-25\,561$  & $\frac{\boldsymbol{769}}{\boldsymbol{2}}$\tabularnewline
\hline 
1  & $-\frac{136\,829}{2}$  & $-18\,622$  & 1419\tabularnewline
\hline 
0  & $-\boldsymbol{43\,335}$  & $\boldsymbol{-11\,683}$  & $\frac{4907}{2}$\tabularnewline
\hline 
\end{tabular}\global\long\def\arraystretch{1}

\bigskip{}
 The figures suggest that the proton/neutron transmutations discussed
in the previous section, by requiring a plateau $L_{m_{0}}=L_{8}=240$,
are the result of nature choosing a minimum (absolute) value of $a_{2}$.
The condition holds that way for the variant $L_{m_{0}}^{+}=L_{8}^{+}=241$,
too. But the figures in the third column have a deviation in store:
The minimum there is $\vert a_{2}\vert=\frac{769}{2}$ and indicates
plateau formation at $m_{0}\!-2=\!6$, a kind of `shadow' of the real
plateau at $m_{0}=8$, if you will! The weight function provides a
rationale for the situation. The incidences $I\!\left(\between_{\,\alpha}^{(n)}=L_{m_{0}+\textrm{offset}}^{-}\right)$
exceed both $I\!\left(\between_{\,\alpha}^{(n)}=L_{m_{0}+\textrm{offset}}^{+}\right)$
and $I\!\left(\between_{\,\alpha}^{(n)}=L_{m_{0}+\textrm{offset}}\right)$
for the offset 0 $-$ with proportions 331 vs. $\!249$ vs. $\!184$.
This reveals how often all $x_{m_{0}}$, not only $L_{m_{0}}^{-}$,
woud be used as \texttt{d} constituent `shadow'-wise $-$ unimpeded
by the real plateau $-$, while the overall $w(m)\equiv I\, L_{m}^{-}$
shows the same thing for all offsets by imitating the shape of the
hyperbolic law $1/x_{m}\,$:

\bigskip{}

\noindent \bigskip{}
 \quad{}\quad{}\qquad{}\includegraphics{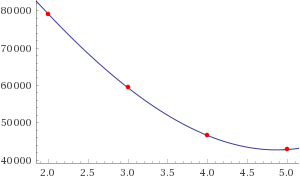}\bigskip{}

\noindent Proton/neutron transmutations are an example of one type
of juxtaposition inducing another. For representing plateau 240, a
boundary with outer nodes of minimally a $3^{18}$-cube complex formed
by $(G_{\rho}^{(31)})$ is needed. The table shows that the ratio
$a_{2}^{+}/a_{2}$ (for: $a_{2}\textrm{ for }I(L_{m}^{+})$ over $a_{2}\textrm{ for }I(L_{m})$)
ranges from 3.63 at $m_{0}-s=8$ to 3.71 at $m_{0}-s=0$, marking
points on an interval $[\Phi^{(31)},\Phi^{(63)}]$ where $\Phi^{(31)}\approx3.43$
and $\Phi^{(63)}\approx3.72$. The quantities $\Phi^{(p_{i})}$, we'll
soon learn, are crucial in that they combine outer parafermionic order
$p_{i}$ with internal or Catalan order $p_{i-2}$. Thus, via $a_{2}^{+}/a_{2}$,
the real plateau implies a $i$ vs. $\!(i-2)$ juxtaposition, while
location of the minima of $a_{2}^{+}$,$a_{2}$ vs. $\!$the minimum
$a_{2}^{-}$ (for $a_{2}\textrm{ for }I(L_{m}^{-})$) in the `shadow'
is accompanied by a $\mu$ vs. $\!(\mu-2)$ juxtaposition.

\subsection{Mass of proton and neutron}

Mass in the pre-geometric setting can only be represented by mass
ratios. With the elimination of $L_{4}=24$ from $\Sigma L_{\textrm{down}}$,
the proton's total content becomes 1836, which, given that the proton-electron
mass ratio is equal to this number, brings the electron mass to a
croton-representable value (needed to justify the use of electric
charge in preons). We can reinterpret this condition as a normalization
of the original proton total content,
\begin{equation}
\textrm{total content }\:1860\:\rightarrow\:\textrm{vetted content }1860/K_{\mathtt{p}}\label{eq:tot-prot-renorm}
\end{equation}
where 
\[
K_{\mathtt{p}}\sim\frac{1860}{1836}(=1.0\overline{1307189542483660}\,).
\]
As the net charge of a neutron is zero, an a priori fixed neutron-electron
mass ratio is no foregone conclusion. If we assume that, on average,
the transmutations $\mathtt{p}\rightarrow\mathtt{n}$ and $\mathtt{n}\rightarrow\mathtt{p}$
combine to a neutron total content $(2359+2356)/2=2357.5$, this content
would by a similar normalization be likely brought to represent the
observed neutron-electron mass ratio, 
\begin{equation}
\textrm{total content }\:2357.5\:\rightarrow\:\textrm{observed content }2357.5/K_{\mathtt{n}}.\label{eq:tot-neut-renorm}
\end{equation}
These innocent looking formul{\ae} in truth point to the effectiveness
of another formal principle besides juxtaposition $x$ vs. $x\!\pm\!1(2)$,
namely the functional 
\begin{equation}
1\divideontimes f^{(a)}\circ(f^{(b)}\divideontimes f^{(c)})\:\textrm{where}\:\divideontimes\in\left\{ \times,\div\right\} .\label{eq:timesonover}
\end{equation}
Substitute $\divideontimes:=\div\:$, $f^{(c)}:=\frac{5}{8}\times\:$,
$f^{(b)}:=\left\lceil \frac{9}{8}\times\right\rceil $ and $f^{(a)}:=\frac{9}{8}\times$
and it yields the expression, interordinal in $p=15,p'=31$, 
\begin{equation}
K_{\mathtt{p}}=\frac{1}{\frac{9}{8}\times\frac{\left\lceil \frac{9}{8}\times15\right\rceil }{\frac{5}{8}\times31}}=\frac{1860}{1836},\label{eq:norm-K_p}
\end{equation}
and the substitution $\divideontimes:=\div\:$, $f^{(c)}:=\left\lceil \frac{9}{8}\times\right\rceil $,
$f^{(b)}:=\frac{5}{8}\times\:$ and $f^{(a)}:=\frac{9}{8}\times$
yields 
\begin{equation}
\frac{1}{K_{\mathtt{n}}}=\frac{1}{\frac{9}{8}\times\frac{\frac{5}{8}\times31}{\left\lceil \frac{9}{8}\times15\right\rceil }}=\frac{1088}{1395},\label{eq:norm-K_n}
\end{equation}

\noindent resulting in a neutron-electron mass ratio $1838.681$.
We will soon learn that two other applications of the functional (\ref{eq:timesonover}),
one interordi\-nal in $p=15,p'=31$, the other interordinal in $p=31$,
$p'=63$, will lead to terms $\kappa$, $\kappa'$ which help approximate
the electromagnetic and weak coupling constants, respectively.

\subsection{\label{GJnum}${\cal GJ}$ bulk behavior and superconductivity}

At the most basic level, numbers involved in the control of electronic
bulk behavior must be among the 11 $G$-numbers (see AppendixA) 
\[
{\cal G}=\{3,5,11,17,19,41,43,113,115,155,429\}
\]
as compliant with the $G$ matrix elements: $\{G_{\alpha\beta}^{(p_{i})}\,\vert\; p_{i}\leq x_{e}^{-1};\, C_{7}\geq G_{\alpha\beta}^{(p_{i})}>1\}$
(a value equal unity would result in zero\,$=$ no control as an
electron is added). Also involved are the 6 $J$-numbers

\[
{\cal J}=\{3,13,15,117,143,149\}
\]
left over using the same restriction criteria: $\{J_{\alpha\beta}^{(p_{i})}\,\vert\; p_{i}\leq x_{e}^{-1};\, C_{7}\geq J_{\alpha\beta}^{(p_{i})}>1\}$.
In both number types, 3 is the minimum element, hence, for a way in
to zero-avoiding control, no more than two electrons may be admitted.
For reasons of symmetry, up to two positrons would be admitted as
well, so the control exerted is the result of applying $x$ vs. $x\pm1(2)$
juxtaposition to ${\cal G}$, \vspace*{-1mm}

\[
{\cal G}^{+}=\{4,6,12,18,20,42,44,114,116,156,430\},
\]
\[
{\cal G}^{++}=\{5,7,13,19,21,43,45,115,117,157,431\},
\]
\[
{\cal G}^{-}=\{2,4,10,16,18,40,42,112,114,154,428\},
\]
\[
{\cal G}^{--}=\{1,3,9,15,17,39,41,111,113,153,427\},
\]
{\em mutatis mutandis} to $J^{+}$, $J^{++}$, $J^{-}$ and $J^{--}$,
and combining the new terms with ${\cal G}$ and ${\cal J}$. Together,
they form the 56 controls, or ${\cal GJ}$-numbers: 
\[
1,{\scriptstyle \ldots},\!7,9,{\scriptstyle \ldots},\!21,39,{\scriptstyle \ldots},\!45,111,{\scriptstyle \ldots},\!119,141,{\scriptstyle \ldots},\!145,147,{\scriptstyle \ldots},\!151,
\]
\[
153,{\scriptstyle \ldots},\!157,427,{\scriptstyle \ldots},\!431.
\]
In Sect. \ref{sec:Pregeometric-categories-relevant}, it was noted
that to the class of cases \vspace*{-2mm}

\[
\between_{\,\alpha}^{(n)}=x_{m}\pm\Delta_{m}\quad x_{m}\!\in\!\left\{ L_{m}^{+},L_{m},L_{m}^{-}\right\} 
\]
belongs the index rule \vspace*{-2mm}

\[
\alpha=2a+\delta\Rightarrow\beta(\beta')=2b+\delta-1>1\qquad(\delta\in\{0,1\})
\]
under the constraint \vspace*{-3mm}

\[
\begin{array}{c}
x_{m}\pm\Delta_{m}\lessgtr x_{m\pm1}\mp\Delta_{m}\end{array}
\]
for $\between_{\beta}^{(n')}=\Delta_{m}$ with $n'=n\pm1$ and $n_{\Delta}=1$.
The rule can be given the more general form \vspace*{-3mm}

\[
\alpha=2a+\delta\Rightarrow\sum\beta_{\nu}=2b+\delta-n_{\Delta}>1\qquad(\delta\in\{0,1\}),
\]
dealing with dispersion/dissipation, \vspace*{-2mm}

\[
\between_{\alpha}^{(n)}+\between_{\beta_{1}}^{(n')}+\between_{\beta_{2}}^{(n')}+\ldots+\between_{\beta_{n_{\Delta}}}^{(n')}\;=\; L_{m}^{,\pm}
\]
where $n_{\Delta}>1$, $\sum_{\nu=1}^{n_{\Delta}}\between_{\beta_{\nu}}^{(n')}=\Delta_{m}$
and order is transitive:

\vspace*{-2mm}

\[
\between_{\alpha}^{(n)}\gg\;\between_{\beta_{1}}^{(n')},\quad\between_{\beta_{1}}^{(n')}>\,\between_{\beta_{2}}^{(n')},\;\ldots,\;\between_{\beta_{-1+n_{\Delta}}}^{(n')}>\,\between_{\beta_{n_{\Delta}}}^{(n')}.
\]

\noindent Dispersion (as related to squashed terms) and dissipation
(as related to stretched ones) may be tracked for $n'$ owing to,
say, Mersennian increments: $n'_{2}=n\pm3\Rightarrow n_{\Delta}=2$,
possibly followed by $n'_{3}=n\pm7\Rightarrow n_{\Delta}=3$, possibly
followed by $n'_{4}=n\pm15\Rightarrow n_{\Delta}=4$, etc.

\noindent An example, with $x_{21}=L_{21}=27720$, $\between_{\,278}^{(14146)}=25927$
and $\Delta_{21}=1793$, is given by 
\[
\between_{\,700}^{(14143)}=1743,\;\between_{\,698}^{(14143)}=50,
\]
\[
\between_{\,423}^{(14153)}=1522,\;\between_{\,270}^{(14153)}=202,\;\between_{\,698}^{(14153)}=69,
\]
\[
\between_{\,529}^{(14131)}=1504,\;\between_{\,108}^{(14131)}=187,\;\between_{\,199}^{(14131)}=76,\;\between_{\,564}^{(14131)}=26,
\]
\vspace*{-2mm}

\[
\between_{\,271}^{(14115)}=1269,\;\between_{\,316}^{(14115)}=223,\;\between_{\,632}^{(14115)}=201,\;\between_{\,24}^{(14115)}=81,\;\between_{\,22}^{(14115)}=19.
\]

\noindent The adaptation to bulk electronic behavior is carried out
in three steps:

\noindent (i) We note: when admitting both electrons and positrons
whose annihilation would produce radiation, Mersennian increments
are a perfect choice $-$ as we have seen in Sect. $\!\!$\ref{sec:Pregeometric-categories-relevant},
radiation is linked to Mersennian time-like refinements (relative
to $n_{0}$, so $n$ takes the role of $n_{0}$). (ii) Since superconductivity
shows up at a temperature decrease past a critical point, we impose
{\em one} direction along which $n'$ varies, starting with $n'_{5}=n_{0}-31$,
with $n'_{4}=n_{0}-15$ next and continuing all along down to $n'_{0}=n_{0}-0$.
(iii) Instead of presuming the obvious condition ${n_{\Delta}=\textrm{l}\textrm{o}\textrm{g}}_{2}(n_{0}-n'+1)$,
we impose $n_{\Delta}=s_{\textrm{c}\textrm{a}\textrm{r}\textrm{d}}-3(5-{\textrm{l}\textrm{o}\textrm{g}}_{2}(n_{0}-n'+1)$,
$s_{\textrm{c}\textrm{a}\textrm{r}\textrm{d}}$ being the sum of the
cardinalities of ${\cal G}$ and ${\cal J}$, 17.

\noindent Because of computability limits, continued fractions for
$n$ as high as 14146 above have too few terms to bear out the envisioned
relationship, so we switch to a more comfortable size such as is afforded
by $x_{21}=L_{21}=27720$, $\between_{\,71}^{(1526)}=28176$ and $\Delta_{21}=456$.

\noindent Before switching, we must first show the equivalence of
the two situations. At $n=14146,$ electron 1 is accelerated toward
electron 2 during $n$ and $n_{\textrm{acc}}=n+15^{2}=14371$, then
exchanging a (virtual) photon with it during $n_{\textrm{acc}}$ and
$n_{\textrm{rad}}=n_{\textrm{acc}}+\ell_{2}=14385$ where $\ell_{2}=14$:
\[
\begin{array}{cc}
\between_{\,34}^{(14371)}=1763, & \between_{\,121}^{(14371)}=30,\\
\between_{\,24}^{(14385)}=1763, & \between_{\,78}^{(14385)}=30
\end{array}\qquad(\Delta_{21}=1793).
\]
Equivalently, at $n=1526,$ electron 1 is accelerated toward electron
2 between $n=1526$ and $n_{\textrm{acc}}=n+15^{2}=1751$, with a
shorter duration of the subsequent photon exchange, $n_{\textrm{rad}}=n_{\textrm{acc}}+\ell_{1}=1757$
where $\ell_{1}=6$: 
\[
\begin{array}{cc}
\between_{\,676}^{(1751)}=453, & \between_{\,677}^{(1751)}=3,\\
\between_{\,674}^{(1757)}=453, & \between_{\,679}^{(1757)}=3
\end{array}\qquad(\Delta_{21}=456).
\]

\noindent The equivalence is due to the fact that both $\between_{\,278}^{(14146)}=25927$
and $\between_{\,71}^{(1526)}=28176$ represent peaks in their respective
Mersenne fluctuations. A {\em complementary} equivalence class
arises if we replace peaks by $x_{m}$-related left-leg $\between_{\,\alpha_{0}-\nu}^{(n)}$
and right-leg partners $\between_{\,\alpha_{0}+\nu}^{(n)}\,(=\,\between_{\,\alpha_{0}-\nu}^{(n)}\!\!+\epsilon$,
$\epsilon\in\left\{ -1,0,1\right\} )$. Those duos imply \emph{two
$\Delta$'}s, one $\Delta_{\textrm{outer}}$, one $\Delta_{\textrm{inner}}$.
The associated processes can be `time-symmetric,' {\em i.e.}, can
be read in normal order or in reverse order and come with a head (tail)
photon exchange having electron 1(2) accelerate.$\!$%
\footnote{An asymmetric process would be the accretion of self-energy in form
of a succession of two acceleration phases, one before interaction
with the environment, one thereafter, like in the following example:
$x_{19}=L_{19}=10668$, $\Delta_{19}=856$ for the outer regime and
$x'_{19}=L_{19}^{-}=10667$, $\Delta'_{19}=855$ for the inner regime,
as realized by the left-legs $\between_{\,559}^{(14751)}(=9812)$,
$\between_{\,314}^{(14887)}(=856)$ and right-legs $\between_{\,567}^{(14753)}(=9812)$,
$\between_{\,307}^{(14881)}(=855)$ so that $(n_{x_{19}}-n_{\Delta_{19}})+(n_{x'_{19}}-n_{\Delta'_{19}})-\ell_{3}=3^{2}+15^{2}$.
The associated correlations between $\between$'s and space-like refinements
(see main text) read 
\[
\begin{array}{l}
\textrm{corr}\!\left(\between_{\,314}^{(14887)}+\between_{\,307}^{(14881)},(559-314)+(567-307)\right)\models\delta\,=2\times3^{2}\times67\\
\textrm{corr}\!\left(\between_{\,559}^{(14751)}+\between_{\,567}^{(14753)},(559-314)+(567-307)\right)\models\delta\,=3\times6373
\end{array}
\]
} We present an example with $x_{21}=L_{21}=27720$, left-leg $\between_{\,34}^{(14375)}=28217$
and right-leg $\between_{\,24}^{(14381)}=28218$, hence $\Delta_{\textrm{outer}}=497$
and $\Delta_{\textrm{inner}}=498$. Under the outer regime, a photon
is exchanged at $n_{\textrm{rad1}}=14375$ and $n_{\textrm{rad2}}=14452$,
respectively, each exchange of length $\ell_{2}=14$, which under
the inner regime has electron 1(2) accelerate for $(n_{\textrm{rad2}}-\ell_{2})-(n_{\textrm{rad1}}+\ell_{2})=7^{2}$:
\[
\begin{array}{cccc}
\quad\between_{\,34}^{(14375)}=28217, &  &  & (\Delta_{\textrm{outer}}=497)\\
\quad\between_{\,24}^{(14381)}=28218, & \between_{\,543}^{(14389)}=498, &  & (\Delta_{\textrm{inner}}=498)\\
\between_{\,141}^{(14438)}=354, & \between_{\,348}^{(14438)}=144, &  & (\Delta_{\textrm{inner}}=498)\\
\!\between_{\,675}^{(14452)}=497 &  &  & (\Delta_{\textrm{outer}}=497)
\end{array}.
\]

\noindent $\between$'s are correlated with space-like refinements
via the differences they model, which are multiples of the previously
mentioned minimal element 3:

\[
\begin{array}{lcl}
\textrm{corr}\!\left(\between_{\,348}^{(14438)},141\right)\models\delta=3, &  & \textrm{corr}\!\left(\between_{\,141}^{(14438)},348\right)\models\delta=6,\\
\textrm{corr}\!\left(\between_{\,675}^{(14452)},543-34\right)\models\delta=-12, &  & \textrm{corr}\!\left(\between_{\,543}^{(14389)},675-24\right)\models\\
 &  & \qquad\,\,\delta=-(3+6+12^{2}).
\end{array}
\]
Moreover, the space-like refinements of $\Delta_{m}$-related $\between$'s
in the peak class are coupled to those of $x_{m}$-related $\between$'s
in the `leggy' class: 
\[
\begin{array}{cc}
\between_{\,34}^{(14371)}=1763, & \between_{\,34}^{(14375)}=28217,\\
\between_{\,24}^{(14385)}=1763, & \between_{\,24}^{(14381)}=28218.
\end{array}
\]
Since we have taken precautions to including electrons that travel
backward in time (positrons), the leggy class is amenable to bulk
behavior examination as well. However, with the above example suffering
from the same computability limits, we turn to the former equivalence
class for its lower-scale $n=1526$.

\noindent Between accelerations, beginning with $n'=1526-31$, we
find seventeen terms $\between_{\beta_{\nu}}^{1495)}$ that match
the ${\cal GJ}$-numbers

\[
117,113,43,42,19,17,16,15,14,13,12,11,10,5,4,3,2
\]
and form the sum 456. They have a dispersion\,:\,dissipation ratio
of 10\,:\,7.

\noindent For $n'=1526-15$, we find fourteen terms $\between_{\beta_{\nu}}^{1511)}$
matching the ${\cal GJ}$-numbers

\[
153,118,43,42,19,18,15,13,12,9,5,4,3,2
\]
and summing up to 456. They have a dispersion\,:\,dissipation ratio
of 8\,:\,6.

\noindent The corresponding figures for $n'=1526-7$ are eleven terms
$\between_{\beta_{\nu}}^{1519)}$ 
\[
117,115,114,43,19,18,15,6,4,3,2,
\]
and a ratio 6\,:\,5. Those for $n'=1526-3$ are eight terms $\between_{\beta_{\nu}}^{1523)}$
\[
155,119,113,18,16,13,12,10,
\]
ratio 4\,:\,4, and those for $n'=1526-1$ five terms $\between_{\beta_{\nu}}^{1525)}$
\[
154,141,114,42,5,
\]
ratio 2\,:\,3. (The smaller terms $\between_{\beta_{\nu}}^{(n')}=2,3,4,5$
can always be chosen such that the adapted index rule is satisfied.)
The closer the critical point $n'$ gets, the more the tipping point
between dispersion and dissipation shifts in favor of the latter.
At the critical point $n'=n_{0}-0=1526$, a duo of ${\cal GJ}$-numbers
matching with 456 would be needed, but no two of them fit the bill.
This suggests that the fate of superconductivity depends on a critical
balance between $x_{m}$ and $\Delta_{m}\,$.

\section{Organizers: boundaries and characteristic quantities}

Much of what we learned about boundaries and controls in the global
context of ideation may be put to use again in the context of organization.
In Sect.\ref{sub:Croton-field-duality}, the part of duality control
was taken by Catalan numbers in conjunction with $5\cdot2^{i}\:(i\in\mathbb{N}_{0})$,
numbers we renamed $M_{\nicefrac{5}{8}}^{+}$ when introducing physically
relevant pregeometric categories (Sect. \ref{sec:Pregeometric-categories-relevant}).
Understood as an expansion factor, the $5\cdot2^{i}$ are the result
of rescaling $\frac{2^{i+4}}{\pi}$ by a constant factor $\frac{5\pi}{16}=1/1.01859\ldots$;
cf. Eqs.(\ref{eq:rows-1}). The `Mersennians' of 5, $5\cdot2$, $5\cdot2^{2},\ldots$
in turn coincide with the numbers $(p+q)/2\quad(q=(p-3)/4)$, which
we distinguished from regular Mersenne numbers $M_{\textrm{reg}}:=p_{i}=1,3,7,\ldots$
by dubbing them $M_{\nicefrac{5}{8}}:=o_{i}=4,9,19,\ldots$ .

\subsection{\label{sub:Duality-controls:-the}Combining Catalan numbers and interordinality:
the electromagnetic coupling constant}

Resuming the discussion of labels on the boundary, we remember that
all deficiencies are removed in the interordinal case: out of 192
distinct crotons ensuing from the enlarged basis 
\[
(G_{\rho}^{(7,15)})=(\underline{1},3,\underline{5},11,17,41,113),
\]
(including a singularity and neglecting sign reversals) \emph{all
}become realizable on the enlarged boundary $\Gamma^{(7,15)}$. As
we will see now, the next interordinality, $p_{i-1}=15$,$p_{i}=31$,
is powerful enough to deal with the electromagnetic coupling constant.
The make-up of this constant, we will learn, leads us in turn to new
boundaries, Bound$(\cdot)$.

\noindent We already mentioned that the Catalan number $C_{q_{i}}$,
where $p_{i}$ and $q_{i}$ are $\left(i,i-2\right)$-juxtaposed via
$q=(p-3)/4$, is central to croton base numbers that are in line with
$C_{q_{i}+1},C_{q_{i}+2},$ $\ldots,C_{2q_{i}}$ $-$ a fundamental
connection leading to the span parameter
\[
\Phi^{(p_{i})}=\left(G_{\textrm{max}}^{(p_{i})}/C_{q_{i}}\right)^{1/q_{i}}.
\]
 Quark constituents are often deemed too artifical to be true elements
of nature, but, as we have seen in Sect. \ref{sec:Application-to-subatomic},
if we concentrate on their role as carriers of fractional electric
charge, $x_{e}^{-1}=31$ turns out to provide a natural framework
for dealing with the electromagnetic force. Indeed, the above set
of Catalan numbers supply up-type and down-type \emph{interordinal}
bounds for the electromagnetic coupling constant, the dimensionless
quantity \textgreek{a}. When normalized with the factor $(C_{q_{i}}-C_{(q_{i}+q_{i-1})/2})^{-1}\;(i=\log_{2}(x_{e}^{-1}+1))$,
the down-type parameter $\Phi^{(p_{i-1})}$ provides a tight upper
bound 
\[
(C_{7}-C_{5})^{-1}\Phi^{(15)}\approx1/136.88
\]
to the currenty measured value  \textgreek{a}$\;=1/137.035999206(11)$. And the
up-type parameter $\Phi^{(p_{i})}$ normalized with the plus-sign
counterpart $(C_{q_{i}}+C_{(q_{i}+q_{i-1})/2})^{-1}$ yields an even
better lower bound: 
\[
(C_{7}+C_{5})^{-1}\Phi^{(31)}\approx1/137.04
\]
\noindent The location where  \textgreek{a}$^{-1}$ interpolates the interval
$\left[136.88,137.04\right]$ can to good approximation be given by
the respective down-type and up-type expressions\\

{\small Down-type form of  \textgreek{a}$^{-1}$:}
\[
\textup{\ensuremath{\textrm{\textgreek{a}}}}^{-1}=(C_{7}-C_{5})/\Phi^{(15)}+\frac{2f_{n}+\kappa}{f_{n+1}+\kappa}\,\Delta_{\textrm{b}}=137.035999547\ldots
\]

{\small Up-type form of \textgreek{a}$^{-1}$:}
\[
\textup{\ensuremath{\textrm{\textgreek{a}}}}^{-1}=(C_{7}+C_{5})/\Phi^{(31)}-\frac{1}{f_{n+1}+\kappa}\,\Delta_{\textrm{b}}=137.035999547\ldots
\]
 where $\Delta_{\textrm{b}}=(C_{7}+C_{5})/\Phi^{(31)}\!-(C_{7}-C_{5})/\Phi^{(15)}$,
$f_{n}=15$, $f_{n+1}=31$ and $\kappa$ is the result of applying
the functional  (\ref{eq:timesonover}) interordinally to $G_{\textrm{max}}^{(15)}/C_{3}$
and $G_{\textrm{max}}^{(31)}/C_{7}$ with $\divideontimes:=\times\:$,
$f^{(c)}:=\:^{{\scriptscriptstyle 7}}\!\!\!\surd\:$, $f^{(b)}:=\:^{{\scriptscriptstyle 3}}\!\!\!\surd\:$
and $f^{(a)}:=\frac{3}{8}+\surd\:$: 
\begin{equation}
\kappa=\frac{3}{8}+\sqrt{\Phi^{(15)}\Phi^{(31)}}.\label{eq:control functional}
\end{equation}
Thus, a surprisingly close approximation of \textgreek{a}$^{-1}$ is achieved
as a result of combining Catalan numbers and interordinality, alongside
twin coincidences, $p_{i-1}=f_{n}=15$, $p_{i}=f_{n+1}=31$, and $\kappa-\sqrt{\Phi^{(15)}\Phi^{(31)}}=\frac{5}{8}p_{n}-o_{n-2}\equiv\frac{3}{8}$,
that arise due to a strange fusion of interordinality and juxtaposition
$n$ vs. $n-2$.

\subsection{\label{sub:Duality-controls:-the-1}Combining Catalan numbers and
juxtaposition: particle-related dimensions }

We have seen combinations of Catalan numbers and $M_{\nicefrac{5}{8}}^{+}$
controlling the boundary conditions in croton amplitudes and phases,
and combinations of Catalan numbers and $M_{\nicefrac{5}{8}}$ regulating
in latent form the electromagnetic coupling constant in both its up-type
and down-type expression. While the former combination may be linked
to the question, given a boundary, how many $L_{\textrm{up}}$'s and
$L_{\textrm{down}}$'s are there that crotons can choose as targets,
the latter combination, with $M_{\nicefrac{5}{8}}$ in manifest form,
would allow to ask the question, how many Euclidean dimensions out
there get involved in particle creation. Key to the approach taken
here is the $(n,n\!-\!2)$-juxtaposed quotient formed by $L_{\nu}-\prod_{i=1}^{n}(\cdot)$
and $\prod_{i=1}^{n-2}(\cdot)$ where $L_{\nu}$ is least among $L_{m}>\prod_{i=1}^{n}(\cdot)$,
$\nu'$ is a natural number and the product arguments are taken from
$M_{\textrm{reg}}$ and $M_{\nicefrac{5}{8}}$: 
\begin{equation}
\left(L_{\nu}-\prod_{i=1}^{n}(\cdot)\right)/\prod_{i=1}^{n-2}(\cdot)=\nu'\qquad(\textrm{natural}).\label{eq:naturalness-1}
\end{equation}
 It's instructive to tabularize the instantiations of Eq. (\ref{eq:naturalness-1})
for choices of $n$ such that $p_{n}$ and $p_{n-1}$ are less or
equal $f_{n+1}(=x_{e}^{-1})$ and $f_{n}(=x_{e}^{-1})$, respectively:
\begin{table}[H]
\caption{\label{tab:Key-particle-creation-related-2}Key particle creation-related
dimensions}
\bigskip{}
\hfill{}%
\begin{tabular}{|>{\centering}p{1cm}|>{\centering}p{2cm}|>{\centering}p{4cm}|>{\centering}p{1.5cm}|>{\centering}p{1cm}|}
\hline 
$n-2$ & $L_{\nu}$ & $\prod_{i=1}^{n}(\cdot)$ & $\prod_{i=1}^{n-2}(\cdot)$ & $\nu'$\tabularnewline
\hline 
\hline 
1 & $L_{4}=24$ & $1\cdot3\cdot7=21$ & $1$ & 3\tabularnewline
\hline 
2 & $L_{10}=336$ & $1\cdot3\cdot7\cdot15=315$ & $1\cdot3$ & 7\tabularnewline
\hline 
\end{tabular}\hfill{}\vspace{-0.5mm}

\hfill{}%
\begin{tabular}{|>{\centering}p{1cm}|>{\centering}p{2cm}|>{\centering}p{4cm}|>{\centering}p{1.5cm}|>{\centering}p{1cm}|}
\hline 
3 & $L_{19}=10668$ & $1\cdot3\cdot7\cdot15\cdot31=9765$ & $1\cdot3\cdot7$ & 43\tabularnewline
\hline 
\hline 
1 & $L_{12}=756$ & $4\cdot9\cdot19=684$ & $4$ & 18\tabularnewline
\hline 
2 & $L_{21}=27720$ & $4\cdot9\cdot19\cdot39=26676$ & $4\cdot9$ & 29\tabularnewline
\hline 
\end{tabular}\hfill{}
\end{table}

\noindent The first observation worth mentioning is that $\nu$-sums
from the respective parts of the table, $\Sigma_{r=n-2}^{n}\nu_{r}$
and $\Sigma_{s=n-2}^{n-1}\nu_{s}$, are invariant: $4+10+19=33$ and
$12+21=33$. This suggests that the $\nu_{r}$ and the $\nu_{s}$
can be combined into a basis, $N_{\textrm{source}}=(4,10,12,19,21)$.
With coefficients $\pm1$ or 0, linear combinations of its elements
with positive result may then be said to span a variety of `source
dimensions:' Out of 66 potentially realizable ones, 11 remain unrepresented
on $\textrm{Bound}(N_{\textrm{source}})$ $-$ the name connoting
a set of labels on the outer nodes of a $T$-cube complex $(T=5)$
constructed in the manner $\Gamma$ and $\chi$ have been $-,$ namely
$49,51,\ldots,65$, and, neglecting sign reversals, two linear combinations
yield singularities (zeros) on $\textrm{Bound}(N_{\textrm{source}})$. 

\noindent The $\nu'$ ensuing from the two table parts can in a similar
manner be combined into a basis, $N_{\textrm{sink}}=(3,7,18,29,43)$,
spanning `sink dimensions.' Out of 100 potentially realizable ones,
41 remain unrepresented on $\textrm{Bound}(N_{\textrm{sink}})$: $2,5,\ldots,99$,
and, neglecting sign reversals, one linear combination yields a singularity
on $\textrm{Bound}(N_{\textrm{sink}})$. The two boundaries can in
turn be combined, to the effect that the number of unrepresentable
dimensions shrinks to 22. 

\noindent The overall picture emerging from these numbers is as follows:
The particle creation-potential of the first 100 Euclidean dimensions
is governed by the number 11. All other features arising \emph{en
route} are identical to or multiples of this number $-$ the $\nu$-invariant
33, the unrepresentable 11 `source dimensions' out of 66 potential
ones, the number of unrepresentable `sink dimensions' plus the number
of singularities under the union $\textrm{Bound}(N_{\textrm{source}})\cup\textrm{Bound}(N_{\textrm{sink}})$,
$41+3$, as well as the number of dimensions staying uninvolved in
particle creation even when the distinction between source and sink
dimensions is dropped, the said 22 unrepresentable dimensions $49,55,\ldots,99$.
On the other hand, the number of `sink dimensions' denying representation
on the single $\textrm{Bound}(N_{\textrm{sink}})$, 41, to which we
may add 1 to account for the one singularity remaining, coincides
with $C_{5}$, the interpolating term that makes $(C_{7}\mp C_{5})^{-1}$
bound the span parameters $\Phi^{(p_{i-1})},\Phi^{(p_{i})}$ so tightly
they approach the electromagnetic coupling constant in the first place.
So it's worthwhile to go into the details of the dimensional branching
process.

\noindent At the heart of it resides a unique link leading from $C_{(q_{i+1}+q_{i\,})/2}$
as starting point to $C_{(q_{i}+q_{i-1\,})/2}$ as end point $(i\!=\!\log_{2}(x_{e}^{-1}\!+\!1))$.
In fact, $C_{(15+7)/2}=C_{11},$ in typical interordinal manner, gets
down to the last member of the Catalan number sequence $(C_{3},C_{4},C_{5},C_{6})$
via 
\begin{equation}
\!\!\!\!\!\!\!\!\!\!\!\!\!\!\!\!\!\!\!\!\!\!\!\!\!\!\!(C_{11}\textrm{B}(11,12))^{-1}=11\cdot12=132=C_{6},\label{eq:one-spine-2}
\end{equation}
and $C_{6}$, in turn, completes the descent intraordinally to the
end point 

\begin{equation}
\;\quad(C_{6}\textrm{B}(6,7))^{-1}\:=\:6\cdot7=42=C_{5}\:(=C_{(7+3)/2}).\label{eq:one-spine-1-1}
\end{equation}
How the numbers $11,12$ and $6,7$ work is sort of like in a double
strand, $\begin{array}{ccc}
\!\!{}_{11} & _{\cdot} & \!\!{}_{12}\\
\!\!\cdot &  & \!\!\cdot\\
\!\!{}^{6} & ^{\cdot} & \!\!{}^{7}
\end{array}$: Horizontally, they serve as upper and lower ties, vertically, as
left and right strands. The number of `source dimensions' evolves
to $11\cdot6$ (left strand) of which $11\cdot6-11$ can be represented
on $\textrm{Bound}(N_{\textrm{source}})$. Panning to the `sink side'
(right strand), branching doesn't end up until at $12\cdot7+(11+C_{3})=100$
`sink dimensions' of which, after all, {12$\cdot7-(11+C_{4})=59$}
can be represented on $\textrm{Bound}(N_{\textrm{sink}})$. The complete
$(G_{\rho})^{(15)}$'s sequence of Catalan controls $(C_{3},C_{4},C_{5},C_{6})$
is exhausted in the process.

\subsection{\label{sub:The-role-of-centralizers}The role of organizers }

\subsection*{\label{sub:Following-the-global-givens}Following the global givens of ideation}

\noindent Does the creation of chunks of space and matter follow
the global givens of ideation? Take, for example, the creation of
a chunk of the eighteenth dimension. A croton $\varphi_{\alpha}^{(n)}=L_{18}(=7398)$
as an ideation of this dimension was found in a Mersenne fluctuation
of Type I (see Tbl. \ref{tab:Specific-fractions-in}), at time-like
refinement level $n\!=\!1144$ and secondary expansion $s\!=\!6$.
We note $1144=11\cdot104$. So we may ask if, with the same fluctuation
type, ideation includes the remaining basic `source' and `sink dimensions.'
What we find is that $\varphi_{\alpha}^{(n)}$ to match their kissing
numbers, if available, do share this global property pairing $11\!\mid\! n$,
$s=6$ $-$ a hint that creation follows the global givens by way
of the left-strand effect presented above: 
\begin{table}[H]
\caption{\label{tab:Key-particle-creation-related-1-1}Key particle dimensions
via$\protect\begin{array}{lc}
\textrm{Mersenne fluctuations Type}\:\textrm{I}: & \!\!\!\left\lfloor \log_{2}C_{63}\right\rfloor \protect\end{array}\!\!\!\left(\frac{2^{n}}{\pi}\right)^{-1}$}
\bigskip{}
\hfill{}%
\begin{tabular}{|>{\centering}m{2cm}|>{\centering}m{2.8cm}|>{\centering}m{2.2cm}|>{\centering}m{3.2cm}|}
\hline 
\phantom{}$L_{\nu}$ & $\varphi_{\alpha}^{(n)}$ & $L_{\nu'}$ & $\varphi_{\alpha'}^{(n')}$\tabularnewline
\hline 
\hline 
\phantom{}$L_{4}=24$ & $\varphi_{315}^{(33=11\cdot3)}=24$ & $L_{3}=12$ & $\varphi_{67}^{(11)}=12$\tabularnewline
\hline 
\phantom{}$L_{10}=336$ & $\varphi_{496}^{(737=11\cdot67)}=336$ & $L_{7}=126$ & $\varphi_{77}^{(55=11\cdot5)}=126$\tabularnewline
\hline 
\end{tabular}\hfill{}\vspace{-0.45mm}

\hfill{}%
\begin{tabular}{|>{\centering}m{2cm}|>{\centering}m{2.8cm}|>{\centering}m{2.2cm}|>{\centering}m{3.2cm}|}
\hline 
\phantom{}$L_{19}=10668$ & n/a & $L_{43}=\:?$ & n/a\tabularnewline
\hline 
\hline 
\phantom{}$L_{12}=756$ & $\varphi_{330}^{(616=11\cdot56)}=756$ & $L_{18}=7398$ & $\varphi_{499}^{(1144=11\cdot104)}=7398$\tabularnewline
\hline 
\phantom{}$L_{21}=27720$ & n/a & $L_{29}=198506$ & n/a\tabularnewline
\hline 
\end{tabular}\hfill{}
\end{table}

\subsection*{\medskip{}
}

\subsection*{\label{sub:Crotonic-implementation-locally}Creating the world container}

For a detailed description of the creation, a type of crotons embracing
the highest standard of precision is required: the natural choice
is organizers 
\[
(\textrm{CFR})\quad2^{-n}\kappa\!\rightarrow\!\left[\varkappa_{0}^{(n)}\!;\varkappa_{\alpha}^{(n)}\right],
\]
where $\kappa=\frac{3}{8}+\sqrt{\Phi^{(15)}\Phi^{(31)}}$. A croton
$\varkappa_{\alpha}^{(n)}$ qualifying as a pivot for $L_{m}$ must
have the special property that the gap between it and the target can
be bridged by integer addition within a collection of organizer co-amplitudes
$\varkappa_{\xi}^{(n)}$ $-$ a procedure that is mirrored on a refinement-dependent,
organizer boundary $\Lambda_{\alpha}^{(n)}$. For the envisioned relationship
we establish the following rules:\\
(1) When $\varkappa_{\alpha}^{(n)}>L_{m}$, $\varkappa_{\alpha}^{(n)}$
qualifies as a pivot with target $L_{m}$ if $\varkappa_{\alpha}^{(n)}/L_{m}$
does not exceed a prespecified range, say, $\sqrt{\frac{\varkappa_{\alpha_{\phantom{}}}^{(n)}}{L_{m}}}\lesssim$
$\frac{16}{5\pi}$ $(=1.01859\ldots)$ as a heuristic, and the organizer
co-amplitudes $\varkappa_{\xi}^{(n)}$ to participate in the collection
are located to the pivot's right; conversely, for $\varkappa_{\alpha}^{(n)}<L_{m},$
$\sqrt{\frac{\varkappa_{\alpha_{\phantom{}}}^{(n)}}{L_{m}}}\gtrsim$
$\frac{5\pi}{16}$ is required, and the participation of $\varkappa_{\xi}^{(n)}$
takes place to the left of the pivot;\\
(2) the prime factorization of $n$ determines how many co-amplitudes
$\varkappa_{\xi}^{(n)}$ are to be included; if it contains at least
one factor $p\in M_{\textrm{reg}}$ ($o\in M_{\textrm{5/8}}$), then
the number of co-amplitudes, in a success-dependent way, is \\
(2a) directly equal to this factor or \\
(2b) interpreted as $p_{j-2}$ $(o_{j-2})$, and $p_{j}$ $(o_{j})$
is assigned to the number of inclusions. \\
 (1) and (2) are only necessary conditions. A $\varkappa_{\alpha}^{(n)}$
obeying them has an entourage of co-amplitudes $\varkappa_{\xi}^{(n)}$
that still contain duplicates. $\theta$ that are distinct enter the
tuple $(\textrm{co-}\varkappa)_{\alpha}$ from which the nodes of
a $\theta$-cube complex  $\Lambda_{\alpha}^{(n)}$ as organizer boundary
get encoded. Either under (2a) or (2b), we get a label $\Delta_{\varkappa}$
that is equal to $\mid\varkappa_{\alpha}^{(n)}-L_{m}\!\mid$; if (2b)
applies, then, additionally, `Catalan tie' and $M_{5/8}$ -properties
have to be deployed to ensure a canonical form of $\Delta_{\varkappa}$
: \\
(2b$'$) If co-amplitudes $\varkappa_{\xi}^{(n)}\in(\textrm{co-}\lambda)_{\alpha}$
are multiples of tie numbers 6,7 or 11,12, they induce a sign divide:
multiples of 6 (7) and co-amplitudes of the form $9+4$ $(9-4)$ have
their signs preserved while all others incur sign inversion; if they
instead are upper-tie type multiples, namely of 11 (12), then only
they and co-amplitudes of the form $39+8$ $(19-4)$ escape sign inversion.\medskip{}

\noindent Let us resume the creation of the eighteen-dimensional
chunk. For $L_{18}$, we reported a match with $\varphi_{499}^{(1144)}=7398$
via Mersenne fluctuations of Type I, $\left\lfloor \log_{2}C_{63}\right\rfloor \left(\frac{2^{n}}{\pi}\right)^{-1}$
. As far as can be told, no matches with any $\varkappa_{\alpha}^{(n)}$
exist for this kissing number, only close pivots. One obeying constraint
(1) is $\varkappa_{78}^{(2016)}=7223$ $-$ a lower-than-target situation.
The factor 7 in $2016=2^{5}\cdot3^{2}\cdot7$ directly gives the number
of inclusions for $\varkappa_{\xi}^{(n)}$ on the pivot's left: $\begin{array}{lr}
83,1,7,1,1,7,84, & \mathbf{\!\!\!\!7223}\\
{\scriptstyle \:71} & {\scriptstyle \!\!\!\!78}
\end{array}$. Cleaned for duplicates, this yields a quadruple $(\textrm{co-}\varkappa)_{78}=(84,\!7,\!1,\!83)$
of distinct co-amplitudes, leading to a label $\Delta_{\varkappa}$
on $\Lambda_{78}^{(2016)}$ encoded by
\begin{equation}
(1,1,1,1)\cdot(\textrm{co-}\varkappa)_{78}^{\vphantom{}^{{\scriptstyle t}}}\:(=175)\label{eq:L_18-1}
\end{equation}
 to mirror the equation $L_{18}=\varkappa_{78}^{(2016)}+\Delta_{\varkappa}$. 

\noindent For the creation of the nineteenth dimensional chunk, the
same online CFR calculator offers a candidate pivot $\varkappa_{341}^{(1736)}=10808$
with target $L_{19}=10668$. It satisfies constraint (1), and, with
$n=1736$ factoring as $2^{3}\cdot7\cdot31$, we encounter a $p_{j-2},p_{j}\in M_{\textrm{reg}}$
situation, hence constraint (2b) applies and we are led to a string
of $31$ co-amplitudes $\varkappa_{\xi}^{(n)}$ to the pivot's (boldface)
right,\\
$\begin{array}{rr}
\mathbf{\!\!\!\!10808}, & \!\!\!\!\!4,\!7,1,\!73,\!18,1,6,1,1,\!39,1,1,2,2,\!34,5,1,1,2,4,2,\!13,3,1,\!9,1,1,3,8,\!16,\!24.\\
\!\!\!\!\!{\scriptstyle 341} & {\scriptstyle \!\!\!\!\!372}
\end{array}$ $\theta=16$ are distinct, so $(\textrm{co-}\varkappa)_{341}=(4,7,1,73,18,6,39,2,34,5,$$13,3,9,8,16,24).$
We recognize that three of them are lower-tie type multiples, namely
of 6, and one is equal to $9+4.$ Then, under constraint (2b$'$),
\vspace{-0.2cm}
\begin{equation}
(-1,-1,-1,-1,1,1,-1,-1,-1,-1,1,-1,-1,-1,-1,1)\cdot(\textrm{co-}\varkappa)_{341}^{\vphantom{}^{{\scriptstyle t}}}\:(=-140)\label{eq:L_19-1}
\end{equation}
 is a node label $\Delta_{\varkappa}$ on the organizer boundary $\Lambda_{341}^{(1736)}$
that mirrors the equation $L_{19}=\varkappa_{341}^{(1736)}+\Delta_{\varkappa}$.

\noindent Similarly works the creation of the twenty-dimensional
chunk with a candidate pivot for $L_{20}=17400$, $\varkappa_{374}^{(1520)}=17949$.
Again, constraint (1) is obeyed. The factorization of $n=1520$ being
$2^{4}\cdot5\cdot19$, it turns out that nineteen co-amplitudes to
the pivot's right fail to fill the gap; instead, 19 is treated as
$o_{j-2}\!\in\! M_{\nicefrac{5}{8}}$ so that $o_{j}=79$ co-amplitudes
to the pivot's (boldface) right are included, \emph{i.e.} $\begin{array}{lr}
\mathbf{17949},\:28,1,1,2,12,\ldots, & \!\!\!\!\!407\\
{\scriptstyle 374} & {\scriptstyle 453}
\end{array},$ $\theta=19$ of them distinct, and a tuple $(\textrm{co-}\varkappa)_{374}=(28,\!1,\!2,\!12,\!8,\!3,5,\!4,\!7,$
$\!45,\!9,32,\!22,\!17,6,15,20,14,407)$ emerges. Again, we recognize
co-amplitudes for whom there is sign preservation under constraint
(2b$'$): three that are lower-tie type multiples, namely of 7, and
one of the form $9-4$, so that
\begin{equation}
\begin{array}{ll}
(1,-1,-1,-1,-1,-1,1,-1,1,-1,-1,-1,-1,-1,-1,-1,-1,1,-1)\\
\cdot\:(\textrm{co-}\varkappa)_{374}^{\vphantom{}^{{\scriptstyle t}}}\:(=-549)
\end{array}\label{eq:L_20-1}
\end{equation}
 is a node label $\Delta_{\varkappa}$ on a new organizer boundary
$\Lambda_{374}^{(1520)}$, mirroring $L_{20}=\varkappa_{374}^{(1520)}+\Delta_{\varkappa}$.

\medskip{}

\noindent The lower-tie situations bear the imprint of duality controls
$-$ the Catalan number part $(C_{6}\textrm{B}(6,7))^{-1}$ being
expressed by sign preservation in multiples of tie numbers 6 and 7;
and the $M_{\nicefrac{5}{8}}$ part by sign preserving in $9\pm4$
respectively: 13 is the minimal residue in the set of residues $\{13,27\}$
$(=(5\cdot2^{n}-1+2^{n+1})\,\textrm{mod}\, C_{5})$ and 5 the minimal
one in $\{5,11,9\}$ $\left(=(5\cdot2^{n}-1-2^{n+1})\,\textrm{mod}\, C_{4}\right)$.
Transposed to their upper-tie counterparts, with signs preserved in
multiples of tie numbers 11 and 12 and in $39+8$,$19-4$ respectively,
the analogous minimal residue in the set of residues $\{(11,23,\!)47,95,59,\!119,\!107,83,(35,\!)71\}$
$(=(5\cdot2^{n}-1+2^{n})$ $\textrm{mod}\, C_{6})$ as restricted
to the interval $]C_{5},C_{6}[$ would be 47, and the analogous minimal
residue in the set of residues $\{(7,)15,31,21(,1)\}$ $\left(=(5\cdot2^{n}-1-2^{n})\,\textrm{mod}\, C_{5}\right)$
as restricted to the interval $]C_{4},C_{5}[$ would be 15. In summary,
we may conclude minimal residue-selection of the above kind, combined
with selection of the tie number with the most multiples in $(\textrm{co-}\lambda)_{\alpha}$,
are part and parcel in coding sign regulations in gap filling-in.
More instances of dimensional creation have to be examined before
one can definitely say a new invariant is looming here $-$ namely
that the number of distinct co-amplitudes having their signs preserved
in gap filling-in equals 4, as suggested by Eqs. (\ref{eq:L_18-1})-(\ref{eq:L_20-1}).

\noindent The case with multiples will be pursued in Sect.\ref{sub:`Field-',-`projection-'-and}
but before that, a peculiar aspect raised by Eqs. (\ref{eq:L_18-1})-(\ref{eq:L_20-1})
concerning the connection between key particle creation-related dimensions
and nature's forces is paid attention.

\medskip{}

\subsection*{Dimensions (4;18,19,20), (24;46,47,48) and forces}

Subtract $G_{\textrm{max}}^{(15)}$ from the fine-tuned approximation
\textgreek{a}$^{-1}=137.035999547\ldots$ and you come close to the kissing
number $L_{4}=24$. Now the same Table \ref{tab:Key-particle-creation-related-2}
that reveals $L_{4}$ as the kissing number of the least key particle
creation-related source dimension also assigns $L_{19}=10668$ to
the highest such dimension (in that table part). Do these assignments
give us a handle on the glueing (strong) force? Denoting the greatest
prime factor of $n$, $\textrm{gpf}(n)$, we see (cf. Table \ref{tab:Prime-factors-of})
that $\textrm{gpf}(L_{19})$ also determines the $L_{\textrm{up}}$
of the charm quark, the heaviest up-type quark allowing bound states
to exist, while $L_{4}$ is the least among the $L_{\textrm{down}}$'s
of the down quark $-$ a set-up which encourages us to formulate the
following CFR-related
\begin{conjecture}
\label{If-a-dimension-1-1}If a dimension $\nu$ is greatest among
particle creation-related source dimensions, it induces the following
pattern: (A) there exists a lower pivot \textup{$\varkappa_{\alpha}^{(n_{\nu-1})}<L_{\nu-1}$
}\textup{\emph{with}}\textup{ $\sqrt{\frac{\varkappa_{\alpha_{\phantom{}}}^{(n_{\nu-1})}}{L_{\nu-1}}}\gtrsim$
$\frac{5\pi}{16}$ }\textup{\emph{such that distinct organizer co-amplitudes}}\textup{
$\left(\textrm{co-}\varkappa\right)_{\alpha}$ }\textup{\emph{are
recruited to the left of the pivot and (A')}}\textup{ $\left(\textrm{gpf}(L_{\nu-1})\right)^{-1}$
}\textup{\emph{represents an upper bound to a physical force constant;
conversely (B) there exists a higher pivot}}\textup{ $\varkappa_{\alpha}^{(n_{\nu+1})}>L_{\nu+1}$
with $\sqrt{\frac{\varkappa_{\gamma_{\phantom{}}}^{(n_{\nu+1})}}{L_{\nu+1}}}\lesssim\frac{16}{5\pi}$
}\textup{\emph{such that its}}\textup{ $\left(\textrm{co-}\varkappa\right)_{\alpha}$
}\textup{\emph{in combination with those of }}\textup{$\varkappa_{\alpha}^{(n_{\nu})}>L_{\nu}$,
$\sqrt{\frac{\varkappa_{\gamma_{\phantom{}}}^{(n_{\nu})}}{L_{\nu}}}\lesssim\frac{16}{5\pi}$,
}\textup{\emph{are recruited to the right}} and (B') \textup{$\textrm{gpf}(L_{\nu-1})\textrm{gpf}(L_{\nu})L_{\nu+1}^{-1}$
}represents a lower bound to a complementary physical force\textup{
}\textup{\emph{constant}}\textup{.}\textup{\emph{ }}
\end{conjecture}
\noindent Parts (A) and (B)  led to Eqs. (\ref{eq:L_18-1})$-$(\ref{eq:L_20-1})
and have proven true for them. As to (A'), we refer the reader to
Table \ref{tab:Prime-factors-of}, where $\textrm{gpf}(L_{18})$ corresponds
to the integer part of \textgreek{a}$^{-1}$, hence leads to an upper bound
for $\alpha$, $1/137>1/137.035999$. As for (B'), the same table
shows that $\textrm{gpf}(L_{18})\textrm{gpf}(L_{19})L_{20}^{-1}=17399/17400<1$.
We are free to interpret 1 as the coupling constant of the glueing
force $-$ and may ask if the weak force fit in in this picture as
well. 

\noindent It is known that the weak and electromagnetic forces are
siblings, so the solution to approach the weak force coupling constant
is to generalize the fine-tuned calculation of \textgreek{a}$^{-1}$ to %
\footnote{\ a tedious {\tt int64} computation yields $G_{\textrm{max}}^{(63)}=3,\!512,\!576,\!820,\!924,\!177$
from which the value of $\Phi^{(63)}$ follows; it was shown in \cite{Merkel}
that $\Phi^{(p_{n})}$ assumes the limit 4 as $n\rightarrow\infty$%
} 
\[
\left(\textrm{\textgreek{a}}'\right)^{-1}=(C_{15}+C_{11})/\Phi^{(63)}-\frac{1}{f_{n+1}+\kappa'}\,\Delta_{\textrm{b}'}=2\,626\,851.2772574808\ldots,
\]
where $\Delta_{\textrm{b}'}=(C_{15}+C_{11})/\Phi^{(63)}\!-(C_{15}-C_{11})/\Phi^{(31)}$,
$f_{n+1}=31$ and where functional (\ref{eq:timesonover}) is applied
to $G_{\textrm{max}}^{(31)}/C_{7}$ and $G_{\textrm{max}}^{(63)}/C_{15}$
with $\divideontimes:=\times\:$, $f^{(c)}:=\:^{{\scriptscriptstyle 15}}\!\!\!\surd\:$,
$f^{(b)}:=\:^{{\scriptscriptstyle 7}}\!\!\!\surd\:$ and $f^{(a)}:=-\frac{1}{8}+\surd\:$
:  %
\footnote{\  note that $-\frac{1}{8}$ arises by way of $\kappa'-\sqrt{\Phi^{(31)}\Phi^{(63)}}=\frac{9}{8}p_{n}-\sigma_{n-1}$.
Here we have a fusion of interordinality and juxtaposition $n$ vs.
$n\!-\!1$ that necessitates a supplementation of Table \ref{tab:Key-particle-creation-related-2}
which will be given in the closing remarks (Sect.\ref{sec:Supplementary-remarks-and})%
}
\begin{equation}
\kappa'=-\frac{1}{8}+\sqrt{\Phi^{(31)}\Phi^{(63)}}.\label{eq:control-functional}
\end{equation}
We note the following: just as the fine-tuning \textgreek{a}$^{-1}-G_{\textrm{max}}^{(15)}$
results in a real number that is only a tiny bit greater than $L_{4}=24$, bringing with it the photon as a massless force carrier,
the analogous fine-tuning $\left(\textrm{\textgreek{a}}'\right)^{-1}-G_{\textrm{max}}^{(31)}$
produces a real a tiny bit greater than the kissing number\label{weakpoin} $L_{24}=196560$.
The 24th dimension then, which is not among the $F_{e}$ and $F_{s}\,$-related
source dimensions mentioned in Table \ref{tab:Key-particle-creation-related-2},
might turn out to be $F_{w}\,$-related and contribute the missing
coupling constant to complement the force trio with \textgreek{a}$_{w}:=\textrm{\textgreek{a}}'\approx3.81\times10^{-7}$ and come along with massive vector bosons as force carriers.
\\
Conjecture \ref{If-a-dimension-1-1} deals with the combination of
a triad of dimensions with a force duo. To show how bounds to all
three couplings could be intertwined, we have to modify that conjecture
such that a dimensional triad covers three couplings. To this end,
we make a series of qualitative statements.\\
 (C) $L_{24}$ is the first kissing number to resume the normal successor
relation $L_{m+1}<2L_{m}$ after the hiatal $L_{23}$ where $L_{24}=2L_{23}+10260$.
We may assume the formation law of the hiatus is $m_{\textrm{hiat}}=p_{k}+2^{k-1}\:(p_{k}\in M_{\textrm{reg}})$;
then $23=15+8$ would be followed by another hiatus at $47=31+16$.\\
 (D) According to Conjecture \ref{If-a-dimension-1-1}, the center
dimension of the dimensional triad is 19. Asssuming a formation law
$m_{\textrm{cent}}=2C_{k-1}+o_{k-2}\:(o_{k-2}\in M_{\nicefrac{5}{8}})$,
the dimension $19=10+9$ would be followed by another center dimension
$47=28+19$.\\
 (E) The left companion in the new triad is dimension 46. A serious
candidate for the kissing number $L{}_{46}$ would be $12\,986\,152$:
Mersenne fluctuations of Type I, $\left\lfloor \log_{2}C_{q_{s}}\right\rfloor \,\left(\frac{2^{n}}{\pi}\right)^{-1}\!,$
at secondary expansion $s=7$ and time-like refinement levels $n$
that are multiples of $12$, do indeed include the global (perhaps
matching) amplitude $\varphi_{93}^{(1152)}=12\,986\,152$ (Table \ref{tab:Specific-fractions-in}).
\\
(F) The right companion is dimension 48, which has a (certified) kissing
number $L_{48}=52\,416\,000$. Because of the assumptions we have
made, the unknown kissing number $L_{47}$ must obey the constraints
$L_{47}<2L_{46}$ and $L_{48}\geq2L_{47}$. We note in passing that
the assumed value $L_{46}$ has the prime factorization $12\,986\,152=2^{3}\cdot1\,623\,269$.
If gpf($L_{46}$) is to give rise to a lower bound to the glueing
force coupling, it must combine with a $L_{47}$ that has 31 among
its prime factors because the closest analog to the previously found
bound $17399/17400$ would be $1\,623\,269\cdot31/52\,416\,000<1$.
The ansatz to assure this is 
\[
L_{47}^{'}:=2L_{46}+(8-31\eta)10260\quad(\eta=1,2,\ldots)
\]
where 10260 is borrowed from the hiatus at $L{}_{23}$. As a rule,
instances of $L_{47}^{'}$ will contain also prime factors $n_{j},n_{k},\ldots>31$.
In order to eliminate them, we have to transform $L_{47}^{'}$ via
$n_{j}\pm1,n_{k}\pm1,\ldots$ into $L_{47}^{''}$. For instance, $L_{47}^{'}:=2L_{46}-85\cdot10260=2^{2}\cdot31\cdot71\cdot2851$
$\mapsto$ $L_{47}^{''}=2^{2}\cdot31\cdot72\cdot2850$ $\equiv$ $2^{6}\cdot3^{3}\cdot5^{2}\cdot19\cdot31$.
Transformation along these lines doesn't always `succeed.' Fortunately
so, since one of the `failures' turns out to be a `lucky take' : $L_{47}^{'}:=2L_{46}-54\cdot10260=2^{3}\cdot17\cdot31\cdot6029$
$\mapsto$ $L_{47}^{''}=2^{3}\cdot17\cdot31\cdot6028$ $\equiv$ $2^{5}\cdot11\cdot17\cdot31\cdot137$,
the latter being the prime factorization of our now preferred candidate
at center dimension 47, $L_{47}:=25\,414\,048$. \\
The combined result of our qualitative statements then is the following
bounds: 

\indent for \textgreek{a}$_{w}$ $\qquad\!\!\!\!\left(\textrm{gpf}(L_{46})\right)^{-1}=1/1\,623\,269>1/2\,626\,851.277257;$\\
\indent for \textgreek{a}, $\qquad\left(\textrm{gpf}(L_{47})\right)^{-1}=1/137>1/137.035999;$
and\\
\indent for \textgreek{a}$_{s},$ $\qquad\!\!\textrm{gpf}(L_{46})\,\textrm{sgpf}(L_{47})L_{48}^{-1}=50\,321\,339/52\,416\,000<1$

\noindent where $\textrm{sgpf}(\cdot)$ connotes `second greatest
prime factor.' Provided the postulated values do indeed embody the
true kissing numbers $L_{46}$ and $L_{47}$, the new triad preserves
the tight bound for the electromagnetic coupling, combining it with
loosened bounds for the weak and the strong couplings, respectively.
To draw a link between the global and local perspective, we first
discuss the tight bounds-only case of the former triad, then proceed
with the loosened bounds-case of the second triad. In what follows
we use the shorthand $(\varphi_{\alpha}^{(n;n_{\textrm{low}}\leq n\leq n_{\textrm{high}})}\mid\varphi_{\alpha_{n_{\textrm{pivot}}}}^{(n_{\textrm{pivot}})}={\bf {pivot}})$
for the Mersenne fluctuation $(\varphi_{\alpha_{n_{\textrm{low}}}}^{n_{\textrm{low}}},\ldots,{\bf {pivot}},\ldots,\varphi_{\alpha_{n_{\textrm{high}}}}^{n_{\textrm{high}}})$.

\subsection*{\label{sub:Link-between-local}Global lead for the dimensional triad
18,19,20}

It was mentioned that if $\nu$ ($\nu'$) is one of the key particle
creation-related `source' (`sink') dimensions (Table \ref{tab:Key-particle-creation-related-1-1}),
there are ideations $\varphi_{\alpha}^{(n)}+\Delta\varphi^{(n)}=L_{\nu}$
$(\varphi_{\alpha'}^{(n')}+\Delta\varphi^{(n')}=L_{\nu'})$ which
are bound to Mersenne fluctuations of Type I, $\left\lfloor \log_{2}C_{q_{s}}\right\rfloor \,\left(\frac{2^{n}}{\pi}\right)^{-1}\!,$
at secondary expansion $s=6$ and time-like refinement levels $n$
that are multiples of $11$; we interpreted this as a hint that creation
follows the global givens by way of the left-strand effect. Here,
we show that global pivots which the organizer pivot $\lambda_{78}^{(2016)}=7223$
snaps into spring from the \emph{same} type of Mersenne fluctuations.
A key feature of this sort of pivots is that they represent the left
strand equivalently by time-like refinement levels $n$ and space-like
refinement offsets $\triangle\alpha(\equiv\alpha_{\psi}-\alpha_{\varphi})$,
where $\varphi$ denotes amplitude, and $\psi$ phase. A organizer
pivot $\varkappa_{\alpha_{w}}^{(n_{w})}$ is said to snap into global
pivots $\varphi_{\alpha_{v}}^{(n_{v})}$ and $\varphi_{\alpha_{u}}^{(n_{u})}$
if it is able to unfreeze them from their frozen-in status they have
since expansion $2{}^{n_{w}}$ went past their ideational counterparts
$2{}^{n_{u}}$ and $2{}^{n_{v}}$. That is, $n_{w}>n_{v},n_{u}$.
For dimension 18, whose kissing number is perfectly matched with $\varphi_{499}^{(1144)}=7398$
but bound up with $\varkappa_{78}^{(2016)}=7223$, we find $\varphi_{\alpha_{u}}^{(n_{u})}:=\varphi_{97}^{(90)}=7569$.
To have this overshoot fit with the organizer pivot, the excess in
amplitude has to be compensated for by a $\varphi_{\alpha_{v}}^{(n_{v})}$
missing its target by $-346$. The global deficient to back this up
is with target $L_{20}$, comes from within the ranks of a Mersenne
fluctuation $(\varphi_{\alpha}^{(n;1681\leq n\leq1701)}\mid\varphi_{35}^{(1690)}=\mathbf{17054})$
and incidentally is as `leggy' a pivot therein as the overshoot is
in $(\varphi_{\alpha}^{(n;79\leq n\leq99)}\mid\varphi_{97}^{(90)}=\mathbf{7569})$:
\begin{table}[H]
\caption{`Leggy' global pivots (boldface), embedded in their respective frozen-in
Mersenne fluctuations $(\varphi_{\alpha}^{(n;79\leq n\leq99)}\mid\varphi_{97}^{(90)}=\mathbf{7569})$,
$(\varphi_{\alpha}^{(n;1681\leq n\leq1701)}\mid\varphi_{35}^{(1690)}=\mathbf{17054})$
of type $\left\lfloor \log_{2}C_{63}\right\rfloor \,\left(\frac{2^{n}}{\pi}\right)^{-1}$}
 \medskip{}

\quad{}\quad{}%
\begin{tabular}{|c|c|c|c|}
\hline 
$n$ & \multicolumn{2}{c|}{$\varphi_{\alpha}$} & $n$\tabularnewline
\hline 
\hline 
89 & \multicolumn{2}{c|}{15139} & 89\tabularnewline
\hline 
88 & 7569 & $\mathbf{7569}$ & 90\tabularnewline
\hline 
87 & 3784 & 3784 & 91\tabularnewline
\hline 
86 & 1892 & 1891 & 92\tabularnewline
\hline 
85 & 946 & 945 & 93\tabularnewline
\hline 
84 & 472 & 472 & 94\tabularnewline
\hline 
83 & 235 & 235 & 95\tabularnewline
\hline 
82 & 117 & 117 & 96\tabularnewline
\hline 
81 & 58 & 58 & 97\tabularnewline
\hline 
80 & 28 & 28 & 98\tabularnewline
\hline 
79 & 13 & 13 & 99\tabularnewline
\hline 
$\vdots$ &  &  & $\vdots$\tabularnewline
\hline 
\end{tabular}\,%
\begin{tabular}{|c|c|c|c|}
\hline 
$n$ & \multicolumn{2}{c|}{$\varphi_{\alpha}\:(\psi_{\alpha})$} & $n$\tabularnewline
\hline 
\hline 
1691 & \multicolumn{2}{c|}{34108} & 1691\tabularnewline
\hline 
1690 & $\mathbf{17054}$ & 17053 (-17054) & 1692\tabularnewline
\hline 
1689 & 8526 & 8526 & 1693\tabularnewline
\hline 
1688 & 4262 & 4263 & 1694\tabularnewline
\hline 
1687 & 2131 & 2131 & 1695\tabularnewline
\hline 
1686 & 1065 & 1065 & 1696\tabularnewline
\hline 
1685 & 532 & 532 & 1697\tabularnewline
\hline 
1684 & 265 & 266 & 1698\tabularnewline
\hline 
1683 & 132 & 132 & 1699\tabularnewline
\hline 
1682 & 66 & 66 & 1700\tabularnewline
\hline 
1681 & 32 & 32 & 1701\tabularnewline
\hline 
$\vdots$ &  &  & $\vdots$\tabularnewline
\hline 
\end{tabular}
\end{table}

\noindent With the phase at $n_{v}=1692$ (in parentheses) looping
into $\psi_{34}^{(1692)}=-17054$, we are lucky to find one value,
$\left|\psi_{34}^{1692}\right|=17054$, that is deficient with respect
to the target $L_{20}=17400$ by the same amount as $\varphi_{35}^{(1690)}$.
With $6\mid n_{v}$ and $5\mid n_{u}$, and also $(6+5)\mid(n_{v}+n_{u})$,
we come upon the time-like refinement realisation of the left strand
depicted in Table \ref{tab:Time-like-refinements mimicking left strand-1}:

\begin{table}[H]
\caption{\label{tab:Time-like-refinements mimicking left strand-1}Time-like
refinement levels for global pivots from fluctuations of type $\left\lfloor \log_{2}C_{63}\right\rfloor \,\left(\frac{2^{n}}{\pi}\right)^{-1}$
the local pivot $\varkappa_{78}^{(2016)}=7223\,$ snaps into }
\bigskip{}

\qquad{}%
\begin{tabular}{|c|c|c|cc|}
\hline 
 & \multicolumn{2}{c|}{overshoot \& deficient} & \multirow{2}{*}{} & \multirow{2}{*}{left strand}\tabularnewline
\multirow{1}{*}{$n$} & \multirow{1}{*}{$\varphi_{u},\left|\psi_{v}\right|$} & \multirow{1}{*}{$\Sigma x\!\mid\!\Sigma n,\, x\!\mid\! n$} &  & \tabularnewline
\hline 
\hline 
 &  & \multirow{2}{*}{$11\!\mid\!1782$} & \multirow{2}{*}{$\longmapsto$} & \multirow{2}{*}{tie number 11}\tabularnewline
 &  &  &  & \tabularnewline
\multirow{3}{*}{90} & \multirow{3}{*}{7569} & \multirow{3}{*}{$5\!\mid\!90$} & \multirow{3}{*}{$\longmapsto$} & \multirow{3}{*}{inter-tie increment 5}\tabularnewline
 &  &  &  & \tabularnewline
 &  &  &  & \tabularnewline
1692 & 17054 & \multirow{1}{*}{$6\!\mid\!1692$} & \multirow{1}{*}{$\longmapsto$} & \multirow{1}{*}{tie number 6}\tabularnewline
\hline 
\end{tabular}
\end{table}
 \noindent The dual realization of the left strand in terms of space-like
refinement offsets works by executing a switch back to the left-leg
$n_{v}=1690$. The effect is as envisioned in the first place: $\varphi_{35}^{(1690)}=17054$
is deficient with respect to $L_{20}$ by $-346$, but now the tale
is told by offsets. 
\begin{table}[H]
\caption{\label{tab:Space-like-refinement-offsets-1}Space-like refinement
offsets for global pivots from fluctuations of type $\left\lfloor \log_{2}C_{63}\right\rfloor \,\left(\frac{2^{n}}{\pi}\right)^{-1}$
the organizer pivot $\varkappa_{78}^{(2016)}=7223\,$ snaps into }
\bigskip{}

\qquad{}%
\begin{tabular}{|c|c|c|c|cc|}
\hline 
 & \multicolumn{3}{c|}{overshoot \& deficient} &  & \multirow{2}{*}{left strand}\tabularnewline
$n$ & \multirow{1}{*}{$\varphi,\psi$} & \multirow{1}{*}{$\alpha_{\varphi},\alpha_{\psi}$} & \multirow{1}{*}{$\Sigma\Delta\alpha,\Delta\alpha$} &  & \tabularnewline
\hline 
\hline 
\multirow{2}{*}{} &  &  & \multirow{2}{*}{11} & \multirow{2}{*}{$\longmapsto$} & \multirow{2}{*}{tie number 11}\tabularnewline
 &  &  &  &  & \tabularnewline
\multirow{2}{*}{90} & 7569, & 97, & \multirow{2}{*}{5} & \multirow{2}{*}{$\longmapsto$} & \multirow{2}{*}{inter-tie increment 5}\tabularnewline
 & -7570 & 102 &  &  & \tabularnewline
\multirow{2}{*}{1690} & 17054, & 35, & \multirow{2}{*}{6} & \multirow{2}{*}{$\longmapsto$} & \multirow{2}{*}{tie number 6}\tabularnewline
 & 17054 & 41 &  &  & \tabularnewline
\hline 
\end{tabular}
\end{table}

\noindent We conclude that for the dimensional triad 18,19,20 global
pivots that are unfrozen by a organizer pivot provide a natural picture
of the left strand in being complementary to the one based on matching
global amplitudes (Table \ref{tab:Key-particle-creation-related-1-1}).
The representations are in terms of time-like refinement levels or
space-like refinement offsets, the former aided by switching the `standing
leg' and substituing absolute value in phase. 

\noindent The unfreezing embraces the entire Mersenne fluctuations,
not only the pivotal amplitudes. This is suggested by 
\[
\begin{array}{c}
\left\lfloor L_{20}/2^{i}\right\rfloor +\varkappa_{\alpha_{w(i)}}^{(w(i))}-\varphi_{\alpha_{v(i)}}^{(v(i))}-\varphi_{\alpha_{u(i)}}^{(u(i))}=\delta\quad(\delta\in\{0,1\}),\end{array}
\]
where $w(i)=2016\mp i,v(i)={}_{{\scriptstyle 1692+i}}^{{\scriptstyle 1690-i}},u(i)=\,_{{\scriptstyle 90+i}}^{{\scriptstyle 88-i}}\;\left(i=\!0,\!1,{\scriptstyle \!\ldots},\!9\right)$,
$(13,\!27,\!55,\!112,\!225$, $451,\!902,\!1805,\!3611,\mathbf{\!7223},3611,\!1805,\!902,\!451,\!225,\!112,\!56,\!27,\!13)$
is the fully traced organizer fluctuation $(\varkappa_{\alpha}^{(n;2007\leq n\leq2025)}\mid\varkappa_{78}^{(2016)}=\mathbf{7223})$
and the expressions $\left\lfloor L_{20}/2^{i}\right\rfloor $ are
the surrogate for the left and right leg of an (auxiliary fourth)
fluctuation.

\subsubsection*{\label{sub:Second-Link-between-local}Organizer lead for the  dimensional
triad 46,47,48}

Within the relatively loose bounds to the weak and strong couplings
that characterize the dimensional triad 46,47,48, organizer amplitudes
defined by 
\[
(\textrm{CFR})\:2^{-n}\kappa'\!\rightarrow\!\left[(\varkappa')_{0}^{(n)}\!;(\varkappa')_{\alpha}^{(n)}\right],
\]
will have properties that are very different from what we know about
$\varkappa{}_{\alpha}^{(n)}$. Instead of a prime factorization of
$n$ going by the rules (2)$-$(2b), the prime factorization of a
key amplitude $(\varkappa')_{\alpha_{z}}^{(n_{z})}$ determines a
co-structure of dispersed amplitudes $(\varkappa')_{\alpha_{y}}^{(n_{y})}$,
$(\varkappa')_{\alpha_{x}}^{(n_{x})},\ldots$ $(n_{z}\geq n_{y}\geq n_{x}\geq\ldots)$
all of which have to be included to achieve a match with a targeted
kissing number. (2c) If this factorization has the form $v\cdot$gpf$(\cdot)$,
and gpf$(\cdot)$ coincides with $p\in M_{\textrm{reg}}$ ($o\in M_{\nicefrac{5}{8}}$),
then $v-1$ amplitudes join the co-structure. And in case $v$ coincides
also with the inter-tie increment 5, the key local amplitude plus
those of the co-structure are found to snap into a global amplitude
matching that kissing number so that the interaction may be caught
in a dual strand representation in terms of time-like refinement levels
and space-like refinement offsets respectively. \\
We pointed out the possibility of $L_{46}:=\varphi_{93}^{(1152)}=12\,986\,152$,
a global amplitude that would satisfy the property pairing $12\mid n$,
$s=7$ and represent an ideation of creation in right-strand mode.
To illustrate how this is realized in organization, let $(\varkappa')_{15}^{(2913)}=95$
take the role of key amplitude, noting that its factorization $95=5\cdot19$
meets the demands. Then the associated space-like refinements $(\alpha_{z})_{\varphi}=15$
and $(\alpha_{z})_{\psi}=17$ are found to work as cross moduli for
the time-like refinement levels and space-like refinement offsets
of the dispersed amplitudes that jointly form the co-structure\\
\indent $\!(\varkappa')_{46}^{(2913)}=209\,406$ with $\alpha_{\varphi}=46$,
$\alpha_{\psi}=51$, \\
\indent $\!(\varkappa')_{410}^{(1784)}=50\,458$ with $\alpha_{\varphi}=410$,
$\alpha_{\psi}=471$, \\
\indent $\!(\varkappa')_{400}^{(1678)}=44\,750$ with $\alpha_{\varphi}=400$,
$\alpha_{\psi}=436$, \\
\indent $\!(\varkappa')_{461}^{(647)}=12\,681\,443$ with $\alpha_{\varphi}=461$,
$\alpha_{\psi}=485$. 

\noindent To begin with the latter, we have (key amplitude in parentheses)\vspace{-0.5cm}
 
\begin{table}[H]
\caption{\label{tab:Space-like-refinement offsets mimicking right strand}Space-like
refinement offsets of organizer pivots $\varkappa'$ snapping into
$\varphi_{93}^{(1152)}=12\,986\,152$}
\bigskip{}

\ \,%
\begin{tabular}{|c|c|c|cc|}
\hline 
 & \multicolumn{2}{c|}{co-structure by cross moduli } & \multirow{2}{*}{} & \multirow{2}{*}{right strand}\tabularnewline
\multirow{1}{*}{$\varkappa'$} & \multirow{1}{*}{$\alpha_{\varphi},\alpha_{\psi}$} & \multirow{1}{*}{$\Sigma\Delta_{\alpha}\textrm{mod}\,\alpha_{z}$} &  & \tabularnewline
\hline 
\hline 
12\,681\,443 & $461,485$ & \multirow{2}{*}{$12\equiv(24+5)(\textrm{mod}\,17)$} & \multirow{2}{*}{$\mapsto$} & \multirow{2}{*}{tie number 12}\tabularnewline
209\,406 & $46,51$ &  &  & \tabularnewline
\multirow{3}{*}{(95)} & \multirow{3}{*}{} & \multirow{3}{*}{$(v=5)$} & \multirow{3}{*}{$\mapsto$} & \multirow{3}{*}{inter-tie increment 5}\tabularnewline
 &  &  &  & \tabularnewline
 &  &  &  & \tabularnewline
50\,458 & $410,471$ & \multirow{2}{*}{$7\equiv(61+36)(\textrm{mod}\,15)$} & \multirow{2}{*}{$\mapsto$} & \multirow{2}{*}{tie number 7}\tabularnewline
44\,750 & $400,436$ &  &  & \tabularnewline
\hline 
\end{tabular}
\end{table}
 \noindent which is flanked by the dual \vspace{-0.5cm}
 
\begin{table}[H]
\caption{\label{tab:Time-like-refinements mimicking lright strand}Time-like
refinement levels of organizer pivots $\varkappa'$ snapping into
$\varphi_{93}^{(1152)}=12\,986\,152$}
\bigskip{}

\begin{tabular}{|c|c|c|cc|}
\hline 
 & \multicolumn{2}{c|}{co-structure by cross moduli } & \multirow{2}{*}{} & \multirow{2}{*}{right strand}\tabularnewline
\multirow{1}{*}{$\varkappa'$} & \multirow{1}{*}{$n$} & \multirow{1}{*}{$\Sigma n\,\textrm{mod}\,\alpha_{z}$} &  & \tabularnewline
\hline 
\hline 
50\,458 & $1784$ & \multirow{2}{*}{$12\equiv(1784\!+\!1678)(\textrm{mod}\,15)$} & \multirow{2}{*}{$\mapsto$} & \multirow{2}{*}{tie number 12}\tabularnewline
44\,750 & $1678$ &  &  & \tabularnewline
\multirow{3}{*}{(95)} & \multirow{3}{*}{} & \multirow{3}{*}{$(v=5)$} & \multirow{3}{*}{$\mapsto$} & \multirow{3}{*}{inter-tie increment 5}\tabularnewline
 &  &  &  & \tabularnewline
 &  &  &  & \tabularnewline
12\,681\,443 & $647$ & \multirow{2}{*}{$7\equiv(647\!+\!2913)(\textrm{mod}\,17)$} & \multirow{2}{*}{$\mapsto$} & \multirow{2}{*}{tie number 7}\tabularnewline
209\,406 & $2913$ &  &  & \tabularnewline
\hline 
\end{tabular}
\end{table}
\noindent Speaking in terms of expansion in the above example, only
$2.35\%$ of global amplitude $\varphi_{93}^{(1152)}\!=\!12\,986\,152$
get `unfrozen' in the usual way by organizer amplitudes: $2{}^{1152}\!<2^{1678}\!<2^{1784}\!<2^{2913}$;
organizer amplitude $(\varkappa')^{(647)}\!=$12\,681\,443 for which
$2{}^{1152}>2^{647}$, by contrast, makes up $97.65\%$. As far as
creation in right-strand mode is concerned, $\varkappa'$ obviously
has all of ideational expansion at its disposal.\\

\subsection{\label{sub:`Field-',-`projection-'-and}`Field-', `projection-' and
`spacetime' simulacrum}

\subsubsection*{From croton base numbers to a simulacrally decorated double strand}

\noindent Constraint (2b$'$) and the comments that followed it suggest
that all entries of the Catalan double strand pattern 
\[
\begin{array}{ccc}
\!\!{}_{11} & _{\cdot} & \!\!{}_{12}\\
\!\!\cdot &  & \!\!\cdot\\
\!\!{}^{6} & ^{\cdot} & \!\!{}^{7}
\end{array}
\]
have equal grounds $-$ a hypothesis which needs to be checked. This
may be done using the subgroups of order 4, $\{1,5,9,13\}$ (powers
of 5) and $\{1,15,33,47\}$ (powers of 15), of the group of units
of the quotient rings $\mathbb{Z}/16\mathbb{Z}$ and $\mathbb{Z}/64\mathbb{Z}\:$,
respectively. Denoting the above quadruples of units by $(\textrm{Sim}_{16})$
and $(\textrm{Sim}_{64})$, we arrive at the following decorated version
of the double strand:
\[
\begin{array}{ll}
11: & 12:\\
\textrm{\hbox{\small characteristic\,\ multiple:}} & \textrm{\textrm{\hbox{\small characteristic\,\ multiple:}}}\\
(1,-1,1,1)\cdot(\textrm{Sim}_{64})^{t}=66(=2n_{c}\cdot11) & (1,1,1,1)\cdot(\textrm{Sim}_{64})^{t}=96(=2n_{c}\cdot12)\\
\textrm{\hbox{\small characteristic\,\ quantity\,39+8:}} & \textrm{\textrm{\hbox{\small characteristic\,\ quantity\,19-4:}}}\\
(0,0,0,1)\cdot(\textrm{Sim}_{64})^{t}=47 & (0,1,0,0)\cdot(\textrm{Sim}_{64})^{t}=15\\
\\
6: & 7:\\
\textrm{\hbox{\small characteristic\,\ multiple:}} & \textrm{\textrm{\hbox{\small characteristic\,\ multiple:}}}\\
(1,-1,1,1)\cdot(\textrm{Sim}_{16})^{t}=18(=n_{c}\cdot6) & (1,1,1,1)\cdot(\textrm{Sim}_{16})^{t}=28(=n_{c}\cdot7)\\
\textrm{\hbox{\small characteristic\,\ quantity\,9+4:}} & \textrm{\textrm{\hbox{\small characteristic\,\ quantity\,9-4:}}}\\
(0,0,0,1)\cdot(\textrm{Sim}_{16})^{t}=13 & (0,1,0,0)\cdot(\textrm{Sim}_{16})^{t}=5\\
\\
\qquad\qquad(n_{c}=3) & \qquad\qquad(n_{c}=4)
\end{array}
\]
Just as the characteristic multiples and projections in the lower
half mark sign-preserving co-amplitudes in the gap fillings of lower-tie
type $-$ 18 and 13 in Eq.(\ref{eq:L_19-1}), 28 and 5 in Eq.(\ref{eq:L_20-1})
$-$, the characteristic multiples and projections in the upper half
$-$ 66 and 47 in case of 11, and 96 and 15 in case of 12 $-$ anticipate
sign-preserving co-amplitudes emerging in gap fillings of the upper-tie
type even if we do not, for the time being, know the kissing number
pivots that evoke these situations. 

\noindent The right strand is marked by characteristic multiples
with signature (1,1,1,1) which are called `field' simulacra here;
its characteristic projections with signature $(0,1,0,0)$, just like
their left-strand counterparts with signature $(0,0,0,1)$, are accordingly
named `projection' simulacra. The characteristic multiples given on
the left strand, in being of signature $(1,-1,1,1)$, may be termed
`spacetime' simulacra. The difference between `field' and `spacetime'
simulacra is one of their associated $n_{c}$, or number of positive
entries in the signature. The label `space-time' perhaps becomes clearer
when the time-like refinements are related to space-like refinements
$-$ much like in the previous subsection. In Eqs. (\ref{eq:L_19-1})
and (\ref{eq:L_20-1}), $n$ $-$ the base 2 logarithm of the time-like
refinement $-$ is a multiple of $p_{j}$ $(o_{j})$ or $p_{j}$ $(o_{j-2})$
and is linked to the space-like refinements $\alpha$ and $\alpha+p_{j}$
$(\alpha+o_{j})$ by the relationship $\left|\alpha-C_{x}\right|=n_{y}\cdot11,$
$\alpha+p_{j}-C_{x}=n_{z}\cdot6$ $(\alpha+o_{j}-C_{x}=n_{z}\cdot6)$,
where $n_{y},n_{z}\in M_{\nicefrac{5}{8}}$ or $M_{\nicefrac{5}{8}}^{+}$
$-$ with interchanged roles for Eq. (\ref{eq:L_19-1}) compared to
Eq. (\ref{eq:L_20-1}): \bigskip{}

\hfill{}%
\begin{tabular}{|c|}
\hline 
\tabularnewline
\hline 
\hline 
$\left|\alpha-C_{x}\right|$\tabularnewline
\hline 
$\begin{array}{c}
\alpha+p_{j}-C_{x}\\
(\alpha+o_{j}-C_{x})
\end{array}$\tabularnewline
\hline 
\end{tabular}\,%
\begin{tabular}{|c|ccc|}
\hline 
Eq. (\ref{eq:L_19-1})  &  & Eq. (\ref{eq:L_20-1}) & \tabularnewline
\hline 
\hline 
$341-C_{6}=19\cdot11$ &  & $C_{7}-374=5\cdot11$ & \tabularnewline
\hline 
$\begin{array}{c}
\vphantom{}\\
\vphantom{}
\end{array}$$372-C_{6}=40\cdot6$ &  & $453-C_{7}=4\cdot6$ & \tabularnewline
\hline 
\end{tabular}\hfill{}

\bigskip{}
\noindent While the second expansion parameter $s$ in Table \ref{tab:Key-particle-creation-related-1-1}
has to `come out of the woodwork' $-$ $C_{63}$ identified as $C_{2^{6}-1}$
$-$ to make the left-strand character of time-like parameter pairing
$11\mid n$, $s=6$ apparent, the left-strand affinity of the pairing
of space-like parameters in the above table $-$ $\alpha$ for the
pivot, $\alpha+p_{j}$ for the rightmost co-amplitude $-$ comes to
the fore by letting the lower-tie character retreat into hiding $-$
into index status $x=6$ and $x=7$ for Eqs. (\ref{eq:L_19-1}) and
(\ref{eq:L_20-1}), respectively. The group-theoretic manifestation
of left-strand affinity is the said `spatio-temporal' simulacrum,
where space-like and time-like aspects combine into one signature.

\noindent Although the distinction of simulacral forms may earn them
merits of their own $-$ left-strand affinity and factorization of
$n$ $\rightarrow$ `spacetime' simulacrum, right-strand affinity
and factorization of key amplitude $\rightarrow$ `field' simulacrum$-$,
with respect to \emph{gap filling-in} they are incomplete because
only one $-$ what is meant by characteristic $-$ multiple is produced
for each tie number and we don't get to know the (number of) other
multiple-type, sign-preserving co-amplitudes contributing to the filling-in.
The group-theoretic background however makes the existence of an invariant
number for them plausible.

\subsubsection*{Bound((Sim$_{16}$)), Bound((Sim$_{64}$)) and M$\phantom{}_{\nicefrac{5}{8}}^{+}$ }

The problem can be narrowed down by the following observations: 

\noindent First, in analogy to the constructions Bound$(N_{\textrm{source}})$,
Bound$(N_{\textrm{sink}})$, we may allow for 4-cubes, Bound$((\textrm{Sim}_{64}))$
and Bound$((\textrm{Sim}_{16}))$, and check which multiples of 11,12
are node labels on the former and multiples of 6,7 node labels on
the latter; it turns out four multiples of 11,12 are represented on
$\textrm{Bound}((\textrm{Sim}_{64}))$ and seven multiples of 6,7
on $\textrm{Bound}((\textrm{Sim}_{16}))$:

I: $\qquad\begin{array}{l}
33,48,66,96;\\
6,7,12,14,18,21,28.
\end{array}$

\noindent (In Bound$((\textrm{Sim}_{64}))$, sixty-four out of the
numbers $1,2,$...,$96$ remain unrepresented: $3,4,$...,12; $20,21$,...,27;
$35,36$,...,45; $50,51$,...,60; $67,68$,...,78; $82,83$,...,$93$.
In Bound$((\textrm{Sim}_{16}))$, out of the numbers $1,2,\ldots,28$,
four remain unrepresented: 11,20,24,25.) 

\noindent The second observation draws on the fact that, while one
`spacetime' simulacrum, $(1,-1,1,1)\cdot(\textrm{Sim}_{64})^{t}(=66)$,
suffices to reproduce the left strand $\begin{array}{c}
\!\!{}_{11}\\
\!\!\cdot\\
\!\!{}^{6}
\end{array}\!$, the binary operation `+' on the two of them is required to reproduce
the right strand, $\begin{array}{c}
\!\!{}_{12}\\
\!\!\cdot\\
\!\!{}^{7}
\end{array}\!$, namely $(1,-1,1,1)\cdot(\textrm{Sim}_{64})^{t}+(1,-1,1,1)\cdot(\textrm{Sim}_{16})^{t}=84$.
The latter being a (Janus-faced) multiple, we can use that operation
to find other multiples, but, similar to the constraint (2) for co-amplitudes
$\varkappa_{\xi}^{(n)}\in(\textrm{co-}\varkappa)_{\alpha}$, need
an $M_{\nicefrac{5}{8}}^{+}$ constraint for them as follows: A multiple
of 6, 7, 11 or 12 assuming the form $(\epsilon_{1},\epsilon_{2},\epsilon_{3},\epsilon_{4})\cdot(\textrm{Sim}_{64})^{t}\pm(\epsilon_{1},\epsilon_{2},\epsilon_{3},\epsilon_{4})\cdot(\textrm{Sim}_{16})^{t}$
is admissible in a gap filling situation only if $[(\epsilon_{1},\epsilon_{2},\epsilon_{3},\epsilon_{4})\cdot(\textrm{Sim}_{64})^{t}\mp(\epsilon_{1},\epsilon_{2},\epsilon_{3},\epsilon_{4})\cdot(\textrm{Sim}_{16})^{t}]\in M_{\nicefrac{5}{8}}^{+}$.

\noindent With only two $\epsilon_{x},\epsilon_{y}\neq0$, we register
five solutions:

IIa: $\qquad\begin{array}{lc}
(1,1,0,0)\cdot(\textrm{Sim}_{64})^{t}\pm(1,1,0,0)\cdot(\textrm{Sim}_{16})^{t} & =16\pm6=\begin{array}{c}
22\\
10
\end{array}\!,\\
(-1,1,0,0)\cdot(\textrm{Sim}_{64})^{t}\pm(-1,1,0,0)\cdot(\textrm{Sim}_{16})^{t} & =14\pm4=\begin{array}{c}
18\\
10
\end{array}\!,\\
(-1,0,1,0)\cdot(\textrm{Sim}_{64})^{t}\mp(-1,0,1,0)\cdot(\textrm{Sim}_{16})^{t} & =32\mp8=\begin{array}{c}
24\\
40
\end{array}\!,\\
(0,1,0,1)\cdot(\textrm{Sim}_{64})^{t}\mp(0,1,0,1)\cdot(\textrm{Sim}_{16})^{t} & =62\mp18=\!\begin{array}{c}
44\\
80
\end{array}\!,\\
(0,0,-1,1)\cdot(\textrm{Sim}_{64})^{t}\pm(0,0,-1,1)\cdot(\textrm{Sim}_{16})^{t} & =14\pm4=\begin{array}{c}
18\\
10
\end{array}\!;
\end{array}$

\noindent and with $\epsilon_{1},\epsilon_{2},\epsilon_{3},\epsilon_{4}\neq0$,
one further:

IIb:$\qquad\begin{array}{lc}
(-1,\!1,-1,\!1)\cdot(\textrm{Sim}_{64})^{t}\pm(-1,\!1,-1,\!1)\cdot(\textrm{Sim}_{16})^{t} & \!\!\!=28\pm8=\begin{array}{c}
36\\
20
\end{array}\!.\end{array}$ 

\noindent (Terms involving $\epsilon_{w}=0$, while $\epsilon_{x},\epsilon_{y},\epsilon_{z}\neq0$,
do not satisfy the condition, nor do terms involving $\epsilon_{w}=1$,
$\epsilon_{x},\epsilon_{y},\epsilon_{z}=0$. That the number of solutions
from I and from IIa+IIb amount to 11 and 6, respectively, while that
of solutions from IIa alone coincides with the inter-tie increment
5, testifies to an all-pervasive left strand affinity.)

\noindent The equal-grounds condition comes closer into focus now.
When we combine contributions from the first bunch (I) with those
from the second (IIa+IIb), we have to be careful: just as for the
non-characteristic multiple of the right strand, 84, we find $84<2n_{c}\cdot12$,
for additional non-characteristic multiples of all tie numbers $t$
of the Catalan double strand the rule $m_{\textrm{add}}<2n_{c}t$
may apply. One actually arrives at an equal number of multiples for
each tie number, as required:
\[
\begin{array}{rlrl}
11: & \textrm{T}=\{33,66\}\cup\{22,44\} & 12: & \textrm{T}=\{48,96\}\cup\{24,36\}\\
6: & \textrm{T}=\{6,12,18\}\cup\{18,24,36\!\!\!\!\!\diagup\} & 7: & \textrm{T}=\{7,14,21,28\}\cup\{\}
\end{array}
\]

\noindent As could be seen from Eqs. (\ref{eq:L_19-1}) and (\ref{eq:L_20-1}),
for lower tie numbers $t=6,7$ the multiples not available to gap
filling-in are $(n_{c}-1)t\,$: 12 and 21 respectively. Generalizing
from that, it appears that, given a suiting gap filling situation,
for the remaining tie numbers $t=11,12$ it's the multiples $n_{c}t$
that would be non-available $-$ 33 and 48 respectively. While this
is in accord with an underlying invariance regarding the number of
sign-preserving co-amplitudes, it seems to violate the equal-grounds
condition. One has to take note of a pecularity, though: The contributions
T for the lower tie numbers contain the `multiples' 6 and 7; for the
upper tie numbers, the contributions T do not contain improper multiples.
So the equal-grounds condition can be reforged formally by requiring
that for all $t=6,7,11,12$, multiple $(n_{c}-\delta_{t\in\textrm{T}})t$
is the one that characterizes an amplitude not partaking in gap filling-in,
where $\delta_{e\in\textrm{S}}=1\;\textrm{if}\enskip e\in\textrm{S},\;0\;\textrm{else}.$

\subsubsection*{In search of more dualities}

For each $t$ we obtain four multiples from the bases $(\textrm{Sim}_{16})$
and $(\textrm{Sim}_{64})$ and one characteristic quantity, a quintuple
from which one entry, the multiple identifiable as $(n_{c}-\delta_{t\in\textrm{T}})t$,
has to be removed to substantiate its non-partaking in gap filling-in.
However, as it seems firmly anchored in the simulacral world, we expect
to see it pop up in related situations: Simulacra may also be extracted
from $N_{\textrm{source}}(=(4,10,12,19,21))$ and $N_{\textrm{sink}}(=(3,7,18,29,43))$
via merger of a pair of basis elements. Going by the $M_{\nicefrac{5}{8}}^{+}$
lead, we may expect a merger to be allowed for $\delta_{x}=\delta_{y}^{'}=1$
$(x\neq y)$ and $\delta_{v,w,y,z},\delta_{v,w,x,z}^{'}=0$ only if
one of the expressions $(\delta_{1},\delta_{2},\delta_{3},\delta_{4},\delta_{5})\cdot N_{\textrm{type}}^{t}\pm(\delta_{1}^{'},\delta_{2}^{'},\delta_{3}^{'},\delta_{4}^{'},\delta_{5}^{'})\cdot N_{\textrm{type}}^{t}\in M_{\nicefrac{5}{8}}^{+}$.
Under this constraint, it turns out there are no solutions of type
`source', but two of type `sink': 
\[
\begin{array}{c}
(0,1,0,0,0)\cdot N_{\textrm{sink}}^{t}\mp(1,0,0,0,0)\cdot N_{\textrm{sink}}^{t}=\begin{array}{c}
4\\
10
\end{array},\\
(0,0,0,0,1)\cdot N_{\textrm{sink}}^{t}\pm(1,0,0,0,0)\cdot N_{\textrm{sink}}^{t}=\begin{array}{c}
46\\
40
\end{array}.
\end{array}
\]
It also turns out that, for each type, two order-4 tuples $-$ one
simulacral, one auxiliary $-$ emerge from this procedure, homogeneous
for type `sink', mixed for type `source,'
\[
\begin{array}{ll}
(\textrm{Sim}\hphantom{}_{\textrm{sink}}^{(2,1)})=(7+3,18,29,43), & (\textrm{Sim}\hphantom{}_{\textrm{source}}^{(5,1)})=(10,12,19,21\!+\!4),\\
\\
\textrm{(Aux}\hphantom{}_{\textrm{sink}}^{(2,1)})=(2,3,4,5), & \textrm{(Aux}\hphantom{}_{\textrm{sink}}^{(5,1)})\:=(20,21,22,23),
\end{array}
\]
where the superscripts in parentheses mark the places of elements
before the confluence, and the auxiliary tuple gives the amplitude
decomposition following the polite partition (staircase Young diagram)
of the sum of the respective outcomes,
\[
\begin{array}{cc}
4+10\:\quad= & 2+3+4+5,\\
46+40\!\quad= & 20+21+22+23.
\end{array}
\]

\noindent The four multiple-type amplitudes 12, 21, 33 and 48 not
involved in gap filling-in now have the following `projecting-out'
representation:
\[
\begin{array}{ll}
11: & 12:\\
\textrm{\ensuremath{}\hbox{\small projection\,\ simulacra:}} & \textrm{\textrm{\hbox{\small projection\,\ simulacra:}}}\\
(0,0,1,0)\cdot(\textrm{Sim}_{\textrm{sink}}^{(2,1)})^{t} & (0,0,0,1)\cdot(\textrm{Sim}_{\textrm{sink}}^{(2,1)})^{t}\\
+(0,0,1,0)\cdot(\textrm{Aux}_{\textrm{sink}}^{(2,1)})^{t}=33 & +(0,0,0,1)\cdot(\textrm{Aux}_{\textrm{sink}}^{(2,1)})^{t}=48\\
(0,1,0,0)\cdot(\textrm{Sim}_{\textrm{source}}^{(5,1)})^{t} & (0,0,0,1)\cdot(\textrm{Sim}_{\textrm{source}}^{(5,1)})^{t}\\
+(0,1,0,0)\cdot(\textrm{Aux}_{\textrm{sink}}^{(5,1)})^{t}=33 & +(0,0,0,1)\cdot(\textrm{Aux}_{\textrm{sink}}^{(5,1)})^{t}=48\\
\\
6: & 7:\\
\textrm{\hbox{\small projection\,\ simulacrum:}} & \textrm{\textrm{\hbox{\small projection\,\ simulacrum:}}}\\
(1,0,0,0)\cdot(\textrm{Sim}_{\textrm{sink}}^{(2,1)})^{t} & (0,1,0,0)\cdot(\textrm{Sim}_{\textrm{sink}}^{(2,1)})^{t}\\
+(1,0,0,0)\cdot(\textrm{Aux}_{\textrm{sink}}^{(2,1)})^{t}=12 & +(0,1,0,0)\cdot(\textrm{Aux}_{\textrm{sink}}^{(2,1)})^{t}=21\\
\\
\end{array}
\]
 Since there are no projection simulacra of type `source' in the lower
tie, we register six simulacral representations in all. We may expect
six simulacral representations to come out as well when, instead of
merger, a deletion scenario is considered: Simple deletion of the
$x$th element again yields quadruples, $(\textrm{Sim}_{\textrm{source}}^{[x]})$
and $(\textrm{Sim}_{\textrm{sink}}^{[x]})$. With these constructions
at hand, we now find that all of the four multiple-type amplitudes
12, 21, 33 and 48 are representable in $\textrm{Bound(}(\textrm{Sim}_{\textrm{source}}^{[2]}))$,
but only 21 and 33 in $\textrm{Bound(}(\textrm{Sim}_{\textrm{sink}}^{[x]}))\;(x=3,4,5)$,
the asymmetry being due to the original non-representability of 12
and 48 in $\textrm{Bound}(\textrm{\emph{N}}_{\textrm{sink}})$ $-$
which does persist after deletion of one of the basis elements. The
`indecorous' double strand, with multiples in $(\textrm{Sim}_{\textrm{source}}^{[2]})$
and, say, $(\textrm{Sim}_{\textrm{sink}}^{[5]})$ representation,
then reads 
\[
\begin{array}{ll}
11: & 12:\\
\textrm{\hbox{\small partial\,\ R-field\,\ simulacrum:}} & \textrm{\textrm{\hbox{\small R-spacetime\,\ simulacrum:}}}\\
(0,1,0,1)\cdot(\textrm{Sim}_{\textrm{source}}^{[2]})^{t}=33 & (-1,1,1,1)\cdot(\textrm{Sim}_{\textrm{source}}^{[2]})^{t}=48\\
\hbox{\small\textrm{partial spacetime\,\ simulacrum}:}\\
(-1,1,0,1)\cdot(\textrm{Sim}_{\textrm{sink}}^{[5]})^{t}=33\\
\\
6: & 7:\\
\textrm{\hbox{\small projection\,\ simulacrum:}} & \textrm{\textrm{\hbox{\small projection\,\ simulacrum:}}}\\
(0,1,0,0)\cdot(\textrm{Sim}_{\textrm{source}}^{[2]})^{t}=12 & (0,0,0,1)\cdot(\textrm{Sim}_{\textrm{source}}^{[2]})^{t}=21\\
 & \hbox{\small\textrm{partial field\,\ simulacrum}:}\\
 & (1,0,1,0)\cdot(\textrm{Sim}_{\textrm{sink}}^{[5]})^{t}=21
\end{array}
\]

\noindent Although the outcome is six simulacral representations
in either scenario, merger and deletion, the multiples represented
differ in detail, the reason being that two Young staircases are used
in the former and only one in the latter $-$ in the implicit form
$x=2,3,4,5$. Apart from that difference, the scenarios are complementary
to one another, in that `tied-to-type-sink' merger and `tied-to-type-source'
deletion reliably identify the quartet 12, 21, 33, and 48. 

\noindent One further complementarity applies lower tie-wise: in
the first decorated double strand of this subsection, `projection'
simulacra apply to characteristic quantities; in the above double
strand, they apply to inexpedient (`projecting-out') multiples. But
there is an unexpected and more profound side to this link between
the two double strand representations, which is why we called the
deletion-based double strand `indecorous' $-$ for the bottom end
of the staircase (at $x=5$), the simulacra of type `sink' are bound
up with the setting familiar from the first decorated double strand
of this subsection: the (partial) `field' simulacrum is right-strand
affine and the (partial) `spacetime' simulacrum left-strand affine.
In contrast, for the top end of the staircase (at $x=2$), laterally
inverted assignment in simulacra of type `source' manifests in the
upper tie: the `spacetime' simulacrum turns right-strand affine and
the (partial) `field' simulacrum left-strand affine (hence the marking
by a prefix R- for reflection). 

\noindent The overall picture and especially the latter unexpected
effect imply that Young staircases play an important part in the very
foundations of particle creation and crotonic implementation involving
gap filling-in.

\section{\label{sec:Supplementary-remarks-and}Closing remarks and outlook}

\vspace{0.2cm}
Table \ref{tab:Key-particle-creation-related-2} is of central ,
regarding information as to which dimensions are partaking in particle
creation and which not. On the face of things, its virtue appears
to be the strict separation of key dimensions of type `source' and
`sink.' On the downside, the numbers of unrepresentable dimensions
it offers seem to differ for these types for no good reason $-$ 11
vs. 41. This however points to an omission: only members of $M_{\textrm{reg}}$
and $M_{\nicefrac{5}{8}}$ have been included, though we have seen
that members of $M_{\nicefrac{9}{8}}$ by necessity enter the stage
as soon as it comes to determining intergenerational quark-mass ratios
or bounds for $\left(\textrm{\textgreek{a}}'\right)^{-1}$ via $\kappa'=-\frac{1}{8}+\sqrt{\Phi^{(31)}\Phi^{(63)}}$.
\footnote{We recall that each member of $M_{\nicefrac{9}{8}}:=\sigma_{n}\!=\!\frac{7}{2},8,17,\ldots$
is induced by $\frac{9}{8}M_{\textrm{reg}}$ such that $\frac{9}{8}p_{2}=\sigma_{1}-\frac{1}{8}$,
$\frac{9}{8}p_{3}=\sigma_{2}-\frac{1}{8},\ldots$, just like $M_{\nicefrac{5}{8}}:=o_{n}\!=\!4,9,19,\ldots$
is induced by $\frac{5}{8}M_{\textrm{reg}}$ such that $\frac{5}{8}p_{3}=o_{1}+\frac{3}{8}$,
$\frac{5}{8}p_{4}=\sigma_{2}+\frac{3}{8},\ldots$.%
} Members of $M_{\nicefrac{9}{8}}$ are not easily incorporated, though.
In fact, one must modify the rule $\left(L_{\nu}-\prod_{i=1}^{n}(\cdot)\right)/\prod_{i=1}^{n-2}(\cdot)=\nu'\:(\textrm{natural})$
to the effect that in the presence of the fraction $\frac{7}{2}$
the upper product bounds are brought to a juxtaposition $n$ vs. $n-1$
instead of $n$ vs. $n-2$. We then get an enhanced version of Table
\ref{tab:Key-particle-creation-related-2},

\bigskip{}

\hspace{0.8cm}%
\begin{tabular}{|c|c|c|c|c|}
\hline 
$n-2$ & $L_{\nu}$ & $\prod_{i=1}^{n}(\cdot)$ & $\prod_{i=1}^{n-2}(\cdot)$ & $\nu'$\tabularnewline
\hline 
\hline 
1 & $L_{4}=24$ & $1\cdot3\cdot7=21$ & $1$ & 3\tabularnewline
\hline 
2 & $L_{10}=336$ & $1\cdot3\cdot7\cdot15=315$ & $1\cdot3$ & 7\tabularnewline
\hline 
3 & $L_{19}=10668$ & $1\cdot3\cdot7\cdot15\cdot31=9765$ & $1\cdot3\cdot7$ & 43\tabularnewline
\hline 
\hline 
1 & $L_{12}=756$ & $4\cdot9\cdot19=684$ & 4 & 18\tabularnewline
\hline 
2 & $L_{21}=27720$ & $4\cdot9\cdot19\cdot39=26676$ & $4\cdot9$ & 29\tabularnewline
\hline 
\hline 
$n-1$ & $L_{\nu}$ & $\prod_{i=1}^{n}(\cdot)$ & $\prod_{i=1}^{n-1}(\cdot)$ & \tabularnewline
\hline 
\hline 
1 & $L_{12}=756$ & $\frac{7}{2}\cdot8\cdot17=476$ & $\frac{7}{2}\cdot8$ & 10\tabularnewline
\hline 
\end{tabular}

\bigskip{}

\noindent which now contains the as-yet-missing ingredient. For one
thing, we see that $\nu=12$ does not contribute anything new to the
basis $N_{\textrm{source}}=(4,10,12,19,21)$, so its associated Bound($N_{\textrm{source}}$)
merely remains a (re-)presenter of 55 out of 66 `source' dimensions.
However, the new two-sidedness of dimension $10$ leads to a radical
change. With $N'_{\textrm{sink}}=(3,7,10,18,29,43)$, the `sink' and
`source' sides are effectively symmetrized in that Bound($N'_{\textrm{sink}}$)
now has the capacity to represent 99 out of 110 `sink' dimensions,
which makes the aforementioned governance of the number 11 a both-sided
and in a way more perspicuous one. Obviously, the apparent loss in
symmetry in the number of basis elements is compensated for by an
increase in symmetry regarding the number of representable dimensions.
Another case in point is the number of basis elements of $J_{\rho}^{(p)}$,
which increases as $\mathbf{4},18,54,\ldots$, while the number of
basis elements of $G_{\rho}^{(p)}$ increases as $6,18,54,\ldots$
($p=15,31,63,\ldots$) (for details see \cite{Merkel}). Again, we
may note how an apparent lack of symmetry (in the number of basis
elements) leads to an increase in symmetry in another form: the first
two elements of the basis $(J_{\rho}^{(15)})=(-\mathbf{5},\mathbf{15},-43,149)$
\emph{uniformly} give rise to order-4 subgroups, namely of the group
of units of the quotient rings $\mathbb{Z}/16\mathbb{Z}$ and $\mathbb{Z}/64\mathbb{Z}\:$
respectively $-$ one formed by the powers of $\left|-5\right|$,
$\{1,5,9,13\}$, and the other by the powers of 15, $\{1,15,33,47\}$;
the remaining elements, which form distinct subgroups of higher order,
\emph{i.e.} powers of $\left|-43\right|$ an order-16 subgroup of
the group of units of the quotient ring $\mathbb{Z}/64\mathbb{Z},$
and powers of 149 an order-64 subgroup of the group of units of the
quotient ring $\mathbb{Z}/256\mathbb{Z}$, will occupy us shortly.
We'll be panning to a synoptic perspective here in order to see how
global and organizer aspects intertwine. As has been mentioned in
Sect.\,\ref{sec:Crotons-on-the}, on $\Gamma$ $-$ the 6-cube complex
ensuing from $(G_{\rho}^{(15)})=(3,5,11,17,41,113)$ $-$ 170 out
of 190 potentially attainable node labels are realized. On $\chi$
$-$ the 4-cube complex ensuing from $J_{\rho}^{(15)}=(-5,15,-43,149)$
$-$, by contrast, out of 212 potentially attainable ones only 40
are realized. The gap between these numbers, $170\!-\!40,$ was shown
(in Sect.\ref{sec:Application-to-subatomic}) to be responsible for
the Magnus equation's remarkable ability to account for the ratio
$F_{e}/F_{g}$. %
\footnote{Viewing the $172\,(=C_{5}+40)$ values that fail with $\chi^{(15)}$
as a partially-veiled-by-$C_{5}$, but otherwise symmetric counterpart
to the values that succeed, would be an alternative option (cf. $C_{5}$'s
role as discussed in subsection \ref{sub:Duality-controls:-the}).%
} In view of this, it comes as no surprise that the respectively realized
and non-realized node labels $\leq96$ on $\Gamma$ and $\chi$ correspond
to multiples (to-be-projected-out ones in square brackets below) and
characteristic numbers of the double strand:\vspace{-0.1cm}
\[
\begin{array}{ll}
11: & 12:\\
\hbox{\small\textrm{multiples}:} & \hbox{\small\textrm{multiples}:}\\
\Gamma\ni x\notin\chi\:(x=22,44,66), & \Gamma\ni x\notin\chi\:(x=24,36),\Gamma\ni96\in\chi,\\
\Gamma\ni[33]\in\chi & \chi\ni[48]\notin\Gamma\\
\hbox{\small\textrm{characteristic number}:} & \hbox{\small\textrm{characteristic number}:}\\
\Gamma\ni47\notin\chi & \Gamma\ni15\in\chi\\
\\
6: & 7:\\
\textrm{\hbox{\small multiples:}} & \textrm{\textrm{\hbox{\small multiples:}}}\\
\Gamma\ni x\notin\chi\:(x=6,18,24), & 7\notin\Gamma\cup\chi,\Gamma\ni14\notin\chi,\Gamma\ni28\in\chi,\\
\Gamma\ni[12]\notin\chi & \Gamma\ni[21]\notin\chi\\
\hbox{\small\textrm{characteristic number}:} & \hbox{\small\textrm{characteristic number}:}\\
\Gamma\ni13\notin\chi & \Gamma\ni5\in\chi
\end{array}
\]

\noindent Thirteen entries can be identified as node labels realized
on $\Gamma$ but not on $\chi$ and five as node labels realized on
both boundaries $-$ numbers that coincide with the lower-tie characteristic
quantities 13 and 5 and whose sum coincides with the number of basis
elements of $G_{\rho}^{(31)}$. The rest are one entry each: One,
not a label on either boundary, and thus coincident with the number
of basis elements of $G_{\rho}^{(7)}$. The other, identifiable as
a node label on $\chi$ but not on $\Gamma$. Of the latter sort there
are three more, 106, 164 and 172. Not being node labels $\leq96$,
they need not be counted. The two left, then, coincide with the number
of basis elements of $J_{\rho}^{(7)}$. If they were counted, their
number would coincide with that of basis elements of $J_{\rho}^{(15)}$.
The 13+5 coincidence brings the bases $(G_{\rho}^{(p)})$ into focus,
whose number of elements is given by the formula 
\[
2\cdot3^{\log_{2}(p+1)-3},\quad M_{\textrm{reg}}\ni p>7.
\]
As was demonstrated in \cite{Merkel}, it can equivalently be expressed
by an ansatz where two triangular matrices %
\footnote{\ $(G_{\xi+(p+1)/4,\zeta})$ and $(G_{\xi+(3p+3)/8,\zeta})$, $\xi=1,2,\ldots(q+1)/2$,
$\zeta=1,\ldots,(q+3)/2-\xi,$%
} with secondary symmetry, each with at most $\frac{(q+1)(q+3)}{8}$
distinct entries, are cleaned for redundant entries: 
\begin{equation}
2\left(\frac{(q+1)(q+3)}{8}-s\right),\quad q=(p-3)/4.\label{eq:subtracts-1}
\end{equation}
For $M_{\textrm{reg}}\ni p>15$, there is a nice twist about this:
As the subtrahends $s$ are keeping company with the units of $\{1,5,9,13\}$,%
\footnote{\ let $\mathcal{L}_{s}$ be the set of numbers $\left\lfloor \log_{2}(C_{q*}C_{2q*+1})\right\rfloor >s$,
where $q*\in M_{\textrm{reg}}$. Then the least element $l_{\textrm{min}}\in\mathcal{L}_{s}$
is found to satisfy $l_{\textrm{min}}-s\equiv u\,\textrm{(mod 16)}$
where $u\in\{1,5,9,13\}$. For $\frac{(q+1)(q+3)}{8}=10$, $s=1$,
we find $\left\lfloor \log_{2}(C_{1}C_{3})\right\rfloor -s\equiv1\,\textrm{(mod 16)};$
also for $\frac{(q+1)(q+3)}{8}=36$, $s=9$, $\left\lfloor \log_{2}(C_{3}C_{7})\right\rfloor -s\equiv1\,\textrm{(mod 16)};$
the next instances are $\frac{(q+1)(q+3)}{8}=136$, $s=55$ with $\left\lfloor \log_{2}(C_{15}C_{31})\right\rfloor -s\equiv5\,\textrm{(mod 16)}$,
$\frac{(q+1)(q+3)}{8}=528$, $s=285$ with $\left\lfloor \log_{2}(C_{63}C_{127})\right\rfloor -s\equiv9$
(mod 16), and so on, %
} the expressions 

\begin{equation}
2\left(\frac{(q+1)(q+3)}{8}-5^{\log_{2}\frac{q+1}{8}}\right)\label{eq:characters-1}
\end{equation}
keep spitting out sums of characteris\-tic quantities per tie $-$
the same that have been dealt with in Sect. \ref{sec:Pregeometric-categories-relevant}.
See the table below for their progression:\bigskip{}

\hfill{}%
\begin{tabular}{|c|}
\hline 
$p$%
\begin{tabular}{c}
\noalign{\vskip0.08cm}
\tabularnewline[0.08cm]
\end{tabular}\tabularnewline
\hline 
\hline 
\multirow{2}{*}{$2\left(\frac{(q+1)(q+3)}{8}-5^{\log_{2}\frac{q+1}{8}}\right)$}\tabularnewline
\tabularnewline
\hline 
\end{tabular}\,%
\begin{tabular}{|c|c|c|c|c|}
\hline 
\negthinspace{}\negthinspace{}\negthinspace{}\negthinspace{}\negthinspace{}\negthinspace{}%
\begin{tabular}{l}
\noalign{\vskip0.08cm}
\tabularnewline[0.08cm]
\end{tabular}\negthinspace{}\negthinspace{}31 & 63 & 127 & 255 & $\ldots$\tabularnewline
\hline 
\hline 
\multirow{2}{*}{\negthinspace{}\negthinspace{}\negthinspace{}\negthinspace{}\negthinspace{}\negthinspace{}%
\begin{tabular}{l}
\noalign{\vskip0.08cm}
\tabularnewline[0.08cm]
\end{tabular}\negthinspace{}\negthinspace{}18} & \multirow{2}{*}{62} & \multirow{2}{*}{222} & \multirow{2}{*}{806} & \multirow{2}{*}{$\ldots$}\tabularnewline
 &  &  &  & \tabularnewline
\hline 
\end{tabular}\hfill{}

\bigskip{}

\noindent So far, we had: $13+5=18$, $47+15=62$ $-$ if we wish
to proceed with the double strand, we know in advance the next tie
should come with the sum of characteris\-tic quantities $\xi+\zeta=222$.
The determination of $\xi$ and $\zeta$ is straightforward. We may
assume they are odd numbers, the left-strand affine lying in the interval
$]C_{6},C_{7}[$ and being the minimum among the remainders $(5\cdot2^{n}-1+2^{n})\:$mod$C_{7}$\enskip{}$(n\in\mathbb{N})$,
and the right-strand affine in the interval $]C_{4},C_{6}[$ and minimal
in $(5\cdot2^{n}-1-2^{n})\:$mod$C_{6}$. The left-strand affine's
interval can be narrowed down further in that 222 serves as an upper
bound: $]C_{6},222[$. Twenty out of sixty remainders are odd-numbered,
and out of the latter seven match the interval $]C_{6},C_{7}[$:
\[
\left(191,(137,149,161,143,167,)215\right).
\]
The parenthesising is put in to make it clear: the minimum remainder
search requires further constraining by 2b$'$. That is, the left-strand
and right-strand characteristic numbers must respectively be of the
form
\[
\begin{array}{c}
\xi=\xi_{0}+2^{m_{\xi}}\\
\zeta=\zeta_{0}-2^{m_{\zeta}}
\end{array}\quad(m_{\xi}>m_{\zeta};\xi_{0},\zeta_{0}\in M_{\nicefrac{5}{8}}).
\]
One possibility, $(79+128)+(79-64)$, fails because 79+128 is not
among the above remainders, so the unique solution under this constraint
is
\[
\begin{array}{c}
\xi=191=159+32,\\
\zeta=31=39-8.
\end{array}
\]
The solution would be none at all if we could not confirm that the
right-strand affine 31 is minimal among the remainders of $(5\cdot2^{n}-1-2^{n})\:$mod$C_{6}$.
Again, we use parenthesising to signal caution: the apparent minimum
$-$ 15 for the interval $]C_{5},C_{6}[$ in the remainders
\[
\left(\,(7,15,)31,63,127,123,115,99,65(,3)\,\right)
\]
$-$ is \emph{not} available in being already in use as right-strand
affine characteristic number, so 31 becomes the true minimum. 

\noindent Next, we have to determine the $t$'s for the new tie.
We have reason to believe that the observed inter-tie increment
\[
t'=t+5
\]
is more than mere heuristics: While the twenty values not realizable
on $\Gamma$ $-$ 7,34,48,51,62,65,79, 106,120,147,161,164,172,175,178,181,183,186,188,189
$-$ are not without formation law, the formation law for the 40 values
realizable on $\chi$ is manifestly one of increase by 5: $5,10,15,20$;
$23,28,\ldots,63$; $86,91,\ldots,126$; $129,134,\ldots,169$; $172,177,\ldots,212$.
So we deduce $16(=11+5)$ and $17(=12+5)$ as new tie numbers. \\
The envisioned continuation furthermore demands (i) specifying the
characteristic multiples, $m_{c}$, for the new tie numbers and (ii)
finding, for them as well as for the new characteristic numbers 191
and 31, equivalents in terms of the units of the order-16 subgroup
formed by powers of $\left|-43\right|$ of the group of units of the
quotient ring $\mathbb{Z}/64\mathbb{Z}$. 

\noindent As regards (i), we simply assume an increase 
\[
m_{c}'=(m+1)n_{c}(t+5),
\]
finding 
\[
\begin{array}{cc}
3n_{c}\cdot16=144 & \textrm{for}\: n_{c}=3,\\
3n_{c}\cdot17=204 & \textrm{for}\: n_{c}=4.
\end{array}
\]
\noindent Regarding (ii), we note that the sixteen units of the subgroup
in question, 
\[
\{\bar{1},43,57,\overline{19},\overline{49},59,41,35,33,\overline{11},25,\overline{51},\overline{17},27,\bar{9},\bar{3}\}\quad(\Sigma^{\,\textrm{sub}}=480),
\]
 can be partitioned into two sequences of order 8, 
\[
\begin{array}{c}
(\textrm{Sim}_{64}^{\textrm{unbar}})=(43,57,59,41,35,33,25,27)\quad(\Sigma^{\,\textrm{unbar}}=320),\\
\\
(\textrm{Sim}_{64}^{\textrm{bar}})=(\bar{1},\overline{19},\overline{49},\overline{11},\overline{51},\overline{17},\bar{9},\bar{3})\quad(\Sigma^{\,\textrm{bar}}=160),
\end{array}
\]
such that the sum of units each realizes $M_{\nicefrac{5}{8}}^{+}$
$-$ the `unbarred' units 320 and the `barred' ones 160, making up 480
in all. In order they form simulacral representations for the new
tie, the very notion of simulacrum asks for generalization: A simulacrum
is now termed `spacetime' or `field' if the number of unbarred units
of $(\textrm{Sim}_{64}^{\textrm{unbar}})$ used occupies the lower
half of the staircase and a `projection' simulacrum if it occupies
the upper half: \\
\vspace{-0.2cm}

\hfill{} %
\begin{tabular}{|c|c|ccccccc}
\cline{1-1} 
{\small `projection'} & \multicolumn{1}{c}{$\:1$} &  & $\overline{1}+1$ &  &  &  &  & \tabularnewline
\cline{1-2} 
\multicolumn{2}{|c|}{{\small "}} & $\:2$ &  & $\overline{2}+1$ &  &  &  & \tabularnewline
\cline{1-3} 
\multicolumn{3}{|c|}{{\small `field'}} & 3 &  &  & $\overline{3}$ &  & \tabularnewline
\cline{1-4} 
\multicolumn{4}{|c|}{{\small `spacetime'}} & 4 &  &  &  & $\overline{4}$\tabularnewline
\cline{1-4} 
\end{tabular}\hfill{}

\medskip{}
\noindent The number of barred units of $(\textrm{Sim}_{64}^{\textrm{bar}})$
used follows a similar pattern but reacts the switch of halves with
unit step lags. The lags are an integral part of the generalization
in that they correspond to the $x=0$ grades of the $2^{x}$-grading
of halves and (halves of) partitions,\\
\hfill{} %
\begin{tabular}{cc}
 & \tabularnewline
 & \tabularnewline
 & \tabularnewline
\end{tabular}%
\begin{tabular}{ccccccc}
 &  & $\qquad\qquad\quad1=(33+47)/l,$ & $(33+47)/r=1,$ &  &  & \tabularnewline
 & \multicolumn{2}{c}{$\qquad2=\Sigma^{\,\textrm{bar}}/l,$} & \multicolumn{2}{c}{$\quad\qquad\Sigma^{\,\textrm{bar}}/r=2,$} &  & \tabularnewline
\multicolumn{3}{c}{$4=\Sigma^{\,\textrm{unbar}}/l,$} & \multicolumn{3}{c}{$\qquad\qquad\Sigma^{\,\textrm{unbar}}/r=4,$} & \tabularnewline
\end{tabular}\hfill{}

\medskip{}

\noindent where $l=\bar{1}+\overline{19}+\overline{49}+\overline{11}$
and $r=\Sigma^{\,\textrm{bar}}-l$. Subtract the $x\!>\!0$ grades
whilst retaining the $x\!=\!0$ ones and you get the original numbers
of units of $(\textrm{Sim}_{64})$ \nopagebreak or $(\textrm{Sim}_{16})$
used:\\

\smallskip{}
\hfill{} %
\begin{tabular}{|c|c|cccc}
\cline{1-1} 
{\small `projection'} & \multicolumn{1}{c}{$1$} & $+\overline{1}\mathbf{+1}\mathbf{-2}=1$ &  &  & \tabularnewline
\cline{1-2} 
\multicolumn{2}{|c|}{{\small "}} & $2$ & $+\overline{2}\mathbf{+1}\mathbf{-4}=1$ &  & \tabularnewline
\cline{1-3} 
\multicolumn{3}{|c|}{{\small `field'}} & $3$ & $+\overline{3}\mathbf{-2}=4$ & \tabularnewline
\cline{1-4} 
\multicolumn{4}{|c|}{{\small `spacetime'}} & $4$ & $+\bar{4}\mathbf{-4}=4$\tabularnewline
\cline{1-4} 
\end{tabular}\hfill{}

\medskip{}

\noindent The enhanced, decorated double strand then reads\hfill{}

\[
\begin{array}{ll}
16: & 17:\\
\textrm{\hbox{\small characteristic\,\ multiple}} & \textrm{\hbox{\small characteristic\,\ multiple}}\\
\textrm{\hbox{\small(`spacetime'\ simulacrum):}} & \textrm{\hbox{\small(`field'\ simulacrum):}}\\
(1,1,-1,1,0,0,0,0)\cdot(\textrm{Sim}_{64}^{\textrm{unbar}})^{t} & (0,0,0,0,1,1,1,0)\cdot(\textrm{Sim}_{64}^{\textrm{unbar}})^{t}\\
+(0,0,0,0,1,1,-1,1)\cdot(\textrm{Sim}_{64}^{\textrm{bar}})^{t} & +(0,0,1,1,1,0,0,0)\cdot(\textrm{Sim}_{64}^{\textrm{bar}})^{t}\\
=144\,(=3n_{c}\cdot16) & =204\,(=3n_{c}\cdot17)\\
\textrm{\hbox{\small characteristic\,\ quantity\,159+32}} & \textrm{\hbox{\small characteristic\,\ quantity\,39-8}}\\
\textrm{\hbox{\small(`projection'\ simulacrum):}} & \textrm{\hbox{\small(`projection'\ simulacrum):}}\\
(0,1,1,0,0,0,0,0)\cdot(\textrm{Sim}_{64}^{\textrm{unbar}})^{t} & (0,0,0,0,0,0,0,1)\cdot(\textrm{Sim}_{64}^{\textrm{unbar}})^{t}\\
+(0,0,1,0,0,1,1,0)\cdot(\textrm{Sim}_{64}^{\textrm{bar}})^{t}=191 & +(1,0,0,0,0,0,0,1)\cdot(\textrm{Sim}_{64}^{\textrm{bar}})^{t}=31\\
11: & 12:\\
\textrm{\hbox{\small characteristic\,\ multiple}} & \textrm{\textrm{\hbox{\small characteristic\,\ multiple}}}\\
\hbox{\small(`\textrm{spacetime}'\ \textrm{simulacrum}):} & \hbox{\small(`\textrm{field}'\ \textrm{simulacrum}):}\\
(1,-1,1,1)\cdot(\textrm{Sim}_{64})^{t}=66(=2n_{c}\cdot11) & (1,1,1,1)\cdot(\textrm{Sim}_{64})^{t}=96(=2n_{c}\cdot12)\\
\textrm{\hbox{\small characteristic\,\ quantity\,39+8}} & \textrm{\textrm{\hbox{\small characteristic\,\ quantity\,19-4}}}\\
\textrm{\hbox{\small(`projection'\ simulacrum):}} & \textrm{\hbox{\small(`projection'\ simulacrum):}}\\
(0,0,0,1)\cdot(\textrm{Sim}_{64})^{t}=47 & (0,1,0,0)\cdot(\textrm{Sim}_{64})^{t}=15\\
6: & 7:\\
\textrm{\hbox{\small characteristic\,\ multiple}} & \textrm{\textrm{\hbox{\small characteristic\,\ multiple}}}\\
\hbox{\small(`\textrm{spacetime}'\ \textrm{simulacrum}):} & \hbox{\small(`\textrm{field}'\ \textrm{simulacrum}):}\\
(1,-1,1,1)\cdot(\textrm{Sim}_{16})^{t}=18(=n_{c}\cdot6) & (1,1,1,1)\cdot(\textrm{Sim}_{16})^{t}=28(=n_{c}\cdot7)\\
\textrm{\hbox{\small characteristic\,\ quantity\,9+4}} & \textrm{\textrm{\hbox{\small characteristic\,\ quantity\,9-4}}}\\
\textrm{\hbox{\small(`projection'\ simulacrum):}} & \textrm{\hbox{\small(`projection'\ simulacrum):}}\\
(0,0,0,1)\cdot(\textrm{Sim}_{16})^{t}=13 & (0,1,0,0)\cdot(\textrm{Sim}_{16})^{t}=5\\
\\
\qquad\qquad(n_{c}=3) & \qquad\qquad(n_{c}=4)
\end{array}
\]

\noindent The construction of one more tie on top of the above, with
hypothetical tie numbers 21 and 22, could proceed along similar lines.
It would involve the units of the order-64 subgroup formed by powers
of 149 of the group of units of the quotient ring $\mathbb{Z}/256\mathbb{Z\,}$.
However, kissing numbers with index > 100, whose pivots of organizer
origin an enhanced double strand might be a suitable study tool for,
are definitely out of reach presently. For kissing numbers with index
$\ll$ 100, the ur-double strand should suffice. 

\noindent  With such prospects, the interplay between pivots of organizer
and global origin and how their carrier Mersenne fluctuations assemble
in qphyla becomes a field of study worth aspiring to. Regarding `sink'
dimension 18, we came upon the specific, alloqphyletic condition under
which an organizer pivot snaps into global pivots (see Tables \ref{tab:Time-like-refinements mimicking left strand-1}
and \ref{tab:Space-like-refinement-offsets-1}). It would be interesting
to find out if this is true of `source' and `sink' dimensions alike. 

\newpage{}

\appendix

\section{\label{sec:Crotons-as-boundary}Crotons on the boundary }

\noindent Bases of order-31 croton base numbers pop up as a by-product
of the matrix constructions $\,\boldsymbol{f}^{(31)}=1\!\!\!\boldsymbol{1}^{\otimes4}\otimes\boldsymbol{b}^{(1)}+(G_{\mu\nu}^{(31)})\otimes c_{3}$
and $\boldsymbol{\, h}^{(31)}=1\!\!\!\boldsymbol{1}^{\otimes4}\otimes\boldsymbol{b}^{(1)}+(J_{\mu\nu}^{(31)})\otimes c_{2}\,$
(for more details of the construction, see \cite{Merkel}). Not all
of the matrix elements $G_{\mu\nu}^{(31)}$ and $J_{\mu\nu}^{(31)}$
need to be considered because the subquadrants UL(LL($\cdot$)) =
LL(UL($\cdot$)) = LL(LR($\cdot$)) just reproduce order-15 croton
base numbers. As shown in Fig. \ref{fig:catal-rep}, order-31 croton
base numbers can be extracted from the non-UR(LL($\cdot$)) parts
of quadrants $\textrm{LL}(G_{\mu\nu}^{(31)})$ and $\textrm{LL}(J_{\mu\nu}^{(31)})$:

\noindent 
\begin{figure}[H]
\caption{\label{fig:catal-rep} Order-31 croton base numbers extracted from
matrices of $\boldsymbol{f}^{(31)}$ and $\boldsymbol{h}^{(31)}$}
\[
\mathrm{LL}\,(G_{\mu\nu}^{(31)})=\qquad\qquad\qquad\qquad
\]
$\quad\left(\begin{array}{c}
\\
\\
\\
\\
\\
\\
\\
\\
\end{array}\right.$\hspace{-0.3cm}%
\begin{tabular}{cccc|cccc}
$\underline{429}$ & 155 & 43 & 19 & $\underline{5}$ & \multicolumn{1}{c|}{3} & $\underline{1}$ & 1\tabularnewline
1275 & $\underline{429}$ & 115 & 43 & 11 & \multicolumn{1}{c|}{$\underline{5}$} & 1 & $\underline{1}$\tabularnewline
\cline{7-8} 
4819 & 1595 & $\underline{429}$ & 155 & 41 & 17 & $\underline{5}$ & 3\tabularnewline
15067 & 4819 & 1275 & $\underline{429}$ & 113 & 41 & 11 & $\underline{5}$\tabularnewline
\cline{5-8} 
58781 & 18627 & 4905 & \multicolumn{1}{c}{1633} & $\underline{429}$ & 155 & 43 & 19\tabularnewline
189371 & 58781 & 15297 & \multicolumn{1}{c}{4905} & 1275 & $\underline{429}$ & 115 & 43\tabularnewline
737953 & 227089 & 58781 & \multicolumn{1}{c}{18627} & 4819 & 1595 & $\underline{429}$ & 155\tabularnewline
2430289 & 737953 & 189371 & \multicolumn{1}{c}{58781} & 15067 & 4819 & 1275 & $\underline{429}$\tabularnewline
\end{tabular}\hspace{-0.2cm}$\left.\begin{array}{c}
\\
\\
\\
\\
\\
\\
\\
\\
\end{array}\right)$\medskip{}
\[
\mathrm{LL}\,(J_{\mu\nu}^{(31)})=\qquad\qquad\qquad\qquad
\]
 $\left(\begin{array}{c}
\\
\\
\\
\\
\\
\\
\\
\\
\end{array}\right.$\hspace{-0.3cm}%
\begin{tabular}{cccc|cccc}
$-\underline{429}$ & 117 & -41 & 13 & $-\underline{5}$ & \multicolumn{1}{c|}{1} & $-\underline{1}$ & 1\tabularnewline
1547 & $-\underline{429}$ & 143 & -41 & 15 & \multicolumn{1}{c|}{$-\underline{5}$} & 3 & $-\underline{1}$\tabularnewline
\cline{7-8} 
-4903 & 1343 & $-\underline{429}$ & 117 & -43 & 15 & $-\underline{5}$ & 1\tabularnewline
18269 & -4903 & 1547 & $-429$ & 149 & -43 & 15 & -$\underline{5}$\tabularnewline
\cline{5-8} 
-58791 & 15547 & -4823 & \multicolumn{1}{c}{1319} & $-429$ & 117 & -41 & 13\tabularnewline
223573 & -58791 & 17989 & \multicolumn{1}{c}{-4823} & 1547 & $-429$ & 143 & -41\tabularnewline
-747765 & 194993 & -58791 & \multicolumn{1}{c}{15547} & -4903 & 1343 & $-\underline{429}$ & 117\tabularnewline
2886235 & -747765 & 223573 & \multicolumn{1}{c}{-58791} & 18269 & -4903 & 1547 & $-\underline{429}$\tabularnewline
\end{tabular}\hspace{-0.2cm}$\left.\begin{array}{c}
\\
\\
\\
\\
\\
\\
\\
\\
\end{array}\right)$
\end{figure}

\noindent Outcomes are the 18-tuples
\begin{eqnarray*}
(G_{\rho}^{(31)}) & = & (19,43,115,155,\underline{429},1275,1595,1633,4819,4905,\\
 &  & \:15067,15297,18627,58781,189371,227089,737953,2430289)
\end{eqnarray*}
and
\begin{eqnarray*}
(J_{\rho}^{(31)}) & = & (13,-41,117,143,-\underline{429},1319,1343,1547,-4823,-4903,\\
 &  & \:15547,17989,18269,-58791,194993,223573,-747765,2886235),
\end{eqnarray*}

\newpage
\noindent from which the outer nodes of 18-cube complexes with boundary
labels $\Gamma_{x}$ and $\chi_{x}$ are respectively formed. 

\vspace*{2mm}
\noindent Just for the sake of completeness, we show an excerpt of the successor
 $16\times16$ matrix

\noindent $\mathrm{LL}\,(G_{\mu\nu}^{(63)})=\begin{array}{cc}
\\
\\
\\
\end{array}$ 

\noindent %
\begin{minipage}[t]{1\columnwidth}%
$\left(\begin{array}{rrrrrrrrr}
\scriptstyle{\underline{9694845}} & \scriptstyle{2926323} & \scriptstyle{747851} & \scriptstyle{230395} & \scriptstyle{58791} & \scriptstyle{18633} & \scriptstyle{4907} & \scriptstyle{1635} & \ldots\\
\scriptstyle{32431347} & \scriptstyle{\underline{9694845}} & \scriptstyle{2461115} & \scriptstyle{747851} & \scriptstyle{189393} & \scriptstyle{58791} & \scriptstyle{15299} & \scriptstyle{4907} & \ldots\\
\scriptstyle{128896939} & \scriptstyle{38331} & \scriptstyle{\underline{9694845}} & \scriptstyle{2926323} & \scriptstyle{738035} & \scriptstyle{227123} & \scriptstyle{58791} & \scriptstyle{18633} & \ldots\\
\scriptstyle{436615771} & \scriptstyle{128896939} & \scriptstyle{32431347} & \scriptstyle{\underline{9694845}} & \scriptstyle{2430515} & \scriptstyle{738035} & \scriptstyle{189393} & \scriptstyle{58791} & \ldots\\
\scriptstyle{1767204399} & \scriptstyle{519693033} & \scriptstyle{130392731} & \scriptstyle{38792251} & \scriptstyle{\underline{9694845}} & \scriptstyle{2926323} & \scriptstyle{747851} & \scriptstyle{230395} & \ldots\\
 &  &  & &  &  \ldots &  &  &   
\end{array}\right)$%
\end{minipage}

\noindent 
\[
\]
 (5 rows, each truncated). The full matrix  yields  the associated  54-tuple,
\[
\begin{array}{c}
(G_{\rho}^{(63)})=\\
\\
\left(\begin{array}[b]{r}
\\
\\
\end{array}\right.
\end{array}\!\!\!\!\!\!\!\!\!\!\!\!\!\!\!\!\!\!\!\!\!\!\!\!\!\!\!\!\!\!\begin{array}[t]{rrr}
1,\!635\,, & 4,\!907\,, & 15,\!299\,,\\
18,\!633\,, & 58,\!791\,, & 189,\!393\,,\\
227,\!123\,, & 230,\!395\,, & 738,\!035\,,\\
747,\!851\,, & 2,\!430,\!515\,, & 2,\!461,\!115\,,\\
2,\!926,\!323\,, & \underline{9,\!694,\!845}\,, & 32,\!431,\!347\,,\\
38,\!331,\!419\,, & 38,\!792,\!251\,, & 38,\!795,\!521\,,\\
128,\!896,\!939\,, & 130,\!392,\!731\,, & 130,\!402,\!545\,,\\
436,\!615,\!771\,, & 441,\!538,\!235\,, & 441,\!568,\!833\,,\\
519,\!693,\!033\,, & 519,7\!30,\!299\,, & 1,\!767,\!204,\!399\,,\\
1,\!767,\!321,\!981\,, & 6,\!044,\!219,\!361\,, & 6,\!044,\!598,\!147\,,\\
7,\!090,\!735,\!179\,, & 7,\!091,\!189,\!425\,, & 7,\!168,\!783,\!827\,,\\
24,\!335,\!131,\!347\,, & 24,\!336,\!607,\!417\,, & 24,\!597,\!422,\!507\,,\\
83,\!908,\!823,\!403\,, & 83,\!913,\!684,\!433\,, & 84,\!796,\!853,\!171\,,\\
99,\!228,\!108,\!067\,, & 343,\!059,\!613,\!221\,, & 1,\!190,\!676,\!037,\!827\,,\\
1,\!390,\!379,\!604,\!203\,, & 1,\!404,\!717,\!639,\!489\,, & 4,\!837,\!348,\!974,\!083\,,\\
4,\!886,\!545,\!335,\!065\,, & 16,\!885,\!007,\!814,\!155\,, & 17,\!054,\!606,\!505,\!569\,,\\
19,\!881,\!172,\!597,\!035\,, & 69,\!531,\!783,\!535,\!237\,, & 243,\!860,\!214,\!616,\!867\,,\\
283,\!858,\!869,\!110,\!417\,, & 997,\!331,\!203,\!563,\!441\,, & 3,\!512,\!576,\!820,\!924,\!177
\end{array}\begin{array}[t]{r}
\\
\\
\\
\\
\\
\\
\\
\\
\\
\\
\\
\\
\\
\\
\\
\left.\!\!\!\!\!\!\begin{array}{r}
\\
\\
\\
\end{array}\right);
\end{array}
\]
and it may be noted that the tuple members can be grouped as follows:

\[
 \begin{array}{l}
\!\!\!\!\!\!\!\!\!\!\!\!\!\!\!\!\!\!(3^{\log_{2}(p+1)-3}+1)/2 \quad\textrm{tuple members}\;
\leq C_{q} \\
\\
\!\!\!\!\!\!\!\!\!\!\!\!\!\!\!\!\!\!(3^{\log_{2}(p+1)-2}-1)/2 \quad\textrm{tuple members}\; 
>  C_{q} \\
\\
\!\!\!\!\!\!\!\!\!\!\!(q=(p-3)/4,\quad M_{\textrm{reg}}\ni p>7). \\
\end{array}
\]
\medskip{}

\newpage

\noindent Now to the different encodings of the CF terms used in Part I; croton
amplitudes and phases in the volume, $\varphi_{\alpha}^{(n)}$, $\psi_{\alpha}^{(n)}\;(n\lesssim3030,\alpha\leq499)$,
corresponding to labels, are given below in the order the croton data occurred
in the text :

\vspace*{3mm}
\noindent Fig. \ref{fig:A-prototype-geometric} Mersenne
fluctuation, $\Gamma$-encoded:\medskip{}

\noindent $\varphi_{336}^{(206)}=13=(0,0,-1,1,-1,-1,-1,-1,0,1,0,0,0,0,0,0,0,0)\cdot(G^{(31)})^{t}$

\noindent $\psi_{384}^{(207)}=-14=(0,-1,1,0,0,0,0,0,1,-1,0,0,0,0,0,0,0,0)\cdot(G^{(31)})^{t}$

\noindent $\varphi_{326}^{(207)}=\psi_{363}^{(207)}=27=(-1,0,-1,0,1,0,1,-1,0,0,1,-1,0,0,0,0,0,0)\cdot(G^{(31)})^{t}$

\noindent $\varphi_{338}^{(208)}=56=(1,1,-1,0,1,1,-1,0,0,0,0,0,0,0,0,0,0,0)\cdot(G^{(31)})^{t}$

\noindent $\psi_{372}^{(208)}=-57=(0,1,-1,0,-1,-1,0,1,-1,1,0,0,0,0,0,0,0,0)\cdot(G^{(31)})^{t}$

\noindent $\varphi_{344}^{(209)}=\psi_{359}^{(209)}=114=(0,1,0,0,1,1,0,-1,0,0,0,0,0,0,0,0,0,0)\cdot(G^{(31)})^{t}$

\noindent $\varphi_{338}^{(210)}=228=(0,0,-1,0,1,0,0,0,1,-1,0,0,0,0,0,0,0,0)\cdot(G^{(31)})^{t}$

\noindent $\psi_{380}^{(210)}=-229=(0,-1,-1,0,-1,-1,0,1,0,0,0,0,0,0,0,0,0,0)\cdot(G^{(31)})^{t}$

\noindent $\varphi_{366}^{(211)}=458=(0,0,1,0,1,0,0,0,1,-1,0,0,0,0,0,0,0,0)\cdot(G^{(31)})^{t}$

\noindent $\psi_{356}^{(211)}=-459=(1,0,-1,1,1,0,0,0,1,0,0,1,0,-1,-1,1,0,0)\cdot(G^{(31)})^{t}$

\noindent $\varphi_{352}^{(212)}=\psi_{371}^{(212)}=918=(0,1,1,0,-1,1,0,0,1,-1,0,0,0,0,0,0,0,0)\cdot(G^{(31)})^{t}$

\noindent $\varphi_{362}^{(213)}=\psi_{375}^{(213)}=459=(-1,0,1,-1,-1,0,0,0,-1,0,0,-1,0,1,1,-1,0,0)\cdot(G^{(31)})^{t}$

\noindent $\varphi_{350}^{(214)}=\psi_{363}^{(214)}=229=(0,1,1,0,1,1,0,-1,0,0,0,0,0,0,0,0,0,0)\cdot(G^{(31)})^{t}$

\noindent $\varphi_{328}^{(215)}=\psi_{371}^{(215)}=114=(0,1,0,0,1,1,0,-1,0,0,0,0,0,0,0,0,0,0)\cdot(G^{(31)})^{t}$

\noindent $\varphi_{336}^{(216)}=56=(1,1,-1,0,1,1,-1,0,0,0,0,0,0,0,0,0,0,0)\cdot(G^{(31)})^{t}$

\noindent $\psi_{351}^{(216)}=57=(0,-1,1,0,1,1,0,-1,1,-1,0,0,0,0,0,0,0,0)\cdot(G^{(31)})^{t}$

\noindent $\varphi_{328}^{(217)}=\psi_{343}^{(217)}=28=(0,-1,0,0,1,1,0,-1,0,0,0,0,0,0,0,0,0,0)\cdot(G^{(31)})^{t}$

\noindent $\varphi_{324}^{(218)}=14=(0,1,-1,0,0,0,0,0,-1,1,0,0,0,0,0,0,0,0)\cdot(G^{(31)})^{t}$

\noindent $\psi_{391}^{(218)}=15=(-1,-1,1,0,0,0,1,-1,0,0,0,0,0,0,0,0,0,0)\cdot(G^{(31)})^{t}$\medskip{}

\noindent Fig. \ref{fig:A-prototype-geometric} Mersenne fluctuation,
$\chi$-encoded:\medskip{}

\noindent $\varphi_{336}^{(206)}=13=(1,0,0,0,0,0,0,0,0,0,0,0,0,0,0,0,0,0)\cdot(J^{(31)})^{t}$

\noindent $\psi_{384}^{(207)}=-14=(1,0,1,0,0,-1,0,-1,0,0,-1,0,1,0,0,0,0,0)\cdot(J^{(31)})^{t}$

\noindent $\varphi_{326}^{(207)}=\varphi_{363}^{(207)}=27=(0,0,-1,0,0,1,0,1,0,0,1,0,-1,0,0,0,0,0)\cdot(J^{(31)})^{t}$

\noindent $\varphi_{338}^{(208)}=56=(0,0,1,1,0,0,1,-1,0,0,0,0,0,0,0,0,0,0)\cdot(J^{(31)})^{t}$

\noindent $\psi_{372}^{(208)}=-57=(-1,1,0,-1,0,-1,-1,0,1,-1,-1,0,1,0,0,0,0,0)\cdot(J^{(31)})^{t}$

\noindent $\varphi_{344}^{(209)}=\varphi_{359}^{(209)}=114=(1,1,1,-1,0,0,1,1,0,0,1,0,-1,0,0,0,0,0)\cdot(J^{(31)})^{t}$

\noindent $\varphi_{338}^{(210)}=228=(1,0,1,1,0,0,-1,0,0,-1,-1,0,-1,-1,1,-1,0,0)\cdot(J^{(31)})^{t}$

\noindent $\psi_{380}^{(210)}=-229=(1,1,0,0,1,-1,0,1,0,0,0,0,0,0,0,0,0,0)\cdot(J^{(31)})^{t}$

\noindent $\varphi_{366}^{(211)}=458=(1,0,1,1,1,-1,-1,-1,-1,0,0,0,0,0,0,0,0,0)\cdot(J^{(31)})^{t}$

\noindent $\psi_{356}^{(211)}=-459=(1,1,1,-1,1,-1,1,0,0,0,0,0,0,0,0,0,0,0)\cdot(J^{(31)})^{t}$

\noindent $\varphi_{352}^{(212)}=\psi_{371}^{(212)}=918=(0,1,0,0,0,1,0,0,-1,1,0,1,-1,0,0,0,0,0)\cdot(J^{(31)})^{t}$

\noindent $\varphi_{362}^{(213)}=\psi_{375}^{(213)}=459=(-1,-1,-1,1,-1,1,-1,0,0,0,0,0,0,0,0,0,0,0)\cdot(J^{(31)})^{t}$

\noindent $\varphi_{350}^{(214)}=\psi_{363}^{(214)}=229=(-1,-1,0,0,-1,1,0,-1,0,0,0,0,0,0,0,0,0,0)\cdot(J^{(31)})^{t}$

\noindent $\varphi_{328}^{(215)}=\psi_{371}^{(215)}=114=(1,1,1,-1,0,0,1,1,0,0,1,0,-1,0,0,0,0,0)\cdot(J^{(31)})^{t}$

\noindent $\varphi_{336}^{(216)}=56=(0,0,1,1,0,0,1,-1,0,0,0,0,0,0,0,0,0,0)\cdot(J^{(31)})^{t}$

\noindent $\psi_{351}^{(216)}=57=(1,-1,0,1,0,1,1,0,-1,1,1,0,-1,0,0,0,0,0)\cdot(J^{(31)})^{t}$

\noindent $\varphi_{328}^{(217)}=\psi_{343}^{(217)}=28=(-1,1,0,-1,-1,0,1,-1,0,0,0,0,0,0,0,0,0,0)\cdot(J^{(31)})^{t}$

\noindent $\varphi_{324}^{(218)}=14=(-1,0,-1,0,0,1,0,1,0,0,1,0,-1,0,0,0,0,0)\cdot(J^{(31)})^{t}$

\noindent $\psi_{391}^{(218)}=15=(0,1,1,1,0,0,1,-1,0,0,0,0,0,0,0,0,0,0)\cdot(J^{(31)})^{t}$

\noindent \medskip{}

\noindent Fig. \ref{fig:Pivotal-amplitude-closing} pivot, $\Gamma$-encoded:

\noindent \medskip{}

\noindent $\varphi_{448}^{(1556)}=\psi_{441}^{(1556)}=1304=(0,0,-1,0,0,1,0,1,-1,0,0,-1,1,0,0,0,0,0)\cdot(G^{(31)})^{t}$

\noindent $\varphi_{442}^{(1557)}=2609=(-1,-1,0,0,-1,1,1,0,0,0,-1,1,0,0,0,0,0,0)\cdot(G^{(31)})^{t}$

\noindent $\psi_{413}^{(1557)}=2610=(0,-1,1,1,0,0,0,0,1,0,0,0,1,-1,-1,1,0,0)\cdot(G^{(31)})^{t}$

\noindent $\varphi_{401}^{(1558)}=5219=(0,0,1,0,1,0,0,0,0,1,1,-1,0,0,0,0,0,0)\cdot(G^{(31)})^{t}$

\noindent $\psi_{430}^{(1558)}=-5220=(-1,0,0,0,0,1,1,1,-1,1,1,1,1,-1,0,0,0,0)\cdot(G^{(31)})^{t}$

\noindent $\varphi_{407}^{(1559)}=2609=(-1,-1,0,0,-1,1,1,0,0,0,-1,1,0,0,0,0,0,0)\cdot(G^{(31)})^{t}$

\noindent $\psi_{430}^{(1559)}=-2610=(0,1,-1,-1,0,0,0,0,-1,0,0,0,-1,1,1,-1,0,0)\cdot(G^{(31)})^{t}$

\noindent \medskip{}

\noindent Fig. \ref{fig:Pivotal-amplitude-closing} pivot, $\chi$-encoded:

\noindent \medskip{}

\noindent $\varphi_{448}^{(1556)}=\psi_{441}^{(1556)}=1304=(0,-1,0,0,0,0,1,0,-1,1,0,0,0,0,0,0,0,0)\cdot(J^{(31)})^{t}$

\noindent $\varphi_{442}^{(1557)}=2609=(-1,0,1,1,-1,0,-1,-1,-1,0,0,0,0,0,0,0,0,0)\cdot(J^{(31)})^{t}$

\noindent $\psi_{413}^{(1557)}=2610=(0,0,0,0,0,0,1,1,0,0,0,1,-1,0,0,0,0,0)\cdot(J^{(31)})^{t}$

\noindent $\varphi_{401}^{(1558)}=5219=(-1,0,0,0,1,0,-1,0,-1,-1,1,0,-1,0,0,0,0,0)\cdot(J^{(31)})^{t}$

\noindent $\psi_{430}^{(1558)}=-5220=(1,0,-1,-1,0,0,-1,-1,0,-1,1,1,1,1,0,0,0,0)\cdot(J^{(31)})^{t}$

\noindent $\varphi_{407}^{(1559)}=2609=(-1,0,1,1,-1,0,-1,-1,-1,0,0,0,0,0,0,0,0,0)\cdot(J^{(31)})^{t}$

\noindent $\psi_{430}^{(1559)}=-2610=(0,0,0,0,0,0,-1,-1,0,0,0,-1,1,0,0,0,0,0)\cdot(J^{(31)})^{t}$\medskip{}

\noindent Fig. \ref{fig:Pivotal-amplitude-closing} residue\,2, $\Gamma$-encoded\medskip{}

\noindent $\varphi_{404}^{(1556)}=102=(-1,0,-1,-1,1,0,1,-1,0,0,0,0,0,0,0,0,0,0)\cdot(G^{(31)})^{t}$

\noindent $\psi_{427}^{(1556)}=103=(0,-1,1,1,0,0,1,-1,1,-1,0,0,0,0,0,0,0,0)\cdot(G^{(31)})^{t}$

\noindent $\varphi_{380}^{(1557)}=\psi_{403}^{(1557)}=205=(1,0,1,0,1,1,0,-1,0,0,0,0,0,0,0,0,0,0)\cdot(G^{(31)})^{t}$

\noindent $\varphi_{390}^{(1558)}=\psi_{423}^{(1558)}=411=(-1,-1,1,0,0,-1,0,1,0,0,0,0,0,0,0,0,0,0)\cdot(G^{(31)})^{t}$

\noindent $\varphi_{398}^{(1559)}=\psi_{419}^{(1559)}=205=(1,0,1,0,1,1,0,-1,0,0,0,0,0,0,0,0,0,0)\cdot(G^{(31)})^{t}$\medskip{}

\noindent Fig. \ref{fig:Pivotal-amplitude-closing} residue\,2, $\chi$-encoded:

\noindent \medskip{}

\noindent $\varphi_{404}^{(1556)}=102=(0,1,0,1,0,0,0,0,0,0,0,0,0,0,0,0,0,0)\cdot(J^{(31)})^{t}$

\noindent $\psi_{427}^{(1556)}=103=(1,1,0,-1,1,1,0,1,-1,0,1,1,1,1,0,0,0,0)\cdot(J^{(31)})^{t}$

\noindent $\varphi_{380}^{(1557)}=\psi_{403}^{(1557)}=205=(0,0,0,0,1,0,0,-1,0,-1,1,0,-1,0,0,0,0,0)\cdot(J^{(31)})^{t}$

\noindent $\varphi_{390}^{(1558)}=\psi_{423}^{(1558)}=411=(0,0,0,-1,0,0,0,-1,-1,0,1,0,-1,0,0,0,0,0)\cdot(J^{(31)})^{t}$

\noindent $\varphi_{398}^{(1559)}=\psi_{419}^{(1559)}=205=(0,0,0,0,1,0,0,-1,0,-1,1,0,-1,0,0,0,0,0)\cdot(J^{(31)})^{t}$\medskip{}

\noindent Fig. \ref{fig:Pivotal-amplitude-closing} residue\,1, $\Gamma$-encoded:

\noindent \medskip{}

\noindent $\varphi_{406}^{(1556)}=100=(1,1,0,0,0,0,-1,1,0,0,0,0,0,0,0,0,0,0)\cdot(G^{(31)})^{t}$

\noindent $\psi_{428}^{(1556)}=-101=(0,1,0,0,0,0,0,-1,1,0,0,1,-1,0,0,0,0,0)\cdot(G^{(31)})^{t}$

\noindent $\varphi_{382}^{(1557)}=49=(0,1,1,0,-1,-1,1,0,0,0,0,0,0,0,0,0,0,0)\cdot(G^{(31)})^{t}$

\noindent $\psi_{404}^{(1557)}=-50=(1,0,0,-1,0,0,0,0,-1,1,0,0,0,0,0,0,0,0)\cdot(G^{(31)})^{t}$

\noindent $\varphi_{392}^{(1558)}=24=(-1,1,0,0,0,0,0,0,0,0,0,0,0,0,0,0,0,0)\cdot(G^{(31)})^{t}$

\noindent $\psi_{424}^{(1558)}=-25=(-1,0,0,1,-1,0,-1,1,0,0,-1,1,0,0,0,0,0,0)\cdot(G^{(31)})^{t}$

\noindent $\varphi_{400}^{(1559)}=\psi_{423}^{(1559)}=12=(0,1,0,0,1,1,1,0,0,0,0,1,-1,0,0,0,0,0)\cdot(G^{(31)})^{t}$

\medskip{}
\noindent Fig. \ref{fig:Pivotal-amplitude-closing} residue\,1,
$\chi$-encoded:

\noindent \medskip{}

\noindent $\varphi_{406}^{(1556)}=100=(1,0,-1,0,0,0,-1,1,0,0,0,0,0,0,0,0,0,0)\cdot(J^{(31)})^{t}$

\noindent $\psi_{428}^{(1556)}=-101=(0,0,0,0,0,-1,0,0,-1,0,-1,0,-1,-1,1,-1,0,0)\cdot(J^{(31)})^{t}$

\noindent $\varphi_{382}^{(1557)}=49=(1,-1,0,1,0,1,0,-1,1,-1,0,0,0,0,0,0,0,0)\cdot(J^{(31)})^{t}$

\noindent $\psi_{404}^{(1557)}=-50=(0,0,1,-1,0,1,-1,0,0,0,0,0,0,0,0,0,0,0)\cdot(J^{(31)})^{t}$

\noindent $\varphi_{392}^{(1558)}=24=(-1,0,1,0,0,0,0,0,-1,1,0,0,0,0,0,0,0,0)\cdot(J^{(31)})^{t}$

\noindent $\psi_{424}^{(1558)}=-25=(1,0,-1,0,1,-1,0,1,0,0,0,-1,1,0,0,0,0,0)\cdot(J^{(31)})^{t}$

\noindent $\varphi_{400}^{(1559)}=\psi_{423}^{(1559)}=12=(-1,-1,0,1,-1,1,0,-1,-1,1,0,1,-1,0,0,0,0,0)\cdot(J^{(31)})^{t}$\medskip{}

\noindent Fig. \ref{fig:Pivotal-amplitude-closing-1}, Table \ref{tab:Pivotal-and-residual-fig4}
pivot, $\Gamma$-encoded:

\noindent \medskip{}

\noindent $\varphi_{448}^{(1003)}=51919=(-1,0,0,0,-1,0,-1,0,-1,0,0,0,0,1,0,0,0,0)\cdot(G^{(31)})^{t}$

\noindent $\psi_{468}^{(1003)}=-51920=(1,-1,0,1,0,0,1,0,0,1,-1,1,0,-1,0,0,0,0)\cdot(G^{(31)})^{t}$

\noindent $\varphi_{442}^{(1004)}=\psi_{464}^{(1004)}=103839=(1,1,-1,1,-1,0,0,-1,-1,-1,-1,0,0,-1,1,0,0,0)\cdot(G^{(31)})^{t}$

\noindent $\varphi_{438}^{(1005)}=\psi_{443}^{(1005)}=207679=(1,0,-1,0,1,1,0,1,0,0,1,0,0,0,1,0,0,0)\cdot(G^{(31)})^{t}$

\noindent $\varphi_{449}^{(1006)}=\psi_{447}^{(1006)}=103839=(1,1,-1,1,-1,0,0,-1,-1,-1,-1,0,0,-1,1,0,0,0)\cdot(G^{(31)})^{t}$

\noindent \medskip{}

\noindent \noindent Fig. \ref{fig:Pivotal-amplitude-closing-1},
Table \ref{tab:Pivotal-and-residual-fig4} pivot, $\chi$-encoded:

\noindent \medskip{}

\noindent $\varphi_{448}^{(1003)}=51919=(-1,1,0,-1,1,0,-1,0,0,1,0,0,0,-1,0,0,0,0)\cdot(J^{(31)})^{t}$

\noindent $\psi_{468}^{(1003)}=-51920=(1,1,-1,-1,-1,0,0,1,0,-1,0,-1,1,1,0,0,0,0)\cdot(J^{(31)})^{t}$

\noindent $\varphi_{442}^{(1004)}=\psi_{464}^{(1004)}=103839=(0,0,-1,1,0,0,0,-1,0,0,0,1,0,-1,-1,1,0,0)\cdot(J^{(31)})^{t}$

\noindent $\varphi_{438}^{(1005)}=\psi_{443}^{(1005)}=207679=(1,0,0,0,-1,0,-1,-1,-1,0,0,0,-1,0,0,1,0,0)\cdot(J^{(31)})^{t}$

\noindent $\varphi_{449}^{(1006)}=\psi_{447}^{(1006)}=103839=(0,0,-1,1,0,0,0,-1,0,0,0,1,0,-1,-1,1,0,0)\cdot(J^{(31)})^{t}$

\noindent \medskip{}

\noindent \noindent Fig. \ref{fig:Pivotal-amplitude-closing-1},
Table \ref{tab:Pivotal-and-residual-fig4} residue\,2, $\Gamma$-encoded:\medskip{}

\noindent $\varphi_{448}^{(1003)}=193=(0,0,0,1,0,0,-1,1,0,0,0,0,0,0,0,0,0,0)\cdot(G^{(31)})^{t}$

\noindent $\psi_{469}^{(1003)}=194=(1,0,0,0,1,1,-1,0,-1,-1,-1,-1,-1,1,0,0,0,0)\cdot(G^{(31)})^{t}$

\noindent $\varphi_{443}^{(1004)}=\psi_{465}^{(1004)}=96=(-1,0,1,0,0,0,0,0,0,0,0,0,0,0,0,0,0,0)\cdot(G^{(31)})^{t}$

\noindent $\varphi_{443}^{(1005)}=\psi_{465}^{(1005)}=48=(0,0,0,0,0,0,1,-1,-1,1,0,0,0,0,0,0,0,0)\cdot(G^{(31)})^{t}$

\noindent $\varphi_{427}^{(1006)}=23=(0,-1,0,0,0,0,0,0,-1,-1,-1,-1,-1,1,0,0,0,0)\cdot(G^{(31)})^{t}$

\noindent $\psi_{449}^{(1006)}=24=(-1,1,0,0,0,0,0,0,0,0,0,0,0,0,0,0,0,0)\cdot(G^{(31)})^{t}$

\noindent \medskip{}

\noindent \noindent Fig. \ref{fig:Pivotal-amplitude-closing-1},
Table \ref{tab:Pivotal-and-residual-fig4} residue\,2, $\chi$-encoded:\medskip{}

\noindent $\varphi_{448}^{(1003)}=193=(-1,-1,0,-1,0,-1,0,1,1,-1,0,0,0,0,0,0,0,0)\cdot(J^{(31)})^{t}$

\noindent $\psi_{469}^{(1003)}=194=(1,0,1,1,-1,1,0,-1,0,0,0,1,-1,0,0,0,0,0)\cdot(J^{(31)})^{t}$

\noindent $\varphi_{443}^{(1004)}=\psi_{465}^{(1004)}=96=(-1,-1,-1,0,1,-1,-1,-1,-1,0,0,0,0,0,0,0,0,0)\cdot(J^{(31)})^{t}$

\noindent $\varphi_{443}^{(1005)}=\psi_{465}^{(1005)}=48=(0,1,-1,-1,-1,0,0,0,-1,1,0,0,0,0,0,0,0,0)\cdot(J^{(31)})^{t}$

\noindent $\varphi_{427}^{(1006)}=23=(1,1,-1,-1,1,0,-1,0,0,1,-1,-1,-1,-1,0,0,0,0)\cdot(J^{(31)})^{t}$

\noindent $\psi_{449}^{(1006)}=24=(0,0,0,0,0,-1,1,0,0,0,0,0,0,0,0,0,0,0)\cdot(J^{(31)})^{t}$

\noindent \medskip{}

\noindent \noindent Fig. \ref{fig:Pivotal-amplitude-closing-1} residue\,1,
$\Gamma$-encoded:

\noindent \medskip{}

\noindent $\varphi_{479}^{(1003)}=58=(-1,0,1,0,0,0,1,-1,0,0,0,0,0,0,0,0,0,0)\cdot(G^{(31)})^{t}$

\noindent $\varphi_{477}^{(1004)}=117=(0,0,0,1,0,0,1,-1,0,0,0,0,0,0,0,0,0,0)\cdot(G^{(31)})^{t}$

\noindent $\varphi_{465}^{(1005)}=\psi_{477}^{(1005)}=58=(-1,0,1,0,0,0,1,-1,0,0,0,0,0,0,0,0,0,0)\cdot(G^{(31)})^{t}$

\noindent $\varphi_{449}^{(1006)}=29=(0,0,1,0,0,0,0,0,1,-1,0,0,0,0,0,0,0,0)\cdot(G^{(31)})^{t}$\\
$\psi_{473}^{(1006)}=30=(0,1,-1,0,0,0,-1,-1,0,0,0,-1,1,0,0,0,0,0)\cdot(G^{(31)})^{t}$

\noindent \medskip{}

\noindent \noindent Fig. \ref{fig:Pivotal-amplitude-closing-1} residue\,1,
$\chi$-encoded:

\noindent \medskip{}

\noindent $\varphi_{479}^{(1003)}=58=(-1,-1,1,-1,0,1,-1,0,1,-1,0,0,0,0,0,0,0,0)\cdot(J^{(31)})^{t}$

\noindent $\varphi_{477}^{(1004)}=117=(0,0,1,0,0,0,0,0,0,0,0,00,0,0,0,0,0)\cdot(J^{(31)})^{t}$

\noindent $\varphi_{465}^{(1005)}=\psi_{477}^{(1005)}=58=(-1,-1,1,-1,0,1,-1,0,1,-1,0,0,0,0,0,0,0,0)\cdot(J^{(31)})^{t}$

\noindent $\varphi_{449}^{(1006)}=29=(1,-1,1,1,1,1,0,1,0,0,1,0,-1,0,0,0,0,0)\cdot(J^{(31)})^{t}$\\
$\psi_{473}^{(1006)}=30=(0,0,1,-1,0,1,-1,0,1,-1,0,0,0,0,0,0,0,0)\cdot(J^{(31)})^{t}$

\noindent \medskip{}

\noindent \noindent Fig. \ref{fig:Pivotal-edge-amplitudes}, Table
\ref{tab:Pivotal-and-residual-fig5} leggy pivot, $\Gamma$-encoded:\medskip{}

\noindent $\varphi_{407}^{(987)}=\psi_{437}^{(987)}=12977=(0,0,-1,0,0,0,-1,1,1,1,0,-1,1,0,0,0,0,0)\cdot(G^{(31)})^{t}$

\noindent $\varphi_{411}^{(988)}=25955=(0,0,-1,-1,-1,0,0,1,1,1,0,1,0,0,0,0,0,0)\cdot(G^{(31)})^{t}$

\noindent $\varphi_{418}^{(988)}=25956=(0,0,0,-1,0,-1,0,-1,0,-1,0,1,1,0,0,0,0,0)\cdot(G^{(31)})^{t}$

\noindent $\varphi_{397}^{(989)}=\psi_{449}^{(989)}=51911=(-1,-1,-1,-1,0,0,0,-1,0,-1,0,0,0,1,0,0,0,0)\cdot(G^{(31)})^{t}$$\varphi_{433}^{(990)}=103823=(-1,-1,-1,1,-1,1,1,0,0,1,0,0,0,1,-1,1,0,0)\cdot(G^{(31)})^{t}$

\noindent $\psi_{441}^{(990)}=103824=(0,-1,-1,1,-1,-1,1,-1,-1,-1,0,-1,0,-1,1,0,0,0)\cdot(G^{(31)})^{t}$

\noindent $\varphi_{417}^{(991)}=207646=(-1,1,0,1,-1,0,1,1,0,0,0,1,0,0,1,0,0,0)\cdot(G^{(31)})^{t}$

\noindent $\psi_{457}^{(991)}=207647=(1,0,0,-1,-1,-1,0,0,1,0,0,1,0,0,1,0,0,0)\cdot(G^{(31)})^{t}$

\noindent $\varphi_{421}^{(992)}=\psi_{469}^{(992)}=415294=(0,0,0,0,1,0,-1,0,0,0,0,0,0,0,1,1,0,0)\cdot(G^{(31)})^{t}$

\noindent $\psi_{423}^{(993)}=\psi_{473}^{(993)}=207647=(1,0,0,-1,-1,-1,0,0,1,0,0,1,0,0,1,0,0,0)\cdot(G^{(31)})^{t}$

\noindent $\varphi_{443}^{(994)}=\psi_{475}^{(994)}=103823=(-1,-1,-1,1,-1,1,1,0,0,1,0,0,0,1,-1,1,0,0)\cdot(G^{(31)})^{t}$

\noindent $\varphi_{425}^{(995)}=51911=(-1,-1,-1,-1,0,0,0,-1,0,-1,0,0,0,1,0,0,0,0)\cdot(G^{(31)})^{t}$

\noindent $\psi_{453}^{(995)}=51912=(1,-1,-1,0,0,0,-1,0,0,-1,1,-1,0,1,0,0,0,0)\cdot(G^{(31)})^{t}$

\noindent $\varphi_{431}^{(996)}=\psi_{477}^{(996)}=25955=(0,0,-1,-1,-1,0,0,1,1,1,0,1,0,0,0,0,0,0)\cdot(G^{(31)})^{t}$

\noindent $\varphi_{439}^{(997)}=12977=(0,0,-1,0,0,0,-1,1,1,1,0,-1,1,0,0,0,0,0)\cdot(G^{(31)})^{t}$

\noindent $\psi_{469}^{(997)}=12978=(1,-1,-1,1,1,-1,0,0,-1,0,0,0,1,0,0,0,0,0)\cdot(G^{(31)})^{t}$

\noindent \medskip{}

\noindent \noindent Fig. \ref{fig:Pivotal-edge-amplitudes}, Table
\ref{tab:Pivotal-and-residual-fig5} leggy pivot, $\chi$-encoded:\medskip{}

\noindent $\varphi_{407}^{(987)}=\psi_{437}^{(987)}=12977=(1,0,-1,0,-1,0,-1,1,-1,-1,-1,0,1,0,0,0,0,0)\cdot(J^{(31)})^{t}$

\noindent $\varphi_{411}^{(988)}=25955=(-1,-1,-1,-1,1,0,1,0,-1,-1,1,0,0,0,0,0,0,0)\cdot(J^{(31)})^{t}$

\noindent $\varphi_{418}^{(988)}=25956=(1,1,0,-1,0,-1,0,-1,1,0,1,0,1,0,0,0,0,0)\cdot(J^{(31)})^{t}$

\noindent $\varphi_{397}^{(989)}=\psi_{449}^{(989)}=51911=(0,-1,0,1,0,1,1,0,1,1,0,0,0,-1,0,0,0,0)\cdot(J^{(31)})^{t}$

\noindent $\varphi_{433}^{(990)}=103823=(1,-1,1,0,0,1,0,0,-1,-1,1,0,1,-1,0,0,0,0)\cdot(J^{(31)})^{t}$

\noindent $\psi_{441}^{(990)}=103824=(-1,0,0,0,0,-1,1,-1,0,0,0,1,0,-1,-1,1,0,0)\cdot(J^{(31)})^{t}$

\noindent $\varphi_{417}^{(991)}=207646=(1,-1,1,1,-1,1,0,0,0,0,0,-1,0,0,0,1,0,0)\cdot(J^{(31)})^{t}$

\noindent $\psi_{457}^{(991)}=207647=(0,0,-1,1,0,-1,0,-1,0,-1,0,-1,0,0,0,1,0,0)\cdot(J^{(31)})^{t}$

\noindent $\varphi_{421}^{(992)}=\psi_{469}^{(992)}=415294=(-1,-1,1,1,0,0,1,0,0,1,0,0,0,0,1,1,0,0)\cdot(J^{(31)})^{t}$

\noindent $\psi_{423}^{(993)}=\psi_{473}^{(993)}=207647=(0,0,-1,1,0,-1,0,-1,0,-1,0,-1,0,0,0,1,0,0)\cdot(J^{(31)})^{t}$

\noindent $\varphi_{443}^{(994)}=\psi_{475}^{(994)}=103823=(1,-1,1,0,0,1,0,0,-1,-1,1,0,1,-1,0,0,0,0)\cdot(J^{(31)})^{t}$

\noindent $\varphi_{425}^{(995)}=51911=(0,-1,0,1,0,1,1,0,1,1,0,0,0,-1,0,0,0,0)\cdot(J^{(31)})^{t}$

\noindent $\psi_{453}^{(995)}=51912=(0,0,0,0,1,0,0,-1,0,1,0,0,0,-1,0,0,0,0)\cdot(J^{(31)})^{t}$

\noindent $\varphi_{431}^{(996)}=\psi_{477}^{(996)}=25955=(-1,-1,-1,-1,1,0,1,0,-1,-1,1,0,0,0,0,0,0,0)\cdot(J^{(31)})^{t}$

\noindent $\varphi_{439}^{(997)}=12977=(1,0,-1,0,-1,0,-1,1,-1,-1,-1,0,1,0,0,0,0,0)\cdot(J^{(31)})^{t}$

\noindent $\psi_{469}^{(997)}=12978=(-1,0,1,-1,1,0,0,0,1,0,0,0,1,0,0,0,0,0)\cdot(J^{(31)})^{t}$

\noindent \medskip{}

\noindent \noindent Fig. \ref{fig:Pivotal-edge-amplitudes}, Table
\ref{tab:Pivotal-and-residual-fig5} residue\,2, $\Gamma$-encoded:\medskip{}

\noindent $\varphi_{397}^{(987)}=9434=(0,-1,1,0,1,-1,1,0,1,0,0,-1,-1,0,-1,1,0,0)\cdot(G^{(31)})^{t}$

\noindent $\psi_{426}^{(987)}=-9435=(-1,1,0,0,1,-1,0,0,-1,0,0,1,1,0,1,-1,0,0)\cdot(G^{(31)})^{t}$

\noindent $\varphi_{399}^{(988)}=4716=(0,0,1,1,0,-1,0,-1,0,0,-1,-1,0,0,-1,1,0,0)\cdot(G^{(31)})^{t}$

\noindent $\psi_{410}^{(988)}=-4717=(0,1,1,0,-1,1,0,1,0,0,1,1,0,0,1,-1,0,0)\cdot(G^{(31)})^{t}$

\noindent $\varphi_{385}^{(989)}=2357=(0,1,0,-1,0,0,0,0,0,1,0,0,1,-1,-1,1,0,0)\cdot(G^{(31)})^{t}$

\noindent $\psi_{438}^{(989)}=-2358=(1,1,1,0,0,0,0,0,1,0,1,1,0,0,1,-1,0,0)\cdot(G^{(31)})^{t}$

\noindent $\varphi_{421}^{(990)}=1178=(0,-1,0,0,1,0,1,1,0,0,0,0,1,-1,-1,1,0,0)\cdot(G^{(31)})^{t}$

\noindent $\psi_{430}^{(990)}=-1179=(0,0,0,0,-1,0,0,1,-1,0,0,0,-1,1,1,-1,0,0)\cdot(G^{(31)})^{t}$

\noindent $\varphi_{409}^{(991)}=\psi_{449}^{(991)}=589=(0,1,0,0,0,-1,-1,0,-1,1,0,-1,1,0,0,0,0,0)\cdot(G^{(31)})^{t}$

\noindent $\varphi_{413}^{(992)}=\psi_{461}^{(992)}=294=(-1,1,-1,1,0,0,0,0,0,0,-1,1,0,0,0,0,0,0)\cdot(G^{(31)})^{t}$

\noindent $\varphi_{405}^{(993)}=\psi_{461}^{(993)}=147=(0,-1,1,-1,0,0,0,0,0,0,-1,1,0,0,0,0,0,0)\cdot(G^{(31)})^{t}$

\noindent $\varphi_{435}^{(994)}=73=(1,-1,0,0,1,0,1,1,0,0,1,0,-1,0,0,0,0,0)\cdot(G^{(31)})^{t}$

\noindent $\psi_{446}^{(994)}=-74=(0,1,0,-1,0,0,-1,1,0,0,0,0,0,0,0,0,0,0)\cdot(G^{(31)})^{t}$

\noindent $\varphi_{413}^{(995)}=\psi_{445}^{(995)}=36=(-1,-1,0,0,0,1,0,0,1,0,1,0,0,-1,-1,1,0,0)\cdot(G^{(31)})^{t}$

\noindent $\varphi_{417}^{(996)}=\psi_{461}^{(996)}=17=(1,0,1,-1,0,0,-1,1,0,0,0,0,0,0,0,0,0,0)\cdot(G^{(31)})^{t}$

\noindent $\varphi_{427}^{(997)}=\psi_{457}^{(997)}=8=(0,1,1,1,-1,0,-1,1,-1,1,0,0,0,0,0,0,0,0)\cdot(G^{(31)})^{t}$

\noindent \medskip{}

\noindent \noindent Fig. \ref{fig:Pivotal-edge-amplitudes}, Table
\ref{tab:Pivotal-and-residual-fig5} residue\,2, $\chi$-encoded:\medskip{}

\noindent $\varphi_{397}^{(987)}=9434=(1,0,0,0,-1,-1,0,0,0,0,0,0,-1,0,-1,1,0,0)\cdot(J^{(31)})^{t}$

\noindent $\psi_{426}^{(987)}=-9435=(-1,0,0,0,-1,0,-1,0,0,1,-1,0,-1,-1,1,-1,0,0)\cdot(J^{(31)})^{t}$

\noindent $\varphi_{399}^{(988)}=4716=(-1,-1,-1,-1,0,0,1,0,0,0,1,0,1,1,-1,1,0,0)\cdot(J^{(31)})^{t}$

\noindent $\psi_{410}^{(988)}=-4717=(0,0,0,-1,0,0,0,-1,-1,-1,1,-1,1,0,1,-1,0,0)\cdot(J^{(31)})^{t}$

\noindent $\varphi_{385}^{(989)}=2357=(0,0,0,0,-1,0,1,0,1,1,0,0,-1,0,-1,1,0,0)\cdot(J^{(31)})^{t}$

\noindent $\psi_{438}^{(989)}=-2358=(-1,0,0,-1,1,0,0,-1,1,1,1,-1,-1,-1,1,-1,0,0)\cdot(J^{(31)})^{t}$

\noindent $\varphi_{421}^{(990)}=1178=(0,0,-1,-1,0,-1,1,1,-1,0,-1,-1,0,0,-1,1,0,0)\cdot(J^{(31)})^{t}$

\noindent $\psi_{430}^{(990)}=-1179=(-1,1,-1,-1,0,0,0,0,-1,-1,0,1,0,0,1,-1,0,0)\cdot(J^{(31)})^{t}$

\noindent $\varphi_{409}^{(991)}=\psi_{449}^{(991)}=589=(0,0,1,1,0,0,0,1,1,0,1,0,1,1,-1,1,0,0)\cdot(J^{(31)})^{t}$

\noindent $\varphi_{413}^{(992)}=\psi_{461}^{(992)}=294=(0,-1,0,0,-1,-1,1,0,1,-1,0,1,-1,0,0,0,0,0)\cdot(J^{(31)})^{t}$

\noindent $\varphi_{405}^{(993)}=\psi_{461}^{(993)}=147=(-1,-1,0,1,0,1,-1,0,0,0,0,0,0,0,0,0,0,0)\cdot(J^{(31)})^{t}$

\noindent $\varphi_{435}^{(994)}=73=(1,-1,0,1,0,0,1,-1,1,-1,0,0,0,0,0,0,0,0)\cdot(J^{(31)})^{t}$

\noindent $\psi_{446}^{(994)}=-74=(1,1,1,0,-1,-1,1,1,-1,0,1,1,1,1,0,0,0,0)\cdot(J^{(31)})^{t}$

\noindent $\varphi_{413}^{(995)}=\psi_{445}^{(995)}=36=(0,-1,1,1,-1,1,1,1,0,1,0,0,0,0,0,0,0,0)\cdot(J^{(31)})^{t}$

\noindent $\varphi_{417}^{(996)}=\psi_{461}^{(996)}=17=(1,-1,-1,0,0,0,0,0,1,-1,0,0,0,0,0,0,0,0)\cdot(J^{(31)})^{t}$

\noindent $\varphi_{427}^{(997)}=\psi_{457}^{(997)}=8=(-1,-1,1,1,0,0,0,0,0,0,0,1,-1,0,0,0,0,0)\cdot(J^{(31)})^{t}$

\noindent \medskip{}

\noindent \noindent Fig. \ref{fig:Pivotal-edge-amplitudes} residue\,1,
$\Gamma$-encoded:\medskip{}

\noindent $\varphi_{485}^{(990)}=\psi_{497}^{(990)}=136=(-1,0,0,1,0,0,0,0,0,0,0,0,0,0,0,0,0,0)\cdot(G^{(31)})^{t}$

\noindent $\varphi_{465}^{(991)}=68=(0,-1,-1,1,1,1,0,-1,0,0,0,0,0,0,0,0,0,0)\cdot(G^{(31)})^{t}$

\noindent $\varphi_{465}^{(992)}=33=(-1,1,-1,0,0,0,-1,1,-1,1,0,0,0,0,0,0,0,0)\cdot(G^{(31)})^{t}$

\noindent $\varphi_{471}^{(993)}=16=(1,-1,-1,1,0,0,0,0,0,0,0,0,0,0,0,0,0,0)\cdot(G^{(31)})^{t}$

\noindent \medskip{}

\noindent \noindent Fig. \ref{fig:Pivotal-edge-amplitudes} residue\,1,
$\chi$-encoded:\medskip{}

\noindent $\varphi_{485}^{(990)}=\psi_{497}^{(990)}=136=(-1,-1,-1,0,-1,0,1,-1,0,0,0,0,0,0,0,0,0,0)\cdot(J^{(31)})^{t}$

\noindent $\varphi_{465}^{(991)}=68=(-1,1,0,-1,1,-1,-1,-1,0,-1,0,0,0,0,0,0,0,0)\cdot(J^{(31)})^{t}$

\noindent $\varphi_{465}^{(992)}=33=(1,1,0,-1,0,0,-1,1,0,0,0,0,0,0,0,0,0,0)\cdot(J^{(31)})^{t}$

\noindent $\varphi_{471}^{(993)}=16=(0,1,0,-1,0,0,0,0,-1,1,0,-1,1,0,0,0,0,0)\cdot(J^{(31)})^{t}$

\noindent  \medskip{}

\noindent \noindent Table \ref{tab:Pivotal-and-residual-fig4-1}
pivot, $\Gamma$-encoded:\medskip{}

\noindent $\varphi_{72}^{(609)}=208430=(0,1,0,1,0,0,0,0,0,0,1,-1,-1,0,0,1,0,0)\cdot(G^{(31)})^{t}$

\noindent $\psi_{72}^{(609)}=-208431=(-1,0,-1,-1,0,-1,1,0,0,0,0,0,1,0,0,-1,0,0)\cdot(G^{(31)})^{t}$

\noindent \medskip{}

\noindent \noindent Table \ref{tab:Pivotal-and-residual-fig4-1}
pivot, $\chi$-encoded:\medskip{}

\noindent $\varphi_{72}^{(609)}=208430=(0,0,1,1,0,1,0,1,0,0,0,0,-1,0,0,1,0,0)\cdot(J^{(31)})^{t}$

\noindent $\psi_{72}^{(609)}=-208431=(0,0,0,0,1,-1,1,0,0,0,1,0,0,0,0,-1,0,0)\cdot(J^{(31)})^{t}$

\noindent \medskip{}

\noindent \noindent Table \ref{tab:Pivotal-and-residual-fig4-1}
residue\,2, $\Gamma$-encoded:\medskip{}

\noindent $\varphi_{78}^{(609)}=\psi_{83}^{(609)}=66=(0,0,0,0,0,0,0,0,-1,-1,-1,-1,-1,1,0,0,0,0)\cdot(G^{(31)})^{t}$

\noindent \medskip{}

\noindent \noindent Table \ref{tab:Pivotal-and-residual-fig4-1}
residue\,2, $\chi$-encoded:\medskip{}

\noindent $\varphi_{78}^{(609)}=\psi_{83}^{(609)}=66=(0,1,0,0,1,0,0,-1,0,1,-1,-1,-1,-1,0,0,0,0)\cdot(J^{(31)})^{t}$

\noindent \medskip{}

\noindent \noindent Table \ref{tab:Pivotal-and-residual-fig4-1}
residue\,1, $\Gamma$-encoded:\medskip{}

\noindent $\varphi_{58}^{(609)}=10=(-1,0,1,0,0,0,0,0,1,-1,0,0,0,0,0,0,0,0)\cdot(G^{(31)})^{t}$

\noindent $\psi_{52}^{(609)}=-11=(0,0,0,-1,0,0,0,1,-1,0,0,-1,1,0,0,0,0,0)\cdot(G^{(31)})^{t}$

\noindent \medskip{}

\noindent \noindent Table \ref{tab:Pivotal-and-residual-fig4-1}
residue\,1, $\chi$-encoded:\medskip{}

\noindent $\varphi_{58}^{(609)}=10=(0,1,-1,0,0,0,1,1,0,0,1,0,-1,0,0,0,0,0)\cdot(J^{(31)})^{t}$

\noindent $\psi_{52}^{(609)}=-11=(1,-1,-1,0,0,1,0,-1,0,0,0,-1,1,0,0,0,0,0)\cdot(J^{(31)})^{t}$\medskip{}

\noindent \noindent Table \ref{tab:Pivotal-and-residual-fig4-1-1}
pivot, $\Gamma$-encoded:\smallskip{}

\noindent $\varphi_{239}^{(1000)}=\psi_{269}^{(1000)}=758=(0,0,-1,0,1,-1,0,1,-1,1,0,0,0,0,0,0,0,0)\cdot(G^{(31)})^{t}$

\noindent \medskip{}

\noindent \noindent Table \ref{tab:Pivotal-and-residual-fig4-1-1}
pivot, $\chi$-encoded:\smallskip{}

\noindent $\varphi_{239}^{(1000)}=\psi_{269}^{(1000)}=758=(1,0,-1,0,1,-1,1,1,0,0,0,1,-1,0,0,0,0,0)\cdot(J^{(31)})^{t}$

\noindent \medskip{}

\noindent \noindent Table \ref{tab:Pivotal-and-residual-fig4-1-1}
residue\,1, $\Gamma$-encoded:\smallskip{}

\noindent $\varphi_{356}^{(1000)}=335=(1,0,0,0,0,0,0,0,-1,1,-1,1,0,0,0,0,0,0)\cdot(G^{(31)})^{t}$

\noindent $\psi_{388}^{(1000)}=-336=(-1,-1,0,1,-1,0,0,0,0,0,0,0,0,0,0,0,0,0)\cdot(G^{(31)})^{t}$

\noindent \medskip{}

\noindent \noindent Table \ref{tab:Pivotal-and-residual-fig4-1-1}
residue\,1, $\chi$-encoded:\smallskip{}

\noindent $\varphi_{356}^{(1000)}=335=(1,-1,0,0,-1,1,0,-1,1,-1,0,0,0,0,0,0,0,0),()\cdot(J^{(31)})^{t}$

\noindent $\psi_{388}^{(1000)}=-336=(0,0,1,0,1,1,-1,0,0,0,0,0,0,0,0,0,0,0)\cdot(J^{(31)})^{t}$

\noindent \medskip{}

\noindent \noindent Table \ref{tab:Pivotal-and-residual-fig4-1-1}
residue\,2, $\Gamma$-encoded:\smallskip{}

\noindent $\varphi_{135}^{(1000)}=5=(0,1,0,0,0,0,1,-1,0,0,0,0,0,0,0,0,0,0)\cdot(G^{(31)})^{t}$

\noindent $\psi_{168}^{(1000)}=-6=(0,0,0,1,-1,0,-1,1,0,0,-1,1,0,0,0,0,0,0)\cdot(G^{(31)})^{t}$

\noindent \medskip{}

\noindent \noindent Table \ref{tab:Pivotal-and-residual-fig4-1-1}
residue\,2, $\chi$-encoded:\smallskip{}

\noindent $\varphi_{135}^{(1000)}=5=(-1,1,1,1,1,-1,0,1,0,0,0,0,0,0,0,0,0,0)\cdot(J^{(31)})^{t}$

\noindent $\psi_{168}^{(1000)}=-6=(0,0,0,1,1,0,0,0,0,0,0,-1,1,0,0,0,0,0)\cdot(J^{(31)})^{t}$

\noindent \medskip{}

\noindent \noindent Both $G^{(31)}$ and $J^{(31)}$ lead to various
singularity assignments. For instance,

\noindent \smallskip{}

\noindent $0=(-1,0,1,0,1,-1,0,-1,1,0,0,0,1,-1,-1,1,0,0)\cdot(G^{(31)})^{t}$

\noindent $0=(0,0,0,0,1,0,1,1,0,1,-1,1,0,0,0,0,0,0)\cdot(J^{(31)})^{t}$

\newpage{}

\section{\label{sec:Crotons-in-the}Crotons in the volume}

\begin{table}[H]
\caption{\label{tab:Specific-fractions-in-1} $L_{m}\!=\! b_{\alpha}^{(n)}$
incidences in (CFR) $\left(\frac{2^{n}}{\sqrt{2}}\right)^{-1}\rightarrow\left[b_{0}^{(n)};b_{\alpha}^{(n)}\right]$
($n\leq3324$, $\alpha\leq499$; $\qquad\qquad\qquad$ match$=$\ \Checkmark{};
closest pivot$=$($\cdot$); largest $b_{\alpha}^{(n)}>L_{31}=$($\cdot$))}
\vspace{0.5cm}
\qquad{}\qquad{}\qquad{}%
\begin{tabular}{|c|c|c|}
\hline 
 & \multicolumn{2}{c|}{$\begin{array}{ccccc}
\\
\\
\\
\end{array}$Type $\sqrt{2}/2^{n}$}\tabularnewline
$L_{m}$ $^{a}$ & $b_{\alpha}^{(n)}$ & $\begin{array}{c}
\textrm{in-}\\
\textrm{cidence}
\end{array}$\tabularnewline
\hline 
\hline 
$\begin{array}{c}
2,6,12,24,\\
40,72,126\\
(m=1,2,{\scriptstyle \ldots,},7)
\end{array}$ & \Checkmark{} & $\begin{array}{c}
\textrm{very}\\
\textrm{high}
\end{array}$ \tabularnewline
\hline 
$240$ & \Checkmark{} & 19\tabularnewline
\hline 
$272$ & \Checkmark{} & 25\tabularnewline
\hline 
$336$ & \Checkmark{} & 9\tabularnewline
\hline 
$438$ & \Checkmark{} & 9\tabularnewline
\hline 
$756$ & \Checkmark{} & 5\tabularnewline
\hline 
$918$ & \Checkmark{} & 1\tabularnewline
\hline 
$1422$ & (1421) & (1)\tabularnewline
\hline 
$2340$ & (2338) & (1)\tabularnewline
\hline 
$4320$ & (4314) & (1)\tabularnewline
\hline 
$5346$ & (5366) & (1)\tabularnewline
\hline 
$7398$ & (7394) & (1)\tabularnewline
\hline 
$10\,668$ & (10\,596) & (1)\tabularnewline
\hline 
$17\,400$ & (17\,502) & (1)\tabularnewline
\hline 
$27\,720$ & (27\,901) & (1)\tabularnewline
\hline 
$49\,896$ & (49\,780) & (1)\tabularnewline
\hline 
$93\,150$ & (94\,869) & (1)\tabularnewline
\hline 
$\vdots$ &  & \tabularnewline
\hline 
$>L_{31}$ & (2\,445\,930) & (1)\tabularnewline
\hline 
\multicolumn{1}{|c}{} & \multicolumn{1}{c}{} & \tabularnewline
\multicolumn{3}{|c|}{$^{a}$ \scriptsize {\tt http://www.math.rwth-aachen.de/Gabriele.Nebe/LATTICES/kiss.html}}\tabularnewline
\hline 
\end{tabular}
\end{table}

\newpage{}

\begin{table}[H]
\caption{\label{tab:Specific-fractions-in} $L_{m}\!=\!\varphi_{\alpha}^{(n)}$
incidences in (CFR) type-I/II/III irrationals $\!\rightarrow\!\left[\varphi_{0}^{(n)}\!;\varphi_{\alpha}^{(n)}\right]$$ $
(Eqs.(\ref{eq:full-mers-1-1})$-$(\ref{eq:trunc-mers-1}); $n\lesssim3330$,
$2\leq s\leq9$, $\alpha\leq499$; \ match$=$\ \Checkmark{}; closest
pivot$=$($\cdot$); largest $b_{\alpha}^{(n)}>L_{31}=$($\cdot$))}
\vspace{1cm}
\begin{tabular}{|c|c|c|c|c|c|c|}
\hline 
 & \multicolumn{2}{c|}{Type I} & \multicolumn{2}{c|}{Type II} & \multicolumn{2}{c|}{Type III}\tabularnewline
$L_{m}$ & $\varphi_{\alpha}^{(n)}$ & $\begin{array}{c}
\textrm{inci-}\\
\!\!\!\!\textrm{dence}\!\!\!\!
\end{array}$ & $\varphi_{\alpha}^{(n)}$ & $\begin{array}{c}
\textrm{inci-}\\
\!\!\!\!\textrm{dence}\!\!\!\!
\end{array}$ & $\varphi_{\alpha}^{(n)}$ & $\begin{array}{c}
\textrm{inci-}\\
\!\!\!\!\textrm{dence}\!\!\!\!
\end{array}$\tabularnewline
\hline 
\hline 
$\begin{array}{c}
\!\!\!\!2,6,\!\!\!\!\\
\!\!\!\!{\scriptstyle \ldots,},126\!\!\!\!
\end{array}$ & \Checkmark{} & $\begin{array}{c}
\textrm{very}\\
\textrm{high}
\end{array}$  & \Checkmark{} & $\begin{array}{c}
\textrm{very}\\
\textrm{high}
\end{array}$  & \Checkmark{} & $\begin{array}{c}
\textrm{very}\\
\textrm{high}
\end{array}$ \tabularnewline
\hline 
$240$ & \Checkmark{} & 131 & \Checkmark{} & 181 & \Checkmark{} & 220\tabularnewline
\hline 
$272$ & \Checkmark{} & 90 & \Checkmark{} & 141 & \Checkmark{} & 163\tabularnewline
\hline 
$336$ & \Checkmark{} & 60 & \Checkmark{} & 93 & \Checkmark{} & 107\tabularnewline
\hline 
$438$ & \Checkmark{} & 47 & \Checkmark{} & 62 & \Checkmark{} & 64\tabularnewline
\hline 
$756$ & \Checkmark{} & 25 & \Checkmark{} & 20 & \Checkmark{} & 25\tabularnewline
\hline 
$918$ & \Checkmark{} & 9 & \Checkmark{} & 21 & \Checkmark{} & 20\tabularnewline
\hline 
$1422$ & \Checkmark{} & 5 & \Checkmark{} & 6 & \Checkmark{} & 4\tabularnewline
\hline 
$2340$ & (2341) & (1) & \Checkmark{} & 1 & \Checkmark{} & 4\tabularnewline
\hline 
$4320$ & (4321) & (1) & \Checkmark{} & 2 & \Checkmark{} & 3\tabularnewline
\hline 
$5346$ & (5344) & (1) & $(5346\pm2)$ & (2) & (5349) & (1)\tabularnewline
\hline 
$7398$ & \Checkmark{} & 1 & (7399) & (2) & $(7398\pm4)$ & (2)\tabularnewline
\hline 
$10668$ & (10\,674) & (1) & (10\,677) & (1) & \Checkmark{} & 1\tabularnewline
\hline 
$17400$ & (17\,409) & (1) & (17\,390) & (1) & (17\,398) & (1)\tabularnewline
\hline 
$27720$ & (27\,738) & (1) & (27\,733) & (1) & (27\,717) & (1)\tabularnewline
\hline 
$49896$ & (49\,679) & (1) & (50\,216) & (1) & (49\,888) & (1)\tabularnewline
\hline 
$93150$ & (92\,646) & (1) & (93\,489) & (1) & (92\,677) & (1)\tabularnewline
\hline 
$\vdots$ &  &  &  &  &  & \tabularnewline
\hline 
$\!\!\!\!207930\!\!\!\!$ & $-$ & $-$ & (207\,679) & (1) & (208\,430) & (1)\tabularnewline
\hline 
$\vdots$ &  &  &  &  &  & \tabularnewline
\hline 
\multirow{2}{*}{$>L_{31}$} & \multirow{2}{*}{$\begin{array}{c}
\!\!\!\!(12\,986\,152)\!\!\!\!\\
=\!\!\!\!\!^{?}\: L_{46}
\end{array}$} & \multirow{2}{*}{(1)} & \multirow{2}{*}{(3\,614\,855)} & \multirow{2}{*}{(1)} & \multirow{2}{*}{(9\,996\,953)} & \multirow{2}{*}{(1)}\tabularnewline
 &  &  &  &  &  & \tabularnewline
\hline 
\end{tabular}
\end{table}

\begin{table}[H]
\caption{\label{tab:CFR-1023}$(\textrm{CFR})\,\left(2^{1023}\log^{2}(2)\log(3)(2^{1023}\log(3)+1)\right)^{-1}\rightarrow\left[\gamma_{0}^{(1023)};\gamma_{\alpha>0}^{(1023)}\right]$}
\medskip{}

{\small  [0;468501467667419229549312377689539844924407660551543874466280005437695581078754471597542961\ 85106148285609224469531623604511871122301327405223890730886112738503306679155403295591866535\ 80930916137792940636148035377851094973289474608114270626188781020096469603973358560240607627\ 52582833446470959519779089630872605772157859088929998686243591882665756030997626912687132482\ 63928840099513485061261996790621980683053499875704951378748164035468315190534758283131147910\ 35077622406873484818364876666081251795168616648329698682949881814965910343984552495077355252\ 630711063949967796135960920246222931797345930397856390079245263791, 2, 1, 113, 1, 4, 1, 2, 1, 7, 7, 1, 2, 3, 1, 1, 2, 21, 30, 1, 9, 1, 12, 6, 1, 1, 1, 4, 2, 19, 1, 2, 1, 1, 466, 1, 5, 63, 2, 6, 1, 1, 1, 15, 29, 1, 2, 3, 12, 1, 9, 1, 12, 9, 1, 24, 7, 1, 6, 4, 56, 1, 2, 1, 1, 1, 2, 3, 9, 6, 9, 1, 5, 1, 1, 15, 12, 2, 2, 1, 4, 8, 1, 22, 11, 1, 1, 9, 1, 2, 1, 9, 14, 1, 1, 1, 3, 17, 3, 1, 1, 3, 4, 2, 20, 1, 1, 1, 10, 2, 1, 21, 1, 11, 1, 3, 2, 3, 1, 1, 2, 1, 1, 1, 1, 1, 7, 61, 12, 2, 1, 1, 3, 1, 4, 1, 1, 1, 4, 1, 6, 1, 38, 16, 1, 1, 2, 2, 3, 3, 1, 62, 2, 11, 1, 7, 8, 1, 2, 1, 1, 4, 23, 13, 2, 1, 1, 1, 4, 1, 1, 70, 1, 1, 1, 5, 1, 27, 11, 2, 4, 1, 8, 2, 1, 2, 1, 2, 1, 1, 1, 3, 1, 2, 2, 2, 2, 1, 11349, 103, 3, 1, 3, 1, 2, 8, 3, 1, 58, 1, 5, 1, 1, 1, 1, 1, 3, 1, 2, 7, 2, 2, 2, 1, 2, 3, 1, 7, 41, 1, 33, 1, 1, 3, 5, 1, 4, 1, 6, 2, 1, 2, 1, 19, 3, 22, 75, 1, 2, 5, 5, 3, 3, 4, 21, 7, 2, 1, 3, 1, 2, 4, 1, 1, 3, 3, 1, 1, 3, 3, 1, 1, 10, 1, 3, 1, 1, 1, 4, 1, 9, 2, 50, 5, 1, 1, 96, 1, 1, 2, 13, 33, 9, 8, 1, 5, 1, 1, 1, 3, 1, 29, 2, 1, 1, 2, 2, 17, 2, 1, 1, 1, 1, 33, 6, 170, 5, 7, 25, 1, 3, 1, 2, 7, 1, 1, 4, 1, 1, 1, 17, 1, 5, 4, 2, 2, 2, 1, 5, 1, 9, 3, 7, 1, 3, 1, 3, 1, 3, 2, 8, 8, 2, 1, 3, 1, 1, 2, 3, 4, 1, 1, 1, 1, 2, 9, 3, 1, 1, 1, 4, 2, 23, 1, 1, 4, 1, 1, 1, 1, 1, 7, 2, 2, 2, 2, 1, 2, 413, 10, 1, 1, 22, 2, 9, 29, 1, 4, 1, 6, 1, 49, 1, 2, 5, 88, 1, 2, 12, 1, 1, 3, 6, 1, 1, 1, 1, 92, 2, 3, 1, 5, 1, 5, 1, 10, 6, 1, 1, 2, 4, 2, 4, 38, 2, 3, 7, 1, 1, 3, 2, 8, 1, 4, 1, 1, 5, 1, 5, 2, 5, 1, 6, 1, 2, 5, 7, 1, 4, 1, 19, 1, 2, 140, 1, 2, 2, 2, 1, 2, 12, 1, 40, 2, 1, 1, 1, 11, 1, 2, 1, 1, 11, 6, 1, 6, 4, 3, 1, 6, 1, 3, 1, 3, 1, 1, 2, 2, 1, 69, 8, 9, 1, 21, 7, 1, 1, 1, 15, 1, 2, 2, 1, 2, 6, 1, 6, 3, 1, 2, 3, 1, 1, 21, 1, 1, 8, 6, 4, 2, 1, 1, 6, 1, 2, 39, 1, 4, 3, 1, 10, 1, 2, 1, 13, 3, 4, 12, 1, 3, 36, 1, 6, 1, 2, 1, 1, 3, 3, 11, 1, 1, 15, 1, 12, 1, 1, 59, 1, 1, 9, 1, 2, 6, 1, 1, 1, 6, 4, 3, 3, 1, 10, 8, 1, 1, 1, 3, 41, 1, 3, 17, 1, 2, 3, 7, 1, 4, 1, 1, 50, 1, 2, 1, 1, 3, 1, 1, 5, 3, 2, 2, 1, 2, 4, 1, 3, 1, 1, 1, 13, 2, 1, 3, 3, 1, 4, 4, 1, 1, 2, 4, 1, 24, 1, 1, 1, 8, 3, 1, 1, 1, 1, 3, 1, 36, 1, 5, 3, 4, 1, 1, 5, 2, 2, 1, 19, 3, 1, 9, 1, 11, 5, 21, 1, 1, 2, 2, 1, 2, 1, 4, 2, 3, 1, 1, 2, 1, 1, 1, 23, 2, 3, 1, 1, 1, 6, 1, 8, 5, 3, 1, 1, 1, 1, 2, 2, 46, 2, 1, 2, 2, 2, 1, 4, 3, 3, 1, 52, 1, 6, 1, 1, 1, 6, 3, 1, 2, 4, 1, 2, 2, 1, 3, 39, 2, 10, 1, 1, 1, 27, 1, 3, 5, 2, 10, 1, 12, 2, 55, 1, 1, 1, 2, 1, 1, 207, 8, 1, 1, 7, 7, 1, 8, 1, 1, 30, 1, 1, 2, 1, 1, 1, 17, 2, 3, 1, 4, 4, 2, 6, 1, 2, 1, 5, 4, 3, 3, 1, 2, 2, 1, 1, 1, 1, 1, 1, 4, 2, 1, 5, 42, 23, 12, 3, 6, 1, 2, 2, 1, 3, 2, 2, 1, 136, 2, 46, 3, 17, 31, 109, 1, 1, 1, 73, 1, 1, 1, 2, 1, 5, 2, 1, 7, 1, 9, 2, 4, 19, 6, 1, 2, 1, 1, 1, 1, 1, 22, 5, 1, 1, 1, 2, 3, 1, 10, 1, 10, 1, 17, 1, 2, 1, 1, 2, 3, 1, 7, 2, 6, 1, 2, 6, 5, 20, 2, 1, 3, 3, 3, 1, 8, 1, 9, 1, 1, 3, 1, 1, 1, 5, 8, 2, 16, 1, 2, 2, 2, 1, 3, 2, 1, 2, 14, 1, 10, 2, 237, 1, 4, 1, 1, 1, 23, 3, 2, 7, 1, 5, 7, 1, 1, 3, 2, 12, 1, 3, 5, 2, 4, 1, 1, 18, 1, 25, 8, 2, 4, 1, 1, 3, 1, 2, 1, 1, 1, 1, 6, 1, 1, 1, 88, 9, 1, 1, 1, 60, 1, 1, 2, 1, 3, 1, 25, 2, 1, 22, 2, 3, 354, 5, 3, 1, 3, 48, 4, 11, 6, 2, 46, 1, 1, 17, 1, 16, 1, 2, 1, 1, 1, 3, 1, 18, 2, 85, 5, 1, 7, 8, 1, 6, 44, 12, 4, 1, 2, 87, 3, 6, 15, 11, 5, 3, 1, 6, 2, 3, 2, 1, 2, 1, 1, 1575, 1, 1, 3, 7, 11, 4, 3, 470, 1, 5, 2, 3, 1, 3, 1, 12, 1, 6, 1, 7, 4, 2, 7, 11, 14, 7, 16, 34, 1, 1, 1, 7, 6, 1, 2, 1, 8, 1, 3, 1, 1, 1, 1, 2, 8, 1, 5, 1, 3, 1, 1, 2, 3, 1, 2, 6, 1, 3, 1, 10, 2, 1, 1, 5, 1, 3, 1, 2, 6, 1, 6, 1, 2, 1, 3, 2, 1, 1, 4, 1, 1, 2, 10, 2, 6, 20, 7, 9, 1, 45, 13, 2, 2, 1, 5, 2, 2, 3, 1, 2, 1365, 2, 7, 4, 1, 4, 2, 1, 1, 16, 2, 1, 2, 2, 1, 16, 2, 1, 2, 2, 1, 6, 1, 3, 5, 3, 1, 1, 1, 2, 1, 1, 1, 1, 2, 5, 8, 1, 8, 66, 1, 57, 1, 14, 1, 1, 1, 2, 1, 9, 1, 1, 2, 1, 101, 2, 6, 10, 1, 36, 1, 5, 4, 1, 1, 5, 3, 1, 2, 146, 2, 41, 2, 394, 41, 3, 1, 2, 3, 16, 22, 1, 4, 1, 1, 2, 1, 3, 3, 2, 11, 1, 7, 7, 1, 4, 10, 9, 21, 1, 3, 1, 4, 1, 3, 1, 12, 1, 4, 1, 1, 1, 19, 1, 3, 3, 4, 1, 2, 1, 2, 3, 5, 2, 2, 1, 1, 1, 1, 2, 5, 7, 6, 4, 19, 1, 7, 3, 1, 1, 1, 2, 2, 18, 1, 1, 9, 1, 4, 1, 1, 55, 4, 5, 3, 6, 1, 3, 2, 1, 16, 7, 5, 1, 1, 1, 1, 7, 1, 1, 2, 3, 3, 2, 7, 1, 1, 4, 1, 11, 146, 1, 14, 1, 19, 4, 5, 1, 1, 3, 1, 2, 1, 1, 1, 5, 1, 1, 1, 8, 1, 21, 5, 2, 1, 1, 2, 9, 1, 1, 1, 2, 1, 2, 3, 1, 1, 1, 1, 3, 3, 1, 4, 1, 1, 2, 1, 7, 2, 1, 4, 1, 3, 1, 1, 2, 2, 1, 1, 1, 3, 1, 1, 6, 1, 1, 1, 9, 1, 1, 5, 2, 1, 2, 2, 7, 1, 6, 1, 1, 1, 1, 1, 81, 1, 100, 21, 1, 4, 1, 5, 2, 4, 1, 2, 8, 1, 1, 3, 1, 1, 5, 1, 1, 5, 1, 1, 14, 3, 1, 3, 1, 3, 37, 22, 6, 1, 1, 26, 1, 2, 5, 1, 1, 1, 5, 102, 1, 1, 1, 1, 2, 1, 1863, 1, 5, 14, 1, 2, 1, 1, 9, 2, 2, 8, 2, 11, 7, 17, 16, 1, 1, 2, 1, 2, 6, 3, 1, 1, 7, 1, 1, 2, 3, 16, ...]}
\end{table}
\begin{table}[H]
\caption{\label{tab:CFR-1013}$(\textrm{CFR})\,\left(2^{1013}\log^{2}(2)\log(3)(2^{1013}\log(3)+1)\right)^{-1}\rightarrow\left[\gamma_{0}^{(1013)};\gamma_{\alpha>0}^{(1013)}\right]$}
\medskip{}

{\small  [0;446797816913050870465576532067813725399406109382194399324684148252196866110567542645972215\ 51042698178872322530299781422149535295773818402503863078008759249213511160998729034034601722\ 53542820108216229091785464647151083920754885299791594148815899868103475192998274383774383189\ 70282395788641891021517839079740008630135896413015936100265961906063964119535356896138152344\ 30571772501916602693436495891143030709703197384256808876096299877501263928380980130737997466\ 86125608979772280262828091947728365371149503468515864646034948474485754382796401019228789558\ 228884483772395515280126760714190895298584447020468622283287, 3, 1, 1, 6, 1, 1, 1, 1, 1, 2, 2, 5, 1, 2, 2, 1, 3, 1, 2, 1, 1, 2, 5, 4, 3, 1, 1, 2, 1, 1, 3, 1, 60, 57, 1, 9, 1, 1, 1, 21, 1, 22, 1, 26, 20, 52, 1, 1, 2, 3, 6, 1, 1, 3, 2, 7, 1, 2, 2, 1, 2, 3, 1, 1, 6, 1, 1, 1, 2, 1, 26, 6, 3, 1\ 2, 2, 1, 3, 1, 1, 30, 4, 1, 4, 4, 2, 6, 1, 5, 2, 166, 100, 1, 3, 3, 1, 1, 1, 3, 3, 10, 1, 10, 1, 12, 2, 2, 1, 13, 13, 2, 4, 5, 1, 9, 1, 1, 1, 1, 3, 1, 8, 2, 2, 1, 5, 1, 1, 4, 8, 1, 5, 1, 1, 1, 1, 22, 1, 2, 34, 1, 3, 2, 3, 6, 1, 1, 1, 1, 1, 126, 2, 10, 1, 2, 1, 5, 32, 3, 2, 1, 65, 15, 2, 3, 1, 7, 1, 15, 1, 87, 18, 4, 3, 8, 1, 2, 2, 2, 5, 44, 1, 55, 1, 1, 11, 1, 6, 2, 1, 12, 1, 3, 6, 1, 2, 1, 1, 9, 6, 3, 5, 1, 2, 1, 1, 10, 1, 1, 1, 1, 11, 1, 1, 4, 6, 1, 3, 231, 3, 1, 1, 18, 1, 1, 1, 7, 1, 2, 23, 2, 2, 6, 1, 1, 4, 1, 2, 3, 8, 1, 1, 1, 1, 4, 38, 6, 2, 1, 1, 2, 2, 4, 2, 2, 21, 2, 55, 1, 1, 7, 27, 120, 1, 2, 1, 1, 1, 2, 9, 2, 2, 1, 1, 22, 15, 1, 3, 1, 1, 12, 195, 1, 1, 2, 1, 4, 6, 1, 4, 8, 13, 4, 1, 8, 2, 2, 1, 48, 4, 5, 2, 2, 2, 1, 10, 27, 1, 1, 19, 1, 2, 1, 1, 1, 4, 2, 7, 1, 3, 1, 1, 1, 15, 2, 1, 3, 17, 2, 2, 2, 1, 1, 2, 56, 1, 1, 3, 1, 1, 1, 11, 2, 2, 13, 1, 1, 6, 1, 1, 1, 1, 1, 7, 1, 2, 4, 2, 2, 1, 43, 4, 1, 14, 3, 6, 84, 2, 2, 4, 4, 2, 2, 3, 1, 5, 1, 3, 1, 89, 1, 1, 1, 104, 1, 2, 2, 1, 38, 1, 3, 1, 1, 2, 2, 1, 1, 3, 1, 1, 1, 6, 1, 1, 6, 2, 1, 3, 1, 2, 19, 1, 5, 1, 1, 3, 2, 5, 1, 2, 1, 3, 1, 8, 3, 6, 5, 10, 17, 1, 7, 1, 1, 1, 1, 1, 1, 2, 2, 1, 1, 2, 4, 1, 7, 1, 1, 1, 10, 1, 3, 4, 2, 1, 1, 1, 10, 2, 7, 1, 3, 1, 1, 1, 1, 2, 24, 1, 1, 1, 1, 2, 1, 38, 1, 4, 1, 2, 3, 1, 4, 3, 67, 1, 10, 3, 4, 6, 1, 2, 1, 1, 9, 4, 2, 1, 36, 1, 1, 34, 3, 1, 2, 1, 6, 1, 9, 33, 33, 6, 5, 17, 1, 2, 3, 1, 1, 1, 9, 6, 2, 2, 9, 3, 1, 1, 1, 6, 1, 2, 3, 8, 18, 1, 7, 2, 15, 1, 5, 3, 3, 1, 2, 3, 1, 1, 121, 1, 2, 4, 5, 1, 62, 4, 1, 1, 2, 1, 14, 11, 6, 21, 39, 1, 11, 1, 27, 2, 4, 6, 1, 2, 2, 9, 1, 1, 1, 19, 2, 1, 1, 11, 1, 2, 1, 3, 1, 5, 1, 1, 4, 1, 1, 1, 5, 1, 3, 7, 7, 1, 7, 6, 21, 1, 1, 1, 9, 1, 5, 3, 1, 1, 2, 1, 1, 22, 2, 1, 7, 13, 7, 4, 1, 109, 1, 3, 2, 54, 2, 2, 3, 2, 1, 3, 1, 12, 1, 6, 168, 1, 1, 1, 3, 1, 1, 1, 3, 10, 1, 9, 64, 1, 4, 1, 3, 2, 5, 2, 1, 12, 2, 3, 1, 6, 1, 7, 1, 2, 4, 18, 1, 1, 1, 5, 3, 1, 15, 20, 20, 3, 2, 22, 10, 1, 1, 1, 3, 2, 2, 4, 1, 2, 10, 1, 4, 2, 1, 2, 2, 1, 1, 39, 1, 13, 1, 12, 3, 1, 8, 1, 2, 1, 2, 3, 8, 2, 1, 3, 2160, 1, 4, 2, 2, 3, 2, 2, 30, 1, 13, 2, 4, 1, 4, 5, 2, 6, 1, 1, 3, 5, 2, 3, 3, 2, 1, 8, 3, 3, 15, 3, 1, 1, 2, 1, 2, 1, 1, 4, 3, 1, 1, 1, 3, 6, 3, 6, 1, 3, 1, 28, 1, 1, 11, 4, 3, 2, 1, 2, 2, 449, 1, 1, 1, 6, 2, 10, 1, 1, 41, 12, 1, 7, 1, 1, 8, 2, 1, 5, 2, 9, 1, 2, 3, 1, 5, 36, 1, 15, 3, 1, 13, 1, 6, 1, 1, 5, 1, 2, 2, 1, 8, 1, 15, 1, 2, 1, 1, 6, 1, 1, 4, 2, 1, 3, 55, 1, 7, 1, 19, 2, 7, 1, 2, 1, 1, 1, 5, 1, 1, 1, 2, 1, 142, 16, 1, 3, 3, 1, 12, 1, 2, 2, 3, 1, 1, 14, 4, 1, 1, 1, 204, 2, 1, 2, 1, 3, 2, 5, 8, 1, 2, 2, 1, 1, 1, 1, 2, 1, 4, 1, 4, 1, 3, 19, 1, 1, 1, 6, 1, 61, 1, 1, 1, 1, 2, 1, 2, 1, 3, 2, 12, 2, 3, 1, 1, 1, 290, 20, 1, 3, 3, 1, 11, 1, 1, 9, 4, 5, 2, 2, 4, 3, 2, 2, 7, 23, 1, 2, 2, 1, 10, 2, 1, 2, 1, 130, 1, 16, 9, 3, 1, 60, 2, 1, 2, 1, 3, 1, 8, 23, 1, 1, 1, 1, 6, 3, 1, 3, 4, 1, 2, 8, 1, 13, 1, 7, 419, 1, 1, 1, 1, 2, 5, 1, 4, 7, 1, 2, 4, 6, 3, 1, 2, 1, 3, 5, 1, 1, 309, 1, 3, 16, 1, 2, 2, 15, 7, 1, 4, 3, 1, 1, 3, 1, 1, 1, 1, 3, 4, 1, 1, 3, 1, 13, 80, 1, 4, 3, 2, 1, 1, 4, 1, 1, 6, 41, 1, 3, 1, 2, 1, 19, 11, 2, 2946, 1, 1, 7, 1, 12, 1, 4, 2, 2, 2, 1, 7, 1, 3, 1, 2, 1, 3, 16, 3, 3, 1, 1, 1, 5, 1, 2, 2, 1, 42, 6, 3, 5, 1, 1, 3, 6, 2, 1, 9, 3, 3, 2, 15, 5, 1, 4, 1, 1, 3, 19, 2, 1, 2, 9, 1, 153, 1, 5, 3, 1, 2, 1, 30, 1, 15, 2, 6, 1, 1, 1, 4, 5, 1, 7, 1, 4, 26, 4, 3, 1, 2, 8, 2, 1, 1, 1, 4, 1, 1, 1, 1, 1, 15, 2, 2, 2, 1, 10, 2, 1, 2, 3, 12, 18, 2, 2, 1, 33, 5, 2, 1, 1, 22, 39, 1, 2, 18, 6, 2, 2, 8, 1, 3, 1, 4, 8, 1, 1, 1, 121, 1, 5, 2, 3, 13, 1, 4, 1, 10, 3, 44, 1, 1, 3, 3, 1, 5, 7, 2, 2, 8, 6, 2, 1, 1, 1, 3, 3, 1, 11, 16, 1, 2, 7, 1, 10, 3, 19, 1, 2, 2, 1, 2, 6, 2, 3, 1, 2, 1, 1, 56, 1, 101, 2, 1, 1, 3, 56, 1, 12, 1, 3, 1, 9, 1, 1, 2, 1, 4, 7, 1, 1, 4, 1, 8, 1, 2, 4, 175, 1, 2, 17, 18, 5, 2, 14, 1, 3, 1, 5, 2, 11, 1, 2, 1, 9, 2, 2, 3, 16, 5, 1, 9, 6, 1, 3, 3, 5, 3, 2, 2, 3, 9, 14, 2, 3, 3, 2, 1, 2, 3, 2, 2, 5, 1, 8, 1, 3, 140, 4, 2, 2, 2, 1, 1, 7, 1, 1, 1, 2, 2, 1, 2, 2, 3, 1, 50, 6, 3, 3, 7, 1, 1, 11, 13, 2, 1, 2, 1, 1, 1, 1, 1, 8, 2, 1, 1, 1, 22, 1, 5, 1, 2, 2, 2, 26, 1, 4, 1, 5, 2, 1, 1, 2, 6, 10, 3, 5, 1, 4, 1, 1, 4, 1, 4, 2, 4, 1, 3, 1, 2, 4, 2, 9, 1, 1, 40, 1, 1, 15, 1, 2, 18, 1, 1, 10, 6, 1, 4, 10, 4, 8, 1, 6, 5, 1, 6, 7, 1, 1, 1, 1, 9, 23, 1, 11, 1, 3, 8, 5, 1, 3, 13, 1, 338, 11, 1, 1, 1, 6, 1, 3, 31, 1, 1, 4, 1, 2, 8, 6, 1, 33, 1, 3858, 4, 1, 4, 101, 2, 4, 4, 3, 1, 3, 4, 2, 5, ... ]}
\end{table}
\medskip{}
\begin{table}[H]
\caption{\label{tab:Specific-fractions-in-1-1-2}Unruh effect for `observer'
$\mathcal{O}$ with $\mathcal{L}^{(1013)}=\{G_{\mu\nu}^{(15)}\setminus113\}\cup\{31,104,419\}$
in $(\textrm{CFR})\,\left(2^{1013}\log^{2}(2)\log(3)(2^{1013}\log(3)+1)\right)^{-1}\rightarrow\left[\gamma_{0}^{(1013)};\gamma_{\alpha>0}^{(1013)}\right]$}
\medskip{}
 \enskip{}%
\begin{tabular}{|c|c|c|c|c|}
\hline 
 & \multicolumn{4}{c|}{`Heat bath' $\mathcal{H}$ in (CFR) $\left(2^{1013}\log^{2}(2)\log(3)(2^{1013}\log(3)+1)\right)^{-1}$}\tabularnewline
$\begin{array}{c}
\\
\\
\end{array}\gamma_{\alpha_{\zeta}}^{(1013)}$ & $o_{i}$ & $o_{j}^{+}$ & $\left\lceil \sigma_{k}\right\rceil $ & $\left\lceil \sigma_{l}^{-}\right\rceil $\tabularnewline
\hline 
\hline 
{[}1{]} &  &  &  & \tabularnewline
\hline 
{[}3{]} &  &  &  & \tabularnewline
\hline 
{[}5{]} &  &  &  & \tabularnewline
\hline 
4 & 4 &  & (4) & \tabularnewline
\hline 
60 &  & 20,40 &  & \tabularnewline
\hline 
57 &  &  &  & 7,16,34\tabularnewline
\hline 
9 & 9 &  &  & \tabularnewline
\hline 
21 &  &  & 4,17 & \tabularnewline
\hline 
22 &  &  & 4,8 & 3,7\tabularnewline
\hline 
26 &  &  &  & 3,7,16\tabularnewline
\hline 
20 &  & 20 &  & \tabularnewline
\hline 
52 & 4,9,39 &  &  & \tabularnewline
\hline 
30 &  & 10,20 &  & \tabularnewline
\hline 
166 &  &  & 8,17,71 & 70\tabularnewline
\hline 
100 &  & 20,80 &  & \tabularnewline
\hline 
10 &  & 10 &  & \tabularnewline
\hline 
8 &  &  & 8 & \tabularnewline
\hline 
34 &  &  &  & 34\tabularnewline
\hline 
32 & 4,9,19 &  &  & \tabularnewline
\hline 
65 &  & 5,20,40 &  & \tabularnewline
\hline 
15 &  & 5,10 &  & \tabularnewline
\hline 
87 &  &  & 71 & 16\tabularnewline
\hline 
18 & 4,9 & 5 & (8) & (3,7)\tabularnewline
\hline 
44 &  &  &  & 3,7,34\tabularnewline
\hline 
55 &  & 5,10,40 &  & \tabularnewline
\hline 
231 &  &  & 17,71,143 & \tabularnewline
\hline 
38 &  &  & 4 & 34\tabularnewline
\hline 
120 &  & 40,80 &  & (16,34,70)\tabularnewline
\hline 
195 &  &  & 17,35,143 & \tabularnewline
\hline 
48 & 9,39 &  &  & \tabularnewline
\hline 
27 &  &  & 17 & 3,7\tabularnewline
\hline 
19 & 19 &  &  & \tabularnewline
\hline 
$\ldots$ &  &  &  & \tabularnewline
\hline 
{[}11{]} &  &  &  & \tabularnewline
\hline 
{[}17{]} &  &  &  & \tabularnewline
\hline 
{[}41{]} &  &  &  & \tabularnewline
\hline 
{[}31{]} &  &  &  & \tabularnewline
\hline 
{[}104{]} &  &  &  & \tabularnewline
\hline 
{[}419{]} &  &  &  & \tabularnewline
\hline 
\end{tabular}
\end{table}
\newpage{}

\begin{table}[H]
\caption{\label{tab:Specific-fractions-in-1-2-1} Coincidences $\between_{\alpha}^{(n)}=L_{m}$$^{,\mp}$$(\: n\leq6262)$
$(\textrm{CFR})\left(2^{n}\log(2)\log(3)\right)^{-1}\!\rightarrow\!\left[\between_{0}^{(n)};\between_{\alpha>0}^{(n)}\right]$
(match$=\:$\Checkmark{}; closest pivot $=$ ($\cdot$); largest $\between_{\alpha}^{(n)}>L_{31}=$($\cdot$))}
\vspace{1cm}
\qquad{}
\begin{tabular}{|c|c|c|c|c|}
\hline 
 & \multicolumn{4}{c|}{$\begin{array}{ccccc}
\\
\\
\\
\end{array}$$(\textrm{CFR})\,\left(2^{n}\log(2)\log(3)\right)^{-1}$}\tabularnewline
$L_{m}$ & $\between_{\alpha}^{(n)}\!=L_{m}$ & $\!\!\begin{array}{c}
\textrm{in-}\\
\textrm{cidence}
\end{array}$ & $\begin{array}{c}
\!\!\between_{\alpha}^{(n)}\!=\qquad\\
L_{m}^{-},L_{m}^{+}
\end{array}$ & $\!\!\begin{array}{c}
\textrm{in-}\\
\textrm{cidence}
\end{array}$\tabularnewline
\hline 
\hline 
$\begin{array}{c}
2,6,12,24,\\
40,72,126
\end{array}$ & \Checkmark{} & $\begin{array}{c}
\textrm{very}\\
\textrm{high}
\end{array}$  & \Checkmark{},\ \Checkmark{} & $\begin{array}{c}
\textrm{very}\\
\textrm{high}
\end{array}$ \tabularnewline
\hline 
$240$ & \Checkmark{} & 184 & \Checkmark{},\ \Checkmark{} & 331,249\tabularnewline
\hline 
$272$ & \Checkmark{} & 161 & \Checkmark{},\ \Checkmark{} & 219,278\tabularnewline
\hline 
$336$ & \Checkmark{} & 101 & \Checkmark{},\ \Checkmark{} & 139,164\tabularnewline
\hline 
$438$ & \Checkmark{} & 65 & \Checkmark{},\ \Checkmark{} & 98,110\tabularnewline
\hline 
$756$ & \Checkmark{} & 24 & \Checkmark{},\ \Checkmark{} & 29,55\tabularnewline
\hline 
$918$ & \Checkmark{} & 8 & \Checkmark{},\ \Checkmark{} & 12,14\tabularnewline
\hline 
$1422$ & \Checkmark{} & 8 & \Checkmark{}~~~~~ & 1,1\tabularnewline
\hline 
$2340$ & \Checkmark{} & 5 & \Checkmark{},\ \Checkmark{} & 2,4\tabularnewline
\hline 
$4320$ &  & 0 & \Checkmark{}~~~~~ & 1,0\tabularnewline
\hline 
$5346$ & \Checkmark{} & 1 &  & 0,0\tabularnewline
\hline 
$7398$ & \multicolumn{2}{c|}{(7406)} & \multicolumn{2}{c|}{(1)}\tabularnewline
\hline 
$10\,668$ & \multicolumn{2}{c|}{(10\,661)} & \multicolumn{2}{c|}{(2)}\tabularnewline
\hline 
$17\,400$ & \multicolumn{2}{c|}{(17\,369)} & \multicolumn{2}{c|}{(1)}\tabularnewline
\hline 
$27\,720$ & \multicolumn{2}{c|}{(27\,762)} & \multicolumn{2}{c|}{(1)}\tabularnewline
\hline 
$49\,896$ & \multicolumn{2}{c|}{$\quad$(49\,872)$\,^{a}$} & \multicolumn{2}{c|}{(1)}\tabularnewline
\hline 
$93\,150$ & \multicolumn{2}{c|}{(93\,532)} & \multicolumn{2}{c|}{(1)}\tabularnewline
\hline 
$\vdots$ & \multicolumn{2}{c|}{} & \multicolumn{2}{c|}{}\tabularnewline
\hline 
$>L_{31}$ & \multicolumn{2}{c|}{(14\,571\,717)} & \multicolumn{2}{c|}{(1)}\tabularnewline
\hline 
\multicolumn{5}{|c|}{}\tabularnewline
\multicolumn{5}{|c|}{$^{a}$ \scriptsize Samples such as $\between^{(3607)}_{\alpha_a}=27448$ $(\simeq L_{21}-L_9)$, $\between^{(3759)}_{\alpha_b}=49872$ $(\simeq L_{22}-L_4)$ or}\tabularnewline
\multicolumn{5}{|c|}{\scriptsize $\between^{(3990)}_{\alpha_c}=2316$ $(\simeq L_{15}-L_4)\,$ suggest that neutrino plateaux aren't the only ones}\tabularnewline
\multicolumn{5}{|c|}{\scriptsize that play a role in $\between^{(n)}_{\alpha}\hspace*{6.9cm}$.}\tabularnewline
\hline 
\end{tabular}
\end{table}

\vspace{1.5cm}

\newpage
\section{\label{sec:A-special-cube}An off-diagonal cube complex}

It was shown in Sect.$\,$\ref{sub:Beta} that various quark and neutrino calculations proceed well when using a bicube complex, plus/minus a term $o_\nu\in M_{5/8}$ in case of a quark (constituent) and  $\sigma_\nu\in M_{9/8}$ in case of a neutrino. The bicube complex is composed  of a $(T_{\scriptscriptstyle\textrm{SD}}=6)$ part ($\textrm{SD}=(227089,15297,1633,113,17,1)$) and a $(T_{\scriptscriptstyle\textrm{MD}}=3)$ part ($\textrm{MD}=(429,5,1)$), respectively taken from the secondary and main diagonals of $\mathrm{LL}\,(G_{\mu\nu}^{(31)})$, $\mathrm{LL}\,(G_{\mu\nu}^{(15)})$ and $\mathrm{LL}\,(G_{\mu\nu}^{(7)})$ (entries $2\,430\,289$ and 1 trimmed off):

\[
\mathrm{LL}\,(G_{\mu\nu}^{(31)})=\qquad\qquad\qquad\qquad
\]
$\quad\left(\begin{array}{c}
\\
\\
\\
\\
\\
\\
\\
\\
\end{array}\right.$\hspace{-0.3cm}%
\begin{tabular}{cccc|cccc}
$\underline{\mathbf{429}}$ & \textcolor{lightgray}{155} & \textcolor{lightgray}{43} & \textcolor{lightgray}{19} & $\underline{\mathbf{5}}$ & \multicolumn{1}{c|}{\textcolor{lightgray}{3}} & $\underline{\textrm{\bf1}}$ & \textcolor{lightgray}{1}\tabularnewline
\textcolor{lightgray}{1275} & $\underline{\mathbf{429}}$ & \textcolor{lightgray}{115} & \textcolor{lightgray}{43} & \textcolor{lightgray}{11} & \multicolumn{1}{c|}{$\underline{\mathbf{5}}$} & {\bf1} & $\underline{\textrm{\bf1}}$\tabularnewline
\cline{7-8} 
\textcolor{lightgray}{4819} & \textcolor{lightgray}{1595} & $\underline{\mathbf{429}}$ & \textcolor{lightgray}{155} & \textcolor{lightgray}{41} & \bf{17} & $\underline{\mathbf{5}}$ & \textcolor{lightgray}{3}\tabularnewline
\textcolor{lightgray}{15067} & \textcolor{lightgray}{4819} & \textcolor{lightgray}{1275} & $\underline{\mathbf{429}}$ & {\bf113} & \textcolor{lightgray}{41} & \textcolor{lightgray}{11} & $\underline{\mathbf{5}}$\tabularnewline
\cline{5-8} 
\textcolor{lightgray}{58781} & \textcolor{lightgray}{18627} & \textcolor{lightgray}{4905} & \multicolumn{1}{c}{\bf1633} & $\underline{\mathbf{429}}$ & \textcolor{lightgray}{155} & \textcolor{lightgray}{43} & \textcolor{lightgray}{19}\tabularnewline
\textcolor{lightgray}{189371} &\textcolor{lightgray}{58781}& {\bf15297} & \multicolumn{1}{c}{\textcolor{lightgray}{4905}} & \textcolor{lightgray}{1275} & $\underline{\mathbf{429}}$ &\textcolor{lightgray}{115} & \textcolor{lightgray}{43}\tabularnewline
\textcolor{lightgray}{737953} & {\bf227089} & \textcolor{lightgray}{58781} & \multicolumn{1}{c}{\textcolor{lightgray}{18627}} & \textcolor{lightgray}{4819} & \textcolor{lightgray}{1595} & $\underline{\mathbf{429}}$ &\textcolor{lightgray}{155}\tabularnewline
\textcolor{lightgray}{2430289} & \textcolor{lightgray}{737953} & \textcolor{lightgray}{189371} & \multicolumn{1}{c}{\textcolor{lightgray}{58781}} & \textcolor{lightgray}{15067} &\textcolor{lightgray}{4819} &\textcolor{lightgray}{1275} & $\underline{\mathbf{429}}$\tabularnewline
\end{tabular}\hspace{-0.2cm}$\left.\begin{array}{c}
\\
\\
\\
\\
\\
\\
\\
\\
\end{array}\right)$
\bigskip{}

\noindent Surprisingly, $L_m$ not involved with quarks can  be calculated by means of a   single $T_{\scriptscriptstyle\textrm{OD}}\,$-cube complex whose  basis consists solely of the off-diagonal elements  (light-gray) of $\mathrm{LL}\,(G_{\mu\nu}^{(31)})$, 
 \[\textrm{OD}=\qquad\qquad\qquad\qquad\qquad\]\[(737953,\!189371,\!58781,18627,\!15067,4905,4819,\!1595,\!1275,\!155,\!115,\!43,\!41,\!19,\!11,\!3),\] 
plus/minus   three auxiliaries to compensate for the missing bicubicity, each  of  its own type: $p_\nu\in M_\textrm{reg}$, $o_\nu\in M_{5/8}$, $\sigma_\nu\in M_{9/8}$:

\vspace*{1mm}
 \noindent  $L_{15}=2340=(0,0,0,0,0,0,0,0,1,0,0,0,-1,0,-1,0)\cdot\textrm{OD}^{t}+p_9+o_7+\sigma_7$,
 
 \noindent $L_{16}=4320=(T_{\scriptscriptstyle\textrm{SD}}$,$\,T_{\scriptscriptstyle\textrm{MD}})$-determined
 
 \noindent $L_{17}=5346=(0,0,0,0,0,0,1,1,1,-1,0,-1,0,1,1,0)\cdot\textrm{OD}^{t}-p_{11}-o_9+\sigma_9$,

 \noindent $L_{18}=7398=(0,0,0,0,0,1,0,1,-1,0,0,-1,1,0,0,0)\cdot\textrm{OD}^{t}+p_{11}+o_9-\sigma_9$,
 
 \noindent $L_{19}=10668=(T_{\scriptscriptstyle\textrm{SD}}$,$\,T_{\scriptscriptstyle\textrm{MD}})$-determined
 
 \noindent $L_{20}=17400=(0,0,0,1,0,1,-1,1,-1,1,1,0,0,1,0,-1)\cdot\textrm{OD}^{t}-p_{11}+o_9-\sigma_9$,
 
 \noindent $L_{21}=27720=(0,0,0,1,\!-1,1,0,1,0,\!-1,\!-1,1,0,\!-1,\!-1,0)\cdot\textrm{OD}^{t}+p_{13}+o_{11}+\sigma_{11}$,
 
 \noindent $L_{22}=49896=(0,0,1,-1,1,1,1,1,1,-1,1,0,1,0,0,-1)\cdot\textrm{OD}^{t}-p_{13}-o_{11}-\sigma_{11}$,

 \noindent $L_{23}=93150=(0,1,0,1,0,1,0,1,1,1,1,0,0,0,-1,-1)\cdot\textrm{OD}^{t}-p_{17}+o_{15}-\sigma_{15}$,

 \noindent $L_{24}=196560=(0,0,1,0,1,0,0,0,0,0,-1,0,-1,0,-1,0)\cdot\textrm{OD}^{t}+p_{17}-o_{15}+\sigma_{15}$,
 
 \noindent $\qquad\;\quad\dots$

\vspace*{0.1mm}
 \noindent While the juxtaposition $x$ vs. $x-2$ lets  dimensions fall in two categories, `source' and `sink' dimensions (Table \ref{tab:Key-particle-creation-related-2}) so that  in quark calculations  contributions cancel each other out and  in neutrino calculations get chirally `sink-fixed,' its effect here is on the dimension of $p$ vs. those of $o$ and $\sigma$. Innately adaptable, OD then conspires to ensure that contributions still cancel each other out $-$ in a novel way  quite unrelated to  quarks or neutrinos.

 \section{\label{sec:Minkowski}CF expansions  in tetrapetalic coordinates}

\begin{table}[H]
\caption{\label{tab:CFR-1422}$(\textrm{CFR})\quad t_\nu= (2\cdot1422^2)^2+\frac{\sqrt{2\cdot1422^2-1}}{2\cdot 1422^2}$}
[16355294812224; 
2011,  83,   1,   3,   1,   4,   6,   1,   3,   1,   2,   6,   2,   1,   1,   7,  15,   1,  73,   3,  12,   6,   1, 133,   1,   1,   1,   1,   1,  12,
   9,  24,   1,   5,   1,   1,   2,  12,   1,  33,   1,   1,   1,   5,   4,   1,  22,   2,  12,   1,  20,   5,   1,  27,   5,   3,  16,   2,   1,   2,
   1,   1,   2,   1,   3,   3,   5,   3,   2,   1,   8,   1,   1,  11,   1,   3,   1,   3,   1,   7,   1,   3,   2,   2,   1,   1,   1,   9,   4,   1, 	
   1,   2,   2,   1,   6,  20,   2,  14,   6,   2,   2,   7,   3,   2,   3,  18,   1,   1,  12,   6,   1,   1,   1,   1,  11,   1,   2,   2,   1,   2,
   2,   2,   4,   1,   1,   2,   5,   1,   5,   2,   4,   1,   7,  19,   7,   1,   1,   1,   1,   1,   1,   2,   1,   1,   1,   7,   2,   6,   1,   1,
   6,   2,   2,   1,   1,   1,   2,   4,   1,  12,   1,  33,   1,   1,   4,   1,  10,  58,   4,  18,   1,  23,   2,   3,  23,   1,   1,   3,   3,   2,	
   3,   1,   1,   3,   2,   3,   2,   3,   5,   1,   3,   1,   2,   2,   2,   1,   3,  38,  22,   5,  47,   1,   1,   1,   6,   2,   2,   1,  27,   5,
   2,   2,   1,   1,   5,   2,   1,   1,   8,  13,   2,   8,   1,   3,   4,  14,   1,   2,   1,   3,   3,   1,   1,   3,   1,   2,   1,   5,   1,   3,
   1,   1,   4,   4,   1,   1,  86,   1,   2,   1, 129,   3,   2,  81,   2,   2,   6,   1,  17,   2,  33,   1,   1,   8,   6,   3,   1,   7,   1,   4,
  22,   1,   1,   1,   1,   1,   1,  10,   1,   1,   2,   4,   1,   1,   1,   1,  16,   1,   1,   1,   4,   8,   1,  35,   2,   1,   1,   1,   3,   8,
   1,   1,   1,   1,   3,  12,   1,   4,   3,   5,   2,   1,   1,   2,   6,   3,  38,   3,   3,   2,   1,   1,   3,   3,   1,   1,   6,   3,   1,   4,
   1,   1,  35,   1,  37,   1,   1,   6,   7,   5,   4,  14,   1,   1, 197,   1,   3,   2,   5,  10,   1,   1,   1,   2,   1,   1,   2,   2,   3,  10, 
   2,   1,   1,   9,   1,   1,   1,   1,  10,   3,   1,   7,   9,  12,   1,  28,   2,  12,   3,   3,   2,   2,   1,  15,   1, \textcolor{red}{286},   1,  25,   1,   3,	
   4,   3,   2,   2,   1,   1,   2,   7,   1,   5,   1,   1,   1,   1,   4,  34,   1,  14,   1,   2,   4,   3,   2,   5,   1,  34,   2,   1,   1,   1,
   3,   1,  88,   1,   2, 219,   2,   1,   1,  11,   2,   1,   6,   2,  19,   1,   1,   1,   9,   1,  47,   1,   1,   1,  17,   1,   1,   4,   1,   9,	
   1,   1,   2,  20,   2,   1,   4,   1,  10,  35,   1,   2,   1,   2,   1,   2,  16,   1,   2,   1,   1,   1,  14,   1,   2,  14,   4,   2,   3,   1,
   2,   1,   1,   1,   1,   7,   2,   6,   1,   1,   1,   7,   1,   3,   4,  40,   1,   6,   1,   7,  13,   8,   1,   1,   6,  11,   1,   1,   5,   1,
  32,   1,   1,   1,   7, 127,   2,  38,   2,   1,   4,   1,   2,   1,   2, 122,   1,   1,  16,  18,   1,   2,  29,   1,   2,   4,   1,   8,  17,  26,	
   2,   2,   2,   1,   4,   6,   1,1391,   1,   2,   1,   1,   1,   1,   2,  10,   5,   2,   5,   2,   9,  10,  26,   4,   1,  29,   1,  56,  10,  14,
   6,  19,   2,  16,   1,   1,  21,   2,   1,   1,   9,   1,   3,   1,   3,   3,   1,   1,   1,  90, 841,   2,   2,   4,   3,   2,   1,   5,   2,   1,
   1,   4,   1,   1,  15,   1,   2,   2,   5,  28,   2,   3,   1,   2,   1,   2,   1,   3,   1,   1,   7,   3,   1,   3,  55,   4,   1,   5,   2,   1,	
   1,   3,   2,   1,   6,   1,   1,   3,   2,   1,   1,   1,   1,   2,   1,   3,   1,   8,   1,  14,   2,   1,   1,   1,   6,  53,   1,  10,   6,   1,
   4,   1,   3,   1,  33,   4,   1,  10,   1,   2,   1,   8,   1,   1,   1,   1,  11,   1,   3,   3,   7,   1,  28,  12,   1,   4,  20,   2,  20,  12,
   1,   1,   1,   5,   3,   3,   3,   1,   8,   1,   5,   1,   2,   2,   1,   1,   1,   1,   1,   1,   1,   1,   1,   2,  11,   1,   2,  12,   1,   4,	
   1,   1,   4,   1,  12,   6,   3,   2,   1,   1,   1,   1,   1,   1,   6,   2,   6,  10,  38,   3,   3,   1,   3,   2,  50,   1,   2,  23,  31,   2,
   2,  12,   3,   1,   1,   3,   3,  12,   7,   2,  85,   1,   1,   4,   2,   4,   2,   2,   3,   3,  75,   1,   1,   2,   8,   2,  18,   1,   7,   6,
   1,   2,   2,   4,   2,   7,   1,   4,   2,   1,   8,   1,  34,   1,   1,   9,   1,   1,   2,   1,   1,   1,   1,   3,   1,   3,   2,   1,   2,   1,	
   1,   1,   1,   1,   3,   5,   1,  10,   4,   1,   6,   1, 110,   1,   2,  10,   1,   3,   5,   1,   1,   1,   7,   2,   1,   1,   4,   4,   1,   1,
   3,   2,   1,   3,   1,  10,   5,   5,   4,  20,   3,   6,  13,   1,   9,   1,   3,   1,   2,   8,   1,  10,   1,   2,   7,   3,   3,   9,   1,   2, 
  77,  38,   1,   2,   2,   1,   1,   3,   2,   2,   1,   8,  18,   1,   1,  14,   6,   2,   2,   4,   4,   3,   3,   6,   2,  97,   1,  12,   1,   4,	
   3,   2,   2,   1,   1,  24,   5,   5,   1,   2,  11,   3,   1,   1,   1,   5,   1,  13,   2,   2,   1,   2,   2,   1,   2,   2,   1,   2,   1,   2,
   1,   4,   2,   1, 139,   6,   4,   1,   6,   1,  14,  20,   2,   1,   1,   9,   1,  13,   1,   6,   2,   9,   1,   2,   1,   2,   1,   2,   4,   6,
   1,   3,   1,   1,   2,  78,   2,  12,   3,   2,   2,   3,   1,   5,   2,   2,   2,   4,   7,   1,   1,   1,   7,   1,  13,   1,   1,   1,   6,   5,	
   1,   2,   6,   1,  15,   1,   1,   1,   2,   2,   7,   1,  44,   7,   1,   1,   6,   2,  35,   2,   4,   1,   1,   3,  10,   1,   2,   1,   1,   4,
   3,   2,   4,   4,   1,   2,  45,   1,   9,   1,   2,   3,   3,   1,  14,   2,   6,   3,   3,   1,  12,  19,   2,   1,   8,   1,   1,   1,   4,   1,
   2,   3,   2,   4,   1,   1,   1, 118,   1,   1,   1,   2,   2,   2,   4,   1,   3,   2,   1,   2,   1,   1,  56,   4,   3,   2,   1,   4,   2,   5,
   1,   4,   1,   4,  10,   2,  19,   1,   5,   1,   4,   5,   1,   1,   1,   1,   1,   1,   1,   5,   7,   1,   1,   1,   1,  26,   1,  30,   1,   2, 
  38,  19,   8,   1,  11,   9,   2,   1, 102,   1,   1,   2,   6,   6,   1,   1,   3,  13,  43,   2,   4,   1,   1,   1,   1,   9,   1,   1,   1,   3,
   1,   1,   2,   2,   1,   2,   2,   1,   1,   4,   3,   1,   7,   2,   1,   1,   2,   1,   3,   7,   2,   2,   6,   9,   2,   3,   1, 125,   1,   1,	
   1,   1,   1,   1,   3,   3,   6,   1,   6,  20,   9,   1,   2,   2,   3,   1,   3,   2,   1,   3,
\end{table}
\newpage
\begin{table}[H]
\caption{\label{tab:CFR-1422}$(\textrm{CFR})\quad t_\nu= (2\cdot1422^2)^2+\frac{\sqrt{2\cdot1422^2-1}}{2\cdot 1422^2}\quad$ ct'd}
   1,   1,   3,   2,   1,   1,   4,   1,   2,   1,
   1,   1,   1,  19,   2,   4,   1,  31,   2,   1,  54,   1,   4,   1,   8,   1,   1,   1,   1,   1,   2,   3,   1,   9,   1, 200,  27,  26, 306,   1,
   1,   1,  24,   1, 121,   2,   1,   3,   2,   2,   1,   2,   1,   6,   5,   1,   4,   1,   1,   4,   1,   1,   1,   1,   1,   4,   9,   1,  38,   1,	
 102,   1,  15,  95,   1,   2, 335,   1,   3,   3,   3,   1,   4,   2,   4, 
  1,   1,  34,   2,   2,   1,  11,   3,   3,   1,   3,   2,   2,   6,   3,
   1,   1,   1,   7, 123,   1,  11,   1,   5,   1,   3,   4,   1,   2,   1,   1,   2,  12,   1,   2,   4,   4,   8,   5,   2,   5,   1,  21,   2,   1,
   1,   4,   1,   1,  10,   5,   7,   2,   1,   4,   2,   1,   2,   9,   1,   2,   1,   5,   2,   6,   3,   1,   8,   2,   6,   3,  17,   1,   3,   1,	
   2,   8,   1,   1,   2,   2,   3,   2,   1,  61,   1,   2,   3,   1,   1,   2,   2,   2,   1,   1,   3,  11,   1,  22,   1,   1,  14,   4,   1,   1, 
  27,   1,   2,  10,  13,  12,   2,   2,   2,   2,   6,   8,   2,   1,   1,   1,   1,   4,   2,   3,   3,   2,   3,   2,   1,  42,   4,   3,   2,   7,
   1,   3,   2,   8, 322,   6,   4,   2,   3,   6,   2,   1,   3,   1,   1,   1,   2,   1,   1,  73,   1,  13,   2,   2,   1,  12,   1,   5,  10,   1,	
   2,   1,   4,   1,  23,   1,  15,   1,   1,   5,   4,   3,   1,   6,   2,   5,   6,   2,   1,   6,   4,   9,   4,   1,   5,   1,   1,   1,   1,   6,
   2,   2,   1,  94,   2,   1,   1,   1,   1,   1,   3,   8,  47,   6,   5,   1,   1,   1,   5,   3,   1, 755,   3,   1,   9,  19,   1,   1,   1,  15,
   2,   1,   5,   3,   1,  12,   2,  14,   4,  16,   1,   2,   1,   1,   1,   1,   1,   5,   1,   2,   5,   4,   1,   1,   1,   2,   1,   1,   1,   4
   1,   1,   1,   3,   1,   1,   1,   2,  12,   4,   1,   3,   1,   6,   4,   1,   2,   3,   2,   2,   2,   1,   1,  32,   1,   5,   4,   1,  13,   1,
   1,   2,   1,   3,   3,   9,   1,   6,  53,   2,   1,   1,   1,   1,   5,   7,   5,   3,   1,   2,1329,   8,   1,   1,  72,   4,  16,   1,   1,   1,
   1,   3,  74,   2, 182,   1,   9,   1,   5,   1,  14,   1,   1,   1,   1,   1,  13,   2,   1,   2,   1,   4,   3,   4,   2,   6,   3,   3,   8,   1,	
   1,   2,   7,   1,   3,   1,  39,   1,   5,   6,   1,   2,   1,   1,   9,   5,   6,   1,   9,   6,   1,   3,   1,   7,  11,   1,1263,  14,   1,   1,
   1,  10,   1,   1,   1,   1,   1,   1,   1,   9,  46,   2,  17,   1,  10,   2,   8,   1,   3,   2,   1,   1,   3,   5,   1,   1,   1,   5,   5,   3, 
  16,   2,   4,   7,   1,   6,   1,   2,   1,   1,   1,   3,   4,   1,   7,   1,   8,   3,   2,   1,   3,   4,   4,   3,   3, 166,   4,   1,   2,   1,
   2,   2,   6,   2,  80,   1,  25,   3,   1,  85,  28,  35,   1,   2,   4,   2,   1,   2,  11,   2,   4,   1,   8,   5,   2,  16,   2,   1,   4,   3,
   1,   2,   1,   3,   3,   1,   2,   1,   3,   1,   1,   7,   1,   1,   5,   1,   3,  13,   2,  10,   1,   1,   8,   1,   1,   1,   1,   1,   3,   3,
   4,   4,   1,   1,   9,   1,   1,   2,   7,   2,   1,   2,   3,  31,   2,   1,  58,   1,   5,   1,   1,   2,   1,   7,   2,  13,   1,  16,   5,  81,	
   1,   6,   1,   2,  15,   1,   1,   3,   1,  40,  17,   1,   4,  22,   5,  15,   1,  40,   1,   1,  28,   2,   1,   2,   1,   2,   1,   2,   1,   2,
   1,   1,   3,   9,   1,   5,   2,   5,  12,   5,   4,   1,   5,   2,   1,   1,   6,   1,  15,   2,   1,   1,   3,   1,   2,   2,  11,   5,   5,   1,
   3,   5,   5,   4,   2,   1,   1,  37,   1,   3,   1,   1,   2,  30,   3,   1,   2,   1,   3,   1,   2,   7,   4,   5,   9,  42,   1,   1,   9,   3,
   1,   2,   2,   1,   1,   3,   1,   1,   5,   1,   1,   1,   1,   1,   1,   1,  11,   3,   1,   4,   2,   1,   2,   1,   1,   5,   2,   1,   8,   1, 
  11,   1,   1,   1,   1,   1,   2,   7,   3,   3,   4,  23,  23,   1,   1,   1,   3,  10,   2,   2,   2,   1,   2,  14,   1,  11,  13,   8, 301,   1,
   2,   1,  28,   3,   4,   3,   1,  19,   1,   3,   5,   1,   3,   8,   2,   7,   2,   3,   3,   1,   1,   1,   1,   1,   6,  36,   1,   4,   1,   2,	
   2,   1,   6,   1,  42,   1,   1,   2,   1,   2,   3,   3,   1,   9,   2,   4,   1,   4,   1,   2,   1,   1,   5,   1,   1,   1,   2,  25,   2,   2, 
 110,  10,   1,   1,   2,   2,   3,   1,   2,   1,  11,   3,   1,   3,   1,  29,   7,   1,   1,   1,  18,   1,   2,   1,   2,   2,   4,   6,   1,   4,
   5,   1,   1,   1,  23,  10,   1,   1,   3,   1,   2,  24,   1,   2,   1,  29,   8,  24,   1,   1,   1,   3,   7,   1,   8,   1,   3,  31,   7,   3,	
   7,   2,   1,  53,   1,   4,   3,  17,   1,  17,   1,   9,   1,  14,   2,   3,   7,   1,  56,   1,   2,   1,  20,   2,   8,   4,   8,   1,   1,   1, 
  87,   1,   2,   2,   3,   1,   1,   2,   1,   3,   8,   7,   1,   3,   1,   1,   1,   9,  28,   1,   7,   1,  16,   1,   3,   1,   1,   2,   1,   6,
   1,   1,   3,   1,  96,  16,   1,   4,   6,   1,   1,   6,   3,   1,   1,   1,   3,   2,   7,   2,  32,   1,   1,   1,   4,   1,   3,   1,   1,   1,
   2,   1,   1,   1,  12,   2, \textcolor{red}{526},   3,   1,   5,   1,   5,   6,   1,   2,   1, 200,   2,  11,   2,   2,  44,   3,   1,   2,   5,  32,   9,   1,   4, 
   1,   2,   1,   2,   2,   1,   2,   1,   6,   4,   2,   4,   2,   19,   1,   3,   5,   71,   8,   2,   3,   2,   1,   14,   1,   3,   2,   14,   1, 
  17,   30,   3,   1,   6,   1,   1,   4,   2,   1,   8,   1,   2,   5,   3,   1,   1,   5,   1,   33,   1,   2,   3,   1,   4,   2,   1,   1,   10,  
 1,   3,   5,   4,   2,   1,   5,   4,   1,   1,   1,   1,   49,   68,   4,   8,   1,   1,   2,   1,   17,   1,   7,   1,   222,   5,   1,   1,   1,   
1,   1,   5,   4,   1,   3,   3,   2,   1,   1,   5,   13,   1,   29,   1,   1,   1,   1,   1,   8,   1,   3,   2,   11,   2,   6,   7,   9,   7,   1,  
 2,   2,   1,   5,   1,   1,   1,   1,   3,   24,   3,   1,   3,   1,   24,   10,   3,   1,   3,   7,   1,   3,   3,   1,   5,   1,   2,   1,   6,   1,  
 2,   3,   1,   1,   2,   1,   1,   8,   11,   8,   6,   9,   2,   4,   2,   22,   44,   1,   39,   1,   44,   1,   62,   2,   3,   4,   3,   17,   1,  
 1,   4,   1,   4,   1,   13,   2,   1,   1,   6,   122,   1,   7,   2,   6,   2,   4,   3,   2,   1,   4,   1,   1,   1,   1,   1,   2,   1,   2,   2,
   2,   4,   15,   3,   2,   1,   1,   1,   1,   1,   6,   1,   1,   1,   1,   12,   1,   1,   3,   17,   7,   2,   1,   1,   5,   1,   5,   2,   3,   4,
   5,   2,   10,   3,   3,   1,   1,   1,   7,   49,   1,   1,   57,   8,   1,   5,   1,   1,   2,   37,   3,   43,   1,   1,   1,   1,   10,   1,   1,
   1,   1,   31,   1,   1,   3,   3,   6,   2,   3,   2,   1,   1,   5,   1,   3,   1,   5,   2,   1,   1,   2,   19,   1,   1,   1,   19,   1,   7,   1,
   1,   1,   19,   1,   1,   1,   3,   2,   1,   3,   2,   5,   1,   2,  ...] 

\end{table}
\newpage
\begin{table}[H]
\caption{\label{tab:CFR-1422}$(\textrm{CFR})\quad1278= (2 x^2)^2-\frac{\sqrt{2 x^2-1}}{2 x^2}$}
{[4; 4, 2, 1, 1, 2, 2, 1, 1, 1, 2, 1, 8, 32, 3, 5, 1, 1, 1, 2, 2, 1, 42, 1, 2, 1, 6, 12, 1, 4, 1, 
    4, 2, 21, 1, 3, 1, 3, 1, 4, 7, 24, 3, 1, 1, 6, 2, 2, 3, 1, 3, 1, 8, 2, 1, 2, 1, 1, 1, 1, 3, 
    1, 1, 1, 17, 2, 1, 3, 2, 1, 1, 1, 1, 6, 1, 1, 6, 2, 1, 1, 2, 1, 6, 2, 1, 3, 4, 3, 1, 1, 2,
    1, 2, 10, 2, 2, 1, 4, 1, 1, 6, 49, 1, 18, 1, 1, 1, 2, 1, 2, 1, 1, 2, 1, 1, 10, 3, 1, 1,178, 6,
    1, 2, 1, 5, 1, 1, 2, 8, 9, 28, 2, 1, 5, 1, 2, 1, 3, 1, 1, 28,248, 50, 19, 2, 1, 1, 2, 1, 2, 6,
    5, 1, 3, 1, 98, 9, 1, 1, 3, 80, 1, 1, 1, 1, 1, 3, 1, 3, 4, 52, 6, 5, 13, 2, 8, 1, 3, 25, 1, 1, 
   54, 1, 3, 2, 4, 1, 8, 3, 1, 1, 4, 2, 1, 3, 40, 1, 1, 16, 1, 2, 3, 1, 9, 4, 1, 9, 7, 2, 2, 9,
    1, 1, 6, 1, 2, 1, 4, 4, 1, 1, 1, 3, 1, 11, 1, 1, 2, 1, 4, 3, 7, 1, 5, 1, 1, 2, 2, 1, 2, 1,
    2, 6, 7, 1, 2, 3, 1, 2, 3, 2, 48, 1, 1, 2, 3, 1, 3, 1,100, 1, 1, 1, 1, 3, 1, 25, 3, 9, 6, 1,
    2, 2, 5, 1, 2, 13, 5, 1, 2, 3, 15, 1, 31, 2, 82, 1, 2, 7, 1, 1, 3, 4, 1, 1, 1, 3, 1, 18, 3, 1,
    1, 12, 11, 1, 14, 1, 3, 1,\textcolor{red}{287}, 22, 6, 6, 5, 6, 15, 1, 3, 3, 5, 12, 1, 3, 44, 3, 1, 1, 9, 3, 1, 6,
   17, 1, 1, 5, 1, 10, 1, 4, 1, 1, 1, 2, 2, 1, 1, 13, 6, 1, 1, 1, 11, 4, 1, 1, 6, 1, 1, 3, 1, 3, 1, 3, 1, 1, 6, 2, 5, 9, 1, 1, 1, 1, 2, 1, 8, 1, 1, 20, 1, 1, 1, 1, 1, 2, 1, 1, 1, 5, 1, 16, 1, 1, 1, 1, 1, 24, 32, 1, 9, 6, 2, 1, 1, 2, 13, 3, 3, 1, 3, 1, 2, 1, 3, 49, 1, 18, 6, 10, 1, 3, 1, 31, 1, 1, 1, 29, 1, 1, 1, 29, 4, 1, 43, 1, 2, 4, 5, 583, 1, 1, 1, 2, 1, 3, 1, 1, 1, 4, 1, 4, 2, 1, 1, 1, 1, 6, 1, 1, 3, 8, 8, 2, 6, 1, 2, 1, 1, 1, 3, 1, 9, 1, 1, 4, 1, 1, 1, 2, 9, 1, 1, 1, 12, 1, 25, 5, 1, 1, 3, 2, 1, 1, 1, 1, 1, 1, 3, 3, 1, ...]
}
\end{table} 
\begin{table}[H]
\caption{\label{tab:CFR-1422}$(\textrm{CFR})\quad1422= (2 x^2)^2-\frac{\sqrt{2 x^2-1}}{2 x^2}$}
 {[4;  
  2, 1, 11, 1, 2, 3, 1, 18, 1, 1, 1,3575, 1, 1,489, 14, 1, 7, 2, 2, 2, 3, 5, 1, 5, 4, 2, 1, 2, 2,
  1, 4, 1, 1, 11, 2, 1, 1, 1, 1, 1, \textcolor{red}{526}, 9, 2, 21, 1, 1, 3, 3, 2, 3, 3, 2, 4, 3, 2, 1, 1, 1, 1, 
 21, 3, 1, 14, 17, 2, 1, 2, 4, 2, 5,  1, 1, 1, 15, 2, 1, 4, 1, 3, 1, 2, 1, 7, 1, 2, 2, 1, 1, 1, 
 36, 1, 2, 1, 1, 10, 4, 1, 57, 2, 1, 2, 11, 1, 1, 1, 1, 2, 4, 1, 1, 1, 23, 2, 1, 5, 7, 1, 10, 64,
  1, 6, 1, 2, 29, 1, 5, 1, 1, 1, 2, 6,  1, 1, 3, 2, 4, 2, 3, 52, 2, 1, 2, 1, 2, 1, 1, 2, 1, 53, 
395, 1, 1, 11, 6, 1, 2, 1, 1, 9, 7, 4,  1, 2, 2, 6, 1, 1, 7, 1, 1, 1, 31, 1, 1, 12, 1, 3, 1, 1,
  1, 6, 47, 2, 1, 30, 1, 15, 1, 1, 6, 18, 23, 4, 5, 1, 1, 1, 1, 2, 2, 1, 16, 1, 11, 11, 1, 1, 1, 8,
  6, 77, 2, 1, 4, 4, 1, 1, 2, 6, 6, 1, 1, 11, 87, 1,112, 10, 6, 1, 1, 2, 1, 4, 1, 3, 1, 6, 2, 3, 
 34, 1, 3, 1, 1, 15, 1, 16, 2, 2, 1,170, 3,  1, 2, 1, 1, 2, 31, 3, 4, 3, 1, 2, 1, 34, 1, 63, 1, 27, 
 59, 28, 1, 2, 1, 19, 29, 1, 6, 1, 6, 5,475,  1, 1, 1, 1, 1, 3, 2, 20, 6, 21, 2, 1, 1, 1, 1, 13, 1,
  1, 1, 2, 1, 4, 1, 1, 4, 2, 6,209, 3, 1,  1, 2, 1, 1,146, 1, 1, 1, 1, 1, 4, 3, 1, 1, 30, 1, 5,
  3, 2, 1, 4, 2, 21, 1, 1, 1, 4, 1, 3, 8,  7, 2, 3, 1, 3, 1, 8, 1, 8, 2, 12, 2, 1, 1, 2, 4, 2,
  1, 1, 3, 1, 4, 45, 3, 29, 2, 1, 2, 1, 3,  9, 1, 1, 3, 9, 2, 1, 1, 1, 3, 1, 6, 7, 1, 4, 2, 1,
  2, 62, 1, 2, 1, 16, 1, 1, 1, 1, 25, 1, 2, 1,  1, 5, 1, 2, 24, 4, 3, 23, 2, 12, 3, 1, 13, 1, 1, 4,
  1, 1, 1, 1, 1, 1, 5, 1, 1, 3, 1, 5, 2, 4,  1,279, 1, 3, 1, 1, 2, 1, 4, 1, 1, 12, 10, 2, 2, 2,
  1, 30, 1, 2, 1, 1, 4, 3, 1, 4, 1, 2, 1, 1,  6, 4, 1, 6, 1, 29, 1, 1, 1, 2, 5, 1, 2, 1, 7, 1,
  7, 1, 3, 4, 1, 3, 7, 1, 1, 1, 1, 5, 1, 4, 49, 1, 2, 9, 2, ...] 
}
\end{table}
\newpage
\begin{table}[H]
\caption{\label{tab:CFR-373}$(\textrm{CFR})\quad t_\nu=(2\cdot 373^2)^2+\frac{\sqrt{2\cdot 373^2-1}}{2\cdot 373^2}$}
{[77427514564; 

\textcolor{red} {527},   1,   1,  95,   2,   2,   3,   1,  17,   1,   2,   1,   3,   1,   2,   1,   2,   1,  14,   3,   6,  10,   2,   1,   4,   1,   7,   1,   6,   1,
   4,   7,   1,   2,   4,   3,   2,   3,   1,   1,   9,   1,   1,  28,   5,   1, 109,   1,   1,   1,   1,   2,   1,   2,  40,   3,   1,   4,   1,  94,
   3,   1,   3,   3,   1,  23,  30,   1,  17,   33,   1,   6,   1,  1,   1,   1,   1, 115,   1,   1,   1,   1,   1,   2,   2,   1,   1,   2,   7,   1,
   2,   1,   1,   6,   1,   4,   1,   4,   1,   1,   3,   1,   1,  30,   7,   1,   2,   3,   1,  12,   4,   7,   1,   5,   1,  12,   3,   3,   2,   1,
   4,   1,  11,   1,   6,   3,   1,   1,   4,   3,   2,   1,   6,  27,   1,   1,   2,   2,   1,   3,  28,   1,   4,   1,   6,   1,   3,   2,   1,   2,
   2,   2,   1,   2,   1,   1,   1,   1,   4,   6,   1,   2,   2,   6,   3,   1,   1,   6,   1, 277,   1,   1,   1,   1,   7,   4,  10,   1,   1,   3,	
   1,   5,   9,   2,   1,   1, 521,   2,   1,   2,  14,   1,   1,  13,   2,  59,   3,   1,  18,   1,  11,   1,   2,   2,   3,   1,   1,   5,   6,  95,
   4,  38,   2,   1,  90,   1,   2,   1,   4,   3,   1,   5,   2,   4,   1,   1,   2,   2,   2,   1,   2,   3,   1,   1,   1,  13,   5,   3,   1,   990
  55,  34,   1,  22,   1,   2,   5,   1,  12,   1,   1,   6,   3,  10,   1,   1,   4,   3,   1,   1,   9,   1,   3,  21,   6,   1,   3,   2,   8,   2,	
  18,   1,   3,   5,   1,   1,   7,   5,   2,   9,   1,   3,   1,   1,   2,   3,   1,   3,   1,   7,   1,   3,   2,   5,   2,  59,   1,   2,   9,   2,
   3,   1,   1,   9,   2,  23,   1,   3,   2,   2,   1,   2,   6,   9,   1,  27, 353,   1,   1,   7,   1,   1,   2,   2,   1,   3,   6, 815,   1,   2,
   1,   2,11088,  1,   3,   3,   5,  11,   2,   1,   1,   3,  20,  20,   1,   1,   2,   3,   3,   1,   1,   1,   1,   4,   2,   2,  15,   1,   1,   4,	
  24,  12,   3,   6,   3,   1,  37,  20,   2,   4,  11,   1,   1,   1,   1,   5,   1,   2,   9,   2,  10,   2,   2, 137,   8,   1,   1,   1,   1,   6,
   3,   3,  11, 612,   1,   1,   3,   1,   1,  11,  20,  14,   1,   9,   1,   1,   4,   4,   1,   1,   7,   1,   3,   4,   1,   9,   6,   1,   2,   3,
   1,   1,   1,   1,   1,   1,  12,   1,   3,   9,   5,   1,   2,   2,   4,   1,   1,   8,  44,   1,   2,   1,   2,   9,   1,   1,   3,   1,   1,   1,	
   1,   1,   1,  11,   9,  15,   3,   2,   6,   6,  16,  12,   1,   1,   2,   6,  16,   3,   1,   1,   8,   1,   1,   1,   4,   4,   3,   1,   2,   1,
   4,   4,   3,   1,   8,   1,   1,   1,   1,   3,   1,   4,   1, 320,   1,   3,   1,   1,   1,  46, 135,   1,   5,   1,   2,   3,   1,   2,   3,   2,
   3,   2,   2,   1,   3,   3,   2,   1,   2,   1,   2,   1,   1,   2,   1,   1,   1,   1,  11,   1,   4,   4,   2,   1,   1,   1,   1,   1,   2,  13,	
   1,   2,  10,  19,   1,   2,   3,   4,   2,   1,   2, 652,   1,   3,   1,  20,   1,   1,  22,   9,   1,   1,   2,   1,   1,   8,   1,   7,   3,   1,
   7,  55,   5,   8,  40,   3,   2,   2,   1,   2,   1,   3,   1,   3,   3,   1,  23,   1,   2,   2,  12,   2,   1,   1,   1,   1,   2,   2,   2,  37,
   1,   5,   4,   2,   1,   2,   1,   5,   1,   7,   1,   1,   4,   1,   1,  10,   1,  15,   3,   2,   2,   1,   1,  64,   3,   3,   1,   1,   5,   6,	  
  21,   1,   1,   1,   1,   1,   2,   1, 105,   5,   1,   2,   1,  24,   2,   2,   1,   1,   1,   1, 105,   1,  49,   2,   5,   2,   2,   1,  26,  43,
   6,   1,   1,   1,   5,   1,   1,  16,  60,   1,   1,  54,   1,   9,   2,   5,   1,  77,   3,  15,   1,   6,   3,   1,   7,  20,   1,   1,   3,   1, 
  15,   3,   4,   2,   1,   3,   1,   7,   2,   4,   1,   9,   5,   1,   3,   1,   1,   2,   1,   2,   6,   3,   1,   5, 161,   1,   1,  30,  15,   1,	
   1,   4,   1,   2,   3,   3,   7,   1,   1,   1, 146,   1,   4,   1,   3,   6,   1,   1,   5,   1,   1,  13,  75,   2,   1,   4,   2,  41,  20,   1,
   3,   1,   3,   2,   1,   1,   3,   2,   2,   1,   1,   1,   1,  24,   1,  26,   1,   2,   2,   1,   2,   1,   1,   1,  10,   1,   1,   2,   3,   1,
   1,   1,   1,   9,   1,  10,   1,   1,   3,  11,  12,   5,   2,   2,  45,   2,   2,  11,  33,   1,   5,   3,   1, 206,   5,  16,   2,   3,   7,   1,	
   1,   3,   1,   1,   5,   1,   3,   3,   1,   1,   6,   2,   5,   7,   1,  10,   1,   4,   1,   2,   2,  40,   4,   1,   5,   1,   1,   2,   1,   1,  
 344,   1,   5,   7,   5,   2,   2,   2,   2,  37,   3,   2,   2,   5,   1,   4,   3,   1,   2,   3,   1,   3,   1,   1,  24,   1,   7,   2,   7,   3,
   1,  16,   2,   1,   9,   4,   1,  30,  22,   2,  11,   5,  21,   2,   1,   2,   1,   1,   2, 235,   1,   7,   1,   9,   2,   7,  13,   4,  12,   1,	
  19,  33,  36,   1,   1,   2,   2,   3,   4,   2,   2,   1,  24,   2,   2,   4,   9,   1,   6,   1,   1,   1,   1,   3,   1,   2,   6,   3,  13,   1, 
  17,   4,   2,   2,   1,   1,   1,  14,   2,   4,   3,   1,   1,   1,  18,   5,   1,   1,   1,   2,   1,   1,   9,   5,   1,   1,   2,   1,   1,   1,
   3,   2, 533,   1,   4,   1,   1,  15,   1,   2,   1,  21,   1,   1,   3,   1,   3,   1,   4,   2,   1,  11,   1,   1,   1,   2,   1,   1,  22,   6,	
   2,   7,   3,   6,   3,   1,   5,   2,   5,   1,   9,   3,   1,   2,   1,   8,   6,   1,   9,   1,   8,  14,  14,   1,   2,   2,   1,   1,   1,  10,
   1,   4,  11,   2,   6,   1,   1,   2,   3,   1,  13,   1,   2,   1,   2,   1,  14,   2,   4,   1,   3,   2,   1,   3,   5,   1,   1,  69,   6,  12,
   1,   3,   2,   1,   1,   2,   2,   4,   1,   2,   2,   1,   6,   2,   1,   1,   1,   5,   1,  10,   2,  20,   3,  18,   1,  12,  14,   2,   1,   1,	
  13,   2,   2,   1,   6,   4,   1,  40,   1,   2,   2,   2,  11,   2,   2,   2,   5,   2,   3,   8,   1, 164,   1,   1,   1,   1,   1,   7,   2,   2, 
  13,   4,   3,   4,   5,  11,   1,   2,  36,   1,   2,   2,  23,   1,  31,  16,   1,   1,   3,   1,   4,   1,   1,   1,   1,   1,   5,   2,   8,   1,
   3,   2, 137,  12,   1,   1,   2,   1,  30,   1,   1,   3,   5,   5,   1,   1,  1, 37,  
 6,   4,  17,   5,   4,   1,   1,   4,   1,   4,   2,   1,  1,   4,   1,  10,  10,   1,   2,   3,   2,  10,   1,   2,   1,   1,   1,   7,   1,   2,  3,   1,   1,   3,   1,   3,
   }
   \end{table}
 \newpage
\begin{table}[H]
\caption{\label{tab:CFR-373c}$(\textrm{CFR})\quad t_\nu= (2\cdot 373^2)^2+\frac{\sqrt{2\cdot 373^2-1}}{2\cdot 373^2}$\quad ct'd}  
{     15,   3,   4,   2,  48,   1,
   4,  10,   2,   1,   3,   1,  23,   1,   1,   8,   1,  17,   1,   2,   5,   5,   1,   3,   3,   1,   6,   1,   2,  25,   3,   1,   6,   3,   6,  18,
   2,  16,   3,   5,   6,   2,   2,   2,   3,   1,   2,   3,  41,  12,   3,  13,   2,   1,   2,   3,   5,   1,   1,   7,   1,   1,   1,   1,   2,   5,	
   1,   1,  25,   2,   1,  10,   1,   6,   5,   5,   2, 250,   2,  24,   1,   1,   1,   5,   1,   4,   4,   5,   2,   2,   1,   1,   6,   1,   4,   2,
   1,  11,   2,  13,   1,   1,   2,   1,   2,   4,   1,  41,   1,   1,   6,  15,  15,   4,   1,   3,   8,   1,   7,  16,   1,   1,  12,  23,   1,   1,
   5,   1,   1,   1,  10,   1,   1,   1,   2,   7,   1,   5,   2,   9,   1,   3,   3,   1,   3,  13,   1,  23,   2,  83,   1,   3,   2,   6,   2,   4,
   7,   2,   2,   1,   1,   5,   1,   1,   2,   1,   1,   2,  22,   1,  28,   3,   2,  10,  25,   1,   2,   1,  41,   4,   4,   1,   8,   1,   4,   1,
   1,   1,   1,  21,   1,   1,   1,   5,   3,  64,   2,   1,   2,  24,   1,   2,   4,   1,  20,   1,   1,   1,   1, 241,   6,   1,   2,   1,   1,  59,
   7,   1,   5,   1,   3,   1,   2,  20,   2, 109,   1,   4,   2,   1,   1,   1,   7,   9,  38,   1,  12,   1,   1,   1,   1,   1,   1,   1,   1,   3,	
   2,   7,   1,   1,  63,   1,   1,   1,   6,   1,   2,   3,   1,  29,   3,   1,   2,   3,   2,   1,   2,   3,   3,   2,   1,   1,   2,   1,   4,   5,
   1,   9,   2,   4,   2,  20,   4,   6,  12,   1,   2,   8,   7, 344,   2,   1,   1,   1,   1,   2,   3,  33,  67,  11,   1,   2,  22,   1,  11,   4,
   2,   1,  11,   1,   3,   2,   1,  47,   2,   8,   2,  10,   1,   2,   3,   1,   3,   1,   5,   4,   1,  14,   2,   1,   1,  25,   1,   5,   1,   8,	
   5,   2,   1,   4,   9,  10,   4,   1,   6,   5,   1,   1,   1,   3,   1,   6,   1,  40,   1,   1,   1,   1,   3,   1,   3,   1,   1,   9,   2,  10,
   2,   3,   2,   3,   5,   1,   3,   1,  15,   1,   3,  36,   1,   2,   2,   1,   1,   5,  32,   1,   2,  35,   3,  43,   1,   1, 175,   3,  19,   1,  
 111,   1,  37,   1,  21,   2,   3,   1, 193,   1,  12,   4,   3,   4,   3,   1,   1,   1,   1,   2,   2,   1,   2,   2,   8,   1,   1,   1,   1,   1,	
   2,   1,   1,   1,   1,   3,   1,  13,   1,   4,   2,  14,   1,   1,   1,   2,  13,   1,   1,   5,   1,   1,   6,   7,   3,  17,   5,   2,   3,   1,
   1,  75,   1,  27,   7,   2,   1,   1,   1,   4,   5,   1,   3,   2,   6,   1,  15,  55,   7,   6,   1,   1,   2,   1,   1,  57,   2,   2,   5,   1,
   1,   3,   1,   1,   1,   1,   1,   1,   2,   1,   2,   3,   1,   6,   1,   1,   1,   1,   1,   4,   1,   2,   1,   6,  17,  13,  22,   1,   8,   2,	
   3,  11,   2,  15,   2,   1,   7,   2,   1,   5,   2,   1,   1,   3,   7,   1,   1,   1,   1,   2,   1,   2,  15,   1,   1,   5,   1,   5,   1,   5,
   1,   1,   2,   1,   1,   1,   4,   1,   6,   1,   1, 397,   1,  43,   1,  14,  46,   1,   2,   3,   1,   7,   3,  19,   2,   2,   1,  15,   3,  18, 
  15,   1,   3,   1,   9,   1,   2,   2,  31,   3,   5,   3,  18,   1,  10,   1,  37,   4,   1,   6,   1,   1,   1,   1,   3,   1,   1,   3,   1,   2,	
   9,   4,   1,   6,  13,   2,   4,   2,   2,   1,   1,   2,   6,   1,   2,   1,   2,   1,   3,   1,   1,   2,   1,   1,   3,   1,   5,   2,   1,  12,
   2,   3,   1,   1,   1,  85,   1,   6,   2,   1,  12,   3,   7,   1,   1,   3,   2, 207,   1,   1,   1,   1,   2,   3,  10,   3,   1,   2,   1,   2,
   1,   2,   1,   4,  18,   1,   1,   1,   1,   6,   1,  11,   1,  15,   1,   1,   1,   7,   3,   1,  43,   7,   1,   8,   1,  10,   5,   5,   2,   2,	
   2,  12,   2,  39,   1,   3,   1,   1,   4,  34,   1,  33,   1,   1,   2,   1,   1,   1,   5,   1,   1, 377,   1,   1,   1,   1,   1,   3,   2,   6,
   1,   3,   2,   4,   3,   6,   6,   1,   7,   1,   1,   1,   1,   7,   1,   1,   5,   1,   1,   3,   1,   4,   1,  11,  33,   2,   4,   1,   5,   1,
   1,   2,   1,   2,   4,   1,   7,   1,   2,   1,   2,   1,   2,   1,   4,   3,   1,   1,  25,   5,   1,  95,   1,   2,   1,   3,   1,   1,   1,   1,	
   1,   1,   2,   6,   5,   1,   1,  16,   1,  30,   1,  13,   5,   1,   5,   1,   1,   1,   2,   2,   6,   5,   5,   3,   4,   9,   2,  11,   1,  15,
   9,   3,   1,   2,  18,   1,   1,   3,   1,   2,   1,   2,   8,   1,   1,   4,   9,   3,   8,   7,   1,   6,   6,   2,  18,   1,   1,   5,   1,   1,
   1,   3,   3,   2,   3,   1,   1,   1,  10,   1,   5,   8,   1,   1,   1,   5,   1,   3,   1,   1,   1,   3,   1,   1,   3,   1, 235,   1,   2,   1,	
   8,   1,  23,   1,   3,   3,   3,   8,   1,   1,   3,   5,   1,   3,   5,   5,  11,   1,   1,   1,   1,   1,  13,   2,   1,   1,   3,   1,  25,   4,
   1,  10,   1,   3,   3,   1,   7,   1,   2,   2,   3,   3,   2,   2,   1,   1,  11,   4,   6,   9,   6,   3,   1, 120,  34,   1,   1,  46,  10,   1,
   4,   2,   6,   1,   1,   1,  10,   1,   4,   3,   4,   1,   1,   3,   1,   1,   3,   3,   2,   1,   1,   2,   1,   2,   4,   6,   2,   2,   1,   6,	
   1,   2,   2,   9,   2,   1, 118,   1,   3,   1,   1,   1,   1,   2,   1,   2,   1,   2,   3,   2,   3,   1,   1,   4,  12,   1,   1,   2,   1,   2,
   1,   1,   3,   8,   1,   9,   3,   1,   5,   1,   3,   4,   3,   1,   8,   1,   6,   4,  54,  12,   9,   3,   1,   6,   1,   2,   1,   1,   1,   1,
   1,   3,   2,   4,   1,  21,   1,   1,   5,   1,   3,   1,   3,   1,   1,  10,   5,   1,   3,   1,   1,   6,  30,   1,  10,   2,   1,   1,   1,   7,	
  12,   1,   2,   1,   3,   2,   1,   1,   3,   1,   1,  10,   1,   1,   4,   1,  17,   1,   1,   1,   1,   1,   2, 110,   1,   3,   2,   1,   5,  42,
   5,   1,   4,  33,   2,   1,   3,   1,   1,   2,   6,   1,   1,  56,   1,   1,   1,   5,   1,  80,   1,   1,  21,   1,   7,   6,   1,   2,   5,   1,
   1,   2,   3,   8,   1,  23,   5,   1,  68,   1,   2,   2,  13,  44,   1,   4,   4,   2,   1,\textcolor{red} {285},20221,  1,  59,   1,   2,   3,   1,   5,   2,  28,	
  18,  20,   1, 400,   1,   1,   1,   1,   1,   5,   1,   1,   1, 145,   2,   2,   2,   2,   1,  18,   4,  10,  10,  12,   2,   5,   1,   4,  50,   1, 
  17,   3,   1,   1,   9,   12,   9,   2,   2,   4,   3,   2,   1,   1,   8,   8,   3,   1,   1,   42,   2,   1,   6,   8,   1,   1,   2,   6,   5,   1,   3,   1,   1,   38,   2,   5,   3,   1,   39,   1,   6,   1,   183,   16,   1,   1,   5,   1,   2,   1,   1,   9,   11,   2,   1,   1,   20,   227,   3,   6,   3,   1,   6,   2,   3,   1,   6,   6,   23,   1,   33,   6,   1,   1,   11,   1,   1,   2,   10,   5,   4,   1,   37,   1,   5,   24,   1,   1,   1,   2,   1,   1,   1,   1,   1,   9,   3,   10,   1,   2,   2,   25,   1,   1,   11,   2,   1,   3,   2,   1,   8,   3,   48,   1, 5;	   4,   5,   1,   5,   1,  10,   2,   1,   3,   1,  96,   1,   2,   1,   5,   2,   1,   1,   5,   1,   5,  13,  ...]
}
\end{table}
\begin{table}[H]
\caption{\label{tab:CFR-3015}$(\textrm{CFR})\quad3015= (2 x^2)^2-\frac{\sqrt{2 x^2-1}}{2 x^2}$}
{[  5;		
   4,   5,   1,   5,   1,  10,   2,   1,   3,   1,  96,   1,   2,   1,   5,   2,   1,   1,   5,   1,   5,  13, 9,   3,   6,   2,   4,   2,   1,   1,
   2,   1,  24,   1,   1,   2,   1,   4,   1,   3,   2,   6,   1,   1,  40,   1,  14,   2,   1,   5,   1,  11,  10,   1,   1,   1,   1,   1,   1,   4,
   1,   1,   2,   1,   3,   1,   1,   1,   6,   2,   3,   1,   2,   4,   3,   1,   7,   5,   2,   2,   1,   3, 105,   2,   3,   7,   1,   1,   2,   4, 
  10,   2,   4,   2,   1,   1,   1,  50,   1,   3,   2,   2,  10,   2,  43,   2,   6,   9,   1,   1,   2,   2,   2,   2,  18,   1,   3,   2,   3,   3,
   1,   2,   7,   1,   3,  55,   4,   1,   1,   6,   4,   1,  12,   1,   5,   3,   2,   1,   2,  36,   1,   1,   5,   1,   1,   1,   3,   3,   1,   1, 
  33,   1,   4,   2,   1,   1,   1,   5,   1,  13,   2,   1,   1,   9,   1,   1,  10,   1,   2,   2,   2,  26,  31,   1,   1,   1,   1,   1,  36,   2, 
  17,   7,   4,   2,   1,   3,   1,   3,   1,   1,   8,  21,   2, 102,   1,  19,   1,   1,   3,   1,   1,   1,   5,   8,   2,   5,   2,   3,   4,   6, 
  14,   1,  12,   2,   2,   4,   1,   2,   1,   1,   3,  42,   3,   8,   5,   1,   1,  14,  23,   1,   1,  72,   1,   1,   1,   2,   7,   1,  12,   1,  
 110,   1,   5,   1,   5,   5,   1,   5,   1,  12,   2,   1,   4,   1,   1,   1,   1,   1,   1,   1,   1,   1,   1,   1,   2,  15,   1,   1,   1,   3,
   2,   4,   1,  20,   1,   1,   5,   1,   1,   2,   5,   2,   1,  29,   1,  16,   1,   1,   5,   4,   6,   1,  10,   1,   3,   9,   1,   2,   1,   5,
   3,   2,   1,   1,   1,   3,   2,   1,   1,   3,   6,   1,   1,   9,   2,   3,  24,   4,   1,   1,  12,   1,   1,   8,   3,   3,   1,   7,   1,   1,
   1,   3,   1,   1,   2,   1,  51,   1,   2,   1,  16,   4,   1,  16,   3,   1,   1,   1,   4,  10,   3,   1,  26,   1,   1,   3,   1,   3,   1,  19, 
  27,   1,   4,  21,   2,   2,   1,   5,   2,   1,   7,   6,   6,   1,   1,   1,  17,   2,   2,  26,  19,   1,   4,   6,   1,  15,   1,  11,   1, 225,
   7,   1,   1,   2,   1,   1,   7,   1,   7,  24,  96,  42,   5,   1,   1,   1,   2,   1,   2,   1,   2,   1,   1, 107,   2,   7,   1,   1,   6,   1,
   3,   1,   7,   7,   2,   1,  27,   1,  24,   1,   2,   4,   3,   1,   1,   2,   6,   2,  13,   1,   1,   5,   1,   1,   1,   1, 394,   6,   2,   7,
   1,   1,   1, \textcolor{red}{286},   2,   6,   1,  16, 169,   1,  22,   1,   8,   7,   2,   1,   4,   8,   3,   2,   4,   2, 139,   1,   2, 279,   4,   1,   3,   2, 
  23,   1,  10,   1,   8, 296,   2,   7,   2,  65,   1,  17,   3,   1,   1,  40,   3,   1,   9, ...]
}
\end{table} 

\medskip {}

\noindent \label{transpoin}Transactional spacetime $\big(\{T\rightarrow X\}_\textrm{cc},\{X\rightarrow T\}_\textrm{pc},\{X\rightarrow T\}_\textrm{pc},\{X\rightarrow T\}_\textrm{pc}\big)$ associated with the Pythagorean triple $373^2=252^2+275^2$ and the combined system $(285,286,287,525,526,527)$:\newline
\[\begin{array}{ll}
\bar\theta^{-1}(1422)[42]=526,&\bar\theta^{-1}(1278)[309]=287,\\
\theta(373)[1]=527,&\bar\theta(372)[2]=525,  \\\theta(373)[2330]=285,& \bar\theta(373)[2331]=285, \\\theta(203)[1]=287, &\bar\theta(203)[2]=286, 
\end{array}\]
\[\begin{array}{c}\big(\bar\alpha_{526}+\bar\alpha_{287}\big)+
\big(\beta_{527}+\bar\beta_{525}+
\beta_{285}+\bar\beta_{285}+
\beta_{287}+\bar\beta_{286}\big)=5018, \\
5018\equiv96
\;(\textrm{mod}\,\,107) \;\quad
96\equiv43
\;(\textrm{mod}\,\,53).\qquad\qquad \end{array}\]

\medskip {}

\noindent Transactional spacetime $\big(\{T\rightarrow X\}_\textrm{cc},\{X\rightarrow T\}_\textrm{pc},\{X\rightarrow T\}_\textrm{pc},\{X\rightarrow T\}_\textrm{pc}\big)$ associated with the Pythagorean triple $3257^2=1232^2+3015^2$ and the combined system $(286,526,4606)$:\newline
\[\begin{array}{ll}
\bar\theta^{-1}(3015)[454]=286,&\bar\theta^{-1}(1422)[42]=526,\\
\theta(1422)[386]=286,&\bar\theta(1422)[387]=286,  \\\theta(1422)[2167]=526,& \bar\theta(1422)[2168]=526, \\\theta(3257)[1]=4606, &\bar\theta(3258)[2]=4606, 
\end{array}\]
\[\begin{array}{c}\big(\bar\alpha_{286}+\bar\alpha_{526}\big)+
\big(\beta_{286}+\bar\beta_{286}+
\beta_{526}+\bar\beta_{526}+
\beta_{4606}+\bar\beta_{4606}\big)=5607, \\
\;\;5607\equiv43
\;(\textrm{mod}\,\,107) \;\;\quad
43\equiv43
\;(\textrm{mod}\,\,53). \qquad\qquad \end{array}\]
\smallskip{}

\noindent Regarding life as consisting of a standby process 
activated by some amplitude $\bar\theta^{-1}(t_\textrm{act})[\,]$, provided a preselected MF emulation's resemblance to the  MF-wise  exposition of $\bar\theta^{-1}(t_\textrm{act})[\,]$ surpasses a  specifiable threshold, a case in point is $t_\textrm{act}=1823$. It is the 281{\small th} prime number, and 281 in turn is a Sophie Germain prime number as well as  hypothenuse of the primitive Pythagorean triangle $281^2=160^2+231^2$. With 1823, we get $\bar\theta^{-1}(1823)[126]=240$ and may build from it the MF-wise exposition  $\,\bar{B}_{t_\textrm{act}}^{(z)}:=2^{-z}\bar\theta^{-1}(1823)\;(z\in\mathbb{Z})$,

\bigskip{}

\label{actpoin}\renewcommand{\arraystretch}{1.2}
\hspace*{1cm}\begin{tabular}{c}
$\bar B_{t_\textrm{act}}^{(1)}[\,]=480$\tabularnewline
$\bar B_{t_\textrm{act}}^{(0)}[\,]=240$\quad
$\bar B_{t_\textrm{act}}^{(2)}[\,]=239$\tabularnewline
$\bar B_{t_\textrm{act}}^{(-1)}[\,]=120$\quad\quad\quad
$\bar B_{t_\textrm{act}}^{(3)}[\,]=119$\tabularnewline
$\bar B_{t_\textrm{act}}^{(-2)}[\,]=59$\quad\quad\quad\quad\quad
$\bar B_{t_\textrm{act}}^{(4)}[\,]=59$\tabularnewline
$\bar B_{t_\textrm{act}}^{(-3)}[\,]=29$\quad\quad\quad\quad\quad\quad\quad
$\bar B_{t_\textrm{act}}^{(5)}[\,]=29$\tabularnewline
$\bar B_{t_\textrm{act}}^{(-4)}[\,]=14$\quad\quad\quad\quad\quad\quad\quad\quad\quad
$\bar B_{t_\textrm{act}}^{(6)}[\,]=14$\tabularnewline
$\cdots$\qquad\quad\quad\quad\quad\quad\quad\quad\quad\quad\quad\quad\quad
$\cdots$\tabularnewline

\end{tabular}

\bigskip{}

\noindent to be met 
by the  $\{X\rightarrow T\}_\textrm{pc}$ emulation

\bigskip{}

\begin{tabular}{c}
$\theta(340)[1 ]=\textcolor{black}{480}\,(0)$   \tabularnewline
$\theta(\textcolor{black}{170})[1]=240\,(0)$\qquad $\bar\theta(170)[2]=\textcolor{black}{239}\,(0)$\tabularnewline
$\theta(85)[1]=120\,(0)$\qquad\quad\quad\quad
$\bar\theta(85)[2 ]=\textcolor{black}{119}\,(0)$\tabularnewline
$\theta(42)[1 ]=59\,(0)$\qquad\qquad\qquad\quad\quad
$\bar\theta(42)[2]=\textcolor{black}{58}\,(-1)$\tabularnewline
$\bar\theta(21)[1 ]=\textcolor{black}{29}\,(0)$\qquad\qquad\qquad\quad\quad\quad\quad
$\bar\theta(21)[2 ]=\textcolor{black}{28}\,(-1)$\tabularnewline
$\theta(\textcolor{black}{10})[1]=14\,(0)$\qquad\qquad\qquad\quad\quad\quad\quad\quad\quad
$\bar\theta(10)[2]=\textcolor{black}{13}\,(-1)$\tabularnewline
$\dots$ \qquad \qquad \qquad \qquad \qquad \qquad \qquad \qquad \qquad \qquad$\dots$
\end{tabular}

\bigskip{}

\noindent With allowance for deviations $\pm1$ across a level,      there's no upper cut-off  and no lower cut-off at least until $x=10$, so the emulation  is  at least 6  storeys high by that measure.
 \newpage {}

\end{document}